\documentclass[10pt,a4paper]{article}
\usepackage{amsmath}
\usepackage{amssymb}
\usepackage{color}

\usepackage{graphicx}
\usepackage{amsfonts,amsmath, amssymb}
\allowdisplaybreaks \allowdisplaybreaks[2]

\numberwithin{equation}{section}

\usepackage{graphics}
\usepackage{epsfig}
\newtheorem{claim}{\bf \t}[part]

\newtheorem{theorem}{Theorem}[section]
\newtheorem{corollary}[theorem]{Corollary}
\newtheorem{lemma}[theorem]{Lemma}
\newtheorem{proposition}[theorem]{Proposition}
\newtheorem{remark}[theorem]{Remark}

\def\v{\varepsilon}

\def\t{\theta}
\def\k{\kappa}
\def\n{\nu}
\def\m{\mu}
\def\a{\alpha}
\def\b{\beta}
\def\g{\gamma}
\def\d{\delta}
\def\l{\lambda}

\def\r{\rho}
\def\s{\sigma}

\def\Om{\Omega}

\def\i{\infty}
\def\f{\frac}
\def\MZ{\mathcal{Z}}

\begin{document}

	\title{ Uniform regularity for the  free surface compressible Navier-Stokes equations with or without surface tension}
	
	\author{ {Yu Mei$^{\ddag}$\footnote{Corresponding author.  \newline \indent
				Email addresses: ymei@math.cuhk.edu.hk(Yu Mei), yongwang@amss.ac.cn(Yong Wang), zpxin@ims.cuhk.edu.hk(Zhouping Xin)},~~~ Yong Wang$^{\dag}$,~~~ Zhouping Xin$^{\ddag}$}
		\\
		\ \\
		{\small \it $^\dag$Institute of Applied Mathematics, AMSS, CAS, Beijing 100190, China}\\
		\  \\
		{\small \it $^\ddag$Institute of Mathematical Science, The Chinese University of Hong Kong, Shatin, Hong Kong } \\
	}
	
	\date{ }
	
	\maketitle

	\begin{abstract}
		
In this paper, we investigate the uniform regularity of solutions to the 3-dimensional isentropic compressible Navier-Stokes system with free surfaces and study the corresponding asymptotic limits of such solutions to that of the compressible Euler system for vanishing viscosity and surface tension. It is shown that there exists an unique strong solution to the free boundary problem for the compressible Navier-Stokes system in a finite time interval which is independent of the viscosity and the surface tension. The solution is uniform bounded both in $W^{1,\infty}$ and a conormal Sobolev space. It is also shown that the boundary layer for the density is weaker than the one for the velocity field.  Based on such uniform estimates, the asymptotic limits to the free boundary problem for the ideal compressible Euler system with or without surface tension as both the viscosity and the surface tension tend to zero, are established by a strong convergence argument.

		\
		
		
		\
		
		AMS: 35Q35, 35B65, 76N10
	\end{abstract}
	
	\tableofcontents
	
\section{Introduction}

Consider the motion of a viscous isentropic compressible fluid with surface tension effect along a free boundary. It is governed by the following isentropic  compressible Navier-Stokes system written in Eulerian coordinate as
\begin{equation}\label{1.1}
\begin{cases}
\r_t^\v+\mbox{div}(\r^\v u^\v)=0,\\		
(\r^\v u^\v)_t+\mbox{div}(\r^\v u^\v\otimes u^\v)+\nabla p^\v=\mu\v \Delta u^\v+(\mu+\l)\v \nabla \text{div} u^\v,
\end{cases}
x\in \Omega_t,~ t>0,
\end{equation}
where $\Omega_t$ is a simply connected  domain in $\mathbb{R}^3$ occupied by the fluid at time $t\geq0$,  $\r^\v, u^\v$, which are unknowns, represent the density and the velocity field respectively. The pressure $p^\v$ is given by the $\g$-law
\begin{equation}
p^\v=(\r^\v)^{\g},\quad\g>1
\end{equation}The constant viscosity coefficients $\m\v, \l\v$ satisfy the physical restrictions
	\begin{equation}
		\mu>0,~~2\m+3\l>0,
	\end{equation}
where the parameter $\v>0$ is the inverse of the Reynolds number.
We assume the  boundary of $\Omega_t$ is given by
	\begin{equation}
		\Sigma_t=\{x \in \mathbb{R}^3| F^\v(x,t)=0\},
	\end{equation}
where $F^\v(x,t)$ is an unknown function which will be uniquely determined by the velocity field.
To study the well-posedness theory of this free boundary problem, the following two boundary conditions are imposed on $\Sigma_t$, $t>0$. On the one hand, the kinetic boundary condition, which states the fluid particles do not cross the free boundary, reads
\begin{equation}\label{1.4}
\partial_tF^\v+u^\v\cdot\nabla_xF^\v=0, ~\mbox{for}~x\in \Sigma_t.
\end{equation}
On the other hand, we also need the dynamic boundary condition to balance the stress tensor on the both side of the free boundary. When the surface tension is taken into consideration, this boundary condition can be written as
	\begin{equation}\label{1.5}
		p^\v\mathbf{n}^\v=(2\m\v Su^\v+\l\v\mbox{div}u^\v) \mathbf{n}^\v+p_e{\bf n}^\v-\sigma H^\v\mathbf{n}^\v, \quad x\in\Sigma_t
	\end{equation}
which describes the stress tensor of the fluid is proportional to the mean curvature of the free boundary $\Sigma_t$. Here $Su^\v=\frac{1}{2}(\nabla u^\v+\nabla^t u^\v)$ and  $\mathbf{n}^\v=\frac{\nabla_xF^\v}{|\nabla_xF^\v|}$ denotes the outward normal vector of $\Omega_t$,  $p_e$ is a given constant external pressure, $\sigma$ denotes the surface tension coefficient, $H$ is the double mean curvature of $\Sigma_t$ which can be expressed in the form
\begin{equation}\label{tension1}
H^\v\mathbf{n}^\v=\Delta_{\Sigma_t}(t)x, \,\, x=(x^1,x^2,x^3),
\end{equation}
	where $\Delta_{\Sigma_t}(t)$ is the Laplace-Beltrami operator on $\Sigma_t$. We also impose the initial data for the compressible Navier-Stokes equations \eqref{1.1} as
	\begin{equation} \label{1.12-2}
	(\rho^\v,u^\v,F^\v)(0,x)=(\rho_0^\v,u_0^\v,F^\v_0)(x),\quad x\in\Om_0
	\end{equation}
such that
	\begin{equation}\label{1.12-1}
	0<\f1{C_0}\leq \rho_0^\v\leq C_0<\infty,
	\end{equation}
where $\Om_0$ is a given initial domain determined by $F_0^\v$ and $C_0$ is a given constant independent of $\v$.

The study of fluid motions with free surfaces  is an important topic in fluid dynamics. For incompressible viscous fluids, we refer to the works by  Beale \cite{Beale}, Tani \cite{Tani}, Solonnikov \cite{Solonnikov} and Guo-Tice \cite{Guo-1} and references therein for  the local well-posedness with or without surface tension. As to compressible viscous fluids, the local-wellposedness theory is established by Secchi-Valli \cite{Secchi} without surface tension, and by Solonnikov-Tani \cite{Solonnikov-1}, Zajaczkowski \cite{Zaja,Zaja-1} and Tanaka-Tani \cite{Tanaka} with surface tension and references therein. For the inviscid fluid, it is much more difficult to get the regularity of free boundary. Wu \cite{Wu,Wu-1} made a big breakthrough for the local well-posedness of  irrotational incompressible Euler system in two and three dimensional, we also refer to \cite{Beale-H-L,Craig,Lannes,Ogawa-Tani,Schneider} and the references therein for related works.
Later, Wu \cite{Wu-2} proved an almost global existence result in the 2-D case and   the global existence results in 3-D case is proved by  Wu \cite{Wu-3} and Germain-Masmoudi-Shatah\cite{Germain,Germain-1}.   We also refer to \cite{Lindblad,Coutand, Lindblad-1,Zeng,Zhang} and the references therein. For the compressible Euler system,  Lindblad \cite{Lindblad-2} proved  a local well-posedness result by using Lagrangian coordinates and Nash-Moser construction. Using the theory of symmetric hyperbolic systems, Trakhinin \cite{Trakhinin} provided a different proof for the local existence of solutions. Both the estimates in \cite{Lindblad-2} and \cite{Trakhinin} had derivatives loss.  Recently, Coutand-Hole-Shkoller \cite{Coutand-1} proved the well-posedness for the motion of a compressible liquid with or without surface tension, and with no derivative loss. The zero surface tension limit is also established in \cite{Coutand-1}.

Another classical and interesting problem in the mathematical theory of fluid mechanics is to study the asymptotic limit of the solutions to the Navier-Stokes equation at high Reynold number which corresponds to small viscosity. There has lots of literature on this problem when the domain has no boundaries, see for instances \cite{Constantin,Constantin-1,Kato-1,Masmoudi}.
However, in the presence of   physical boundaries, the problems become much more complicated. When no-slip boundary condition is imposed on the incompressible fluid in a fixed domain, the vanishing viscosity limit of the incompressible Navier-Stokes is one of the major open problems due to the possible appearance of boundary layers, as illustrated by Prandtl's theory. In \cite{S-C-1,S-C-2}, the authors  proved the (local in time) convergence of the incompressible Navier-Stokes flows to the Euler flows outside the boundary layer and to the prandtl flows in the boundary layer at the inviscid limit for the analytic initial data. Recently, Y. Maekawa \cite{Maekawa} proved this limit when the initial vorticity is located away from the boundary  in 2-D half plane. While, for the incompressible Navier-Stokes system with Navier-slip boundary condition in a fixed domain, considerable progress has been made on this problem. Indeed, the uniform $H^3$ bound and a uniform existence time interval as $\v$ tends to zero are obtained by Xiao-Xin in \cite{Xiao-Xin-1} for flat boundaries, which are generalized to $W^{k,p}$ in \cite{Beiro-1,Beiro-2}. However, such results can not be expected for general curved boundaries since boundary layer may appear due to non-trivial curvature as pointed out in \cite{Iftimie}. In such a case, Iftimie and Sueur have proved the convergence of the viscous solutions to the inviscid Euler solutions in $L^\infty(0,T,; L^2)$-space by a careful construction of boundary layer expansions and energy estimates. However, to identify precisely the asymptotic structure and get the convergence in stronger norms such as $L^\infty(0,T; H^s)(s>0)$, further  a priori estimates and analysis are needed. Recently, Masmoudi-Rousset \cite{Masmoudi-R} established conormal uniform estimates for 3-D general smooth domains with the Naiver-slip boundary condition, which, in particular, implies the uniform boundedness of the normal first order derivatives of the velocity field. This allows the authors to obtain the convergence of the viscous solutions to the inviscid ones by a compact argument. Based on the uniform estimates in \cite{Masmoudi-R}, better convergence with rates have been studied in \cite{Gie-Kelliher} and \cite{Xiao-Xin-2}. In particular, Xiao-Xin  \cite{Xiao-Xin-2} has proved the convergence in $L^\infty(0,T; H^1)$ with an rate of convergence.

For the compressible Navier-Stokes equations in fixed domain, Xin-Yanagisawa \cite{Xin-Y} studied the vanishing viscosity limit of the linearized compressible Navier-Stokes system  with the no-slip boundary condition in the 2-D half plane. Recently, Wang-Williams \cite{Wang-Williams} constructed a boundary layer solution of the compressible Navier-Stokes equations with Navier-slip boundary conditions in 2-D half plane. The layers constructed in \cite{Wang-Williams} are of width $O(\sqrt\v)$ as the Prandtl boundary layer, but are of amplitude $O(\sqrt\v)$ which is similar to the one \cite{Iftimie} for the incompressible case. So, in general, it is impossible to obtain the $H^3$ or $W^{2,p}(p>3)$ estimates for the compressible Navier-Stokes system \eqref{1.1} with the generalized Navier-slip boundary condition.  Later, Paddick \cite{Paddick} obtained an existence and conormal Sobolev regularity of strong solutions to the 3-D compressible isentropic Navier-Stokes system on the half-space with a Navier boundary condition. Recently, Wang-Xin-Yong \cite{Wang-Xin-Yong} also obtained an uniform regularity for the solutions of the compressible Navier-Stokes
with general Navier-slip boundary conditions in 3-D domains with curvature, especially, the vanishing viscosity limit of viscous solution to the corresponding inviscid one was also obtained with  rate of convergence in  $L^\infty$. In \cite{Wang-Xin-Yong}, it is also shown that the boundary layer for  density is weaker than the one for velocity fields.

As to the vanishing viscosity limit problem of fluid motion with free surfaces, an interesting result was achieved recently, by Masmoudi-Rousset \cite{Masmoudi-R-1}, that the local existence of solutions to the incompressible Navier-Stokes system with the gravity field but without surface tension in an uniform in $\v$ time interval by using a suitable functional framework based on conormal Sobolev spaces which minimizes the needed amount of normal regularity but which gives a control of the Lipschitz norm of the solution. These regularities of solution are necessary and reasonable for the free boundary problem. Actually, the Lipschitz regularity guarantee that the free surface moves along the partical path into outside of the fluid. And, in the vicinity of the free boundary, the expected behavior of the solution to the incompressible Navier-Stokes equation is $u^\v(x,t)\sim u(x,t)+\sqrt\v U(t,y,\f{z}{\sqrt{\v}})$ so that it is hard to obtain an uniform  $H^{k}$-norm($k\geq2$) estimate for $u^\v$ in a time interval independent of $\v$. In
particular, a new existence result for the incompressible Euler equations can be obtained by strong convergence argument for the Navier-Stokes equations. Later, Elgindi-Lee has tried to study the similar problem in the presence of the surface tension in \cite{Elgindi-Lee}. It should be noted that in the case surface tension be a positive constant, the pressure term in the Euler system becomes less regular. On the other hand, it is also interesting to investigate the zero surface tension limit for free boundary problems. This is nontrivial since that the surface tension coefficient $\sigma$ is connected to the mean curvature of free surface which is a second order derivative term of the boundary function. It is a subtle issue to obtain the uniform estimate in $\sigma$. Ambrose-Masmoudi \cite{Ambrose-M,Ambrose-M-1} studied the zero surface tension limit of irrotational water waves. The zero surface tension limit of incompressible Navier-Stokes equation has been established by Tan-Wang in \cite{Tan-Wang} for small initial data. For the zero surface tension limit of compressible Euler with free surfaces, we can refer Coutand-Hole-Shkoller \cite{Coutand-1}.

In this paper, we are interested in the existence of strong solution to the compressible Navier-Stokes system \eqref{1.1} with free surfaces in a finite interval of time independent of viscosity $\v$ and the surface tension $\sigma$ and the corresponding asymptotic limits problem of vanishing viscosity and zero surface tension. Formally, when the $\v$ tends to zero, the limit of \eqref{1.1} is compressible Euler system with surface tension
\begin{eqnarray}\label{FET}
			\begin{cases}
				\partial_t\rho+\mbox{div}(\rho u)=0,\\
			\rho \partial_tu+\rho (u\cdot\nabla)u+\nabla p=0,
			\end{cases}
			x\in\Om_t
\end{eqnarray}
with the boundary conditions
		\begin{equation}\label{1.3}
			\partial_tF+u\cdot \nabla F=0,~~~\mbox{and}~~p=p_e-\sigma H,~~x\in\Sigma_t,
		\end{equation}
Moreover, when the surface tension coefficient $\s$ also goes to zero, it can be reduce to the compressible Euler system \eqref{FET}
with the boundary conditions
		\begin{equation}
			\partial_tF+u\cdot \nabla F=0,~~~\mbox{and}~~p=p_e,~~x\in\Sigma_t,
		\end{equation} 
Motivated by Masmoudi-Rosset\cite{Masmoudi-R-1}, we aim to obtain the uniform in $\v$ and $\s$ regularity in the anisotropic conormal Sobolev spaces and a control of the Lipschitz norm for solutions to the compressible Navier-Stokes equations \eqref{1.1} with free surfaces. Based on this uniform regularity, we can study the asymptotic limits to the ideal compressible Euler system with free surfaces as both the viscosity and the surface tension tend to zero by a strong compactness argument. Of course, the Taylor sign condition for the pressure on the boundary   is needed to have local well-posedness for the  Euler system. Usually, one expects that the boundary layer for density is weaker than the one for the velocity fields, we also aim to prove this facts in this paper. It should be mentioned that around same time as this research project, similar theory has been just established by Wang-Xin\cite{Wang-Xin} for the incompressible fluids with rigid bottom below. Although there are some common difficulties caused by the appearance of surface tension for both cases, yet, as we shall explain later, there are major differences in these two cases. Different techniques are needed in the compressible case due to the lack of divergence free conditions which are essential in \cite{Masmoudi-R-1,Wang-Xin}.

	\subsection{Reformulation the problem into local coordinates}

 Due to the possible appearance of boundary layers, we shall consider this problem in the conormal Sobolev space as \cite{Masmoudi-R-1}.
 In order to define the conormal derivatives for the free boundary problem,  we first assume that the initial domain $\Omega_0$ has a covering such that
	\begin{equation}\label{domain}
		\Om_0 \Subset\tilde{\Om}_0\cup_{k=1}^{n}\tilde{\Om}_k
	\end{equation}
	with each $\tilde{\Om}_k,k=0,\cdots,n$ being a convex domain. Here, $\tilde{\Om}_0 \Subset\Om_0$ is chosen to satisfy
	\begin{equation}
		\mbox{dist}(\Sigma_0,\tilde{\Om}_0)=2d_0>0.
	\end{equation}
where $\Sigma_0$ is the boundary of $\Om_0$ and $d_0$ is a given positive constant.
For each fixed $k$ and $a>0$, define $\tilde{\Om}_{k,a}$ as
	\begin{equation}
		 \tilde\Om_{k,a}=\{x\in\mathbb{R}^3~|~\mbox{dist}(x,\tilde{\Om}_k)< ad_0\}.
	\end{equation}
Clearly,
	\begin{equation}\label{domain-a}
		\Om_0\Subset\tilde{\Om}_{0,1}\cup_{k=1}^{n}\tilde{\Om}_{k,1}.
	\end{equation}
$\tilde{\Omega}_{k,a}, k=1,\cdots,n$ will be called the boundary covering. Without loss of generality, we  assume that
	\begin{equation*}
		 \Omega_0\cap\tilde\Om_{k,1}={\{x|x_3<\tilde{h}^\v_k(x_1,x_2)\}}\cap\tilde{\Omega}_{k,1},~~\partial\Om_0\cap\tilde{\Om}_{k,1}=\{x|x_3=\tilde{h}^\v_k(x_1,x_2)\}\cap \tilde{\Omega}_{k,1},
	\end{equation*}
since the other cases can be handled similarly.  Define
	\begin{equation*}
		 \Sigma^{k,a}=\{(x_1,x_2)~|~x\in\tilde{\Om}_{k,a}\cap\Om_0,~~ \mbox{for some}~~x=(x_1,x_2,x_3)\in\mathbb{R}^3 \},~~k=1,\cdots,n.
	\end{equation*}
Then, $\Sigma^{k,a}$ is the projection of $\tilde{\Om}_{k,a}\cap\Om_0$ onto the hyperplane $\mathbb{R}_{x_1}\times\mathbb{R}_{x_2}$.

Similarly, for $\Om_t$, we can define $\tilde{h}^\v_k(t,x_1,x_2)$ and $\Sigma^{k,a}_t$. 
\begin{remark}
Under the assumption of the velocity is bounded, it is easy to know that $\tilde{\Omega}_{0,1},~\tilde{\Omega}_{k,\f18},~k=1,\cdots,n$ is still an effective covering of $\Omega_t$ at least in a short time $t\in[0,T_0]$, i.e.
 	$$\Om_t\Subset\tilde{\Om}_{0,1}\cup_{k=1}^{n}\tilde{\Om}_{k,\f18},~~\forall t\in[0,T_0],~\mbox{with}~T_0\ll1,$$
where $T_0$ depends on the upper bound of velocity and $d_0$.
 \end{remark}
	
As in \cite{Masmoudi-R-1}, to deal with the free boundary problem, one can reduce the problem in each boundary covering into a fixed domain and use the local coordinate to define the conormal derivatives. Instead of using Lagrangian variables, we define a family of differmorphism $\Phi_k^\v(t,\cdot)$ by the following process. Since $\tilde{h}_k^\v$ is defined locally, we first  extend  $\tilde{h}_k^\v$ to $\mathbb{R}^2$ by multiplying a cut-off function
\begin{equation}
h_k^\v(t,y)=\tilde\psi_k(y)\tilde{h}^\v_k(t,y),
\end{equation}
where $\tilde\psi_k(y)$ is a smooth function satisfying
\begin{equation}
\tilde\psi_k(y)=
\begin{cases}
1,~~y\in \Sigma_{k,\f34}(0),\\
0,~~y\in \Sigma_{k,\f78}(0).
\end{cases}
\end{equation}
Then,  define $\eta_k^\v$ as 
\begin{equation}\label{1.16}
\hat{\eta}^\v_k(t,\xi,z)=\kappa(z\langle\xi\rangle)\hat{h}^\v_{k}(t,\xi),
\end{equation}
where $\hat{\cdot}$ denotes the Fourier transform with respect to $y$ and $\kappa$ is defined by
\begin{equation}\label{1.17}
\kappa(z\langle\xi\rangle)=e^{-z^2{\langle\xi\rangle}^2}=e^{-z^2(1+|\xi|^2)},
\end{equation}
which implies that
\begin{equation}
\eta_k^\v(t,y,z)=\frac{1}{2}e^{-z^2}\int_{\mathbb{R}^2}\frac{1}{z^2}e^\frac{|y-\tilde{y}|^2}{4z^2}h^\v_k(t,\tilde{y})d\tilde{y}.
\end{equation}
As in \cite{Masmoudi-R-1}, $\varphi_k^\v$ can be chosen as
\begin{equation}\label{1.18-0}
\varphi^\v_k(t,y,z)=Az+\eta^\v_k(t,y,z),
\end{equation}
where $A>0$ is a constant chosen such that
\begin{equation}
\partial_z\varphi^\v_k(0,y,z)\geq 1, ~(y,z)\in\mathcal{S}.
\end{equation}	
Therefore, the differmorphism $\Phi_k^\v$ can be defined by 
\begin{align}\label{1.14}
\Phi_k^\v(t,\cdot): \mathcal {S}&=\mathbb{R}^2\times (-\infty,0)\rightarrow D_{k}(t)\supset \Om_t\cap\tilde{\Om}_{k,\f58}\nonumber\\
&(y,z)\mapsto x:=(y,\varphi_k^\v(t,y,z)),
\end{align}
where $D_{k}(t)$ is defined by 
\begin{align}\label{1.27}
D_{k}(t):=\{x=(y,x_3)~|~y\in \mathbb{R}^2,~x_3< h^\v_k(t,y) \}.
\end{align}
By the above parameterizations, the free surface becomes $z=0$ locally. 
Although we have defined an extension of $\tilde{h}_k^\v$ in $\mathbb{R}^2\times(-\infty,0)$, but it will be more convenient to work in a vicinity of  the boundary, i.e.
	\begin{equation}
		\mathcal{S}_{k,a}(t)=\{(y,z)~|(y,z)\in\mathbb{R}^3~~ \mbox{such that}~x=(y,\varphi_k^\v(t,y,z))\in \tilde{\Om}_{k,a}\}\subset \mathcal{S},~~~k=1,\cdots,n.
	\end{equation}
with  $0\leq a\leq \f34$.  Since the velocity will be uniform bounded, one has that 
\begin{align}\label{10.1}
\mathcal{S}_{k,\f18}(t)\subset\mathcal{S}_{k,\f14}(0)\subset\mathcal{S}_{k,\f38}(t)\subset\mathcal{S}_{k,\f12}(0)\subset\mathcal{S}_{k,\f58}(t)\subset\mathcal{S}_{k,\f34}(0),~~\forall t\in[0,T_0],~~T_0\ll1.
\end{align}

Since $\Sigma_t$ is given locally by $x_3=h(t,y)$ (Henceforth the subscript $k$ and the supscript $\v$ in $\rho^\v, u^\v,h^\v_k$ will be omitted for notational convenience), it is convenient to use the local coordinates $(y,z)$ defined in \eqref{1.14} at least locally.
Once the choice of $\varphi$ is made, based on the observation of \eqref{10.1}, one can reduce locally(neighborhood near the boundary $\Sigma_t$) the problem into the fixed domain $\mathcal{S}_{\f12}(0)$  by setting
\begin{equation}\label{20.1}
 \varrho(t,y,z)=\rho(t,\Phi(t,y,z)),~~~ v(t,y,z)=u(t,\Phi(t,y,z)),~~~(y,z)\in \mathcal{S}_{\f12}(0).
\end{equation}
For simplicity, we will denote also  $(\r,u)$ as $(\varrho,v)$ in the interior domain.
	
Now, we rewrite the system \eqref{1.1} to the one for $(\varrho,v,h)$ in the new coordinate $(y,z)$ of the fixed domain $\mathcal{S}_{\f12}(0)$. First, one introduces the operators $\partial{_i^\varphi}, i=t,1,2,3$ such that
	\begin{equation}
		\partial{_i^\varphi}\varrho=(\partial_i\r)\circ(t,\Phi(t,\cdot)),~~\partial{_i^\varphi}v=(\partial_iu)\circ(t,\Phi(t,\cdot)), 
	\end{equation}
	which directly yields that
	\begin{equation}
		\begin{cases}
			 \partial{_i^\varphi}=\partial_i-\frac{\partial_i\varphi}{\partial_z\varphi}\partial_z,\,\,i=t,1,2,\\
			\partial{_3^\varphi}=\frac{1}{\partial_z\varphi}\partial_z.
		\end{cases}
	\end{equation}
	Then, the compressible Navier-Stokes system \eqref{1.1}  is locally reduced to the system for $(\varrho, v)$  as
	\begin{equation}\label{1.22}
		\begin{cases}
			\partial{_t^\varphi}\varrho+\mbox{div}^\varphi(\varrho v)=0,\\
			\varrho\partial{_t^\varphi}v+\varrho v\cdot\nabla^\varphi v+\nabla^\varphi p(\varrho)=2\m\v\mbox{div}^\varphi(S^\varphi v)+\l\v\nabla^\varphi{div}^\varphi v,
		\end{cases}
		(y,z)\in \mathcal{S}_{\f12}(0),~t>0,
	\end{equation}
	where $\mbox{div}^\varphi v=\sum_{i=1}^{3}\partial{_i^\varphi}v_i$, $\nabla^\varphi=(\partial{_1^\varphi},\partial{_2^\varphi},\partial{_3^\varphi})^t$, $v\cdot\nabla^\varphi=\sum_{i=1}^{3}v_i\partial{_i^\varphi}$ and $S^\varphi v=\frac{1}{2}(\nabla^\varphi v+{\nabla^\varphi v}^t)$.
	And the two boundary conditions read as
	\begin{equation}\label{1.23}
		\partial_t h=v\cdot \mathbf{N}=-v_y(t,y,0)\cdot\nabla_yh+v_3(t,y,0),~(y,z)\in \mathcal{S}_{\f12}(0)
	\end{equation}
	and
	\begin{equation}\label{1.23-1}
		p\mathbf{N}=(2\m\v S^\varphi v+\l\v\mbox{div}^\varphi v)\mathbf{N}+p_e\mathbf{N}-\sigma\nabla_y\cdot\frac{\nabla_y h}{\sqrt{1+|\nabla_y h|^2}}\mathbf{N},~~
		(y,0)\in\mathcal{S}_{\f12}(0),
	\end{equation}
	where $v=(v_y,v_3)^t:=(v_1,v_2,v_3)^t$,  $\nabla=(\nabla_y,\partial_z)^t:=(\partial_1,\partial_2,\partial_z)^t$ and $\mathbf{N}=(-\partial_1h,-\partial_2h,1)^t$. For later use, we define $\mathbf{n}=\f{\mathbf{N}}{|\mathbf{N}|}$.
\begin{remark}
	By using the covering of $\Om_t$, one can always assume that each vector field $(\varrho,v,h)$ is supported in either  $\mathcal{S}_{k,\f12}(0)$ in $(y,z)$ coordinate, or in $\tilde{\Om}_{0,1}$, since this assumption can be achieved by multiplying a  cut-off function $\psi_k\in C_0^{\infty}(\mathbb{R}^3_+)$  such that $\psi_k \equiv1,~(y,z)\in \mathcal{S}_{k,\f14}(0)$ and $\psi_k\equiv0,~(y,z)\in \mathbb{R}^3_+-\mathcal{S}_{k,\f12}(0)$, and $\psi_0\in C_0^{\infty}(\mathbb{R}^3)$ such that $\psi_0(x)\equiv1$ for $x\in\tilde{\Om}_{0,\f34} $ and $\psi_0(x)\equiv0$ for $x\in \mathbb{R}^3-\tilde{\Om}_{0,1} $, respectively. Therefore, using this localization arguments, one can also assume that the problem  \eqref{1.22} holds in   $\mathcal{S}$ or $\mathbb{R}^3$  and \eqref{1.23}-\eqref{1.23-1} hold on $z=0,~y\in\mathbb{R}^2$.
\end{remark}

To measure the regularity of functions defined in $\mathcal{S}$, we shall use the Sobolev conormal spaces as in \cite{Masmoudi-R-1}. Introduce the vector fields
	\begin{equation}
	Z_i=\partial_i,\, i=1,2,~~~Z_3=\frac{z}{1-z}\partial_z.
	\end{equation}
The Sobolev conormal spaces $H{_{co}^m}$ is defined as
 \begin{equation}
 H{_{co}^m}(\mathcal{S})=\{f\in L^2(\mathcal{S})~|~Z^\a f\in L^2(\mathcal{S}),~~|\a|\leq m\},
  \end{equation}
 where $Z^\a=Z_1^{\a_1}Z_2^{\a_2}Z_3^{\a_3}$ with norm defined as
  \begin{equation*}
 \|f\|{_m^2}=\sum_{|\alpha|\leq m}\|Z^\alpha f\|^2.
  \end{equation*}
Similarly, set
  \begin{equation}
  W{_{co}^{m.\infty}}(\mathcal{S})=\{f\in L^\infty(\mathcal{S})~|~Z^\a f\in L^\infty(\mathcal{S}),~~~|\a|\leq m\},
 \end{equation}
  and
  \begin{equation*}
  \|f\|_{m,\infty}=\sum_{|\a|\leq m}\|Z^\a f\|_{L^\infty}.
  \end{equation*}

However, for the compressible Navier-Stokes, the time derivative should be involved in the energy functional. So, we also set
	\begin{equation}
		 Z_0=\partial_0=\partial_t,~~~\mathcal{Z}^\alpha=Z_0^{\a_0}Z{_1^{\a_1}}Z{_2^{\a_2}}Z{_3^{\a_3}},
	\end{equation}
	and, for smooth space-time function $f(t,y,z)$ with $t\geq 0$, $(y,z)\in S$, set the following notations
	\begin{equation*}
		\|f(t)\|{_{\mathcal{H}^m}^2}=\sum_{|\alpha|\leq m}\|\mathcal{Z}^\alpha f(t)\|^2,~~\|f(t)\|{_{\mathcal{H}^{k,\infty}}^2}=\sum_{|\alpha|\leq k}\|\mathcal{Z}^\alpha f(t)\|^2_{L^\i},
	\end{equation*}
	\begin{equation*}
	\|\nabla \MZ^mf(t)\|^2=\sum_{|\a|\leq m} \|\nabla \MZ^\a f(t)\|^2,
	\end{equation*}
and
\begin{equation*}
|f(t)|^2_{\mathcal{H}^m}=\sum_{|\a|\leq m} | \MZ^\a f(t,\cdot,z=0)|^2_{L^2},~~|\MZ^m f(t)|^2_s=\sum_{|\a|\leq m}|\MZ^\a f(t,\cdot,z=0)|^2_s.
\end{equation*}
\begin{remark}
	It should emphasized that the norms defined above can be used in each covering although we have not written this fact explicitly.
\end{remark}
\begin{remark}
In the present paper, we will focus on the estimates near the boundary, i.e. in $\mathcal{S}_{\f12}(0)$. The interior estimate is easy to obtain in a similar way by a cut-off approach such that no boundary terms involves, since the conormal Sobolev space is equivalent to the standard Sobolev space in the interior domain. The global estimates follow from collecting all the local boundary ones and the interior one. 
\end{remark}
\textbf{Notations:} Throughout this paper, the positive
generic constants that are independent of $\v,\s$ are denoted by
$c,C$. And $\|\cdot\|$ denotes the standard
$L^2$ norm, and $\|\cdot\|_{H^m}~(m=1,2,3,\cdots)$ to
denote the Sobolev $H^m$ norm. The notation $|\cdot|_s$ will be used for the standard Sobolev $H^s$ norm of functions defined on boundary, and note that this norm involves only tangential  derivatives. $a\lesssim b$ used here denotes $a\leq Cb$ for some positive constant independent of $\v,\s$. $\Lambda(\cdot)$ denotes a polynomial function independent of $\v$ and $\sigma$, and may change from line to line.

\
	
	\subsection{Main results}

The aim of this paper is to get a local well-posedness result for strong solutions to \eqref{1.1}  in an interval of time independent of $\v$ and $\s$ for $\v\in(0,1], ~\s\in[0,1]$. Note that such a result will also imply the local existence of strong solutions for the Euler equations with or without surface tension. As it is well-known  \cite{Coutand,Masmoudi-R-1,Zeng,Wu-1}, in the absence of viscosity and surface tension, we need to require the initial pressure $p_0$ for the compressible Euler flow satisfies the Taylor sign condition (see, for example \cite{Taylor} and \cite{Rayleigh}), i.e.
\begin{equation}\label{Talyor}
0<c_0\leq-\frac{\partial p_0^\v}{\partial \mathbf{N}},~~x\in \Sigma_t.
\end{equation}
In the local coordinate $(y,z)$,   the initial data \eqref{1.12-2} is reduced to
\begin{equation} \label{1.12-2-1}
(\varrho^\v,v^\v,h^\v)(y,z,0)=(\varrho_0^\v,v_0^\v,h^\v_0)(y,z),
\end{equation}
such that
\begin{equation}\label{1.12-1-1}
0<\f1{C_0}\leq \varrho_0^\v\leq C_0<\infty,
\end{equation}
and the Taylor sign condition becomes
\begin{equation}\label{Talyor1}
0<c_0\leq-\partial_z^{\varphi}p_0^\v,~~z=0.
\end{equation}

Since boundary layers may appear in the vicinity of physical boundaries, in order to obtain uniform estimates for solutions to the compressible Navier-Stokes system with free surfaces, one needs to find a suitable functional space. Here, we define the functional space $X_m^\v(T)$ for  $(p,v,h)=(p,v,h)(t)$  as follows:
	\begin{equation}
		X^\v_m(T)=\Big\{(p,v,h)\in L^\infty([0,T],L^2);~~
		\mbox{esssup}_{0\leq t\leq T}\|(p,v,h)(t)\|_{X_m^\v}<+\infty\Big\},
	\end{equation}
	where the norm $\|(\cdot,\cdot,\cdot)\|_{X_m^\v}$ is given by
	\begin{align}
		&\|(p,v,h)(t)\|^2_{X^\v_m}=\|( p, v)(t)\|^2_{\mathcal{H}^m}+\|\nabla(p, v)(t)\|^2_{\mathcal{H}^{m-2}}+\|\Delta p(t)\|^2_{\mathcal{H}^1}+|(h,\sqrt{\s}\nabla_yh)(t)|^2_{\mathcal{H}^m}\nonumber\\
		&~~~~~~~~~~~~~~~ +\|\nabla(p, v)(t)\|^2_{\mathcal{H}^{1,\infty}}+\v\|\nabla(p,v)(t)\|^2_{\mathcal{H}^{m-1}}	 +\v\|\Delta p(t)\|^2_{\mathcal{H}^2}+\v\|\nabla^2v(t)\|^2_{L^\infty}+\v|\MZ^mh|_\frac{1}{2}^2.
	\end{align}
	
Note that the {\it a priori} estimates in Theorem \ref{thm3.1} below is obtained in the case that the approximate solution is sufficient smooth up to the boundary. In order to obtain a selfcontained result, one needs  to assume that the approximate initial data  satisfies the boundary compatibility conditions.
Therefore, we set
\begin{align}\label{Initial}
&X_{NS,ap}^{\v,m}=\Big\{(\r,u)\in C^{3m}(\Om_0),F(x(s_1,s_2))\in C^{3m}(U)\Big|~\mbox{The Taylor sign condition }~\eqref{Talyor}~\mbox{holds}; \nonumber\\
&~~~~~~~~~~~~~~\partial_t^k \rho,~\partial_t^k u~\mbox{and}~\partial^k_t F,k=1,\cdots,m~\mbox{are defined   through the Navier-Stokes equations}\nonumber\\
&~~~~~~~~~~~~~~~\eqref{1.1}~ \mbox{and}~\eqref{1.4},~\mbox{respectively};~~\partial_t^k(\rho,u,F),~k=0,1,\cdots,m-1~\mbox{satisfy the } \nonumber\\
&~~~~~~~~~~~~~~~~~~~~~ \mbox{boundary compatibility conditions}~\eqref{1.5}~\Big\},
\end{align}
where the domain  $\Om_0$ is defined through $F$ as previous. Note that
the definition of $X^{\v,m}_{NS,ap}$ is given in the Euler coordinate, but it will be changed automatically  into local coordinate when one evaluates it by some norm, and we should keep this in mind all the time.  Define
\begin{equation}\label{Initial space}
X^{\v,m}_{NS}=\mbox{The closure of}~ X^{\v,m}_{NS,ap}~ \mbox{in the norm }~ \|(\cdot,\cdot,\cdot)\|_{X^\v_m}.
\end{equation}
We assume that the initial data $(\varrho_0^\v,v_0^\v,h_0^\v)\in X^{\v,m}_{NS}$ and satisfies
\begin{equation}\label{1.13-1}
		\sup_{\v\in(0,1],~\s\in[0,1]}\Big\{\|(p_0^\v,v_0^\v,h^\v_0)\|^2_{X^\v_m}+\|\nabla v_0^\v\|^2_{\mathcal{H}^{m-1}}\Big\}\leq \tilde{C}_0,
	\end{equation}
	where $p_0^\v=p(\varrho_0^\v)$, $C_0>0$, $\tilde{C}_0>0$ are positive constants independent of $\v\in(0,1]$ and $\s\in[0,1]$. Thus,  the  initial data $(\varrho^\v_0,u^\v_0)$ is assumed to have higher space regularity and  compatibilities. It should be noted that the initial data may also depends on $\s$, but, for simplicity,  this will not be written explicitly.

Our main result   is as follows:
\begin{theorem}[Uniform Regularity]\label{thm1.1}
Assume that $m\geq 6$,  $\v\in(0,1]$ and $\s\in[0, 1]$. Consider the initial data $(\r^{\v}_0, u^{\v}_0, F_0^{\v})\in X^{\v,m}_{NS}$ given in \eqref{1.12-2} satisfying \eqref{1.12-1-1}, \eqref{Talyor1} and \eqref{1.13-1}.
Then there exists a time $T_0>0$ and $\tilde{C}_1>0$  independent of  $\v$ and $\s$, such that there exists an unique solution $( \rho^\v, u^\v, F^\v)$ to the free surface problem \eqref{1.1}, \eqref{1.4} and \eqref{1.5} which is defined on $[0,T_0]$  and satisfies the following estimates:
          \begin{align}\label{1.18}
           &\sup_{0\leq t\leq T_0}\bigg\{ \|( p^{\v}, v^{\v})(t)\|^2_{\mathcal{H}^m}+\|\nabla(p^{\v}, v^{\v})(t)\|^2_{\mathcal{H}^{m-2}}+\|\Delta p^{\v}(t)\|^2_{\mathcal{H}^1}+|h^{\v}|^2_{\mathcal{H}^m}+\s|\nabla_yh^{\v}|^2_{\mathcal{H}^m}\nonumber\\
           &~~~~~~~~+\|\nabla(p^{\v}, v^{\v})(t)\|^2_{\mathcal{H}^{1,\infty}}+\v\|\nabla(p^{\v},v^{\v})(t)\|^2_{\mathcal{H}^{m-1}}	 +\v\|\Delta p^{\v}(t)\|^2_{\mathcal{H}^2}+\v\|\nabla^2v^{\v}\|^2_{L^\infty}+\v|\MZ^mh|_\frac{1}{2}^2\bigg\}\nonumber\\	 
           &~~~~~~~~+\int_{0}^{T_0}\|\nabla p^{\v}(t)\|^2_{\mathcal{H}^{m-1}}+\|\Delta p^{\v}(t)\|^2_{\mathcal{H}^2}dt
		   +\int_0^{T_0}\|\nabla v^{\v}(t)\|^4_{\mathcal{H}^{m-1}}dt\nonumber\\
           &~~~~~~~~
           +\v\int_{0}^{T_0}\|\nabla v^{\v}(t)\|^2_{\mathcal{H}^m}+\|\nabla^2 v^{\v}(t)\|^2_{\mathcal{H}^{m-2}}dt
		   +\v^2\int_{0}^{T_0}\|\nabla^2v^{\v}(t)\|^2_{\mathcal{H}^{m-1}}dt
		   \leq \tilde{C}_1,
		\end{align}
		and
		\begin{equation}\label{1.19}
			\f1{2C_0}\leq  \varrho^{\v}(t)\leq 2C_0~~\forall t\in [0,T_0],
		\end{equation}
		where $\tilde{C}_1$ only depends on $C_0,~\tilde{C}_0$ and $d_0$.

	\end{theorem}
	
Note that in the above result, time derivatives shall be involved in the energy functional. This is not necessary for incompressible flow without surface tension studied by Masmoudi-Rousset\cite{Masmoudi-R-1}, since the pressure $p$ there solves an elliptic equation with a Dirichlet boundary condition so that elliptic regularity of $p$ in the spatial conormal spaces can be used to obtain the uniform estimates. But, as to the incompressible flow with surface tension studied by Wang-Xin \cite{Wang-Xin}, time derivatives are also needed since the Neumann boundary condition for $p$, which was derived from the momentum equation, shall be used. It should be noted that the main difficulties in \cite{Wang-Xin} for the incompressible flow with surface tension are that there are no any estimates of the highest order time derivatives $\partial_t^mp$. However, the difficulties for the compressible fluid are different, since the pressure $p$ here satisfies a transport equation so that the regularity of $\partial_t^mp$ can be derived by an energy method. As will be seen later, the main difficulties here are to derive delicate estimates of $\text{div}^\varphi v$ and uniform in $\sigma$ estimates (when $\sigma\neq 0$) of the boundary term as \eqref{4.2.15} involving surface tension. We also give some remarks on our uniform regularity in Theorem 1.5 as follows.   
\begin{remark}
Note that  $(\r^\v,u^\v,F^\v)$ in Theorem \ref{thm1.1} is the solution to \eqref{1.1} \eqref{1.4} and \eqref{1.5} given in Eulerian coordinate of $\Om_t$. While, $(\varrho^\v,v^\v,h^\v)$ defined locally in the fixed domain $\mathcal{S}$ is used in all the a priori estimates. Since  the solution has   enough regularity, they can be changed to each other,  equivalently.  These two notations will be used  without confusion throughout this paper.
\end{remark}	

\begin{remark}
The Taylor sign condition \eqref{Talyor} or \eqref{Talyor1} is necessary when one studies both vanishing viscosity and zero surface tension limits to the ideal compressible Euler system. However, when surface tension is considered(for fixed $\sigma>0$),this is not necessary since the surface tension prevent the Rayleigh-Taylor instability。
\end{remark}

\begin{remark}
One can construct a class of initial data to satisfy the compatibility conditions\eqref{Initial}. e.g. Let $\Om_0=B_{R_0}$ be a ball with radius $R_0$. When $\sigma>0$ is fixed, let $(\r_0^\v,u_0^\v)$ be sufficiently smooth functions, and $\r_0^\v$ be a positive constant $\r_{0c}$ and $u_0^\v$ vanishes in a vicinity of the boundary. By choosing $p(\r_{0c})=p_e+\f{2\s}{R_0}$, then it is obvious that $(\r_0^\v,u_0^\v)$ satisfies the boundary compatibility conditions. When $\s=0$, it becomes a little complicated due to the Taylor sign condition. In this case, one can assume that the initial data is radial symmetric (near the boundary of ball), which are polynomial functions of radius $R_0-r$, so that the boundary compatibility conditions can be reduce to algebraic equations for coefficients of polynomials. We omit the details here for simplicity.
\end{remark}

\begin{remark}
	It follows from the uniform estimate of $\|\Delta p^\v(t)\|^2_{\mathcal{H}^1}$  that the boundary layer for the density $\r^\v$ is weaker than the one for the velocity $u^\v$ as expected.
\end{remark}

Based on the uniform estimates given in Theorem \ref{thm1.1}, we can justify the vanishing viscosity limit, zero surface tension limit and  the local existence of solution to the free surface problem for the compressible Euler system.
	\begin{theorem}[Inviscid Limit]\label{thm1.2} For any fixed $\s\geq0$,
		under the assumptions in Theorem \ref{thm1.1}, and in addition that there exists $(\varrho_0, u_0,h_0)$ such that
		\begin{equation}
			\lim_{\v\rightarrow 0}\|(\varrho_0^\v-\varrho_0, v_0^\v-v_0)\|_{L^2} +|h^\v_0-h_0|_{L^2}=0.
		\end{equation}
		Then, there exists  $(\rho,u,F)(x,t)$  on the time interval $[0,T_0]$, and satisfies
		\begin{align}\label{1.18-1}
			&\sup_{0\leq t\leq T_0}\bigg\{\|(p,  v)(t)\|^2_{\mathcal{H}^m}+\| \nabla(p,v) (t)\|^2_{\mathcal{H}^{m-2}}+ |(h,\sqrt{\s}\nabla_yh)(t)|^2_{\mathcal{H}^m}+\|\Delta p(t)\|^2_{\mathcal{H}^1}+\|\nabla(p, v)(t)\|^2_{\mathcal{H}^{1,\infty}}\bigg\}\nonumber\\
			&~~~~~~~+\int_0^{T_0}\|\nabla p(t)\|^2_{\mathcal{H}^{m-1}}dt
			+\int_{0}^{T_0}\|\Delta p(t)\|^2_{\mathcal{H}^2}dt+\int_0^{T_0}\|\nabla v (t)\|^4_{\mathcal{H}^{m-1}}dt\leq \tilde{C}_1<\infty,
		\end{align}
		\begin{equation}\label{1.19-1}
			0<\f1{2C_0}\leq \r(t)\leq 2C_0<\infty,~~\forall t\in[0,T_0],
		\end{equation}
		and
		\begin{equation}\label{1.20}
			\lim_{\v\rightarrow 0} \sup_{t\in[0,T_0]}
			\|(\varrho^\v-\varrho, v^\v-v)(t)\|_{L^\infty}+|h^\v-h|_{W^{1,\i}}
			=0,
		\end{equation}
where $(\varrho,v,h)$ is the localized representation of $(\r,u,F)$. Furthermore, $(\r,u,F)$ is the unique solution to the free surface Euler equations
		\begin{eqnarray}\label{FE1}
			\begin{cases}
				\partial_t\rho+\mbox{div}(\rho u)=0,\\
			\rho \partial_tu+\rho (u\cdot\nabla)u+\nabla p=0,
			\end{cases}
			x\in\Om_t
		\end{eqnarray}
		with the boundary conditions
		\begin{equation}\label{FE2}
			\partial_tF+u\cdot \nabla F=0,~~~\mbox{and}~~p=p_e-\sigma H,~~x\in\Sigma_t,
		\end{equation}
	where $\Om_t$ is the domain occupied by the fluid on time  $t\geq0$ and the boundary of $\Om_t$ is given by
$$\Sigma_t=\{x\in\mathbb{R}^3|~F(x,t)=0 \}.$$

	\end{theorem}

	Similarly, by the strong compactness argument in Theorem \ref{thm1.2}, as $\v$ and $\s$ tend to zero independently, one can obtain that
	\begin{theorem}[Inviscid  and Zero Surface Tension Limit]\label{thm1.2-1}
		Under the assumptions of Theorem \ref{thm1.1}, if we assume in addition that there exists $(\varrho_0, u_0, h_0)$ such that
		\begin{equation}
			\lim_{\v\rightarrow 0,\s\rightarrow0}\|(\varrho_0^{\v,\s}-\varrho_0, v_0^{\v,\s}-v_0)\|_{L^2} +|h^{\v,\s}_0-h_0|_{L^2}=0.
		\end{equation}
		Then, there exists  $(\rho,u,F)(t)$  on the time interval $[0,T_0]$, and satisfies
		\begin{align}\label{1.18-1-1}
			&\sup_{0\leq t\leq T_0}\bigg\{\|(p,  v)(t)\|^2_{\mathcal{H}^m}+\| \nabla(p,v) (t)\|^2_{\mathcal{H}^{m-2}}+ |h(t)|^2_{\mathcal{H}^m}+\|\Delta p(t)\|^2_{\mathcal{H}^1}+\|\nabla(p,v) (t)\|^2_{\mathcal{H}^{1,\infty}}\bigg\}\nonumber\\
			&~~~~~~~+\int_0^{T_0}\|\nabla p(t)\|^2_{\mathcal{H}^{m-1}}dt
			+\int_{0}^{T_0}\|\Delta p(t)\|^2_{\mathcal{H}^2}dt+\int_0^{T_0}\|\nabla v (t)\|^4_{\mathcal{H}^{m-1}}dt\leq \tilde{C}_1<\infty,
		\end{align}
		\begin{equation}\label{1.19-1-1}
			0<\f1{2C_0}\leq \r(t)\leq 2C_0<\infty,~~\forall t\in[0,T_0],
		\end{equation}
		and
		\begin{equation}\label{1.20-1}
			\lim_{\v\rightarrow 0} \sup_{t\in[0,T_0]}
			\|(\varrho^{\v,\s}-\varrho, v^{\v,\s}-v)(t)\|_{L^\infty}+|h^{\v,\s}-h|_{W^{1,\i}}
			=0,
		\end{equation}
where $(\varrho,v,h)$ is the localized version of $(\r,u,F)$. Furthermore, $(\r,u,F)$ is the unique solution to the free surface Euler equations
		\begin{eqnarray}\label{E1}
		\begin{cases}
		\partial_t\rho+\mbox{div}(\rho u)=0,\\
		\rho \partial_tu+\rho (u\cdot\nabla)u+\nabla p=0,
		\end{cases}
		x\in\Om_t
		\end{eqnarray}
		with the boundary conditions
		\begin{equation}\label{E2}
		\partial_tF+u\cdot \nabla F=0,~~~\mbox{and}~~p=p_e,~~x\in\Sigma_t,
		\end{equation}
		where $\Om_t$ is the domain occupied by the fluid on time  $t\geq0$ and the boundary of $\Om_t$ is given by
		$$\Sigma_t=\{x\in\mathbb{R}^3|~F(x,t)=0 \}.$$
	\end{theorem}
\begin{remark}
Based on the uniform regularity \eqref{1.18}-\eqref{1.19}, for any fixed $\v>0$, one can also obtain the zero surface tension limit of free surface compressible Navier-Stokes system with surface tension by the similar strong compactness argument above.
\end{remark}

\subsection{Sketch of the proof}
	
Since the classical local existence results to smooth solutions of \eqref{1.1} are available \cite{Secchi,Zaja}, the main difficulty in the proof of Theorem \ref{thm1.1} is to get the {\it a priori } estimates of the solution in a small time interval independent of $\v$ and $\s$. For notational convenient, we drop the superscript $\v$ in the a priori estimates. Define
\begin{align}\label{1.68}
&\varTheta_{m}(p,v,h)(T)=\sup_{0\leq t\leq T}(1+\|(p,v,h)(t)\|^2_{X^\v_m})+\int_{0}^{T}\|\nabla p(t)\|^2_{\mathcal{H}^{m-1}}+\|\Delta p(t)\|^2_{\mathcal{H}^2}+\|\nabla v(t)\|^4_{\mathcal{H}^{m-1}}dt
\nonumber\\
&~~~~~~~~~~~+\v\int_{0}^{T}\|\nabla v(t)\|^2_{\mathcal{H}^m}+\|\nabla^2 v(t)\|^2_{\mathcal{H}^{m-2}}dt
 +\v^2\int_{0}^{T}\|\nabla^2v(t)\|^2_{\mathcal{H}^{m-1}}dt<+\infty.
\end{align}
The major goal of this paper will be to derive a uniform bound of the above quantity in a finite time interval. We outline the main steps and ideas for the a priori estimates as follows:

\

{\it Step 1: Conormal energy estimates for $(p,v,h)$.} The first step will be to estimate $\mathcal{Z}^\a(p,v,h)$ for $|\a|\leq m$. The basic energy estimate($\a=0$) of $(p,v,h)$ is easy to obtain from the total energy identity of compressible Navier-Stokes system. In order to get the estimates for higher order conormal derivatives, one can perform the $L^2$ energy estimates on the Alinhac good unknowns $Q^\a=\mathcal{Z}^\a p-\partial_z^{\varphi}p\mathcal{Z}^\a\eta,~V^\a=\mathcal{Z}^\a v-\partial_z^{\varphi}v\mathcal{Z}^\a\eta$ for $|\a|\neq0$ to overcome the loss of $\frac{1}{2}$ derivative for $h$ on the free surfaces as in \cite{Masmoudi-R-1}, when $\sigma=0$. However, one subtle difficulty arises on the uniform estimate of boundary terms for $V^\alpha$ and $Q^\alpha$. This is because the following boundary estimate
\begin{equation}\label{20.120}
|(\partial_z v)^b|_s\lesssim \Lambda(\|v\|_{1,\i}+|h|_{2,\i})(|v|_{s+1}+|h|_{s+1}),
\end{equation}
which plays an key role in the analysis of \cite{Masmoudi-R-1}, is invalid for the compressible case due to $\text{div}^\varphi v$ is not free. Here and henceforth, $()^b$ denotes the function taking value on the boundary $z=0$. So, in general, the term $\v^2\int_{0}^{t}\|\nabla^2 v\|{_{\mathcal{H}^{m-1}}^2}d\tau$ will appear on the right hand side of the estimate for $V^\alpha$ and $Q^\alpha$, when one applies directly the trace estimates to terms involving $\v^2\int_{0}^{t}|(\partial_z v)^b|_{\mathcal{H}^{m-1}}^2d\tau$(cf. Lemma \ref{lem3.3}). Fortunately, with a $\v^2$ there, this term can be bounded by the dissipation terms in the estimates for normal derivatives and $\text{div}^\varphi v$ in the next steps. When the effect of surface tension($\sigma\neq0$) is taken into consideration, the boundary term involving the mean curvature $\sigma\nabla_y\cdot\frac{\nabla_yh}{\sqrt{1+|\nabla_yh|^2}}$ of the free surface seems difficult to be bounded uniformly , since 
\begin{equation}\label{hvl}
\sigma\int\nabla_y\cdot\MZ^\a\frac{\nabla_yh}{\sqrt{1+|\nabla_yh|^2}}(\partial_t^{m-1}v_y^b\cdot\partial_t\nabla_yh)dy\cong\s|\nabla_y^2h|_{\mathcal{H}^{m-1}}^2,~ \text{or} ~\sigma|v^b|_{\mathcal{H}^{m}}^2,
\end{equation}
which both terms are difficult to be bounded uniformly (cf. Lemma \ref{lem3.4}). To deal with this boundary term, we will take full advantage of the third component of dynamic boundary condition (\ref{4.34}) and momentum equations (\ref{4.5}) to turn it into the volume integral (cf. \eqref{tr0}-\eqref{tr6}). The key observation is that the worst boundary terms appearing in \eqref{tr1} and \eqref{tr6} by this approach will be canceled with each other. It should be reminded that, for any fixed $\s>0$, one can derive the uniform in $\v$ and $\sigma$ estimates without using the Alinhac good unknowns since the boundary regularity can be improved due to the surface tension effect. But, when one studies the zero surface tension limit problem, this improved regularity of boundary disappear so that the Alinhac good unknowns should be used here. Precisely, one can obtain that
\begin{align*}
&\|(V^m,Q^m)(t)\|+|(h,\sqrt{\s}\nabla_yh)(t)|_{\mathcal{H}^m}^2+\v\int_{0}^{t}\|\nabla V^m\|^2d\tau\\
&\leq \Lambda(\tilde{C}_0)+\d\int_{0}^{t}\|\nabla p\|_{\mathcal{H}^{m-1}}^2+\v^2\|\nabla^2v\|_{\mathcal{H}^{m-1}}^2d\tau+\Lambda(\varTheta_m)\int_0^t\|\nabla v\|^2_{\mathcal{H}^{m-1}}d\tau+C_{\d}t\Lambda(\varTheta_m).
\end{align*} 
provided the Taylor sign condition \eqref{Talyor1} holds true.

\

{\it Step 2: Estimates for $(\nabla^{\varphi}p,\mbox{div}^{\varphi}v)$.} Next, one needs to bound $\nabla p$ to close the argument. We shall use an energy method to estimate $\nabla p$, since the pressure $p=\varrho^\g$, shown in $\eqref{1.22}_2$ for compressible flows, satisfies a transport equation due to$\eqref{1.22}_1$. Thus, one also needs to bound $\text{div}^\varphi v$. This is different from the estimate of pressure in \cite{Masmoudi-R-1,Wang-Xin} for incompressible flows, which satisfies an elliptic equation due to the divergence free condition.  Furthermore, since
\begin{equation*}
\partial_z v\cong \mbox{div}^{\varphi}v+\Pi(\partial_zv)+\nabla_yv,
\end{equation*} 
the estimate of $\text{div}^\varphi v$ is also helpful to the estimates of normal derivative $\partial_z v$. It follows from the standard energy method that
\begin{align*}
&\|\nabla v(t)\|{_{\mathcal{H}^{m-2}}^2}+\v\|(\text{div}^\varphi v,\nabla^\varphi p)(t)\|{_{\mathcal{H}^{m-1}}^2}
+\int_{0}^{t}\v\|\nabla^{\varphi} \MZ^{m-2}\text{div}^\varphi v\|^2+\v^2\|\nabla^\varphi \MZ^{m-1}\text{div}^\varphi v\|^2d\tau\nonumber\\
&\qquad\qquad\leq\Lambda(\tilde{C}_0)+\d\v^2\int_{0}^{t}\|\nabla^2 v\|^2_{\mathcal{H}^{m-1}}d\tau
+\Lambda(\varTheta_m)\int_0^t\|\nabla v\|^2_{\mathcal{H}^{m-1}}d\tau+C_{\d}t\Lambda(\varTheta_m),\\
&\|\nabla p(t)\|{_{\mathcal{H}^{m-2}}^2}+\int_{0}^{t}\|\nabla p\|{_{\mathcal{H}^{m-1}}^2}d\tau\leq \Lambda(\tilde{C}_0) +\v^2\int_0^t\|\nabla^2v\|^2_{\mathcal{H}^{m-1}}d\tau+t\Lambda(\varTheta_m),\\
&\|\text{div}^\varphi v\|_{\mathcal{H}^{m-1}}^2\leq \Lambda(\varTheta_m).
\end{align*}
Here the estimate of $\|(\nabla v,\nabla p)(t)\|{_{\mathcal{H}^{m-2}}^2}$ is a consequence of the fundamental theorem of calculus. It should be reminded that the estimate of
\begin{equation}\label{1.63}
\v\|\text{div}^\varphi v(t)\|{_{\mathcal{H}^{m-1}}^2}+\int_{0}^{t}\v\|\nabla^{\varphi} \MZ^{m-2}\text{div}^\varphi v\|^2+\v^2\|\nabla^\varphi \MZ^{m-1}\text{div}^\varphi v\|^2d\tau
\end{equation} 
is the key to bound $\v\|\nabla v(t)\|{_{\mathcal{H}^{m-1}}^2}+\int_{0}^{t}\v\|\nabla^2v\|_{\mathcal{H}^{m-2}}^2+\v^2\|\nabla^2 v\|^2_{\mathcal{H}^{m-1}}d\tau$. This is not necessary for incompressible case, since $\int_{0}^{t}\v\|\nabla^2v\|_{\mathcal{H}^{m-2}}^2+\v^2\|\nabla^2 v\|^2_{\mathcal{H}^{m-1}}d\tau$ does not appear due to \eqref{20.120}.

\

{\it Step 3: Normal derivative estimates Part I.}  Since $\nabla_yv$ and  $\mbox{div}^{\varphi}v$ have been bounded in the above steps, in order to bound $\partial_zv$, it remains to estimate $\Pi(\partial_zv)$, where $\Pi=Id-\mathbf{n}\otimes\mathbf{n}$ denotes the tangential vector field. One can follow the arguments in \cite{Masmoudi-R-1} to derive the estimate for
$$S_{\mathbf{n}}=\Pi(S^\varphi v \mathbf{N}),$$
which is equivalent to $\Pi(\partial_zv)$ but vanishes on the boundary. The main conclusion of this step is
\begin{align}
&\|S_\mathbf{n}(t)\|{_{\mathcal{H}^{m-2}}^2}+\v\|S_\mathbf{n}(t)\|{_{\mathcal{H}^{m-1}}^2}+\int_{0}^{t}\v\|\nabla^\varphi\MZ^{m-2} S_\mathbf{n}\|^2+\v^2\|\nabla^\varphi\MZ^{m-1} S_\mathbf{n}\|^2d\tau\nonumber\\
&\leq\Lambda(\tilde{C}_0)
+\Lambda(\varTheta_m)\int_0^t\|\nabla v\|^2_{\mathcal{H}^{m-1}}d\tau+C_{\d}t\Lambda(\varTheta_m).\nonumber
\end{align}
This, together with \eqref{1.63}, yields immediately that $\v\|\nabla v(t)\|{_{\mathcal{H}^{m-1}}^2}+\int_{0}^{t}\v\|\nabla^2 v\|_{\mathcal{H}^{m-2}}^2+\v^2\|\nabla^2v\|_{\mathcal{H}^{m-1}}^2d\tau,$
which does not appear in the arguments of \cite{Masmoudi-R-1} due to the divergence free condition. It should be noted that $\v\|(\nabla v,\nabla p)(t)\|{_{\mathcal{H}^{m-1}}^2}$ is necessary to close the estimate of $\int_{0}^{t}\v^2\|\nabla^2v\|_{\mathcal{H}^{m-1}}^2d\tau.$

\

{\it Step 4: Normal derivative estimates Part II.}  In order to close the a priori estimate, it remains to bound the $m-1$ order conormal derivatives of $\partial_z v$. As indicated in \cite{Masmoudi-R-1}, one can only expect the estimate of $\int_0^t\|\nabla v\|^4_{\mathcal{H}^{m-1}}d\tau$. Based on the estimate of $\text{div}^\varphi v$ in Step 2 and the fact $\partial_z v\cong \mbox{div}^{\varphi}v+\omega\times\mathbf{N}+\nabla_yv$, it suffices to bound $\int_0^t\|\omega\times\mathbf{N}\|^4_{\mathcal{H}^{m-1}}d\tau.$ Due to the less regularity of $\nabla p$ and $\text{div}^\varphi v$, one chooses to bound $\omega\times\mathbf{N}$ as in \cite{Masmoudi-R-1} instead of $S_\mathbf{n}$, since  $\nabla^\varphi p$ and $\nabla^\varphi\text{div}^\varphi v$ can be eliminated in the equation of the vorticity $\omega$. We will follow the approach of microlocal symmetrizer in \cite{Masmoudi-R-1} to derive this kind of estimate. However, if one applies directly this argument to $\int_0^t\|\omega\|^4_{\mathcal{H}^{m-1}}d\tau$ as \cite{Masmoudi-R-1} so that
\begin{align*}
\left(\int_0^t\|\omega\|^4_{\mathcal{H}^{m-1}}d\tau\right)^{\f12}
\lesssim \Lambda(\tilde{C}_0)+\sqrt{\v}\int_0^t|\omega^b(\tau)|^2_{\mathcal{H}^{m-1}}d\tau+t^{\f12}\Lambda(\varTheta_m).
\end{align*}
Then it seems difficult to bound $\sqrt{\v}\int_0^t|\omega^b(\tau)|^2_{\mathcal{H}^{m-1}}d\tau$, since it involves $|(\partial_zv)^b|_{\mathcal{H}^{m-1}}$ but in general \eqref{20.120} is not valid for the compressible flows. It should be noted that even though it follows from trace theorem that
\begin{equation*}
\sqrt{\v}\int_0^t|\omega(\tau)|^2_{\mathcal{H}^{m-1}}d\tau
\lesssim\sqrt{\v}\int_{0}^{t}\|\nabla\omega\|_{\mathcal{H}^{m-1}}\|\omega\|_{\mathcal{H}^{m-1}}d\tau+t^{\f12}\Lambda(\varTheta_m)\lesssim\v\int_{0}^{t}\|\nabla\omega\|^2_{\mathcal{H}^{m-1}}d\tau+t^{\f12}\Lambda(\varTheta_m),
\end{equation*}
but, $\v\int_{0}^{t}\|\nabla\omega\|_{\mathcal{H}^{m-1}}^2d\tau$ may not be expected to be bounded. However, direct calculations yield(c.f. \eqref{4.4.5} and \eqref{8.17}) that 
$$(\omega\times\mathbf{N})^b=-2\Pi\{(\partial_1v\cdot\mathbf{N},\partial_2v\cdot\mathbf{N},0)^t\}.$$
So, one can expect to bound $\int_{0}^{t}\|\omega\times\mathbf{N}\|_{\mathcal{H}^{m-1}}^4d\tau$. To avoid involving too much regularity of $h$(due to the commutator for $\MZ^{m-1}\Delta^\varphi \omega\times\mathbf{N}$ ), it is equivalent to bound $\int_0^t\|\MZ^{m-1}\omega\times \mathbf{N}\|^4d\tau$ base on the fact 
\begin{equation*}
\int_0^t\|\omega\times \mathbf{N}\|^4_{\mathcal{H}^{m-1}}d\tau\lesssim \int_0^t\|\MZ^{m-1}\omega\times \mathbf{N}\|^4d\tau
+t\Lambda(\varTheta_m).
\end{equation*}
Note that the boundary condition \eqref{1.23-1} yields 
\begin{equation}\nonumber
(\MZ^{m-1}\omega\times\mathbf{N})^b\cong \MZ^{m-1}\nabla_yh+(\MZ^{m-1}\nabla_yv)^b+(\MZ^{m-2}\omega)^b+\mbox{l.o.t}.
\end{equation}
Then, it follows from the argument of the microlocal symmetrizer in \cite{Masmoudi-R-1} that
\begin{align*}
\int_0^t\|\omega\times \mathbf{N}\|^4_{\mathcal{H}^{m-1}}d\tau
&\lesssim \Lambda(\tilde{C}_0)
+\sqrt{\v}\Lambda(\varTheta_m)\int_0^t|(\MZ^{m-1}\omega\times\mathbf{N})|^2_{L^2}d\tau+t\Lambda(\varTheta_m)\\
&\lesssim \Lambda(\tilde{C}_0)+\v\int_0^t\|\nabla V^m\|^2+\|\nabla^2v\|_{\mathcal{H}^{m-2}}^2d\tau+t\Lambda(\varTheta_m).
\end{align*}
Therefore, the conclusion of this step is
\begin{equation}\nonumber
\Big(\int_0^t\|\nabla v\|^4_{\mathcal{H}^{m-1}}\Big)^{\f12}\lesssim \Lambda(\tilde{C}_0)+ \v\int_{0}^{t}\|\nabla V^m\|^2+\|\nabla^2v\|_{\mathcal{H}^{m-2}}^2d\tau +t\Lambda(\varTheta_m).
\end{equation}

{\it Step 5: $L^\infty$-estimates.}  Finally, we will bound $L^\infty$-norm of $\|(\nabla p,\nabla v)\|_{\mathcal{H}^{1,\infty}}$ and $\sqrt{\v}\|\nabla^2v\|_{L^\infty}$ to close the priori estimate. In order to estimate $\|\nabla p\|_{\mathcal{H}^{1,\i}}$, using the anisotropic Sobolev inequality, one only needs to estimate $\|\Delta^\varphi p\|_{\mathcal{H}^1}$. Actually, one can obtain, for $m\geq 6$, that
\begin{equation*}
\Big(\|\Delta^\varphi p(t)\|^2_{\mathcal{H}^1}+\v\|\Delta^\varphi p(t)\|^2_{\mathcal{H}^2}\Big)+\int_{0}^{t}\|\MZ^\a\Delta^\varphi p\|^2d\tau
\lesssim \Lambda(\tilde{C}_0)+t\Lambda(\varTheta_m).
\end{equation*}
This estimate also implies that the boundary layer for density is weaker than the one for velocity.

In order to estimate  $\|\nabla v\|_{\mathcal{H}^{1,\infty}}$ and $\sqrt{\v}\|\nabla^2v\|_{L^\infty}$, it suffices to to estimate $\|S_\mathbf{n}\|_{\mathcal{H}^{1,\i}}^2$ and $\v\|\nabla S_\mathbf{n}\|^2_{L^\infty}$, since
\begin{equation*}
\|\nabla v\|^2_{\mathcal{H}^{1,\infty}}+\v\|\nabla^2v\|^2_{L^\infty}\lesssim \Lambda(\tilde{C}_0)+\|\Delta^\varphi p\|^2_{\mathcal{H}^1}+\|S_\mathbf{n}\|_{\mathcal{H}^{1,\i}}^2+\v\|\nabla S_\mathbf{n}\|^2_{L^\infty}.
\end{equation*}
However, $\|S_\mathbf{n}\|_{\mathcal{H}^{1,\i}}^2$ and $\v\|\nabla S_\mathbf{n}\|^2_{L^\infty}$ seems difficult to be bounded if one performs $L^\i$ estimates directly on the convection diffusion equation solved by $S_\mathbf{n}$ due to the appearance of the term $\|\nabla\text{div}^\varphi v\|_{\mathcal{H}^{2,\i}}^2$, which is difficult to controll. Therefore, we try to estimate
$$\zeta_{\mathbf{n}}=\Pi(\omega\times\mathbf{N})-2\Pi\{(\nabla_y,0)^t\partial_t\eta-(\nabla^{\varphi}\mathbf{N})^tv \},$$
which also solves a convection diffusion equation but involving $\|\nabla\text{div}^\varphi v\|_{\mathcal{H}^{1,\i}}$ only.
By using the similar argument in \cite{Masmoudi-R-1}, one can obtain, for $m\geq 6$, that
\begin{align*}\nonumber
\|\nabla v\|^2_{\mathcal{H}^{1,\infty}}+\v\|\nabla^2v\|^2_{L^\infty} \lesssim \Lambda(\tilde{C}_0)+t^{\f12}\Lambda(\varTheta_m).
\end{align*}

The estimate of $\varTheta_m$ on some uniform time follows by collecting the estimates of the five steps. Note that at the end, we need to check that the Taylor sign condition and $\Phi$ being a differmorphism  remain true.

\

The rest of the paper is organized as follows: In the next section, we collect some elementary facts and inequalities to be used later.  In section 3, we study the equations satisfied by $(\MZ^\a p, \MZ^\a v,\MZ^\a h)$ and prove the estimates for the lower order commutators and boundary values. In section 4,  the a priori estimates Theorem \ref{thm3.1} will be established, which is the main part of this paper. Using the a priori estimates,  Theorem \ref{thm1.1} and Theorems \ref{thm1.2}, \ref{thm1.2-1} are proved  in section   5.  The Appendix collects a generalization of Lemma 14 and Lemma 15 for $L^\i$ estimates in \cite{Masmoudi-R} so that it can be applied to compressible Navier-Stokes equations.

	
	\section{Preliminaries }
		In this section, we will list some basic properties of conormal function spaces and some elementary facts which will be used frequently in the process of a priori estimates.
		For later use, define
		\begin{equation*}
		\mathcal{W}^m([0,T]\times{\mathcal{S}})=\{f(t,y,z)\in L^2([0,T]\times{\mathcal{S}})|~\mathcal{Z}^\alpha f\in L^2([0,T]\times{\mathcal{S}}),~~~|\alpha|\leq m\},
		\end{equation*}
		\begin{equation*}
		\mathcal{W}^{m,\infty}([0,T]\times{\mathcal{S}})=\{f(t,y,z)\in L^\infty([0,T]\times{\mathcal{S}})|~\mathcal{Z}^\alpha f\in L^\infty([0,T]\times{\mathcal{S}}),~|\alpha|\leq m\},
		\end{equation*}
		and the following notations
		\begin{equation*}
		\|f\|^2_{\mathcal{W}_t^{m,\infty}}=\sup_{0\leq\tau\leq t}\sum_{|\alpha|\leq m}\|\mathcal{Z}^\alpha f(\tau)\|^2_{L^\infty(\mathcal{S})}, ~m\geq 1,~\|f\|_{\infty,t}=\sup_{0\leq\tau\leq t}\|f(\tau)\|_{L^\infty(\mathcal{S})}.
		\end{equation*}

	\subsection{General Inequalities}
	The following Gagliardo-Nirenberg-Morser type inequality will be used repeatedly whose proof can be found in \cite{G-1}.
	\begin{lemma}\label{l2.2}
The following products and commutator estimates hold:
		\item[(1)]For $ u,v\in L^\infty([0,T]\times\mathcal{S})\cap \mathcal{W}^k([0,T]\times\mathcal{S})$ with  $k\in \mathbb{N}$ be an integer. It holds for any $0\leq t\leq T$ that
		\begin{equation}\label{2.6}
			 \int_{0}^{t}\|(\mathcal{Z}^{\b}u\mathcal{Z}^{\n}v)(\tau)\|^2d\tau
			\lesssim \|u\|^2_{\infty,t}\int_{0}^{t}\|v(\tau)\|^2_{\mathcal{H}^k}d\tau
			 +\|v\|^2_{\infty,t}\int_{0}^{t}\|u(\tau)\|^2_{\mathcal{H}^k}d\tau,~~|\b|+|\n|=k
		\end{equation}
		\item[(2)]For $1\leq|\alpha|\leq k$, $0\leq t\leq T$, $g\in\mathcal{W}^{k-1}([0,T]\times\mathcal{S})\cap L^\infty([0,T]\times\mathcal{S})$, $f\in \mathcal{W}^k([0,T]\times\mathcal{S})$ such that $\mathcal{Z}f\in L^\infty([0,T]\times\mathcal{S})$, we have
		\begin{equation}\label{2.7}
			\int_{0}^{t}\|[\mathcal{Z}^\alpha,f]g\|^2d\tau\lesssim \|g\|^2_{\infty,t}\int_{0}^{t}\|\mathcal{Z}f\|{_{\mathcal{H}^{k-1}}^2}d\tau+\|\mathcal{Z}f\|^2_{\infty,t}\int_{0}^{t}\|g\|{_{\mathcal{H}^{k-1}}^2}d\tau.
		\end{equation}
		\item[(3)]For $|\alpha|=k\geq 2$,  $0\leq t\leq T$, and $f,~g\in\mathcal{W}^{k-1}([0,T]\times\mathcal{S})\cap \mathcal{W}^{1,\infty}([0,T]\times\mathcal{S})$. Define the symmetric commutator as $[\mathcal{Z},f,g]=\mathcal{Z}^\alpha(fg)-\mathcal{Z}^\alpha fg-f\mathcal{Z}^\alpha g$. Then, it holds for any $0<t<T$ that
		\begin{equation}\label{2.8}
			 \int_{0}^{t}\|[\mathcal{Z}^\alpha,f,g]\|^2d\tau\lesssim\|\mathcal{Z}f\|^2_{\infty,t}\int_{0}^{t}\|\mathcal{Z}g\|{_{\mathcal{H}^{k-2}}^2}d\tau+\|\mathcal{Z}g\|^2_{\infty,t}\int_{0}^{t}\|\mathcal{Z}f\|{_{\mathcal{H}^{k-2}}^2}d\tau.
		\end{equation}
	\end{lemma}

We also need the following anisotropic Sobolev embedding and trace estimates in $\mathcal{S}$ whose proof can be found in \cite{Masmoudi-R-1}:
	\begin{lemma}\label{l2.3}
		Let  $m_1\geq 0,~m_2\geq 0$ be integers,  $f\in H^{m_1}_{co}(\mathcal{S})\cap H^{m_2}_{co}(\mathcal{S})$ and  $\nabla f\in H^{m_2}_{co}(\mathcal{S})$.\\
		1) The following the anisotropic Sobolev embedding holds:
		\begin{equation}\label{2.9}
			\|f\|_{L^\infty}^2\leq C (\|\nabla f\|_{m_2}+\|f\|_{m_2})\|f\|_{m_1}.
		\end{equation}
provided $m_1+m_2\geq 3$.\\
2)	 The following  trace estimate holds:
		\begin{equation}\label{2.10}
			|f|^2{_{H^{s}(\mathbb{R}^2)}}\leq C (\|\nabla f\|_{m_2}+\|f\|_{m_2})\|f\|_{m_1}.
		\end{equation}
		with $m_1+m_2\geq 2s\geq 0$.
	\end{lemma}

The following classical Sobolev inequalities and commutator estimates in $\mathbb{R}^2$ hold:
	\begin{lemma}\label{lem2.4}
		For $s\in\mathbb{R}$, $s\geq 0$, it holds that
		\begin{align}
			&|fg|_{H^s(\mathbb{R}^2)}\leq C_s(|f|_{L^\infty(\mathbb{R}^2)}|g|_{H^s(\mathbb{R}^2)}+|g|_{L^\infty(\mathbb{R}^2)}|f|_{H^s(\mathbb{R}^2)}),\nonumber\\
			&|[\mathcal{F}^s,f]\nabla g|\leq C_s(|\nabla f|_{L^\infty(\mathbb{R}^2)}|g|_{H^s(\mathbb{R}^2)}+|\nabla g|_{L^\infty(\mathbb{R}^2)}|f|_{H^s(\mathbb{R}^2)}),\nonumber\\
			&|uv|_\frac{1}{2}\lesssim|u|_{1,\infty}|v|_\frac{1}{2},\nonumber
		\end{align}
where $\mathcal{F}^s$ is a Fourier multiplier of $(1+|\xi|^2)^{\f s2}$.
	\end{lemma}

\subsection{Estimates involving $\varphi$}
We shall derive some estimates of $\eta$ chosen by \eqref{1.16} and \eqref{1.17} in terms of $h$ and $v$. It shows that $\eta$ has $\frac{1}{2}$ higher order standard Sobolev regularity in $\mathcal{S}$ than $h$.
Note that the different choice of cut-off function $\kappa$ in \eqref{1.17} from the one in \cite{Masmoudi-R-1} allows one to control the $L^2$ norm of $\eta$ itself.
Assume that $A$ is chosen such that $\partial_z\varphi_0(y,z)\geq 1$ at the initial time and
\begin{equation}
\partial_z\varphi(t,y,z)\geq c_0,\quad \forall t\in[0,T^\v],
\end{equation}
	for some $c_0>0$ and small time $T^\v$. With this a priori assumption, one obtains that
	\begin{lemma}\label{l2.4}
		It holds for $\eta$, defined in \eqref{1.16},\eqref{1.17}, that:
		\begin{align}
			&\forall k\in \mathbb{N},\quad \|\eta\|_{H^k(\mathcal{S})}\leq C_k|h|_{k-\frac{1}{2}},\nonumber\\
			&\forall k,l\in \mathbb{N},\quad \|\partial{_t^l}\eta\|_{H^k(\mathcal{S})}\leq C_k|\partial{_t^l}h|_{k-\frac{1}{2}},\nonumber
		\end{align}
		and moreover, we also have the $L^\infty$ estimates
		\begin{align}
			&\forall k\in \mathbb{N},\quad \|\eta\|_{W^{k,\infty}}\leq C_k|h|_{k,\infty},\nonumber\\
			&\forall k,l\in \mathbb{N},\quad \|\partial{_t^l}\eta\|_{W^{k,\infty}}\leq C_k|\partial{_t^l}h|_{k,\infty}.\nonumber
		\end{align}
		
	\end{lemma}
	
	\noindent \textbf{Proof}. It follows from the expressions \eqref{1.16} and \eqref{1.17} that
	\begin{align}
		\|\eta\|^2&
		=\int_{-\infty}^{0}dz\int_{\mathbb{R}^2}|\hat{\eta}(t,\xi,z)|^2d\xi
		 \leq\int_{\mathbb{R}^2}|\hat{h}(t,\xi)|^2d\xi\int_{-\infty}^{0}|\k(z\langle\xi\rangle)|^2dz
		\nonumber\\
		 &\lesssim\int_{\mathbb{R}^2}\langle\xi\rangle^{-1}|\hat{h}(t,\xi)|^2d\xi=|h|_{-\frac{1}{2}}.
	\end{align}
	The proof of remaining parts is similar to \cite{Masmoudi-R}, so we omit the details here. $\hfill\Box$
	
Lemma \ref{l2.2} and \ref{l2.4} yield immediately the following estimates which will be used repeatly, whose proof can be found in \cite{Masmoudi-R-1}.
	\begin{lemma}\label{l2.5}
		For any $m\in\mathbb{N}$, it holds that
		\begin{equation}\label{2.20}
			 \int_{0}^{t}\|\frac{f}{\partial_z\varphi}\|{_{\mathcal{H}^m}^2}d\tau\lesssim\Lambda(\frac{1}{c_0},\|f\|_{\infty,t}+|h|_{\mathcal{W}{^{1,\infty}_t}})\int_{0}^{t}\|f\|{_{\mathcal{H}^m}^2}+|\mathcal{Z}^mh|{_\frac{1}{2}^2}d\tau,
		\end{equation}
		and moreover, for standard Sobolev norms, we also have that, for $0\leq k\leq m$,
		\begin{equation}\label{2.21}
			 \int_{0}^{t}\|\partial{_t^k}(\frac{f}{\partial_z\varphi})\|{_{H^{m-k}}^2}d\tau\leq\Lambda(\frac{1}{c_0},\|f\|_{\infty,t}+|h|_{\mathcal{W}{^{1,\infty}_t}})\int_{0}^{t}(\sum_{k=0}^{m}\|\partial{_t^k}f\|{_{H^{m-k}}^2}+|\partial{_t^k}h|{_{m-k+\frac{1}{2}}^2})d\tau.
		\end{equation}
	\end{lemma}
	
	As a consequence of Lemma \ref{l2.5}, one has
	\begin{corollary}\label{cor2.8}
		For any $m\geq 1$ and any sufficiently smooth function $f$, it holds that
		\begin{align}
			&\int_{0}^{t}\|\nabla^\varphi f\|{_{\mathcal{H}^m}^2}d\tau\leq\Lambda(\frac{1}{c_0}, |h|_{\mathcal{W}{^{1,\i}_t}})\int_{0}^{t}\|\nabla f\|{_{\mathcal{H}^m}^2}d\tau\nonumber\\
			&~~~~~~~~~~~~~~~~~~~~~~~+\Lambda(\frac{1}{c_0},\|\partial_zf\|_{\infty,t}+|h|_{\mathcal{W}{^{2,\i}_t}})\int_{0}^{t}\|\partial_zf\|{_{\mathcal{H}^{m-1}}^2}+|\MZ^mh|{_\frac{1}{2}^2}d\tau,\label{2.23}\\
			&\int_{0}^{t}\|\nabla f\|{_{\mathcal{H}^m}^2}d\tau\leq\Lambda(\frac{1}{c_0}, |h|_{\mathcal{W}{^{1,\i}_t}})\int_{0}^{t}\|\nabla^\varphi f\|{_{\mathcal{H}^m}^2}d\tau\nonumber\\
			&~~~~~~~~~~~~~~~~~~~~~~+\Lambda(\frac{1}{c_0},\|\partial_zf\|_{\infty,t}+|h|_{\mathcal{W}{^{2,\i}_t}})\int_{0}^{t}\|\partial_zf\|{_{\mathcal{H}^{m-1}}^2}+|\MZ^mh|{_\frac{1}{2}^2}d\tau,\label{2.23-1}
		\end{align}
		and
		\begin{align}\label{2.24}
			&\int_{0}^{t}\|\nabla^\varphi\nabla^\varphi f\|{_{\mathcal{H}^m}^2}d\tau\leq\Lambda(\frac{1}{c_0}, |h|_{\mathcal{W}{^{2,\i}_t}})\Big\{\int_{0}^{t}\|\nabla^2 f\|^2_{\mathcal{H}^m}d\tau+\int_{0}^{t}\|\nabla f\|^2_{\mathcal{H}^m}d\tau\nonumber\\
			&~~~~~~~~~~~~~~~~~~~~~~~~~~~+\|\nabla^2 f\|^2_{\infty,t}\int_{0}^{t}|\MZ^mh|{_\frac{1}{2}^2}d\tau
			+\|\nabla f\|^2_{\infty,t}\int_0^t|\nabla_y \mathcal{Z}^m h|^2_{\f12}d\tau\Big\}.
		\end{align}
		
	\end{corollary}
\noindent \textbf{Proof}.
These estimates follow from \eqref{2.7}, \eqref{2.20} and the identity, for $i=1,2,3$,
\begin{equation*}
\MZ^\a(\partial{_i^\varphi}f)=\MZ^\a\partial_if-\frac{\partial_i\varphi}{\partial_z\varphi}\MZ^\a\partial_zf-[\MZ^\a,\frac{\partial_i\varphi}{\partial_z\varphi}]\partial_zf,
\end{equation*}	
where $\partial_i\varphi$ should be replaced by $-1$ when $i=3$.$\hfill\Box$
	\begin{corollary}
		For $1\leq|\a|\leq k$, $0<t<T$, and any sufficiently smooth function $f$, it holds that
		\begin{equation}\label{2.25}
			 \int_{0}^{t}\|[\partial{_i^\varphi},\MZ^\a]f\|^2d\tau\leq\Lambda(\frac{1}{c_0}, \|\nabla f\|_{\infty,t}+|h|_{\mathcal{W}{^{2,\i}_t}})\int_{0}^{t}\|\nabla f\|{_{\mathcal{H}^{k-1}}^2}+|\MZ^k h|{_\frac{1}{2}^2}d\tau.
		\end{equation}
	\end{corollary}
\noindent \textbf{Proof}. This corollary follows from \eqref{2.7}, \eqref{2.20} and the identity, for $i=1,2,3$,
\begin{equation*}
[\partial_i^\varphi,\MZ^\a]f=[\MZ^\a,\frac{\partial_i\varphi}{\partial_z\varphi}](\partial_z f)+\frac{\partial_i^\varphi}{\partial_z\varphi}[\MZ^\a,\partial_z],
\end{equation*}	
where $\partial_i\varphi$ should be replaced by $-1$ when $i=3$.$\hfill\Box$

\
	
The following lemma gives the regularity of free surfaces once we gain the corresponding one of the velocity.
\begin{lemma}\label{vh}
For any $m\in\mathbb{N}$, $\v\in(0,1]$, it holds that
\begin{equation}
\v|\MZ^mh|_\frac{1}{2}^2d\tau\leq \v|\MZ^mh_0|_\frac{1}{2}^2+\v\int_{0}^{t}\|\nabla v\|_{\mathcal{H}^m}^2d\tau+\Lambda(|\nabla_yh|_{\i,t}+\|v\|_{\mathcal{W}^{1,\i}_t})\int_{0}^{t}(\|v\|_{\mathcal{H}^m}^2+\v|\MZ^mh|_\frac{1}{2}^2)d\tau.
\end{equation}
\end{lemma}	
\noindent \textbf{Proof}. The proof of this lemma is similar to Proposition 3.4 in \cite{Masmoudi-R-1}	except for the time derivative involving. By taking time derivative in the boundary condition \eqref{1.23} and following the same argument in \cite{Masmoudi-R-1}, we can prove this lemma. The details are omitted here for brevity. $\hfill\Box$

\

	Since the Jacobian of the change of variable \eqref{1.14} is $\partial_z\varphi$, it is natural to use on $\mathcal{S}$ the following weighted $L^2$ scalar products:
	\begin{equation}
		 \int_\mathcal{S}fgd\mathcal{V}_t,~\mbox{with}~~d\mathcal{V}_t=\partial_z\varphi(t,y,z)dydz.
	\end{equation}
	With this notation, one has the following integration by parts identities for the operators $\partial{_i^\varphi}$:
	\begin{lemma}\label{l3.1}For any sufficient smooth function $f$ and $g$, it holds that
		\begin{equation}
			 \int_\mathcal{S}\partial{_i^\varphi}fgd\mathcal{V}_t=-\int_\mathcal{S}f\partial{_i^\varphi}gd\mathcal{V}_t+\int_{z=0}fg\mathbf{N}_idy,~~~
			i=1,2,3,
		\end{equation}
		\begin{equation}
			 \int_\mathcal{S}\partial{_t^\varphi}fgd\mathcal{V}_t=\partial_t\int_\mathcal{S}fgd\mathcal{V}_t-\int_\mathcal{S}f\partial{_t^\varphi}gd\mathcal{V}_t-\int_{z=0}fg\partial_thdy,
		\end{equation}
		where   $\mathbf{N}=(-\partial_1h,-\partial_2h,1)^t$.
	\end{lemma}

The  following Korn's inequality is also needed(see \cite{Masmoudi-R-1} for the proof).
\begin{lemma}\label{lem20.2}
Assume that $\partial_z\varphi\geq c_0$ and $\|\nabla\varphi\|_{L^\infty}+\|\nabla^2\varphi\|_{L^\infty}\leq \f1{c_0}$ for some $c_0>0$, then there exists $\Lambda_0:=\Lambda(\f1{c_0})$, such that for every $v\in H^1(\mathcal{S})$, it holds that
\begin{equation}
\|\nabla v\|^2_{L^2(\mathcal{S})}\leq \Lambda_0\int_{\mathcal{S}}|\nabla^{\varphi}v|^2d\mathcal{V}_t,
\end{equation}
and
\begin{equation}
\|\nabla v\|^2_{L^2(\mathcal{S})}\leq \Lambda_0\Big(\int_{\mathcal{S}}|S^{\varphi}v|^2d\mathcal{V}_t
+\|v\|^2_{L^2(\mathcal{S})}\Big),
\end{equation}
where $S^\varphi v=\f12(\nabla^{\varphi}v+(\nabla^{\varphi}v)^t)$.
\end{lemma}

	
\section{Equations for higher order conormal derivative terms}
In this section, we derive the equations for $(\MZ^\a v,\MZ^\a p, \MZ^\a \eta)$ involving high order conormal derivatives. As explained in \cite{Masmoudi-R-1}, the commutators between $\MZ^\a$ and $\partial_t^\varphi+v\cdot\nabla^\varphi$ or $\nabla^\varphi p$ in the equations yields a loss of $\frac{1}{2}$ derivative due to the lower regularity of $\varphi$ when $\s$ is not a given constant independent of $\v$. We shall introduce the Alinhac good unknowns\cite{Alinhac} to yields some nonlinear cancellation in the equations so that the commutators are all the lower order terms.  
	\subsection{A commutator estimate}
	In order to perform higher order conormal estimates, we first derive the equations satisfied by $(\mathcal{Z}^\alpha v, \mathcal{Z}^\alpha p)$. We thus need to commute the conoraml vector field $\mathcal{Z}^\alpha$ with each term in the equation \eqref{1.22}. The following commutation relations are needed, which has been shown in \cite{Masmoudi-R-1} except for $Z_0=\partial_t$.
	
	It is easy to check that for $i=0,1,2,3$, and a smooth function $f$, one has
	\begin{equation}\label{4.1}
		 \mathcal{Z}^\alpha\partial{_i^\varphi}f=\partial{_i^\varphi}\mathcal{Z}^\alpha f-\partial{_z^\varphi}f\partial{_i^\varphi}\mathcal{Z}^\alpha\eta+\mathcal{C}{_i^\alpha}(f),
	\end{equation}
	where the commutator $\mathcal{C}{_i^\alpha}(f)$ is given for $\alpha\neq 0$ and $i\neq 3$ by
	\begin{equation}\label{4.2}
		 \mathcal{C}{_i^\alpha}(f)=\mathcal{C}{_{i,1}^\alpha}(f)+\mathcal{C}{_{i,2}^\alpha}(f)+\mathcal{C}{_{i,3}^\alpha}(f),
	\end{equation}
	with
	\begin{align*}
		 &\mathcal{C}{_{i,1}^\alpha}(f)=-[\mathcal{Z}^\alpha, \frac{\partial_i\varphi}{\partial_z\varphi}, \partial_zf],\\
		 &\mathcal{C}{_{i,2}^\alpha}(f)=-\partial_zf[\mathcal{Z}^\alpha,\partial_i\varphi,\frac{1}{\partial_z\varphi}]-\partial_i\varphi(\mathcal{Z}^\alpha(\frac{1}{\partial_z\varphi})+\frac{\mathcal{Z}^\alpha\partial_z\eta}{(\partial_z\varphi)^2})\partial_zf,\\
		 &\mathcal{C}{_{i,3}^\alpha}(f)=-\frac{\partial_i\varphi}{\partial_z\varphi}[\mathcal{Z}^\alpha,\partial_z]f+\frac{\partial_i\varphi}{(\partial_z\varphi)^2}\partial_zf[\mathcal{Z}^\alpha,\partial_z]\eta.
	\end{align*}
	For $i=3$, we only need to replace $\partial_i\varphi$ by $-1$ in the above expression.
The following commutator estimates hold:
	\begin{lemma}\label{l4.1}
		For $1\leq|\alpha|\leq m$, $i=0,1,2,3$, one has that
		\begin{equation}\label{4.3}
			 \int_{0}^{t}\|\mathcal{C}{_i^\alpha}(f)\|^2d\tau\lesssim\Lambda(\frac{1}{c_0},|h|_{\mathcal{W}{^{2,\infty}_t}}+\|\nabla f\|_{\mathcal{W}{^{1,\infty}_t}})\int_{0}^{t}\|\nabla f\|{_{\mathcal{H}^{m-1}}^2}+|\mathcal{Z}^{m-1}h|{_\frac{1}{2}^2}d\tau.
		\end{equation}
	\end{lemma}
	\noindent\textbf{Proof}. The proof is similar to \cite{Masmoudi-R-1} which  involves only spacial norms, since the role of $Z_0=\partial_t$ is the same to tangential derivatives. The details are omitted here. $\hfill\Box$
	
	\
	
	\subsection{Equations satisfied by ($\mathcal{Z}^\alpha v$, $\mathcal{Z}^\alpha p$, $\mathcal{Z}^\alpha\eta$)}
	We now derive the system satisfied by $\mathcal{Z}^\alpha v, \mathcal{Z}^\alpha p, \mathcal{Z}^\alpha \eta$. Alihnac's good unknown are defined as 
	\begin{align}\label{10.2}
V^\alpha=\mathcal{Z}^\alpha v-\partial{_z^\varphi v\mathcal{Z}^\a\eta},~\mbox{and}~~Q^\alpha=\mathcal{Z}^\alpha p-\partial{_z^\varphi p}\mathcal{Z}^\alpha\eta.
	\end{align} 
It should be noted that the following computation work only in a neighborhood of the boundary, but, for notational convenience,  we write it in a general form without confusion.
	\begin{lemma}\label{lem3.2}
		For $1\leq |\alpha|\leq m$, $V^\alpha$ and $Q^\alpha$ satisfy
		\begin{align}\label{4.4}
			&\mbox{div}^\varphi V^\alpha+\frac{1}{\g p}\partial{_t^\varphi}Q^\alpha+\frac{1}{\g p}v\cdot\nabla^\varphi Q^\alpha\nonumber\\
			&=\mathcal{Z}^\alpha\eta(\partial{_z^\varphi}(\frac{1}{\g p})\partial{_t^\varphi}p+\partial{_z^\varphi}(\frac{v}{\g p})\cdot\nabla^\varphi p)
			-\mathcal{C}^\alpha(d)-\frac{1}{\g p}\mathcal{C}{_0^\alpha}(p)-\frac{v}{\g p}\cdot\mathcal{C}^\alpha(p)\nonumber\\
			&~~-[\mathcal{Z}^\alpha,\frac{1}{\g p}]\partial{_t^\varphi}p-[\mathcal{Z}^\alpha,\frac{v}{\g p}]\cdot\nabla^\varphi p\nonumber\\
			&=:\mathcal{R}_d^\a,
		\end{align}
		and
		\begin{align}\label{4.5}
			 &\varrho(\partial{_t^\varphi}+v\cdot\nabla^\varphi)V^\alpha+\nabla^\varphi Q^\alpha-2\m\v\mbox{div}^\varphi S^\varphi V^\alpha-\l\v\nabla^\varphi \mbox{div}^\varphi V^\alpha\nonumber\\
			 &=\mathcal{Z}^\alpha\eta\partial{_z^\varphi}\varrho(\partial{_t^\varphi}v+v\cdot\nabla^\varphi v)+\varrho\mathcal{Z}^\alpha\eta(\partial{_z^\varphi}v\cdot\nabla^\varphi)v-[\mathcal{Z}^\alpha,\varrho](\partial{_t^\varphi}v+v\cdot\nabla^\varphi v)-\mathcal{C}^\alpha(p)\nonumber\\
			 &~~~~-\varrho\mathcal{C}^\alpha(T)+2\m\v{div}^\varphi\mathcal{E}^\alpha(v)+2\m\v\mathcal{D}^\alpha(S^\varphi v)+\l\v\nabla^\varphi\mathcal{C}^\alpha(d)+\l\v\mathcal{C}^\alpha(\mbox{div}^\varphi v)\nonumber\\
			 &=:\mathcal{R}_{M}^\a,
		\end{align}
		where the commutators $\mathcal{C}_0^\alpha(p),~\mathcal{C}^\alpha(p)$, $\mathcal{C}^\alpha(d)$, $\mathcal{C}^\alpha(T)$ and $\mathcal{E}^\alpha(v)$ admit the estimates:
		\begin{align}\label{4.6}
			 &\int_{0}^{t}\|(\mathcal{C}_0^\alpha(p),\mathcal{C}^\a(p))\|^2d\tau\leq\Lambda(\frac{1}{c_0},|h|_{\mathcal{W}{^{2,\infty}_t}}+\|\nabla p\|_{\mathcal{W}{^{1,\infty}_t}})\int_{0}^{t}\|\nabla p\|{_{\mathcal{H}^{m-1}}^2}+|\mathcal{Z}^{m-1}h|{_\frac{1}{2}^2}d\tau,\\\label{4.7}
			 &\int_{0}^{t}\|\mathcal{C}^\alpha(d)\|^2+\|\mathcal{E}^\alpha(v)\|^2d\tau\leq\Lambda(\frac{1}{c_0},|h|_{\mathcal{W}{^{2,\infty}_t}}+\|\nabla v\|_{\mathcal{W}{^{1,\infty}_t}})\int_{0}^{t}\|\nabla v\|{_{\mathcal{H}^{m-1}}^2}+|\mathcal{Z}^{m-1}h|{_\frac{1}{2}^2}d\tau,\\\label{4.8}
			 &\int_{0}^{t}\|\mathcal{C}^\alpha(T)\|^2d\tau\leq\Lambda(\frac{1}{c_0},\|(v,\nabla v)\|_{\mathcal{W}{^{1,\infty}_t}}+|h|_{\mathcal{W}{^{2,\infty}_t}})\int_{0}^{t}\|v\|{_{\mathcal{H}^m}^2}+\|\nabla v\|{_{\mathcal{H}^{m-1}}^2}+|\mathcal{Z}^{m}h|^2_{L^2}d\tau,
		\end{align}
		and 
\begin{align}\label{4.8-1}
&\int_0^t\|(\mbox{div}^{\varphi}\mathcal{E}^\alpha(v),\mathcal{D}^\alpha(S^{\varphi}v),\nabla^\alpha\mathcal{C}^\alpha(d),\mathcal{C}^\alpha(\mbox{div}^\varphi v))\|^2d\tau\\
&\leq \Lambda(\f1{c_0},|h|_{\mathcal{W}{^{3,\infty}_t}}+\|\nabla v\|_{\mathcal{W}{^{1,\infty}_t}}+\sqrt\v\|\nabla^2v\|_{L^{\infty}_t})
\int_0^t\v^2\|\nabla^2 v\|^2_{\mathcal{H}^{m-1}}+\|\nabla v\|^2_{\mathcal{H}^{m-1}}+\v^2|\nabla_y\MZ^{m-1}h|^2_{\f12}d\tau.\nonumber
\end{align}		
	\end{lemma}
	\noindent\textbf{Proof}. It follows from \eqref{4.1} that
	\begin{align}
		&\mathcal{Z}^\alpha\mbox{div}^\varphi v
		=\mbox{div}^\varphi\mathcal{Z}^\alpha v-\partial{_z^\varphi}v\cdot\nabla^\varphi\mathcal{Z}^\alpha\eta+\mathcal{C}^\alpha(d),\label{4.9}\\
		&\mathcal{Z}^\alpha(\frac{1}{\g p}\partial{_t^\varphi}p)
		=\frac{1}{\g p}(\partial{_t^\varphi}\mathcal{Z}^\alpha p-\partial{_z^\varphi}p\partial{_t^\varphi}\mathcal{Z}^\alpha\eta+\mathcal{C}{_0^\alpha}(p))+[\mathcal{Z}^\alpha,\frac{1}{\g p}]\partial{_t^\varphi}p,\label{4.10}\\
		&\mathcal{Z}^\alpha(\frac{v}{\g p}\cdot\nabla^\varphi p)
		=\frac{v}{\g p}\cdot(\nabla^\varphi\mathcal{Z}^\alpha p-\partial{_z^\varphi}p\nabla^\varphi\mathcal{Z}^\alpha\eta+\mathcal{C}^\alpha(p))+[\mathcal{Z}^\alpha,\frac{v}{\g p}]\cdot\nabla^\varphi p,\label{4.11}
	\end{align}
	where $\mathcal{C}^\alpha(d)=\sum_{i=1}^{3}\mathcal{C}{_i^\alpha}(v_i)$, $\mathcal{C}^\alpha(p)=(\mathcal{C}{_1^\alpha}(p),\mathcal{C}{_2^\alpha}(p),\mathcal{C}{_3^\alpha}(p))^t$.  It follows from Lemma \ref{l4.1} that
	\begin{align}
		 &\int_{0}^{t}\|\mathcal{C}^\alpha(d)\|^2d\tau\leq\Lambda(\frac{1}{c_0},|h|_{\mathcal{W}{^{2,\infty}_t}}+\|\nabla v\|_{\mathcal{W}{^{1,\infty}_t}})\int_{0}^{t}\|\nabla v\|{_{\mathcal{H}^{m-1}}^2}+|\mathcal{Z}^{m-1}h|{_\frac{1}{2}^2}d\tau,\nonumber\\
		 &\int_{0}^{t}\|(\mathcal{C}_0^\alpha(p),\mathcal{C}^\a(p))\|^2d\tau\leq\Lambda(\frac{1}{c_0},|h|_{\mathcal{W}{^{2,\infty}_t}}+\|\nabla p\|_{\mathcal{W}{^{1,\infty}_t}})\int_{0}^{t}\|\nabla p\|{_{\mathcal{H}^{m-1}}^2}+|\mathcal{Z}^{m-1}h|{_\frac{1}{2}^2}d\tau,\nonumber
	\end{align}
	which proves \eqref{4.6} and a part of \eqref{4.7}. Adding \eqref{4.9}-\eqref{4.11} and using $\eqref{1.22}_1$,  one obtains  \eqref{4.4}.
	
	Note that
	\begin{equation}\label{4.14}
		 \mathcal{Z}^\alpha(\varrho(\partial{_t^\varphi}v+v\cdot\nabla^\varphi v))=\varrho\mathcal{Z}^\alpha(\partial{_t^\varphi}v+v\cdot\nabla^\varphi v)+[\mathcal{Z}^\alpha,\varrho](\partial{_t^\varphi}v+v\cdot\nabla^\varphi v).
	\end{equation}
Rewrite the transport operator $\partial{_t^\varphi}+v\cdot\nabla^\varphi$ as
	\begin{equation}\label{4.15}
		 \partial{_t^\varphi}+v\cdot\nabla^\varphi=\partial_t+v_y\nabla_y+V_z\partial_z,
	\end{equation}
	where
	\begin{equation}\label{4.16}
		 V_z=\frac{1}{\partial_z\varphi}v_z=\frac{1}{\partial_z\varphi}(v\cdot\mathbf{N}-\partial_t\varphi)=\frac{1}{\partial_z\varphi}(v\cdot\mathbf{N}-\partial_t\eta),
	\end{equation}
	and $\mathbf{N}(t,y,z)$ defined as 
	\begin{equation*}
		 \mathbf{N}(t,y,z)=(-\partial_1\eta(t,y,z),-\partial_2\eta(t,y,z),1)^t.
	\end{equation*}
	is an extension of the outward normal vector to the boundary. Henceforth, we will use it and the outward normal vector to the boundary without confusion.\\
	It follows that
	\begin{equation}\label{4.17}
		 \mathcal{Z}^\alpha(\partial{_t^\varphi}+v\cdot\nabla^\varphi)v=(\partial{_t^\varphi}+v\cdot\nabla^\varphi)\mathcal{Z}^\alpha v-\partial{_z^\varphi}v(\partial{_t^\varphi}+v\cdot\nabla^\varphi)\mathcal{Z}^\alpha\eta+\mathcal{C}^\alpha(T),
	\end{equation}
	where the commutator $\mathcal{C}^\alpha(T)$ is given by
	\begin{equation*}
		\mathcal{C}^\alpha(T)=\sum_{i=1}^{6}T{_i^\alpha},
	\end{equation*}
	and
	\begin{align*}
		& T{_1^\alpha}=[\mathcal{Z}^\alpha,v_y]\nabla_yv,~ T{_2^\alpha}=[\mathcal{Z}^\alpha,V_z,\partial_zv], ~ T{_3^\alpha}=\frac{\partial_z v}{\partial_z\varphi}([\mathcal{Z}^\alpha,v_y]\nabla_y\eta+\mathcal{Z}^\alpha v_z),\\
		& T{_4^\alpha}=(\mathcal{Z}^\alpha(\frac{1}{\partial_z\varphi})+\frac{\mathcal{Z}^\alpha\partial_z\eta}{(\partial_z\varphi)^2})v_z\partial_zv,~ T{_5^\alpha}=-\frac{v_z\partial_zv}{(\partial_z\varphi)^2}[\mathcal{Z}^\alpha,\partial_z]\eta+V_z[\mathcal{Z}^\alpha,\partial_z]v,~  T{_6^\alpha}=[\mathcal{Z}^\alpha,\frac{1}{\partial_z\varphi},v_z]\partial_zv.
	\end{align*}
	Finally, for the viscous term, it follows from \eqref{4.1} that
	\begin{equation}\label{4.19-1}
		\mathcal{Z}^\alpha\mbox{div}^\varphi S^\varphi v=\mbox{div}^\alpha\mathcal{Z}^\alpha S^\varphi v-\partial{_z^\varphi}S^\varphi v\cdot\nabla^\varphi\mathcal{Z}^\alpha\eta+\mathcal{D}^\alpha(S^\varphi v),
	\end{equation}
	where $\mathcal{D}^\alpha(S^\varphi v)_i=\mathcal{C}{_j^\alpha}(S^\varphi v)_{i,j}$.
Note that
	\begin{equation}\label{4.19}
		\mathcal{Z}^\alpha S^\varphi v=S^\varphi\mathcal{Z}^\alpha v-\partial{_z^\varphi}v\otimes\nabla^\varphi\mathcal{Z}^\alpha\eta-\nabla^\varphi\mathcal{Z}^\alpha\eta\otimes\partial{_z^\varphi}v+\mathcal{E}^\alpha(v),
	\end{equation}
	with $(\mathcal{E}^\alpha v)_{ij}=\mathcal{C}{_i^\alpha}(v_j)+\mathcal{C}{_j^\alpha}(v_i)$.
Therefore, it follows from \eqref{4.19-1} and \eqref{4.19} that
	\begin{align}\label{4.20}
	\mathcal{Z}^\alpha\mbox{div}^\varphi S^\varphi v=& \mbox{div}^\varphi S^\varphi\mathcal{Z}^\alpha v-\mbox{div}^\varphi(\partial{_z^\varphi}v\otimes\nabla^\varphi\mathcal{Z}^\alpha\eta+\nabla^\varphi\mathcal{Z}^\alpha\eta\otimes\partial{_z^\varphi}v)\nonumber\\
		 &+\mbox{div}^\varphi\mathcal{E}^\alpha(v)-\partial{_z^\varphi}S^\varphi v\cdot\nabla^\varphi\mathcal{Z}^\alpha\eta+2\m\v\mathcal{D}^\alpha(S^\varphi v).
	\end{align}
	Similarly, one has
	\begin{align}\label{4.22}
		\mathcal{Z}^\alpha\nabla^\varphi\mbox{div}^\varphi v&=\nabla^\varphi\mathcal{Z}^\alpha\mbox{div}^\varphi v-\partial{_z^\varphi}\mbox{div}^\varphi v\nabla^\varphi\mathcal{Z}^\alpha\eta+\mathcal{C}^\alpha(\mbox{div}^\varphi v)\\
		&=\nabla^\varphi\mbox{div}^\varphi\mathcal{Z}^\alpha v-\nabla^\varphi(\partial{_z^\varphi}v\cdot\nabla^\varphi\mathcal{Z}^\alpha\eta)-\partial{_z^\varphi}\mbox{div}^\varphi v\nabla^\varphi\mathcal{Z}^\alpha\eta+\nabla^\alpha\mathcal{C}^\alpha(d)+\mathcal{C}^\alpha(\mbox{div}^\varphi v),\nonumber
	\end{align}
	where \eqref{4.9} has been used in the last equality.
	Collecting \eqref{4.14}, \eqref{4.17}, \eqref{4.20}, \eqref{4.22} and using $\eqref{1.22}_2$ lead to \eqref{4.5}.  Then Lemma \ref{l4.1} yields  that
	\begin{equation}
	\int_{0}^{t}\|\mathcal{E}^\alpha(v)\|^2d\tau\lesssim\Lambda(\frac{1}{c_0},|h|_{\mathcal{W}{^{2,\infty}_t}}+\|\nabla v\|_{\mathcal{W}{^{1,\infty}_t}})\int_{0}^{t}\|\nabla v\|{_{\mathcal{H}^{m-1}}^2}+|\mathcal{Z}^{m-1}h|{_\frac{1}{2}^2}d\tau,\nonumber
	\end{equation}
		and 
		\begin{align}
		&\int_0^t\|(\mbox{div}^{\varphi}\mathcal{E}^\alpha(v),\mathcal{D}^\alpha(S^{\varphi}v),\nabla^\alpha\mathcal{C}^\alpha(d),\mathcal{C}^\alpha(\mbox{div}^\varphi v))\|^2d\tau\nonumber\\
		&\leq \Lambda(\f1{c_0},|h|_{\mathcal{W}{^{3,\infty}_t}}+\|\nabla v\|_{\mathcal{W}{^{1,\infty}_t}}+\sqrt\v\|\nabla^2v\|_{L^{\infty}_t})
		\int_0^t\v^2\|\nabla^2 v\|^2_{\mathcal{H}^{m-1}}+\|\nabla v\|^2_{\mathcal{H}^{m-1}}+\v^2|\nabla_y\MZ^{m-1}h|^2_{\f12}d\tau.\nonumber
		\end{align}	
			
	It remains to estimate the commutator $\mathcal{C}^\alpha(T)$. One first needs some estimates of $v_z$ and $V_z$ which will be used later. It follows from the definition  of $v_z$ and $V_z$ in \eqref{4.16} that
	\begin{align}\label{4.23}
		 &\|v_z\|_{\infty,t}+\|V_z\|_{\infty,t}\leq\Lambda(\frac{1}{c_0},\|v\|_{\infty,t}+\|\nabla\eta\|_{\infty,t}+\|\partial_t\eta\|_{\infty,t})\leq\Lambda(\frac{1}{c_0},\|v\|_{\infty,t}+|h|_{\mathcal{W}{^{1,\infty}_t}}),
		\\ \label{4.24}
		 &\|\mathcal{Z}v_z\|_{\infty,t}+\|\mathcal{Z}V_z\|_{\infty,t}\leq\Lambda(\frac{1}{c_0},\|\mathcal{Z}v\cdot\mathbf{N}+v\cdot\mathcal{Z}\mathbf{N}\|_{\infty,t}+\|\partial_t\mathcal{Z}\eta\|_{\infty,t}+\|v_z\|_{\infty,t}+\|\mathcal{Z}\partial_z\varphi\|_{\infty,t})\nonumber\\
		 &\qquad\qquad\qquad\qquad~~~~~\leq\Lambda(\frac{1}{c_0},\|v\|_{\mathcal{W}{^{1,\infty}_t}}+|h|_{\mathcal{W}{^{2,\infty}_t}}).
	\end{align}
Using \eqref{2.6} and Lemma \ref{l2.4}, one has that
	\begin{align}\label{4.25}
		 &\int_{0}^{t}\|\mathcal{Z}v_z\|{_{\mathcal{H}^{m-2}}^2}d\tau\lesssim\int_{0}^{t}\|\eta\|{_{\mathcal{H}^{m}}^2}+\|v\cdot\mathbf{N}\|{_{\mathcal{H}^{m-1}}^2}d\tau\nonumber\\
		 &\leq\Lambda(\frac{1}{c_0},\|v\|_{\infty,t}+|h|_{\mathcal{W}^{1,\infty}_t})\int_{0}^{t}\|v\|{_{\mathcal{H}^{m-1}}^2}+|h|{_{\mathcal{H}^m}^2}d\tau.
	\end{align}
	It follows from Lemma \ref{l2.5} that
	\begin{align}\label{4.26}
		 &\int_{0}^{t}\|\mathcal{Z}V_z\|{_{\mathcal{H}^{m-2}}^2}d\tau\leq\int_{0}^{t}\|\frac{1}{\partial_z\varphi}\mathcal{Z}v_z\|{_{\mathcal{H}^{m-2}}^2}+\|v_z\mathcal{Z}(\frac{1}{\partial_z\varphi})\|{_{\mathcal{H}^{m-2}}^2}d\tau\nonumber\\
		 &\leq\Lambda(\frac{1}{c_0},\|v_z\|_{\mathcal{W}{^{1,\infty}_t}}+|h|_{\mathcal{W}{^{2,\infty}_t}})\int_{0}^{t}\|(v_z,\mathcal{Z}v_z)\|{_{\mathcal{H}^{m-2}}^2}+|\mathcal{Z}^{m-2}h|{_\frac{1}{2}^2}d\tau+\|\frac{\mathcal{Z}\partial_z\eta}{(\partial_z\varphi)^2}\|{_{\mathcal{H}^{m-2}}^2}d\tau\nonumber\\
		 &\leq\Lambda(\frac{1}{c_0},\|v\|_{\mathcal{W}{^{1,\infty}_t}}+|h|_{\mathcal{W}{^{2,\infty}_t}})\int_{0}^{t}\|v\|{_{\mathcal{H}^{m-1}}^2}+|h|_{\mathcal{H}^m}^2d\tau.
	\end{align}
As a consequence of Lemma \ref{l2.2} and the above estimates for $v_z,V_z$, one gets that
	\begin{equation}\nonumber
		 \int_{0}^{t}\|\mathcal{C}^\alpha(T)\|^2d\tau\leq\Lambda(\frac{1}{c_0},\|(v,\nabla v)\|_{\mathcal{W}{^{1,\infty}_t}}+|h|_{\mathcal{W}{^{2,\infty}_t}})\int_{0}^{t}\|v\|{_{\mathcal{H}^m}^2}+\|\nabla v\|{_{\mathcal{H}^{m-1}}^2}+|h|_{\mathcal{H}^m}^2d\tau.
	\end{equation}
	Therefore, the proof of  Lemma \ref{lem3.2} is completed. $\hfill\Box$

\subsection{Boundary conditions for ($\mathcal{Z}^\alpha v$, $\mathcal{Z}^\alpha p$,  $\mathcal{Z}^\alpha\eta$)}
In the following,	we derive the boundary conditions for $(\mathcal{Z}^\alpha v,\mathcal{Z}^\alpha p,\mathcal{Z}^\alpha\eta)$. Note that the only interesting case occurs when $\alpha_3=0$, since $\mathcal{Z}^\alpha v=\mathcal{Z}^\alpha\eta=\mathcal{Z}^\alpha p=0$ with $\a_3\neq0$ on the boundary. We start with the dynamic boundary condition.
	\begin{lemma}\label{lem3.3}
		Let $1 \leq |\alpha|\leq m$ with $\alpha_3=0$, on the boundary $\{z=0\}$, it holds that
		\begin{align}\label{4.34}
			&(2\m\v S^\varphi V^\alpha+\l\v\mbox{div}^\varphi V^\alpha)\mathbf{N}+\mathcal{Z}^\alpha h(2\m\v\partial{_z^\varphi}S^\varphi v+\l\v\partial{_z^\varphi}\mbox{div}^\varphi v)\mathbf{N}-\sigma\nabla_y\cdot\mathcal{Z}^\alpha\left(\frac{\nabla_yh}{\sqrt{1+|\nabla_yh|^2}}\right)\mathbf{N}\nonumber\\
			&-\mathcal{Z}^\alpha p\mathbf{N}+\left(2\m\v S^\varphi v+\l\v\mbox{div}^\varphi v-\sigma\nabla_y\cdot\frac{\nabla_yh}{\sqrt{1+|\nabla_yh|^2}}-(p-p_e)\right)\mathcal{Z}^\alpha\mathbf{N}=\mathcal{C}^\alpha(B),
		\end{align}
		where the commutator $\mathcal{C}^\alpha(B)$ satisfies
		\begin{equation}\label{4.35}
			 \int_{0}^{t}|\mathcal{C}^\alpha(B)|^2_{L^2}d\tau\leq\Lambda(\frac{1}{c_0},|h|_{\mathcal{W}{^{2,\infty}_t}}+\|\nabla v\|_{\mathcal{W}{^{1,\infty}_t}})\v^2\int_{0}^{t}|(\nabla v)^b|{_{\mathcal{H}^{m-1}}^2}+|\nabla_yh|{_{\mathcal{H}^{m-1}}^2}d\tau.
		\end{equation}
    Here and henceforth, $()^b$ denotes the value of function on the boundary ${z=0}$. 
	\end{lemma}
\begin{remark}
It should be emphasized that the \eqref{4.34} holds only locally on the boundary. But, this is enough since one can assume that the solution $(\varrho,v)$ is supported compactly on a neighborhood of the boundary. While, for the interior domain the boundary condition will not be involved in the conormal estimates.
\end{remark}	
	
	\noindent\textbf{Proof}. Applying the operator $\mathcal{Z}^\alpha$ to $\eqref{1.23-1}$ and using \eqref{4.9}, \eqref{4.19}, and the Alinhac good unknowns, one obtains  \eqref{4.34} with the commutator $\mathcal{C}^\alpha(B)$ of the form
	\begin{equation*}
		 \mathcal{C}^\alpha(B)=[\mathcal{Z}^\alpha,p-p_e+\sigma\nabla_y\cdot\frac{\nabla_yh}{\sqrt{1+|\nabla_yh|^2}}-2\m\v S^\varphi v-\l\v\mbox{div}^\varphi v,\mathbf{N}]-2\m\v\mathcal{E}^\alpha(v)-\l\v\mathcal{C}^\alpha(d).
	\end{equation*}
In order to control  $\mathcal{C}^\alpha(B)$, we first notice,  from \eqref{1.23-1},  that
	\begin{equation*}
		p-p_e+\sigma\nabla_y\cdot\frac{\nabla_y h}{\sqrt{1+|\nabla_y h|^2}}=2\m\v (S^\varphi v\mathbf{n})\cdot\mathbf{n}+\l\v{div}^\varphi v.
	\end{equation*}
	which implies that
	\begin{equation}\label{4.36}
		 \int_{0}^{t}|\mathcal{C}^\alpha(B)|{_{L^2}^2}d\tau\lesssim\v^2\int_{0}^{t}|[\MZ^\alpha,\Pi(S^\varphi v),\mathbf{N}]|{_{L^2(\mathbb{R}^2)}^2}+|\mathcal{E}^\alpha(v)|{_{L^2(\mathbb{R}^2)}^2}+|\mathcal{C}^\alpha(d)|{_{L^2(\mathbb{R}^2)}^2}d\tau.
	\end{equation}
	Here and henceforth, $\Pi:=\text{Id}-\mathbf{n}\otimes\mathbf{n}$ denotes the tangential vector field. It follows from Lemma \ref{l2.2} that
	\begin{align}\label{4.37}
		\int_{0}^{t}|[\MZ^\alpha,\Pi(S^\varphi v),\mathbf{N}]|{_{L^2}^2}d\tau&\lesssim\int_{0}^{t}\|\MZ\Pi(S^\varphi v)\|{_{\infty,t}^2}|\MZ\mathbf{N}|{_{\mathcal{H}^{m-2}}^2}+\|\MZ\mathbf{N}\|{_{\infty,t}^2}|\MZ\Pi(S^\varphi v)|{_{\mathcal{H}^{m-2}}^2}d\tau\nonumber\\
		&\lesssim\Lambda(\frac{1}{c_0},|h|_{\mathcal{W}{^{2,\infty}_t}}+\|\nabla v\|_{\mathcal{W}{^{1,\infty}_t}})\int_{0}^{t}|\nabla_yh|{_{\mathcal{H}^{m-1}}^2}+|\nabla v^b|{_{\mathcal{H}^{m-1}}^2}d\tau.
	\end{align}
	To estimate $\mathcal{E}^\alpha(v)$ and $\mathcal{C}^\alpha(d)$ on the boundary, one needs only to control $\int_{0}^{t}|\mathcal{C}{_i^\alpha}(v_j)|{_{L^2(\mathbb{R}^2)}^2}d\tau$. Note that $\mathcal{C}{_{i,3}^\alpha}(v_j)=0$ since we only consider the case that $\alpha_3=0$. Then, using the similar argument as in the proof of \eqref{4.3}, one gets that
	\begin{equation}\label{4.38}
		 \int_{0}^{t}|\mathcal{C}{_i^\alpha}(v_j)|^2d\tau=\Lambda(\frac{1}{c_0},|h|_{\mathcal{W}{^{2,\infty}_t}}+\|\nabla v\|_{\mathcal{W}{^{1,\infty}_t}})\int_{0}^{t}|(\nabla v)^b|{_{\mathcal{H}^{m-1}}^2}+|\nabla_yh|{_{\mathcal{H}^{m-1}}^2}d\tau.
	\end{equation}
Substituting \eqref{4.37}, \eqref{4.38} into \eqref{4.36}, we prove \eqref{4.35}.  Therefore, the proof of Lemma \ref{lem3.3} is completed. $\hfill\Box$

	\begin{lemma}\label{lem3.4}
		For any $1 \leq |\alpha|\leq m$ such that $\alpha_3=0$,   on the boundary $\{z=0\}$, it holds that 
		\begin{equation}\label{4.41}
			\partial_t\mathcal{Z}^\alpha h-v^b\cdot\mathcal{Z}^\alpha\mathbf{N}-V^\alpha\cdot\mathbf{N}=\mathcal{C}^\alpha(h),
		\end{equation}
		where the commutator $\mathcal{C}^\alpha(h)$ satisfies the estimate
		\begin{align}\label{4.42}
			 \int_{0}^{t}|\mathcal{C}^\alpha(h)|^2_{L^2}d\tau\leq\Lambda(\frac{1}{c_0},\|v\|_{\mathcal{W}{^{1,\infty}_t}}+\|\nabla v\|_{\infty,t}+|h|_{\mathcal{W}{^{2,\infty}_t}})\int_{0}^{t}|h|{_{\mathcal{H}^{m}}^2}+|v{_y^b}|{_{\mathcal{H}^{m-1}}^2}.
		\end{align}
	\end{lemma}
	\noindent\textbf{Proof}. Applying the operator $\mathcal{Z}^\alpha$ to the kinematic boundary condition $\eqref{1.23}_1$ shows that
	\begin{equation}
		\partial_t\mathcal{Z}^\alpha h-v^b\cdot\mathcal{Z}^\alpha\mathbf{N}-V^\alpha\cdot\mathbf{N}=-[\mathcal{Z}^\alpha,v{_y^b},\nabla_yh]+\frac{(\partial_zv)^b}{\partial_z\varphi}\mathcal{Z}^\alpha h\mathbf{N}:=\mathcal{C}^\alpha(h).
	\end{equation}
It follows from Lemma \ref{l2.2}  that
	\begin{equation*}
		 \int_{0}^{t}|\mathcal{C}^\alpha(h)|^2_{L^2}d\tau\leq\Lambda(\frac{1}{c_0},\|v\|_{\mathcal{W}{^{1,\infty}_t}}+\|\nabla v\|_{\infty,t}+|h|_{\mathcal{W}{^{2,\infty}_t}})\int_{0}^{t}|h|{_{\mathcal{H}^{m}}^2}+|v{_y^b}|{_{\mathcal{H}^{m-1}}^2},
	\end{equation*}
which complete the proof of this lemma. $\hfill\Box$

	
\section{A priori Estimates}
	 \renewcommand{\theequation}{\arabic{section}.\arabic{subsection}.\arabic{equation}}

	The aim of this section is to derive the   {\it a priori} estimates in Theorem \ref{thm3.1}, which is a crucial step to prove Theorem \ref{thm1.1}.  We drop the superscript $\v$ throughout this section for notation convenience.

	\begin{theorem}[A priori estimates]\label{thm3.1}
		Let $m$ be an integer satisfying $m\geq 6$, $\v\in(0,1]$, $\s\in[0,1]$, and $(\rho, u, F)$ be sufficiently smooth solution, defined on $[0,T^\v]$, to \eqref{1.1}, \eqref{1.4}, \eqref{1.5}, \eqref{1.12-2}, \eqref{1.12-1}  with the initial data satisfying \eqref{Talyor} and  \eqref{1.13-1}. Then it holds that
		\begin{equation}\label{3.0-3}
		|\varrho_0(y,z)|\exp(-\int_0^t\|\mbox{div}^\varphi v(\tau)\|_{L^\infty}d\tau)\leq\varrho(y,z,t)\leq |\varrho_0(y,z)|\exp(\int_0^t\|\mbox{div}^\varphi v(\tau)\|_{L^\infty}d\tau), \forall t\in[0,T^\v].
		\end{equation}
		Furthermore, there exists a time $T_a>0$  independent of $\s$ and $\v$ such that the following {\it a priori} estimate holds:
		\begin{align}\label{3.0-1}
		&\varTheta_{m}(t)=\sup_{0\leq \tau\leq t}(1+\|(p,v,h)(\tau)\|^2_{X^\v_m})+\int_{0}^{t}\|\nabla p(\tau)\|^2_{\mathcal{H}^{m-1}}+\|\Delta p(\tau)\|^2_{\mathcal{H}^2}+\|\nabla v(\tau)\|^4_{\mathcal{H}^{m-1}}d\tau
		\nonumber\\
		&~~~~~~~~~~~~~~~~~~~~~~~~~~+\v\int_{0}^{t}\|\nabla v(\tau)\|^2_{\mathcal{H}^m}+\|\nabla^2 v(\tau)\|^2_{\mathcal{H}^{m-2}}d\tau
		 +\v^2\int_{0}^{t}\|\nabla^2v(\tau)\|^2_{\mathcal{H}^{m-1}}d\tau\nonumber\\
		&\qquad\qquad\leq \Lambda(C_0, \tilde{C}_0),~~\forall ~t\in[0,\min(T_a,T_\v)]
		\end{align}
		where $T_a$  depends only on $C_0$, $d_0$ and $\tilde{C}_0$.

	\end{theorem}

	\
	
The main difficulties in the proof of Theorem \ref{thm3.1} are to  derive higher order conormal energy estimates for $(p,v,h)$.  To overcome these difficulties, we will apply the Alinhac good unknowns
	\begin{equation}\label{Alin}
	V^\alpha=\mathcal{Z}^\alpha v-\partial{_z^\varphi}v\mathcal{Z}^\alpha\eta,~ ~ ~Q^\alpha=\mathcal{Z}^\alpha p-\partial{_z^\varphi}p\mathcal{Z}^\alpha\eta,~~\alpha\neq 0,
	\end{equation}
and for $\alpha=0$, we denote $V^0=v$ and $Q^0=p$. As emphasized before, the good unknowns are chosen so that the corresponding commutators in the equations for $V^\alpha$ and $Q^\alpha$, which are computed in the previous section, are all lower order terms. The control of these good unknowns and $\mathcal{Z}^\alpha h$, which can be obtained by the standard energy method, yields a control of   $v$, $p$ and $h$.
	Set
	\begin{equation*}
	\|V^m(t)\|^2=\sum_{|\alpha|\leq m}\|V^\alpha(t)\|^2,\quad \|\nabla V^m(t)\|^2=\sum_{|\alpha|\leq m}\|\nabla V^\alpha(t)\|^2,\quad\|Q^m(t)\|^2=\sum_{|\alpha|\leq m}\|Q^\alpha(t)\|^2.
	\end{equation*}
	
	Throughout this section, we assume the following   a priori assumptions hold:
	\begin{equation}\label{5.1}
	0<\f{1}{4C_0}\leq\varrho(t)\leq  4C_0,~ \partial_z\varphi\geq c_0,~  |h|_{\mathcal{H}^{3,\infty}}+|\nabla_y h|_{\mathcal{H}^{[\f{m}2]+1}}\leq\frac{1}{c_0},~-\partial{_z^\varphi}p|_{z=0}\geq\frac{c_0}{2},~\forall t\in[0,T^\v].
	\end{equation}
	Thus, by using the assumption \eqref{5.1}, one has that
	\begin{align}
	&\|v(t)\|{_{\mathcal{H}^m}^2}\leq \|V^m(t)\|^2+\sum_{|\alpha|\leq m}\|\partial{_z^\varphi}v(t)\mathcal{Z}^\alpha \eta(t)\|^2\leq \|V^m(t)\|^2+\Lambda(\frac{1}{c_0},\|\nabla v\|_{L^\infty})|h|_{\mathcal{H}^m}^2,\label{5.2}\\
	&\|p(t)\|_{\mathcal{H}^m}^2\leq \|Q^m(t)\|^2+\sum_{|\alpha|\leq m}\|\partial{_z^\varphi}p(t)\mathcal{Z}^\alpha \eta(t)\|^2\leq \|Q^m(t)\|^2+\Lambda(\frac{1}{c_0},\|\nabla p\|_{L^\infty})|h|_{\mathcal{H}^m}^2\label{5.3}.
	\end{align}

	\begin{remark}
		As explained above, throughout this section, it is always assumed that $(\varrho,v,h)$ has compact support in $\mathbb{R}^3$ or in a vicinity of $\partial\mathcal{S}$,  and the system \eqref{1.22} holds in $\mathbb{R}^3$ or $\mathcal{S}$  because one can realize this assumption by multiplying a cut-off function and at the cost of some source terms(which are easy to estimate in the conormal Sobolev space because at least one derivative is applied to the cut-off function) appearing in the right hand side of \eqref{1.22}.
	\end{remark}

	Throughout this section, we shall work on the interval of time $[0,T^\v]$ such that \eqref{5.1} holds. And we point out that the generic constant $C$ may depend on $\f1{c_0}, ~\mu$ and $\l$ in this section.  Since the proof of Theorem \ref{thm3.1} is very complicated, we divide the proof into the following subsections.

\subsection{Basic energy estimate}

\setcounter{equation}{0}

First, we give the basic energy estimate which corresponds to the physical energy. Since it is not necessary to use the local coordinates to prove the basic energy estimate for \eqref{1.1}, we use the Euler coordinates here. Precisely, the following lemma holds:
	\begin{lemma}\label{lem4.1}
		For a smooth solution to \eqref{1.1}, \eqref{1.4} and \eqref{1.5}, it holds that
		\begin{align}\label{3.1-01}
		& \int_{\Om_t}\frac{1}{2}\rho u^2+\frac{p(\rho)}{\g-1}dx+p_e|\Om_t|+\sigma|\Sigma_t|+\int_0^t\int_{\Om_\tau}2\m\v|Su|^2+\l\v|\text{div} u|^2dxd\tau\nonumber\\
		&\leq \int_{\Om_0}\frac{1}{2}\rho_0 u_0^2+\frac{p(\rho_0)}{\g-1}dx+p_e|\Om_0|+\sigma|\Sigma_0|.
		\end{align}
			where $|\Om_t|$, $|\Sigma_t|$ denote the volume of  $\Om_t$  and surface area of $\Sigma_t$, respectively.
	\end{lemma}
	\noindent\textbf{Proof}. Multiplying $\eqref{1.1}_2$ by $u$ and integrating over $\Om_t$, one gets from integration by parts and the kinetic boundary condition \eqref{1.4} that
	\begin{align}\label{3.2}
	0=&\frac{d}{dt}\int_{\Om_t}\frac{1}{2}\rho u^2dx-\int_{\Om_t}p\text{div}udx+2\mu\v\int_{\Om_t}|Su|^2dx+\l\v\int_{\Om_t}|\text{div}u|^2dx\nonumber\\
	&+\int_{\Sigma_t}(p\mathbf{n}-2\m\v Su\mathbf{n}-\l\v\text{div}u\mathbf{n})\cdot udS.
	\end{align}
	where we have used the mass equation $\eqref{1.1}_1$ and the transport formula. The dynamic boundary condition \eqref{1.5} yields that
	\begin{equation}\label{3.3}
	\int_{\Sigma_t}(p\mathbf{n}-2\m\v Su\mathbf{n}-\l\v\text{div}u\mathbf{n})\cdot udS=\int_{\Sigma_t}(p_e\mathbf{n}-\sigma H\mathbf{n})\cdot udS=\int_{\Om_t}p_e\text{div}udx-\sigma\int_{\Sigma_t}H \mathbf{n}\cdot udS.
	\end{equation}
	Using the mass equation $\eqref{1.1}_1$ and the transport formula, one obtains that
	\begin{equation}\label{3.4}
	 -\int_{\Om_t}p\text{div}udx+\int_{\Om_t}p_e\text{div}udx=\frac{d}{dt}\int_{\Om_t}\frac{\rho^\g}{\g-1}dx+\frac{d}{dt}\int_{\Om_t}p_edx.
	\end{equation}
	Note that the mean curvature of the free surface is give by \eqref{tension1}. Let $\Sigma_t$ be determined by a local coordinate system $x=x(s_1,s_2,t),~ (s_1,s_2)\in U\subset\mathbb{R}^2$ so that its Riemannian metric given by
		\begin{equation*}
			ds=g_{\a\beta} dx^\a dx^\beta,
		\end{equation*}
		where $U$ is an open set, $g_{\a\beta}=x_\a\cdot x_\beta$ with $x_\a=\partial_{s^\a}x$ and the convention summation over the repeated indices has been used. It is well known that the Laplace-Beltrami operator can be written as
		\begin{equation}\label{tension2}
			 \Delta_{\Sigma_t}(t)=\frac{1}{\sqrt{g}}\partial_{s^\a}(\sqrt{g}g^{\a\beta}\partial_{s^\beta}),\quad\alpha, \beta=1,2,
		\end{equation}
		with $g=det{g_{\a\beta}}$ and $(g^{\a\beta})=(g_{\a\beta})^{-1}$.
	Then, one gets that
	\begin{align}\label{3.5}
	&\quad-\sigma\int_{\Sigma_t}H\mathbf{n}\cdot udS=-\sigma\int_{\Sigma_t}g^{-\frac{1}{2}}\partial_{s^\a}(\sqrt{g}g^{\a\beta}\partial_{s^\beta}x)\cdot udS\nonumber\\
	 &=-\sigma\int_U\partial_{s^\a}(\sqrt{g}g^{\a\beta}\partial_{s^\beta}x)\cdot x_tds_1ds_2=\frac{\sigma}{2}\int_U\sqrt{g}g^{\a\beta}\partial_tg_{\a\beta}ds_1ds_2\nonumber\\
	 &=\sigma\frac{d}{dt}\int_U\sqrt{g}ds_1ds_2=\sigma\frac{d}{dt}\int_{\partial\Om_t}dS.
	\end{align}
	Substituting \eqref{3.3}-\eqref{3.5} into \eqref{3.2}, one gets that
	\begin{align}\nonumber
	&\frac{d}{dt}(\int_{\Om_t}\frac{1}{2}\rho u^2+\frac{p(\rho)}{\g-1}dx+p_e|\Om_t|+\sigma|\Sigma_t|)+\int_{\Om_t}2\m\v|Su|^2+\l\v|\text{div} u|^2dx=0.
	\end{align}
	Integrating the above equation, we proved \eqref{3.1-01}. Therefore,
	 the proof of this lemma is completed. $\hfill\Box$

\subsection{Higher order conormal energy estimates}

			\setcounter{equation}{0}

Next, we perform the standard $L^2$ energy estimates on the Alinhac good unknowns $V^\a$ and $Q^\a$ solved by \eqref{4.4}, \eqref{4.5}. However, a new difficulty shall be overcame when the effect of surface tension($\sigma\neq0$) is taken into consideration. It seems difficult to bound the boundary term involving the highest order derivatives of $\sigma\nabla_y\cdot\frac{\nabla_yh}{\sqrt{1+|\nabla_yh|^2}}$ due to the less regularity of $h,v$(cf.\eqref{hvl}). We will use the third equation in the boundary condition (\ref{4.34}) and the equation (\ref{4.5}) to turn this boundary term back into the volume integral. On the other hand, since \eqref{20.120} is invalid for compressible flows, $\v^2\int_{0}^{t}\|\nabla^2 v\|{_{\mathcal{H}^{m-1}}^2}d\tau$ will appear in the estimate of $V^\alpha$ and $Q^\alpha$, when one directly apply the trace estimates to terms involving $\v^2\int_{0}^{t}|(\partial_z v)^b|_{\mathcal{H}^{m-1}}^2d\tau$(cf. Lemma \ref{lem3.3}). Fortunately, with a $\v^2$ there, this term can be bounded by the dissipation terms in the estimates for normal derivatives and $\text{div}^\varphi v$ in the next steps. Precisely, one has the following lemma. 

	\begin{lemma}\label{p5.1}
		For any $t\in[0,T^\v]$ and $m\geq 5$, it holds that
		\begin{align}\label{4.1-1}
			 &\|V^m(t)\|^2+\|Q^m(t)\|^2+|h(t)|{_{\mathcal{H}^m}^2}+\sigma|\nabla_yh(t)|{_{\mathcal{H}^m}^2}+\v\int_{0}^{t}\|\nabla V^m\|^2d\tau\nonumber\\
			 &\leq\Lambda(\f1{c_0},\|V^m(0)\|^2+\|Q^m(0)\|^2+\|v(0)\|^2_{\mathcal{H}^m}+\sigma|\nabla_yh(0)|{_{\mathcal{H}^m}^2}+|h(0)|{_{\mathcal{H}^m}^2})+\Lambda_0\delta\int_{0}^{t}\|\nabla p\|{_{\mathcal{H}^{m-1}}^2}d\tau\nonumber\\
			&\quad +\Lambda_0\delta\v^2\int_{0}^{t}\|\nabla^2 v\|{_{\mathcal{H}^{m-1}}^2}d\tau+C_\d\Lambda_\infty(t)\int_{0}^{t}\Lambda(Y_m(\tau))+\|\nabla v\|{_{\mathcal{H}^{m-1}}^2}+\v|\MZ^mh|_{\frac{1}{2}}^2d\tau,
		\end{align}
		where $\Lambda_0=\Lambda(\frac{1}{c_0}, C_0)$, $\Lambda_\infty(t)=\Lambda(\frac{1}{c_0},C_0,\mathcal{Q}(t))$ and
			\begin{align}
			&\mathcal{Q}(t):=\|(v,\nabla v)\|_{\mathcal{W}{^{1,\infty}_t}}+\|(p,\nabla p)\|_{\mathcal{W}{^{1,\infty}_t}}+|h|_{\mathcal{W}{^{3,\infty}_t}}+\v^\frac{1}{2}\|\partial_{zz}v\|_{\infty,t},\label{5.5}\\[1mm]
			&Y_m(t):= \|(V^m,Q^m)(t)\|^2+\|\nabla(p,v)(t)\|^2_{\mathcal{H}^{m-2}}+|(h,\sqrt{\s}\nabla_yh)(t)|^2_{\mathcal{H}^m}+\|\Delta^{\varphi}p(t)\|^2_{\mathcal{H}^1}\nonumber\\
			&~~~~~~~~~~+\v\|\nabla(p,v)(t)\|^2_{\mathcal{H}^{m-1}}+\v\|\Delta^{\varphi}p(t)\|^2_{\mathcal{H}^2}.\label{5.5-1}
	\end{align}
	\end{lemma}
\noindent\textbf{Proof}. The case for $m=0$ is proved already in Lemma \ref{lem4.1}. Assume that \eqref{4.1-1} is proved for $k\leq m-1$. We shall prove that it holds for $k=m\geq1$.  Multiplying the equation \eqref{4.5} by $V^\alpha$ and integrating over the domain $\mathcal{S}$ yield
	\begin{align}\label{5.4}
		 &\frac{d}{dt}\int_{\mathcal{S}}\frac{1}{2}\varrho|V^\alpha|d\mathcal{V}_t+\int_{\mathcal{S}}2\m\v|S^\varphi V^\alpha|^2+\l\v|{div}^\varphi V^\alpha|^2d\mathcal{V}_t\nonumber\\
		&=\int_{z=0}(2\m\v S^\varphi V^\alpha\mathbf{N}+\l\v{div}^\varphi V^\alpha\mathbf{N}-Q^\alpha\mathbf{N})V^\alpha dy+\int_{\mathcal{S}}Q^\alpha{div}^\varphi V^\alpha d\mathcal{V}_t+\int_{\mathcal{S}}\mathcal{R}_M^\a V^\a d\mathcal{V}_t\nonumber\\
	    &=:\mathcal{B}+J+\int_{\mathcal{S}}\mathcal{R}_M^\a V^\a d\mathcal{V}_t,
	\end{align}
	where we have used the integration by parts and the boundary condition $\eqref{1.23}$.
	
For the boundary term $\mathcal{B}$ in \eqref{5.4}, it follows from \eqref{4.34} that
	\begin{align}
		&\int_{0}^{t}\mathcal{B}d\tau=\sigma\int_{0}^{t}\int_{z=0}\nabla_y\cdot\mathcal{Z}^\alpha(\frac{\nabla_y h}{\sqrt{1+|\nabla_yh|^2}})\mathbf{N}\cdot V^\alpha dyd\tau+\int_{0}^{t}\int_{z=0}\partial{_z^\varphi}p\mathcal{Z}^\alpha h\mathbf{N}\cdot V^\alpha dyd\tau\nonumber\\
		&-\int_{0}^{t}\int_{z=0}\left(2\m\v S^\varphi v+\l\v{div}^\varphi v-\sigma\nabla_y\cdot\frac{\nabla_yh}{\sqrt{1+|\nabla_yh|^2}}-(p-p_e)\right)\mathcal{Z}^\alpha\mathbf{N}\cdot V^\alpha dyd\tau\nonumber\\
		&-\int_{0}^{t}\int_{z=0}\mathcal{Z}^\alpha h\left(2\m\v\partial{_z^\varphi} S^\varphi v+\l\v\partial{_z^\varphi}{div}^\varphi v\right)\mathbf{N}V^\alpha dyd\tau+\int_{0}^{t}\int_{z=0}\mathcal{C}^\alpha(B)V^\alpha dyd\tau=:\sum_{i=1}^{5}K_i,
	\end{align}
	which will be dealt term by term. First, we deal with the most difficult one $K_1$. The kinetic boundary condition \eqref{4.41} yields
	\begin{equation}\label{k1}
		 K_1=\sigma\int_{0}^{t}\int_{z=0}\nabla_y\cdot\MZ^\alpha(\frac{\nabla_yh}{\sqrt{1+|\nabla_yh|^2}})(\partial_t\MZ^\alpha h-v^b\cdot\MZ^\alpha\mathbf{N}-\mathcal{C}^\alpha(h))dyd\tau.
	\end{equation}Note that
	\begin{equation*}
		 \mathcal{Z}^\alpha\frac{\nabla_yh}{\sqrt{1+|\nabla_yh|^2}}=\frac{\mathcal{Z}^\alpha\nabla_yh}{\sqrt{1+|\nabla_yh|^2}}-\frac{\nabla_yh\langle\nabla_yh,\mathcal{Z}^\alpha\nabla_yh\rangle}{(1+|\nabla_yh|^2)^\frac{3}{2}}
		+\mathcal{C}^\a(S),
	\end{equation*}
	where $\mathcal{C}^\alpha(S)=[\MZ^\alpha,\nabla_yh,\frac{1}{\sqrt{1+|\nabla_yh|^2}}]$. It follows from the  integration by parts  that
	\begin{align}\label{pt1}
		&\quad\sigma\int_{0}^{t}\int_{z=0}\nabla_y\cdot(\frac{\mathcal{Z}^\alpha\nabla_yh}{\sqrt{1+|\nabla_yh|^2}}-\frac{\nabla_yh\langle\nabla_yh,\mathcal{Z}^\alpha \nabla_yh\rangle}{(1+|\nabla_yh|^2)^\frac{3}{2}})\partial_t\mathcal{Z}^\alpha hdyd\tau\nonumber\\
		&=-\sigma\int_{0}^{t}\int_{z=0}(\frac{\mathcal{Z}^\alpha \nabla_yh}{\sqrt{1+|\nabla_yh|^2}}-\frac{\nabla_yh\langle\nabla_yh,\mathcal{Z}^\alpha \nabla_yh\rangle}{(1+|\nabla_yh|^2)^\frac{3}{2}})\cdot\partial_t\mathcal{Z}^\alpha \nabla_yhdyd\tau\nonumber\\
		&=-\sigma\int_{0}^{t}\frac{d}{d\tau}\int_{z=0}(\frac{|\mathcal{Z}^\alpha \nabla_yh|^2}{2\sqrt{1+|\nabla_yh|^2}}-\frac{|\langle\nabla_yh,\mathcal{Z}^\alpha \nabla_yh\rangle|^2}{2(1+|\nabla_yh|^2)^\frac{3}{2}})dyd\tau\nonumber\\
		&\quad+\sigma\int_{0}^{t}\int_{z=0}(-\frac{|\mathcal{Z}^\alpha \nabla_yh|^2\langle\nabla_yh,\partial_t\nabla_yh\rangle}{2(1+|\nabla_yh|^2)^\frac{3}{2}}+\frac{3|\langle\nabla_yh,\mathcal{Z}^\alpha \nabla_yh\rangle|^2\langle\nabla_yh,\partial_t\nabla_yh\rangle}{2(1+|\nabla_yh|^2)^\frac{5}{2}})dyd\tau\nonumber\\
		&\quad-\sigma\int_{0}^{t}\int_{z=0}\frac{\langle\nabla_yh,\mathcal{Z}^\a \nabla_yh\rangle\langle\partial_t\nabla_yh,\MZ^\alpha \nabla_yh\rangle}{(1+|\nabla_yh|^2)^\frac{3}{2}}dyd\tau\nonumber\\
		&\leq-\sigma\int_{z=0}\frac{|\mathcal{Z}^\alpha \nabla_yh(t)|^2}{2(1+|\nabla_yh(t)|^2)^\frac{3}{2}}dy+\sigma\int_{z=0}(\frac{|\mathcal{Z}^\alpha \nabla_yh_0|^2}{2\sqrt{1+|\nabla_yh_0|^2}}-\frac{|\langle\nabla_yh_0,\mathcal{Z}^\alpha \nabla_yh_0\rangle|^2}{2(1+|\nabla_yh_0|^2)^\frac{3}{2}})dy\nonumber\\
		&\quad+\Lambda_\infty(t)\sigma\int_{0}^{t}|\nabla_yh|{_{\mathcal{H}^m}^2}d\tau,
	\end{align}
	and
	\begin{align}\label{pt2}
		&\quad-\sigma\int_{0}^{t}\int_{z=0}\nabla_y\cdot(\frac{\mathcal{Z}^\alpha \nabla_yh}{\sqrt{1+|\nabla_yh|^2}}-\frac{\nabla_yh\langle\nabla_yh,\mathcal{Z}^\alpha \nabla_yh\rangle}{(1+|\nabla_yh|^2)^\frac{3}{2}})v^b\cdot\mathcal{Z}^\alpha\mathbf{N}\nonumber\\
		&=\sigma\int_{0}^{t}\int_{z=0}(\frac{\mathcal{Z}^\alpha \nabla_yh}{\sqrt{1+|\nabla_yh|^2}}-\frac{\nabla_yh\langle\nabla_yh,\mathcal{Z}^\alpha \nabla_yh\rangle}{(1+|\nabla_yh|^2)^\frac{3}{2}})\cdot\nabla_y(v^b\cdot\mathcal{Z}^\alpha\mathbf{N})dyd\tau\nonumber\\
		&\leq\sigma\int_{0}^{t}\int_{z=0}(\frac{\mathcal{Z}^\alpha \nabla_yh}{\sqrt{1+|\nabla_yh|^2}}-\frac{\nabla_yh\langle\nabla_yh,\mathcal{Z}^\alpha \nabla_yh\rangle}{(1+|\nabla_yh|^2)^\frac{3}{2}})\cdot(v{_y^b}\cdot\nabla_y)\mathcal{Z}^\alpha\nabla_yhdyd\tau\nonumber\\
		&\quad+\sigma\int_{0}^{t}\int_{z=0}(\frac{\mathcal{Z}^\alpha \nabla_yh}{\sqrt{1+|\nabla_yh|^2}}-\frac{\nabla_yh\langle\nabla_yh,\mathcal{Z}^\alpha \nabla_yh\rangle}{(1+|\nabla_yh|^2)^\frac{3}{2}})\cdot(\nabla_yv^b\cdot\mathcal{Z}^\alpha\mathbf{N})dyd\tau\nonumber\\
		&\leq\sigma\int_{0}^{t}\int_{z=0}\frac{v{_y^b}\cdot\nabla_y(|\mathcal{Z}^\alpha\nabla_y h|^2)}{2\sqrt{1+|\nabla_yh|^2}}-\frac{v{_y^b}\cdot\nabla_y(|\langle\nabla_yh,\mathcal{Z}^\alpha \nabla_yh\rangle|^2)}{2(1+|\nabla_yh|^2)^\frac{3}{2}}dyd\tau\nonumber\\
		&\quad+\sigma\int_{0}^{t}\int_{z=0}\frac{\mathcal{Z}^\alpha \nabla_yh\langle\nabla_yh,\mathcal{Z}^\alpha \nabla_yh\rangle}{(1+|\nabla_yh|^2)^\frac{3}{2}}(v{_y^b}\cdot\nabla_y)\nabla_y hdyd\tau\nonumber\\
		&\quad+\sigma\int_{0}^{t}\int_{z=0}(\frac{\mathcal{Z}^\alpha \nabla_yh}{\sqrt{1+|\nabla_yh|^2}}-\frac{\nabla_yh\langle\nabla_yh,\mathcal{Z}^\alpha \nabla_yh\rangle}{(1+|\nabla_yh|^2)^\frac{3}{2}})\cdot(\nabla_yv{_y^b}\cdot\mathcal{Z}^\alpha \nabla_yh)dyd\tau\nonumber\\
		&\leq\Lambda_\infty(t)\sigma\int_{0}^{t}|\nabla_yh|{_{\mathcal{H}^m}^2}d\tau.
	\end{align}
Using  Lemma \ref{l2.2}, one has that
	\begin{equation}
		 \int_{0}^{t}|\mathcal{C}^\alpha(S)|^2_{L^2}d\tau\leq\Lambda_\infty(t)\int_{0}^{t}|\nabla_yh|{_{\mathcal{H}^{m-1}}^2}d\tau,
	\end{equation}
	and
	\begin{equation}\label{5.10}
		 \int_{0}^{t}|\partial_t\mathcal{C}^\alpha(S)|^2_{L^2}+|\nabla_y\mathcal{C}^\alpha(S)|^2_{L^2}d\tau\leq\Lambda_\infty(t)\int_{0}^{t}|\nabla_yh|{_{\mathcal{H}^{m}}^2}d\tau.
	\end{equation}
	Therefore, from \eqref{5.10}, \eqref{4.42} and trace estimate \eqref{2.10}, one obtains  that
	\begin{align}\label{cs1}
		&\left|\sigma\int_{0}^{t}\int_{z=0}\nabla_y\cdot\mathcal{C}^\alpha(S)(-v^b\cdot\MZ^\alpha\mathbf{N}-\mathcal{C}^\alpha(h))dyd\tau\right|\nonumber\\
		&\leq\sigma\Lambda_\infty(t)\left\{\int_{0}^{t}|\nabla_y\cdot\mathcal{C}^\alpha(S)|^2_{L^2}\right\}^\frac{1}{2}\left\{ \int_{0}^{t}|\nabla_yh|{_{\mathcal{H}^{m}}^2}+|\mathcal{C}^\alpha(h)|^2d\tau\right\}^\frac{1}{2}\nonumber\\
		&\leq  \Lambda_\infty(t)\int_{0}^{t}|h|{_{\mathcal{H}^m}^2}+\|v\|{_{\mathcal{H}^{m-1}}^2}+\sigma|\nabla_yh|{_{\mathcal{H}^m}^2}+\|\nabla v\|{_{\mathcal{H}^{m-1}}^2}d\tau
	\end{align}
	For $\int_{0}^{t}\int_{z=0}\nabla_y\cdot\mathcal{C}^\alpha(S)\partial_t\MZ^\alpha hdyd\tau$, if $\MZ^\alpha\neq\partial{_t^m}$, then \eqref{5.10} yields immediately that
	\begin{equation}\label{cs2}
		 \left|\sigma\int_{0}^{t}\int_{z=0}\nabla_y\cdot\mathcal{C}^\alpha(S)\partial_t\MZ^\alpha h\right|\leq \Lambda_\infty(t)\sigma\int_{0}^{t}|\nabla_yh|{_{\mathcal{H}^m}^2}d\tau.
	\end{equation}
For $\MZ^\alpha=\partial{_t^m}$, it follows from integration by parts and \eqref{5.10}   that
	\begin{align}\label{cs3}
		 &\left|\sigma\int_{0}^{t}\int_{z=0}\nabla_y\cdot\mathcal{C}^\alpha(S)\partial{_t^{m+1}}hdyd\tau\right|=\left|\sigma\int_{0}^{t}\int_{z=0}\mathcal{C}^\alpha(S)\partial{_t^{m+1}}\nabla_yhdyd\tau\right|\nonumber\\
		 &\leq\left|\int_{z=0}\sigma\mathcal{C}^\alpha(S)\partial{_t^m}\nabla_yhdy\right|+\left|\s\int_{z=0}\mathcal{C}^\alpha(S_0)\partial{_t^m}\nabla_yh_0dy\right|+\left|\sigma\int_{0}^{t}\int_{z=0}\partial_t\mathcal{C}^\alpha(S)\partial{_t^m}\nabla_yhdyd\tau\right|\nonumber\\
		&\leq \Lambda_0[|\nabla_yh(0)|{_{\mathcal{H}^{m-1}}^2}+\sigma|\nabla_yh_0|{_{\mathcal{H}^m}^2}] +\f\sigma8\int_{z=0}\frac{|\mathcal{Z}^\alpha \nabla_yh(t)|^2}{2(1+|\nabla_yh(t)|^2)^\frac{3}{2}}dy+\Lambda_0\s|\nabla_yh|{_{\mathcal{H}^{m-1}}^2}\nonumber\\
		&~~~~~+\Lambda_{\infty}(t)\sigma\int_{0}^{t}|\nabla_yh|{_{\mathcal{H}^m}^2}d\tau,
	\end{align}
	where one has used that, for $m\geq 5$, 
	\begin{align}
		 |\mathcal{C}^\alpha(S)(t)|^2_{L^2}&\lesssim\sum_{a=1}^{m-1}|\partial{_t^{m-a}}\nabla_yh\partial{_t^a}\frac{1}{\sqrt{1+|\nabla_yh|^2}}|^2_{L^2}\nonumber\\
		 &\lesssim\sum_{a=1}^{[\frac{m}{2}]}|\partial{_t^a}\frac{1}{\sqrt{1+|\nabla_yh|^2}}|{_{t,\infty}^2}|\partial{_t^{m-a}}\nabla_yh|^2_{L^2}+\sum_{a=[\frac{m}{2}]+1}^{m-1}|\partial{_t^{m-a}}\nabla_yh|^2_{t,\infty}|\partial{_t^a}\frac{1}{\sqrt{1+|\nabla_yh|^2}}|^2_{L^2}\nonumber\\
		 &\leq\Lambda(\f1{c_0},|\nabla_y h|^2_{\mathcal{H}^{[\f{m}2]+1}}) |\nabla_yh|{_{\mathcal{H}^{m-1}}^2}\leq \Lambda_0|\nabla_yh|{_{\mathcal{H}^{m-1}}^2},
	\end{align}
	It should be reminded that the estimate of
	\begin{equation}\label{4.2.15}
		-\sigma\int_{0}^{t}\int_{z=0}\nabla_y\cdot(\frac{\mathcal{Z}^\alpha \nabla_yh}{\sqrt{1+|\nabla_yh|^2}}-\frac{\nabla_yh\langle\nabla_yh,\mathcal{Z}^\alpha \nabla_yh\rangle}{(1+|\nabla_yh|^2)^\frac{3}{2}})\mathcal{C}^\alpha(h)dyd\tau,
	\end{equation}
	is a subtle issue. If one uses directly the integration by parts, the term $\int_{0}^{t}|\nabla_yv^b|{_{\mathcal{H}^{m-1}}^2}$,  which is difficult to estimate, must appear due to the boundary term $\MZ^\beta v{_y^b}\cdot\MZ\nabla_yh$ with $|\beta|=m-1$ in $\mathcal{C}^\alpha(h)$. Therefore, to deal with this boundary term, we will take full advantage of the boundary condition (\ref{4.34}) and  equation (\ref{4.5}) to turn it into a volume integral. In fact, the dynamic boundary condition \eqref{4.34} yields that
	\begin{align}\label{tr0}
		&\quad-\sigma\nabla_y\cdot\MZ^\alpha(\frac{\nabla_yh}{\sqrt{1+|\nabla_yh|^2}})-\partial{_z^\varphi}p\MZ^\alpha h\nonumber\\
		&=Q^\alpha-\left((2\mu\v S^\varphi V^\alpha+\l\v\text{div}^\varphi V^\alpha)\mathbf{N}\right)_3-(\MZ^\alpha h(2\mu\v\partial{_z^\varphi}S^\varphi v+\l\v\partial{_z^\varphi}\text{div}^\varphi v)\mathbf{N})_3+\mathcal{C}^\alpha(B)_3.
	\end{align}
	where $(\cdot)_3$ denotes the third component of a vector. Thus, it follows from \eqref{4.35} and trace estimate \eqref{2.10} that
	\begin{align}\label{tr1}
		 &-\sigma\int_{0}^{t}\int_{z=0}\nabla_y\cdot\left(\frac{\mathcal{Z}^\alpha \nabla_yh}{\sqrt{1+|\nabla_yh|^2}}-\frac{\nabla_yh\langle\nabla_yh,\mathcal{Z}^\alpha \nabla_yh\rangle}{(1+|\nabla_yh|^2)^\frac{3}{2}}\right)\MZ^\beta v{_y^b}\cdot\MZ\nabla_yhdyd\tau\nonumber\\
		&\leq\int_{0}^{t}\int_{z=0}Q^\alpha\MZ^\beta v{_y^b}\cdot\MZ\nabla_yh dyd\tau-\int_{0}^{t}\int_{z=0}\Big((2\mu\v S^\varphi V^\alpha+\l\v\text{div}^\varphi V^\alpha)\mathbf{N}\Big)_3\MZ^\beta v{_y^b}\cdot\MZ\nabla_yhdyd\tau\nonumber\\
		& ~~~~~+\Lambda_{\infty}(t)\left(\int_0^t|v_y^b|^2_{\mathcal{H}^{m-1}}d\tau\right)^{\f12}
		 \left(\int_0^t|h|^2_{\mathcal{H}^{m}}+|\mathcal{C}^\alpha(B)_3|^2_{L^2}+\sigma|\nabla_y\cdot\mathcal{C}^\alpha(S)|^2_{L^2}d\tau\right)^{\f12}\nonumber\\
		&\leq\int_{0}^{t}\int_{z=0}Q^\alpha\MZ^\beta v{_y^b}\cdot\MZ\nabla_yh dyd\tau-\int_{0}^{t}\int_{z=0}\Big((2\mu\v S^\varphi V^\alpha+\l\v\text{div}^\varphi V^\alpha)\mathbf{N}\Big)_3\MZ^\beta v{_y^b}\cdot\MZ\nabla_yhdyd\tau\nonumber\\
		& ~~~~~+\d\int_0^t\v^2\|\nabla^2v\|^2_{\mathcal{H}^{m-1}}d\tau+C_\d\Lambda_{\infty}(t)
		\int_0^t|h|^2_{\mathcal{H}^{m}}+\|v\|_{\mathcal{H}^{m-1}}^2+\|\nabla v\|^2_{\mathcal{H}^{m-1}}+\sigma|\nabla_yh|_{\mathcal{H}^m}^2d\tau,
	\end{align}
	where $|\beta|=m-1$. To estimate the boundary term $\int_{0}^{t}\int_{z=0}Q^\alpha\MZ^\beta v{_y^b}\cdot\MZ\nabla_yh dyd\tau$ in \eqref{tr1},  using  integration by parts, we note  that
	\begin{align}\label{tr2}
		&\int_{0}^{t}\int_{z=0}Q^\alpha\MZ^\beta v{_y^b}\cdot\MZ\nabla_yhdyd\tau\nonumber\\
		&=\int_{0}^{t}\int_{\mathcal{S}}\partial{_z^\varphi}Q^\alpha\MZ^\beta v_y\cdot\MZ\nabla_y\eta d\mathcal{V}_\tau d\tau+\int_{0}^{t}\int_{\mathcal{S}}Q^\alpha\partial{_z^\varphi}(\MZ^\beta v_y\cdot\MZ\nabla_y\eta) d\mathcal{V}_\tau d\tau\\
		&\leq \int_{0}^{t}\int_{\mathcal{S}}\partial{_z^\varphi}Q^\alpha\MZ^\beta v_y\cdot\MZ\nabla_y\eta d\mathcal{V}_\tau d\tau
		+\int_{0}^{t}\|\nabla v\|{_{\mathcal{H}^{m-1}}^2}d\tau+\Lambda_\infty(t)\int_{0}^{t}\|v\|_{\mathcal{H}^{m-1}}^2+\|Q^m\|^2d\tau.\nonumber
	\end{align}
	Using the equation \eqref{4.5}, one obtains that
	\begin{align}\label{tr4}
		 &\int_{0}^{t}\int_{\mathcal{S}}\partial{_z^\varphi}Q^\alpha\MZ^\beta v_y\cdot\MZ\nabla_y\eta d\mathcal{V}_\tau d\tau\nonumber\\
		 &\leq-\int_{0}^{t}\int_{\mathcal{S}}\varrho(\partial_t+v_y\cdot\nabla_y+V_z\partial_z)(V^\alpha)_3\MZ^\beta v_y\cdot\MZ\nabla_y\eta d\mathcal{V}_\tau d\tau\\
		 &\quad+\int_{0}^{t}\int_{\mathcal{S}}(2\mu\v\text{div}^\varphi S^\varphi V^\alpha+\l\v\nabla^\varphi\text{div}^\varphi V^\alpha)_3\MZ^\beta v_y\cdot\MZ\nabla_y\eta d\mathcal{V}_\tau d\tau+\Lambda_{\infty}(t)\int_0^t\|\mathcal{R}^\a_M\|\|v\|_{\mathcal{H}^{m-1}}d\tau,\nonumber
	\end{align}
	where $|\beta|=m-1$ and $\|\mathcal{R}^\a_M\|$ will be estimated later. Integration by parts, and using  $V_z=0$ on $z=0$, one gets that
	\begin{align}\label{tr5}
		 &\quad-\int_{0}^{t}\int_{\mathcal{S}}\varrho(\partial_t+v_y\cdot\nabla_y+V_z\partial_z)(V^\alpha)_3\MZ^\beta v_y\cdot\MZ\nabla_y\eta d\mathcal{V}_\tau d\tau\nonumber\\
		 &=-\int_{\mathcal{S}}\varrho(V^\alpha)_3\MZ^\beta v_y\cdot\nabla_y\eta d\mathcal{V}_t+\int_{\mathcal{S}}\varrho_0(V{^\alpha_0})_3\MZ^\beta v_y(0)\cdot\nabla_y\eta_0d\mathcal{V}_0\nonumber\\
		 &~~~~~~~~~~~~~~~~~~~~~~~~~~~~~~~~~~~~~~~~~~~~~~~~~~~~~~~~+\Lambda_\infty(t)\int_{0}^{t}\|V^\a(\tau)\|^2+\|v(\tau)\|{_{\mathcal{H}^m}^2}d\tau\nonumber\\
		 &\leq\f{1}{16}\int\varrho|V^\a(t)|^2d\mathcal{V}_t+\Lambda_0[\|V^\a(0)\|^2+\|v(0)\|{_{\mathcal{H}^{m-1}}^2}]+\Lambda_\infty(t)\int_{0}^{t}\|V^m\|^2+\|v\|_{\mathcal{H}^m}^2+|h|{_{\mathcal{H}^m}^2}d\tau.
	\end{align}
	For the second term in \eqref{tr4}, it follows from the integrating by parts that
	\begin{align}\label{tr6}
		&\quad\int_{0}^{t}\int_{\mathcal{S}}(2\m\v\text{div}^\varphi S^\varphi V^\alpha+\l\v\nabla^\varphi \text{div}^\varphi V^\alpha)_3\MZ^\beta v_y\cdot\MZ\nabla_y\eta d\mathcal{V}_\tau d\tau\nonumber
		\\
		&=-\int_{0}^{t}\int_{\mathcal{S}}\m \v(\nabla^\varphi V{_3^\alpha}+\partial{_z^\varphi}V^\alpha)\cdot\nabla^\varphi(\MZ^\beta v_y\cdot\MZ\nabla_y\eta)+\l\v \text{div}^\varphi V^\alpha\partial{_z^\varphi}(\MZ^\beta v_y\cdot\MZ\nabla_y\eta)d\mathcal{V}_\tau d\tau\nonumber\\
		&\quad+\int_{0}^{t}\int_{z=0}\Big((2\mu\v S^\varphi V^\alpha+\l\v\text{div}^\varphi V^\alpha)\mathbf{N}\Big)_3\MZ^\beta v{_y^b}\cdot\MZ\nabla_yhdyd\tau\nonumber\\
		&\leq \d\v\int_{0}^{t}\|\nabla V^m\|^2d\tau+C_\d\Lambda_\infty(t)\int_{0}^{t}\v\|\nabla v\|^2_{\mathcal{H}^{m-1}}\nonumber\\
		&~~~~~~~~~~~~~~+\int_{0}^{t}\int_{z=0}\Big((2\mu\v S^\varphi V^\alpha+\l\v\text{div}^\varphi V^\alpha)\mathbf{N}\Big)_3\MZ^\beta v{_y^b}\cdot\MZ\nabla_yhdyd\tau.
	\end{align}
	Combining \eqref{tr1}-\eqref{tr6} and noticing the cancellation between the worst boundary terms in \eqref{tr1} and \eqref{tr6}, it follows from \eqref{5.2}, for $|\beta|=m-1$,  that
	\begin{align}\label{tr-7}
		 &-\sigma\int_{0}^{t}\int_{z=0}\nabla_y\cdot\left(\frac{\mathcal{Z}^\alpha \nabla_yh}{\sqrt{1+|\nabla_yh|^2}}-\frac{\nabla_yh\langle\nabla_yh,\mathcal{Z}^\alpha \nabla_yh\rangle}{(1+|\nabla_yh|^2)^\frac{3}{2}}\right)\MZ^\beta v{_y^b}\cdot\MZ\nabla_yhdyd\tau\nonumber\\
		&\leq \f{1}{16}\int\varrho|V^\a(t)|^2d\mathcal{V}_t+C\d\int_0^t\v^2\|\nabla^2v\|^2_{\mathcal{H}^{m-1}}+\v\|\nabla V^m\|^2_{\mathcal{H}^{m-1}}d\tau\nonumber\\
		&~+\Lambda_0[\|V^m(0)\|^2+\|v(0)\|^2_{\mathcal{H}^{m-1}}]+C_\d\Lambda_{\infty}(t)
		\int_0^tY_m+\|\nabla v\|^2_{\mathcal{H}^{m-1}} +\|\mathcal{R}^\a_M\|Y_md\tau.
	\end{align}

For the terms of the form $\MZ^\beta v_y\MZ^\gamma\nabla_yh$ with $|\beta|+|\gamma|\leq m,~1\leq |\beta|\leq m-2$ and the term $\dfrac{(\partial_zv)^b}{\partial_z\varphi}\MZ^\alpha h\mathbf{N}$ in $\mathcal{C}^\alpha(h)$, it follows from the integration by parts, Lemma \ref{l2.2}, trace estimate \eqref{2.10} and \eqref{5.2}   that
	\begin{align}\label{gt1}
		&\quad-\sigma\int_{0}^{t}\int_{z=0}\nabla_y\cdot\left(\frac{\mathcal{Z}^\alpha \nabla_yh}{\sqrt{1+|\nabla_yh|^2}}-\frac{\nabla_yh\langle\nabla_yh,\mathcal{Z}^\alpha \nabla_yh\rangle}{(1+|\nabla_yh|^2)^\frac{3}{2}}\right)\left(\MZ^\beta v_y\MZ^\gamma\nabla_yh+\dfrac{(\partial_zv)^b}{\partial_z\varphi}\MZ^\alpha h\mathbf{N}\right)dyd\tau\nonumber\\
		&=\sigma\int_{0}^{t}\int_{z=0}\left(\frac{\mathcal{Z}^\alpha \nabla_yh}{\sqrt{1+|\nabla_yh|^2}}-\frac{\nabla_yh\langle\nabla_yh,\mathcal{Z}^\alpha \nabla_yh\rangle}{(1+|\nabla_yh|^2)^\frac{3}{2}}\right)\nabla_y\left(\MZ^\beta v_y\MZ^\gamma\nabla_yh+\dfrac{(\partial_zv)^b}{\partial_z\varphi}\MZ^\alpha h\mathbf{N}\right)dyd\tau\nonumber\\
		&\leq\Lambda_\infty(t)\int_{0}^{t}\|\nabla v\|{_{\mathcal{H}^{m-1}}^2}+\|V^m\|^2+|h|{_{\mathcal{H}^m}^2}+\sigma|\nabla_yh|{_{\mathcal{H}^m}^2}d\tau.
	\end{align}
Substituting \eqref{pt1}-\eqref{cs3}, \eqref{tr-7} and \eqref{gt1} into \eqref{k1}, one obtains, for  $\mathcal{Z}^\a\neq\partial_t^m$, that
	\begin{align}\label{k1e}
		K_1&\leq \Lambda_0\Big(\sigma|\nabla_yh_0|{_{\mathcal{H}^m}^2}+\|V^m(0)\|^2+\|v(0)\|^2_{\mathcal{H}^{m-1}}\Big)-\sigma\int_{z=0}\frac{|\mathcal{Z}^\alpha \nabla_yh(t)|^2}{4(1+|\nabla_yh(t)|^2)^\frac{3}{2}}dy\nonumber\\
		& ~~+\f{1}{16}\int\varrho|V^\a(t)|^2d\mathcal{V}_t+C\d\int_0^t\v^2\|\nabla^2v\|^2_{\mathcal{H}^{m-1}}+\v\|\nabla V^m\|^2_{\mathcal{H}^{m-1}}d\tau\nonumber\\
		&~~+C_\d\Lambda_{\infty}(t)
		\int_0^tY_m+\|\nabla v\|^2_{\mathcal{H}^{m-1}} +\|\mathcal{R}^\a_M\|Y_md\tau,
	\end{align}
	and for $\mathcal{Z}^\a=\partial_t^m$, it holds that
	\begin{align}\label{k1e-1}
		K_1&\leq \Lambda_0[|\nabla_yh(0)|{_{\mathcal{H}^{m-1}}^2}+\sigma|\nabla_yh_0|{_{\mathcal{H}^m}^2}+\|V^m(0)\|^2+\|v(0)\|^2_{\mathcal{H}^{m-1}}]+\Lambda_0|\nabla_yh|{_{\mathcal{H}^{m-1}}^2}\nonumber\\
		& ~~-\sigma\int_{z=0}\frac{|\mathcal{Z}^\alpha \nabla_yh(t)|^2}{4(1+|\nabla_yh(t)|^2)^\frac{3}{2}}dy+\Lambda_0\d\int_0^t\v^2\|\nabla^2v\|^2_{\mathcal{H}^{m-1}}+\v\|\nabla V^m\|^2_{\mathcal{H}^{m-1}}d\tau\nonumber\\
		&~~+\f{1}{16}\int\varrho|V^\a(t)|^2d\mathcal{V}_t+C_\d\Lambda_{\infty}(t)
		\int_0^tY_m+\|\nabla v\|^2_{\mathcal{H}^{m-1}} +\|\mathcal{R}^\a_M\|Y_md\tau,
	\end{align}
	where the remaining term $\|\mathcal{R}^\a_M\|$ in \eqref{k1e} and \eqref{k1e-1} will be estimated later.

It follows from \eqref{4.41}, \eqref{4.42} and \eqref{2.10}   that
	\begin{align}\label{k2e}
		K_2&=\int_{0}^{t}\int_{z=0}\partial{_z^\varphi}p\mathcal{Z}^\alpha h(\partial_t\mathcal{Z}^\alpha h-v^b\cdot\mathcal{Z}^\alpha\mathbf{N}-\mathcal{C}^\alpha(h))dyd\tau\nonumber\\
		&=\frac{1}{2}\int_{z=0}\partial{_z^\varphi}p(t)|\mathcal{Z}^\alpha h(t)|^2dy-\frac{1}{2}\int_{z=0}\partial{_z^\varphi}p(0)|\mathcal{Z}^\alpha h(0)|^2dy\nonumber\\
		&\quad-\frac{1}{2}\int_{0}^{t}\int_{z=0}\partial_t\partial{_z^\varphi}p|\mathcal{Z}^\alpha h|^2dyd\tau-\int_{0}^{t}\int_{z=0}\partial{_z^\varphi}p\mathcal{Z}^\alpha h(v^b\cdot\mathcal{Z}^\alpha\mathbf{N}+\mathcal{C}^\alpha(h))dyd\tau\nonumber\\
		&\leq\frac{1}{2}\int_{z=0}\partial{_z^\varphi}p(t)|\mathcal{Z}^\alpha h(t)|^2dy-\frac{1}{2}\int_{z=0}\partial{_z^\varphi}p(0)|\mathcal{Z}^\alpha h(0)|^2dy+\Lambda_\infty(t)\int_{0}^{t}|h|{_{\mathcal{H}^m}^2}d\tau\nonumber\\
		&\quad-\int_{0}^{t}\int_{z=0}\partial{_z^\varphi}pv^b\cdot\nabla_y(\frac{|\mathcal{Z}^\alpha h|^2}{2})+\partial{_z^\varphi}p\mathcal{Z}^\alpha h\mathcal{C}^\alpha(h)dyd\tau\nonumber\\
		&\leq\frac{1}{2}\int_{z=0}\partial{_z^\varphi}p(t)|\mathcal{Z}^\alpha h(t)|^2dy-\frac{1}{2}\int_{z=0}\partial{_z^\varphi}p(0)|\mathcal{Z}^\alpha h(0)|^2dy\nonumber\\
	    &\quad+\Lambda_\infty(t)\int_{0}^{t}\|\nabla v\|{_{\mathcal{H}^{m-1}}^2}+|h|{_{\mathcal{H}^{m}}^2}+\|v\|{_{\mathcal{H}^{m-1}}^2}d\tau.
	\end{align}
Since the boundary condition $\eqref{1.23-1}$ yields immediately that
	\begin{equation*}
		p-p_e=(2\m\v S^\varphi v+\l\v {div}^\varphi v)\mathbf{n}\cdot\mathbf{n}-\sigma\nabla_y\cdot\frac{\nabla_y h}{\sqrt{1+|\nabla_y h|^2}}.
	\end{equation*}
This, together with  Lemma \ref{lem2.4} and trace theorem, implies that
	\begin{align}\label{k3e}
		K_3&\leq\Lambda_\infty(t)\v\left(\int_{0}^{t}|\mathcal{Z}^\a\mathbf{N}|_{-\f{1}{2}}^2\right)^\frac{1}{2}\left(\int_{0}^{t}|(V^\alpha)^b|^2_{\f12}d\tau\right)^\frac{1}{2}\nonumber\\
		&\leq\delta\v\int_{0}^{t}\|\nabla V^m\|^2d\tau+C_\d\Lambda_\i(t)\int_{0}^{t}\v|\mathcal{Z}^m h|{_{\f12}^2}d\tau.
	\end{align}
	Direct calculations yield immediately that
	\begin{align}\label{k4e}
		 K_4&\leq\Lambda(\frac{1}{c_0},\|v\|_{\mathcal{W}^{2,\infty},t}+\|\nabla v\|_{\mathcal{W}^{1,\infty},t})\v^\frac{1}{2}\|\partial_{zz}v\|_{\infty,t}\v^\frac{1}{2}\int_{0}^{t}|h|_{\mathcal{H}^m}|V^\alpha|d\tau\nonumber\\
		 &\leq\delta\v\int_{0}^{t}\|\nabla V^m\|^2d\tau+C_\d\Lambda_\infty(t)\int_{0}^{t}(|h|{_{\mathcal{H}^{m}}^2}+\|V^m\|^2)d\tau.
	\end{align}
It follows from   \eqref{4.35} that
	\begin{align}\label{k5e}
		 &K_5\leq\int_{0}^{t}|\mathcal{C}^\alpha(B)|_{L^2}|V^\alpha|_{L^2}d\tau
		 \leq\left(\int_{0}^{t}|\mathcal{C}^\alpha(\mathcal{B})|^2_{L^2}d\tau\right)^\frac{1}{2}\left(\int_{0}^{t}\|V^\a\|\|\nabla V^\a\|\right)^\frac{1}{2}\nonumber\\
		 &\leq\Lambda_\i(t)\Big\{\v^2\int_{0}^{t}|(\nabla v)^b|{_{\mathcal{H}^{m-1}}^2}+|\nabla_yh|{_{\mathcal{H}^{m-1}}^2}d\tau\Big\}^\frac{1}{2}\Big\{\int_{0}^{t}\|\nabla V^\a\|\|V^\a\|+\|V^\a\|^2\Big\}^\frac{1}{2}
	\nonumber	\\
		&\leq\delta\v\int_{0}^{t}\|\nabla V^m\|^2d\tau+\delta\v^2\int_{0}^{t}\|\nabla^2 v\|{_{\mathcal{H}^{m-1}}^2}d\tau +C_\d\Lambda_\infty(t)\int_{0}^{t}Y_m+\|\nabla v\|{_{\mathcal{H}^{m-1}}^2}d\tau.
	\end{align}
Therefore, from \eqref{k1e}-\eqref{k5e}, one  obtains the boundary estimate for $\mathcal{Z}^\a\neq\partial_t^m$ that
	\begin{align}\label{5.15}
	\int_{0}^{t}\mathcal{B}d\tau\leq& \Lambda_0[|h_0|{_{\mathcal{H}^{m}}^2}+\sigma|\nabla_yh_0|{_{\mathcal{H}^m}^2}+\|V^m(0)\|^2+\|v(0)\|^2_{\mathcal{H}^{m-1}}]+\frac{1}{2}\int_{z=0}\partial{_z^\varphi}p(t)|\mathcal{Z}^\alpha h(t)|^2dy\nonumber\\
		 &-\sigma\int_{z=0}\frac{|\nabla_y\mathcal{Z}^\alpha h(t)|^2}{4(1+|\nabla_yh(t)|^2)^\frac{3}{2}}dy+C\d\int_0^t\v^2\|\nabla^2v\|^2_{\mathcal{H}^{m-1}}+\v\|\nabla V^m\|^2d\tau\nonumber\\
		 &+C_\d\Lambda_{\infty}(t)
		\int_0^tY_m+\|\nabla v\|^2_{\mathcal{H}^{m-1}} +\|\mathcal{R}^\a_M\|Y_md\tau+\f{1}{16}\int\varrho|V^\a(t)|^2d\mathcal{V}_t,
	\end{align}
	and for $\mathcal{Z}^\a=\partial_t^m$, it holds that
	\begin{align}\label{5.15-1}
		\int_{0}^{t}\mathcal{B}d\tau\leq& \Lambda_0[ |h_0|{_{\mathcal{H}^{m}}^2}+\sigma|\nabla_yh_0|{_{\mathcal{H}^m}^2}+\|V^m(0)\|^2+\|v(0)\|^2_{\mathcal{H}^{m-1}}]+\Lambda_0|\nabla_yh|{_{\mathcal{H}^{m-1}}^2}\nonumber\\
		 &+\f{1}{16}\int\varrho|V^\a(t)|^2d\mathcal{V}_t+\frac{1}{2}\int_{z=0}\partial{_z^\varphi}p(t)|\mathcal{Z}^\alpha h(t)|^2dy-\sigma\int_{z=0}\frac{|\nabla_y\mathcal{Z}^\alpha h(t)|^2}{4(1+|\nabla_yh(t)|^2)^\frac{3}{2}}dy\\
		 &+C\d\int_0^t\v^2\|\nabla^2v\|^2_{\mathcal{H}^{m-1}}+\v\|\nabla V^m\|^2d\tau+C_\d\Lambda_{\infty}(t)
		\int_0^tY_m+\|\nabla v\|^2_{\mathcal{H}^{m-1}} +\|\mathcal{R}^\a_M\|Y_md\tau.\nonumber
	\end{align}	
Using \eqref{4.4}  and integration by parts, one obtains that
	\begin{align}\label{5.16}
		&\int_{0}^{t}Jd\tau=\int_{0}^{t}\int_{\mathcal{S}}Q^\alpha(-\frac{1}{\g p}\partial{_t^\varphi}Q^\alpha-\frac{1}{\g p}v\cdot\nabla^\varphi Q^\alpha+\mathcal{R}{_d^\alpha})d\mathcal{V}_{\tau}d\tau\nonumber\\
		&=-\int_{\mathcal{S}}\frac{|Q^\alpha(t)|^2}{2\g p(t)}d\mathcal{V}_t+\int_{\mathcal{S}}\frac{|Q^\alpha(0)|^2}{2\g p(0)}d\mathcal{V}_0\nonumber\\
		&\quad+\int_{0}^{t}\int_{\mathcal{S}}(\partial_t(\frac{1}{2\g p})+{div}^\varphi(\frac{v}{2\g p}))|Q^\alpha|^2d\mathcal{V}_{\tau}+\int_{0}^{t}\int_{\mathcal{S}}Q^\alpha\mathcal{R}{_d^\alpha}d\mathcal{V}_{\tau}\nonumber\\
		&\leq-\int_{\mathcal{S}}\frac{|Q^\alpha(t)|^2}{2\g p(t)}d\mathcal{V}_t+\int_{\mathcal{S}}\frac{|Q^\alpha(0)|^2}{2\g p(0)}d\mathcal{V}_0+\Lambda_\infty(t)\int_{0}^{t}\|Q^m(t)\|^2d\tau+\int_{0}^{t}\int_{\mathcal{S}}Q^\alpha\mathcal{R}{_d^\alpha}d\mathcal{V}_\tau d\tau,
	\end{align}
	where
	\begin{equation*}
		 \mathcal{R}{_d^\alpha}=\mathcal{Z}^\alpha\eta(\partial{_z^\varphi}(\frac{1}{\g p})\partial{_t^\varphi}p+\partial{_z^\varphi}(\frac{v}{\g p})\cdot\nabla^\varphi p)-\mathcal{C}^\alpha(d)-\frac{1}{\g p}\mathcal{C}{_0^\alpha}(p)-\frac{v}{\g p}\cdot\mathcal{C}^\alpha(p)-[\mathcal{Z}^\alpha,\frac{1}{\g p}]\partial{_t^\varphi}p-[\mathcal{Z}^\alpha,\frac{v}{\g p}]\cdot\nabla^\varphi p.
	\end{equation*}
	Furthermore, it is easy to obtain that
	\begin{equation}\label{5.18}
		 \int_{0}^{t}\int_{\mathcal{S}}Q^\alpha\mathcal{Z}^\alpha\eta(\partial{_z^\varphi}(\frac{1}{\g p})\partial{_t^\varphi}p+\partial{_z^\varphi}(\frac{v}{\g p})\cdot\nabla^\varphi p)d\mathcal{V}_\tau d\tau\leq\Lambda_\infty(t)\int_{0}^{t}\|Q^m(t)\|^2+|h|{_{\mathcal{H}^{m}}^2}d\tau.
	\end{equation}
	It follows from \eqref{4.6}-\eqref{4.7} that
	\begin{align}\label{5.17}
		 &\int_{0}^{t}\int_{\mathcal{S}}Q^\alpha(-\mathcal{C}^\alpha(d)-\frac{1}{\g p}\mathcal{C}{_0^\alpha}(p)-\frac{v}{\g p}\cdot\mathcal{C}^\alpha(p))d\mathcal{V}_\tau d\tau\nonumber\\
		 &\leq\Lambda(\frac{1}{c_0},\|v\|_{\infty,t})\left(\int_{0}^{t}\|Q^\alpha\|^2d\tau\right)^\frac{1}{2}\left(\int_{0}^{t}(\|\mathcal{C}^\alpha(d)\|+\|\mathcal{C}{_0^\alpha}(p)\|+C^\alpha(p)\|)^2d\tau\right)^\frac{1}{2}\nonumber\\
		 &\leq\delta\int_{0}^{t}\|\nabla p\|{_{\mathcal{H}^{m-1}}^2}d\tau+C_\d\Lambda_\infty(t)\int_{0}^{t}\|\nabla v\|{_{\mathcal{H}^{m-1}}^2}+\|Q^m(t)\|^2+|h|{_{\mathcal{H}^{m}}^2}d\tau.
	\end{align}
Using \eqref{2.7}, \eqref{5.2} and \eqref{5.3}, one gets that
	\begin{align}\label{5.19}
		 &\int_{0}^{t}\int_{\mathcal{S}}Q^\alpha(-[\mathcal{Z}^\alpha,\frac{1}{\g p}]\partial{_t^\varphi}p-[\mathcal{Z}^\alpha,\frac{v}{\g p}]\cdot\nabla^\varphi p)d\mathcal{V}_\tau d\tau\nonumber\\
		 &\leq\int_{0}^{t}\|Q^\alpha\|(\|[\mathcal{Z}^\alpha,\frac{1}{\g p}]\partial{_t^\varphi}p\|+\|[\mathcal{Z}^\alpha,\frac{v}{\g p}]\cdot\nabla^\varphi p\|)d\tau\nonumber\\
		 &\leq\Lambda_\infty(t)\left(\int_{0}^{t}\|Q^\alpha\|^2d\tau\right)^\frac{1}{2}\left(\int_{0}^{t}(\|p\|{_{\mathcal{H}^{m}}^2}+\|\nabla p\|{_{\mathcal{H}^{m-1}}^2}+\|v\|{_{\mathcal{H}^{m}}^2}+|h|^2_{\mathcal{H}^m}d\tau\right)^\frac{1}{2}\nonumber\\
		 &\leq\delta\int_{0}^{t}\|\nabla p\|{_{\mathcal{H}^{m-1}}^2}d\tau+C_\d\Lambda_\infty(t)\int_{0}^{t}\|v\|{_{\mathcal{H}^{m}}^2}+\|p\|{_{\mathcal{H}^{m}}^2}+\|Q^\alpha\|^2+|h|^2_{\mathcal{H}^m}d\tau\nonumber\\
		 &\leq\delta\int_{0}^{t}\|\nabla p\|{_{\mathcal{H}^{m-1}}^2}d\tau+C_\d\Lambda_\infty(t)\int_{0}^{t}\|V^m(t)\|^2+\|Q^m(t)\|^2+|h|^2_{\mathcal{H}^m}d\tau.
	\end{align}
	Combining \eqref{5.16}-\eqref{5.19} leads to that
	\begin{align}\label{5.20}
		 \int_{0}^{t}Jd\tau&\leq-\int_{\mathcal{S}}\frac{|Q^\alpha(t)|^2}{2\g p(t)}d\mathcal{V}_t+\int_{\mathcal{S}}\frac{|Q^\alpha(0)|^2}{2\g p(0)}d\mathcal{V}_t+\delta\int_{0}^{t}\|\nabla p\|{_{\mathcal{H}^{m-1}}^2}d\tau\nonumber\\
		 &\quad+C_\d\Lambda_\infty(t)\int_{0}^{t}\|\nabla v\|{_{\mathcal{H}^{m-1}}^2}+\|V^m(t)\|^2+\|Q^m(t)\|^2+|h|{_{\mathcal{H}^{m}}^2}d\tau.
	\end{align}
	Now we estimate the terms involve $\mathcal{R}_M^\a$. Using \eqref{4.5}-\eqref{4.8-1}, it is easy to obtain  that
	\begin{equation*}
		\int_0^t\|\mathcal{R}^\alpha_M\|^2d\tau\leq \Lambda_{\infty}(t)\int_0^t\v^2\|\nabla^2v\|^2_{\mathcal{H}^{m-1}}+Y_m+\v|\MZ^{m} h|^2_{\f12} d\tau.
	\end{equation*}
	Therefore, using holder inequality,  one gets that
	\begin{align}\label{5.21}
		&\left|\int_{0}^{t}\int_{\mathcal{S}}\mathcal{R}_M^\a V^\a d\mathcal{V}_\tau d\tau\right|\leq \int_{0}^{t}\|\mathcal{R}^\alpha_M\|\|V^\alpha\|d\tau  \nonumber\\
		 &\leq\d\int_0^t\v^2\|\nabla^2v\|^2_{\mathcal{H}^{m-1}}d\tau+C_\d\Lambda_{\infty}(t)\int_0^tY_m+\|\nabla v\|^2_{\mathcal{H}^{m-1}}+\v|\MZ^{m} h|^2_{\f12}d\tau,
	\end{align}
	and
	\begin{equation}\label{5.21-1}
\int_0^t\|\mathcal{R}^\a_M\|Y_md\tau \leq\d\int_0^t\v^2\|\nabla^2v\|^2_{\mathcal{H}^{m-1}}d\tau+C_\d\Lambda_{\infty}(t)\int_0^tY_m+\|\nabla v\|^2_{\mathcal{H}^{m-1}}+\v|\MZ^{m} h|^2_{\f12}d\tau.
	\end{equation}
Combining \eqref{5.4}, \eqref{5.15}, \eqref{5.20}, \eqref{5.21} and \eqref{5.21-1}, we obtain, for $\MZ^\a\neq\partial^m_t$, that
	\begin{align}\label{4.1.40}
		 &\int_{\mathcal{S}}\frac{\varrho}{4}|V^\alpha(t)|+\frac{|Q^\alpha(t)|^2}{2\g p(t)}d\mathcal{V}_t+\frac{1}{2}\int_{z=0}-\partial{_z^\varphi}p(t)|\mathcal{Z}^\alpha h(t)|^2+\sigma\frac{|\nabla_y\mathcal{Z}^\alpha h(t)|^2}{(1+|\nabla_yh(t)|^2)^\frac{3}{2}}dy\nonumber\\
		 &\quad+\int_{0}^{t}\int_{\mathcal{S}}2\m\v|S^\varphi V^\alpha|^2+\l\v|{div}^\varphi V^\alpha|^2d\mathcal{V}_\tau d\tau\nonumber\\
		 &\leq\Lambda_0[|(h,\sqrt\sigma\nabla_yh)(0)|^2_{\mathcal{H}^m}+\|(V^m,Q^m)(0)\|^2+\|v(0)\|^2_{\mathcal{H}^m}]+\delta\int_{0}^{t}\|\nabla p\|{_{\mathcal{H}^{m-1}}^2}d\tau\nonumber\\
		 &~~+\Lambda_0\delta\int_{0}^{t}\v\|\nabla V^m\|^2+\v^2\|\nabla^2 v\|{_{\mathcal{H}^{m-1}}^2}d\tau+\Lambda_\infty(t)\int_{0}^{t}\Lambda(Y_m)+\|\nabla v\|{_{\mathcal{H}^{m-1}}^2}+\v|\MZ^mh|_{\frac{1}{2}}^2d\tau.
	\end{align}
	Therefore, summing over $\alpha$, using the a priori assumptions \eqref{5.1} on the Taylor sign condition and  Lemma \ref{lem20.2}, and taking $\d$ suitably small,  one proves \eqref{4.1-1} for the case $\MZ^\a\neq\partial^m_t$.
	
Collecting \eqref{5.4}, \eqref{5.15-1}, \eqref{5.20}, \eqref{5.21}, \eqref{5.21-1},  and using \eqref{4.1-1} for the case of $\MZ^\a\neq\partial^m_t$ whose proof  has already been closed in \eqref{4.1.40}, one yields  immediately \eqref{4.1-1} for the case 	 $\MZ^\a=\partial^m_t$.  Therefore, the proof of Lemma \ref{p5.1} is completed. $\hfill\Box$

	
\subsection{Estimates for $\nabla^\varphi p$ and ${div}^\varphi v$}
	
	\setcounter{equation}{0}
	
To deal with the compressibility, one needs to  obtain some uniform  estimates for the pressure and the divergence of the velocity field. For incompressible flows, the divergence of the velocity field is zero and the pressure plays a role of Lagrangian multiplier so that the pressure can be estimated by elliptic regularities, see \cite{Masmoudi-R-1}. While, the pressure $p(\varrho)=\varrho^\g$, shown in $\eqref{1.22}_2$ for compressible flows, satisfies a transport equation due to$\eqref{1.22}_1$ so that it can be estimated by energy method. However, it seems that one can only expect uniform estimates of $m-2$ order conormal derivatives of $\nabla p$, which is one order lower than the one in incompressible flows in \cite{Masmoudi-R-1}, since the divergence of the velocity field is not free but depends on the regularity of the free surfaces in compressible flows.  

Precisely, one can obtain the following  estimates for $\nabla^\varphi p$ and $\text{div}^\varphi v$.
\begin{proposition}\label{p6.6}
For any $t\in[0,T^\v]$ and $m\geq2$, it holds that
\begin{align}\label{10.3}
&\|\nabla v(t)\|^2_{\mathcal{H}^{m-2}}+\v\|(\mbox{div}^\varphi v, \nabla^\varphi p)(t)\|^2_{\mathcal{H}^{m-1}}+\v^2\int_{0}^{t}\|\nabla^\varphi\MZ^{m-1}\mbox{div}^\varphi v\|^2 d\tau\nonumber\\
&\leq\|\nabla v_0\|^2_{\mathcal{H}^{m-2}}+\Lambda_0\v\|(\mbox{div}^{\varphi_0} v_0,\nabla^{\varphi_0} p_0)\|^2_{\mathcal{H}^{m-1}}+(\d_1+\d\Lambda_0)\v^2\int_{0}^{t}\|\nabla^2 v\|^2_{\mathcal{H}^{m-1}}d\tau\nonumber\\
&\quad+\frac{\Lambda_0}{\d_1}\v\int_{0}^{t}\|\nabla V^m\|^2d\tau+C_\d\Lambda_\i(t)(1+\|\nabla\text{div}^\varphi v\|_{\i,t}^2)\int_0^t\Lambda(Y_m)+\|\nabla v\|^2_{\mathcal{H}^{m-1}}+\v|\MZ^mh|_\frac{1}{2}^2d\tau,
\end{align}
and 
\begin{align}\label{10.4}
\|\nabla p(t)\|^2_{\mathcal{H}^{m-2}}+\int_0^t\|\nabla p\|^2_{\mathcal{H}^{m-1}}d\tau\leq \Lambda_0\v^2\int_{0}^{t}\|\nabla^2 v\|^2_{\mathcal{H}^{m-1}}d\tau+\Lambda_\i(t)\int_0^t\Lambda(Y_m)d\tau,
\end{align}
where $\d_1>0$ and $\d>0$ are small constants which will be chosen later.
\end{proposition}

It should be noted that the estimate of $\v^2\int_{0}^{t}\|\nabla^\varphi\MZ^{m-1}\mbox{div}^\varphi v\|^2 d\tau$ is the key to control $\v^2\int_{0}^{t}\|\nabla^2 v\|^2_{\mathcal{H}^{m-1}}d\tau$. This is not necessary for incompressible case, since $\v^2\int_{0}^{t}\|\nabla^2 v\|^2_{\mathcal{H}^{m-1}}d\tau$ does not appear due to \eqref{20.120}.

The proof of this proposition is a consequence of the following lemmas. Note that 
\begin{eqnarray}\label{4.3.2}
\begin{cases}
\|\nabla v(t)\|^2_{\mathcal{H}^{m-2}}
\leq \|\nabla v_0\|^2_{\mathcal{H}^{m-2}}+\int_0^t\|\partial_t\nabla v\|^2_{\mathcal{H}^{m-2}}d\tau\leq \|\nabla v_0\|^2_{\mathcal{H}^{m-2}}+\int_0^t\|\nabla v\|^2_{\mathcal{H}^{m-1}}d\tau,\\
\|\nabla  p(t)\|^2_{\mathcal{H}^{m-2}}
\leq \|\nabla  p_0\|^2_{\mathcal{H}^{m-2}}+\int_0^t\|\partial_t\nabla p\|^2_{\mathcal{H}^{m-2}}d\tau\leq\|\nabla  p_0\|^2_{\mathcal{H}^{m-2}}+\int_0^t\|\nabla p\|^2_{\mathcal{H}^{m-1}}d\tau.
\end{cases}
\end{eqnarray}
We will estimate $\v\|(\nabla^\varphi p, \text{div}^\varphi v)\|_{\mathcal{H}^{m-1}}$, $\v^2\int_{0}^{t}\|\nabla^\varphi\MZ^{m-1}\mbox{div}^\varphi v\|^2 d\tau$  and $\int_0^t\|\nabla  p\|^2_{\mathcal{H}^{m-1}}d\tau$ in the remaining part of this subsection. The estimate of $\int_{0}^{t}\|\nabla v\|_{\mathcal{H}^{m-1}}^2d\tau$ will be left to later subsections.

	The first lemma gives the equations satisfied by $\MZ^\a v$ and $\MZ^\a\nabla^\varphi p$.
	\begin{lemma}
		For $|\alpha|\leq m$,  the equations for $\mathcal{Z}^\alpha\nabla^\varphi p$ and $\mathcal{Z}^\alpha\mbox{div}^\varphi v$ read as
		\begin{eqnarray}
			\begin{cases}
				\mathcal{Z}^\alpha\mbox{div}^\varphi v+\frac{1}{\g p}(\partial{_t^\varphi}+v\cdot\nabla^\varphi)\mathcal{Z}^\alpha p=\mathcal{C}{_T^\alpha}(p),\label{6.1}\\
				 \varrho(\partial{_t^\varphi}+v\cdot\nabla^\varphi)\mathcal{Z}^\alpha v+\mathcal{Z}^\alpha\nabla^\varphi p-\mathcal{Z}^\alpha(2\m\v\mbox{div}^\varphi S^\varphi v+\l\v\nabla^\varphi\mbox{div}^\varphi v)=\mathcal{C}{_M^\alpha}(p),\label{6.2}
			\end{cases}
		\end{eqnarray}
		where the commutators are given by
		\begin{align}
			&\mathcal{C}{_T^\alpha}(p):=-\frac{1}{\g p}([\mathcal{Z}^\alpha, v_y]\nabla_yp+[\mathcal{Z}^\alpha,V_z]\partial_zp+V_z[\mathcal{Z}^\alpha,\partial_z]p)-[\mathcal{Z}^\alpha,\frac{1}{\g p}](\partial{_t^\varphi}p+v\cdot\nabla^\varphi p), \nonumber\\
			&\mathcal{C}{_M^\alpha}(p):=-\varrho([\mathcal{Z}^\alpha,v_y]\nabla_yv+[\mathcal{Z}^\alpha,V_z]\partial_zv+V_z[\mathcal{Z}^\alpha,\partial_z]v)-[\mathcal{Z}^\alpha,\varrho](\partial{_t^\varphi}v+v\cdot\nabla^\varphi v).\nonumber
		\end{align}
	Moreover, these commutators vanish when $|\a|=0$.
	\end{lemma}
	\noindent\textbf{Proof}. The mass equation $\eqref{1.22}_1$ yields  that
	\begin{equation}\label{6.3}
		\mbox{div}^\varphi v+\frac{1}{\g p}(\partial{_t^\varphi}p+v\cdot\nabla^\varphi p)=0.
	\end{equation}
	Applying $\mathcal{Z}^\alpha$ to \eqref{6.3} and $\eqref{1.22}_2$ shows the lemma by trivial calculations.$\hfill\Box$
	
	\

\begin{lemma}\label{lem6.5}
For any $t\in[0,T^\v]$, $m\geq 1$, it holds that
\begin{align}\label{6.50}
&\v\|(\mbox{div}^\varphi v, \nabla^\varphi p)(t)\|^2_{\mathcal{H}^{m-1}}+\v^2\int_{0}^{t}\|\nabla^\varphi\MZ^{m-1}\mbox{div}^\varphi v\|^2 d\tau\nonumber\\
&\leq\Lambda_0\v\|(\mbox{div}^{\varphi_0} v_0,\nabla^{\varphi_0} p_0)\|^2_{\mathcal{H}^{m-1}}+\d_1\int_{0}^{t}\v^2\|\nabla^2 v\|^2_{\mathcal{H}^{m-1}}+\d\int_{0}^{t}\|\nabla p\|^2_{\mathcal{H}^{m-1}}d\tau\nonumber\\
&~~~~~~~~~~~~+\frac{\Lambda_0}{\d_1}\v\int_{0}^{t}\|\nabla V^m\|^2d\tau+C_\d\Lambda_\i(t)\v^2\int_{0}^{t}\|\nabla^2 v\|^2_{\mathcal{H}^{m-2}}d\tau\nonumber\\
&~~~~~~~~~~~~~~~~~~+C_{\d_1}\Lambda_\i(t)(1+\|\nabla\mbox{div}^\varphi v\|^2_{L{_{t,x}^\infty}})\int_0^t\Lambda(Y_m(\tau))+\v|\MZ^mh|_\frac{1}{2}^2d\tau
\end{align}
where $\d$ and $\d_1$ will be chosen later.
	\end{lemma}
\noindent\textbf{Proof}. 
Multiplying $\eqref{6.2}_2$ by $\v\nabla^\varphi\MZ^\a\mbox{div}^\varphi v$ and integrating over $\mathcal{S}\times[0,t]$ lead to that
	\begin{align}\label{6.11}
		 &\v\int_{0}^{t}\int_{\mathcal{S}}\varrho(\partial{_t^\varphi}+v\cdot\nabla^\varphi)\mathcal{Z}^\alpha v\cdot\nabla^\varphi\MZ^\a\mbox{div}^\varphi vd\mathcal{V}_\tau d\tau+\v\int_{0}^{t}\int_{\mathcal{S}}\mathcal{Z}^\alpha\nabla^\varphi p\cdot\nabla^\varphi\MZ^\a\mbox{div}^\varphi vd\mathcal{V}_{\tau}d\tau\nonumber\\
		 &=(2\m+\l)\v^2\int_{0}^{t}\int_{\mathcal{S}}\mathcal{Z}^\alpha\nabla^\varphi\mbox{div}^\varphi v\cdot\nabla^\varphi\MZ^\a\mbox{div}^\varphi vd\mathcal{V}_{\tau}d\tau\nonumber\\
		 &\quad-\m\v^2\int_{0}^{t}\int_{\mathcal{S}}\MZ^\a(\nabla^\varphi\times(\nabla^\varphi\times v))\cdot\nabla^\varphi\MZ^\a\mbox{div}^\varphi vd\mathcal{V}_{\tau}d\tau\nonumber\\
		 &\quad+\v\int_{0}^{t}\int_{\mathcal{S}}\mathcal{C}{_M^\alpha}(p)\cdot\nabla^\varphi\MZ^\a\mbox{div}^\varphi vd\mathcal{V}_{\tau}d\tau=\sum_{i=1}^{3}J_i.
	\end{align}
	Integration by parts shows that
	\begin{align}\label{6.12}
		 &\v\int_{0}^{t}\int_{\mathcal{S}}\varrho(\partial{_t^\varphi}+v\cdot\nabla^\varphi)\mathcal{Z}^\alpha v\cdot\nabla^\varphi\MZ^\a\mbox{div}^\varphi vd\mathcal{V}_{\tau}d\tau\nonumber\\
		 &=-\v\int_{0}^{t}\int_{\mathcal{S}}\varrho(\partial{_t^\varphi}+v\cdot\nabla^\varphi)\mbox{div}^\varphi\mathcal{Z}^\alpha v\MZ^\a\mbox{div}^\varphi vd\mathcal{V}_{\tau}d\tau\nonumber\\
		 &\quad-\v\int_{0}^{t}\int_{\mathcal{S}}(\nabla^\varphi\varrho(\partial_t+v_y\cdot\nabla_y+V_z\partial_z)\mathcal{Z}^\alpha v+\varrho (\nabla^\varphi v)^t\nabla^\varphi\mathcal{Z}^\alpha v)\MZ^\a\mbox{div}^\varphi vd\mathcal{V}_{\tau}d\tau\\
		 &\quad+\v\int_{0}^{t}\int_{z=0}\varrho(\partial_t+v_y\cdot\nabla_y+V_z\partial_z)\mathcal{Z}^\alpha v\cdot\mathbf{N}\MZ^\a\mbox{div}^\varphi vd\mathcal{V}_{\tau}d\tau=\sum_{1}^{3}K_i.\nonumber
	\end{align}
Note that
	\begin{equation*}
		 \MZ[\partial{_i^\varphi},\MZ^\a]f=-\frac{\partial_i\varphi}{\partial_z\varphi}[\MZ,\partial_z]\MZ^\alpha f-\MZ(\frac{\partial_i\varphi}{\partial_z\varphi})\partial_z\MZ^\alpha f+[\partial{_i^\varphi},\MZ\MZ^\a]f, \quad[Z_3,\partial_z]=-\frac{1}{(1-z)^2}\partial_z
	\end{equation*}
It follows from \eqref{2.6}, \eqref{2.7} and \eqref{2.25} that
	\begin{align}\label{6.13}
		 K_1&=-\frac{1}{2}\v\int_{\mathcal{S}}\varrho|\MZ^\a\mbox{div}^\varphi v(t)|^2d\mathcal{V}_t+\frac{1}{2}\v\int_{\mathcal{S}}\varrho_0|\MZ^\a\mbox{div}^{\varphi_0} v(0)|^2d\mathcal{V}_0\nonumber\\
		 &\quad-\v\int_{0}^{t}\int_{\mathcal{S}}\varrho(\partial_t+v_y\cdot\nabla_y+V_z\partial_z)[\mbox{div}^\varphi,\mathcal{Z}^\alpha] v\MZ^\a\mbox{div}^{\varphi_\tau} vd\mathcal{V}_\tau d\tau\nonumber\\
		 &\leq-\frac{1}{2}\v\int_{\mathcal{S}}\varrho|\MZ^\a\mbox{div}^\varphi v(t)|^2d\mathcal{V}_t+\frac{1}{2}\v\int_{\mathcal{S}}\varrho_0|\MZ^\a\mbox{div}^{\varphi_0} v(0)|^2d\mathcal{V}_0\nonumber\\
		 &\quad+\Lambda_\i(t)\v\int_{0}^{t}\|(\partial_t+v_y\cdot\nabla_y+\frac{(1-z)V_z}{z}Z_3)[\mbox{div}^\varphi,\mathcal{Z}^\alpha] v\|\|\MZ^\a\mbox{div}^\varphi v\|d\tau\nonumber\\
		 &\leq-\frac{1}{2}\v\int_{\mathcal{S}}\varrho|\MZ^\a\mbox{div}^{\varphi} v(t)|^2d\mathcal{V}_t+\frac{1}{2}\v\int_{\mathcal{S}}\varrho_0|\MZ^\a\mbox{div}^{\varphi_0} v(0)|^2d\mathcal{V}_0\nonumber\\
		 &\quad~~~~~~+\Lambda_\i(t)\int_{0}^{t}\v\|\nabla v\|{_{\mathcal{H}^{m-1}}^2}+\v|\mathcal{Z}^{m}h|{_\frac{1}{2}^2}d\tau,
	\end{align}
	where one has the fact $V_z|_{z=0}=0$ and
	\begin{equation}\label{6.14}
		\|\frac{V_z}{z}\|_{\infty,t}\leq\|V_z\|_{\infty,t}+\|\partial_z V_z\|_{\infty,t}\leq\Lambda(\frac{1}{c_0},\|(v,\partial_z v)\|_{\infty,t}+|h|_{\mathcal{W}{_t^{2,\i}}}).
	\end{equation}
Similarly, one gets that
	\begin{equation}\label{6.15}
		 K_2\leq\Lambda_\i(t)\int_{0}^{t}\|v\|{_{\mathcal{H}^{m}}^2}+\|\nabla v\|{_{\mathcal{H}^{m-1}}^2}+\v|\MZ^{m-1}h|{_{\frac{1}{2}}^2} d\tau,
	\end{equation}
As to the boundary term $K_3$, one can obtain from the trace estimates \eqref{2.10} that
	\begin{align}\label{6.16}
		K_3&\leq\Lambda_\i(t)\v\int_{0}^{t}|\MZ\MZ^\a v|_{L^2(\mathbb{R}^2)}|\MZ^\a\mbox{div}^\varphi v|_{L^2(\mathbb{R}^2)}d\tau\nonumber\\
		&\leq\Lambda_\i(t)\v\int_{0}^{t}\Big(\|\nabla \MZ\MZ^\a v\|^\frac{1}{2}\|\MZ\MZ^\a v\|^\frac{1}{2}+\|\MZ\MZ^\a v\|\Big)\nonumber\\
		&~~~~~~~~~~~~~~~~~~~~~~~~~~~\cdot\Big(\|\nabla\MZ^\a\text{div}^\varphi v\|^\frac{1}{2}\|\MZ^\a\text{div}^\varphi v\|^\frac{1}{2}+\|\MZ^\a\text{div}^\varphi v\|\Big)d\tau\nonumber\\
		&\leq\d_1\v^2\int_{0}^{t}\|\nabla^2v\|_{\mathcal{H}^{m-1}}^2d\tau+\d_1\v\int_0^t\|\nabla V^m\|^2d\tau\nonumber\\
		&~~~~~~~~~+C_{\d_1}\Lambda_\infty(t)\int_{0}^{t}\|v\|^2_{\mathcal{H}^{m}}+\|\nabla v\|^2_{\mathcal{H}^{m-1}}+\v|\MZ^mh|^2 d\tau.
	\end{align}
where \eqref{6.14} and Corollary \ref{cor2.8} have been used. 

Substituting \eqref{6.13}, \eqref{6.15} and \eqref{6.16} into \eqref{6.12}, one has that
	\begin{align}\label{6.17}
		&\v\int_{0}^{t}\int_{\mathcal{S}}\varrho(\partial{_t^\varphi}+v\cdot\nabla^\varphi)\mathcal{Z}^\alpha v\cdot\nabla^\varphi\MZ^\a\mbox{div}^\varphi vd\mathcal{V}_\tau d\tau\nonumber\\
		&\leq-\frac{1}{2}\v\int_{\mathcal{S}}\varrho|\MZ^\a\mbox{div}^\varphi v(t)|^2d\mathcal{V}_t+\frac{1}{2}\v\int_{\mathcal{S}}\varrho_0|\MZ^\a\mbox{div}^{\varphi_0} v(0)|^2d\mathcal{V}_0+\d_1\v^2\int_{0}^{t}\|\nabla^2v\|_{\mathcal{H}^{m-1}}^2d\tau\nonumber\\
		&~~~~~~~~~~~+\d_1\v\int_0^t\|\nabla V^m\|^2d\tau+C_{\d_1}\Lambda_\infty(t)\int_{0}^{t}\|v\|^2_{\mathcal{H}^{m}}+\|\nabla v\|^2_{\mathcal{H}^{m-1}}+|\MZ^mh|^2 d\tau.
	\end{align}
	For the pressure term, it follows from   \eqref{6.3} that
	\begin{align}\label{6.18}
		&\v\int_{0}^{t}\int_{\mathcal{S}}\MZ^\a\nabla^\varphi p\cdot\nabla^\varphi\MZ^\alpha\mbox{div}^\varphi vd\mathcal{V}_{\tau}d\tau\nonumber\\
		&=\v\int_{0}^{t}\int_{\mathcal{S}}\MZ^\a\nabla^\varphi p\cdot\MZ^\alpha\nabla^\varphi\mbox{div}^\varphi vd\mathcal{V}_{\tau}d\tau+\v\int_{0}^{t}\int_{\mathcal{S}}\MZ^\a\nabla^\varphi p[\nabla^\varphi,\MZ^\alpha]\mbox{div}^\varphi vd\mathcal{V}_{\tau}d\tau\nonumber\\
		&=\v\int_{0}^{t}\int_{\mathcal{S}}\MZ^\a\nabla^\varphi p\cdot\Big([\nabla^\varphi,\MZ^\alpha]\mbox{div}^\varphi v+\MZ^\a(\nabla^\varphi(\frac{1}{\g p})\partial{_t^\varphi}p+\nabla^\varphi(\frac{v}{\g p})\cdot\nabla^\varphi p)\Big)d\mathcal{V}_{\tau}d\tau\nonumber\\
		&\qquad-\v\int_{0}^{t}\int_{\mathcal{S}}\MZ^\a\nabla^\varphi p\cdot\MZ^\a(\frac{1}{\g p}(\partial{_t^\varphi}+v\cdot\nabla^\varphi)\nabla^\varphi p)d\mathcal{V}_{\tau}d\tau=:L_1+L_2.
	\end{align}
Using \eqref{2.6} and \eqref{2.25}, one obtains immediately  that
	\begin{align}\label{6.19}
L_1&\leq\d\int_0^t\|\nabla p\|^2_{\mathcal{H}^{m-1}}d\tau+\Lambda_\i(t)(1+\|\nabla\mbox{div}^\varphi v\|^2_{\infty,t})\int_{0}^{t}\Lambda(Y_m) d\tau\nonumber\\
&~~~~~~~~~~~~~~~~~ +C_\d\Lambda_\i(t)\v^2\int_{0}^{t}\|\nabla^2 v\|^2_{\mathcal{H}^{m-2}}d\tau
	\end{align}
As to $L_2$, note that
	\begin{align}
		&\MZ^\a(\frac{1}{\g p}(\partial{_t^\varphi}+v\cdot\nabla^\varphi)\nabla^\varphi p)=\frac{1}{\g p}(\partial{_t^\varphi}+v\cdot\nabla^\varphi)\MZ^\a\nabla^\varphi p+\frac{V_z}{\g p}[\MZ^\a,\partial_z]\nabla^\varphi p\nonumber\\
		&\quad+\sum_{\substack{\beta+\nu=\a\\|\beta|\geq 1}}c_{\beta,\nu}\Big(\MZ^\beta(\frac{1}{\g p})\MZ^\nu\partial_t\nabla^\varphi p+\MZ^\beta(\frac{v_y}{\g p})\MZ^\nu\nabla_y\nabla^\varphi p+\MZ^\beta(\frac{V_z}{\g p})\MZ^\nu\partial_z\nabla^\varphi p\Big),
	\end{align}
	then, it follows from  integration by parts  and \eqref{2.23} that
	\begin{align}\label{6.21}
		&-\v\int_{0}^{t}\int_{\mathcal{S}}\MZ^\a\nabla^\varphi p\cdot \frac{1}{\g p}(\partial{_t^\varphi}+v\cdot\nabla^\varphi)\MZ^\a\nabla^\varphi pd\mathcal{V}_\tau d\tau\nonumber\\
		&=-\v\int_{\mathcal{S}}\frac{1}{2\g p}|\MZ^\a\nabla^\varphi p|^2d\mathcal{V}_t+\v\int_{\mathcal{S}}\frac{1}{2\g p_0}|\MZ^\a\nabla^{\varphi_0} p_0|^2d\mathcal{V}_0\nonumber\\
		&\quad~~~~~~~~~~~~~~+\v\int_{0}^{t}\int_{\mathcal{S}}(\partial{_t^\varphi}(\frac{1}{2\g p})+\mbox{div}^\varphi(\frac{v}{2\g p}))|\MZ^\a\nabla^\varphi p|^2d\mathcal{V}_td\tau\nonumber\\
		&\leq-\v\int_{\mathcal{S}}\frac{1}{2\g p}|\MZ^\a\nabla^\varphi p|^2d\mathcal{V}_t+\v\int_{\mathcal{S}}\frac{1}{2\g p_0}|\MZ^\a\nabla^{\varphi_0} p_0|^2d\mathcal{V}_0+\Lambda_\i(t)\int_{0}^{t}\v\|\nabla p\|_{\mathcal{H}^{m-1}}^2+\v|\MZ^{m-1}h|_\frac{1}{2}^2d\tau.
	\end{align}
	Since $V_z$ vanishes on the boundary, one has that
	\begin{equation}
		-\v\int_{0}^{t}\int_{\mathcal{S}}\MZ^\a\nabla^\varphi p\cdot \frac{V_z}{\g p}[\MZ^\alpha,\partial_z]\nabla^\varphi pd\mathcal{V}_\tau d\tau\leq\Lambda_\infty(t)\int_{0}^{t}\v\|\nabla p\|{_{\mathcal{H}^{m-1}}^2}+\v|\MZ^{m-1}h|{_\frac{1}{2}^2}d\tau.
	\end{equation}
	For the term like
	\begin{equation*}
		-\v\int_{0}^{t}\int_{\mathcal{S}}c_{\beta,\nu}\MZ^\a\nabla^\varphi p\MZ^\beta(\frac{1}{\g p})\MZ^\nu\partial_t\nabla^\varphi p,
	\end{equation*}
	where $\beta$ and $\nu$ satisfy that $\beta\neq 0$ and $\beta+\nu=\alpha$, it follows from  \eqref{2.6} that
	\begin{align}\label{6.22}
		&-\v\int_{0}^{t}\int_{\mathcal{S}}c_{\beta,\nu}\MZ^\a\nabla^\varphi p\MZ^\beta(\frac{1}{\g p})\MZ^\nu Z_0\nabla^\varphi p\leq\int_{0}^{t}\v\|\nabla^\varphi p\|{_{\mathcal{H}^{m-1}}^2}d\tau+\v\int_{0}^{t}\|\MZ^{\beta-1}\MZ(\frac{1}{\g p})\MZ^\nu Z_0\nabla^\varphi p\|^2d\tau\nonumber\\
		&~~~~~~~~~~~~~~~~~~~~~~~~~~~~~~\leq\Lambda_\i(t)\int_{0}^{t}\v\|\nabla p\|{_{\mathcal{H}^{m-1}}^2}+\|p\|{_{\mathcal{H}^{m-1}}^2}+\v|\mathcal{Z}^{m-1}h|{_\frac{1}{2}^2}d\tau.
	\end{align}
	In a similar way, one gets that
	\begin{equation}\label{6.23}
		-\v\int_{0}^{t}\int_{\mathcal{S}}c_{\beta,\nu}\MZ^\a\nabla^\varphi p\MZ^\beta(\frac{v_y}{\g p})\MZ^\nu\nabla_y\nabla^\varphi p\leq\Lambda_\i(t)\int_0^tY_m(\tau)d\tau.
	\end{equation}
	To estimate the term like
	\begin{equation*}
		-\v\int_{0}^{t}\int_{\mathcal{S}}c_{\beta,\nu}\MZ^\a\nabla^\varphi p\MZ^\beta(\frac{V_z}{\g p})\MZ^\nu\partial_z\nabla^\varphi p,
	\end{equation*}
	one can rewrite it as
	\begin{equation*}
		 -\v\int_{0}^{t}\int_{\mathcal{S}}c_{\tilde{\beta},\tilde{\nu}}\MZ^\a\nabla^\varphi p\MZ^{\tilde{\beta}}(\frac{(1-z)V_z}{\g pz})Z_3\MZ^{\tilde{\nu}}\nabla^\varphi p,
	\end{equation*}
	where $|\tilde{\beta}|+|\tilde{\nu}|\leq m-1$, $\tilde{\nu}\leq m-2$ and $c_{\tilde{\beta},\tilde{\nu}}$ are smooth bounded functions depend only on $z$. Indeed, this follows from the fact that $Z_3((1-z)/z)=\tilde{c}(1-z)/z$ for some smooth bounded function $\tilde{c}$.
	
	If $\tilde{\beta}=0$, $|\tilde{\nu}|\leq m-2$, then \eqref{6.14} implies  that
	\begin{equation}\label{6.24}
		-\v\int_{0}^{t}\int_{\mathcal{S}}\MZ^\a\nabla^\varphi p(\frac{(1-z)V_z}{\g pz})Z_3\MZ^{\tilde{\nu}}\nabla^\varphi p\leq\Lambda_\i(t)\int_{0}^{t}\v\|\nabla p\|{_{\mathcal{H}^{m-1}}^2}+\v|\MZ^{m-1}h|{_\frac{1}{2}^2}d\tau.
	\end{equation}
	When $\tilde{\beta}\neq 0$, it follows from \eqref{2.6}   that
	\begin{align}\label{6.25}
	    &\quad-\v\int_{0}^{t}\int_{\mathcal{S}}\MZ^\a\nabla^\varphi p\MZ^{\tilde{\beta}}(\frac{(1-z)V_z}{\g pz})Z_3\MZ^{\tilde{\nu}}\nabla^\varphi p\nonumber\\
		&\leq\v\Big\{\int_{0}^{t}\|\MZ^\a\nabla^\varphi p\|^2d\tau\Big\}^\frac{1}{2}\Big\{\int_{0}^{t}\|\MZ^{\tilde{\beta}}(\frac{(1-z)V_z}{\g pz})Z_3\MZ^{\tilde{\nu}}\nabla^\varphi p\|^2d\tau\Big\}^\frac{1}{2}\\
		&\leq\Lambda_{\infty}(t)\v\Big\{\int_{0}^{t}\|\nabla^\varphi p\|{_{\mathcal{H}^{m-1}}^2}d\tau\Big\}^\frac{1}{2}\Big\{\int_{0}^{t}\|\nabla^\varphi p\|{_{\mathcal{H}^{m-1}}^2}+\|\frac{(1-z)V_z}{\g pz}\|{_{\mathcal{H}^{m-1}}^2}d\tau\Big\}^\frac{1}{2},\nonumber
	\end{align}
	where one has used \eqref{4.23}, \eqref{4.24} and \eqref{6.14}.
Using \eqref{2.6} and the Hardy inequality (Lemma 8.4 in \cite{Masmoudi-R-1}), one gets that
	\begin{align}\label{6.27}
		\int_{0}^{t}\|\frac{(1-z)V_z}{\g pz}\|{_{\mathcal{H}^{m-1}}^2}d\tau&\leq\Lambda_\i(t)\int_{0}^{t}\|p\|{_{\mathcal{H}^{m-1}}^2}+\|\frac{1-z}{z}V_z\|{_{\mathcal{H}^{m-1}}^2}d\tau\nonumber\\
		&\leq\Lambda_\i(t)\int_{0}^{t}\|p\|{_{\mathcal{H}^{m-1}}^2}+\|V_z\|{_{\mathcal{H}^{m-1}}^2}+\|\nabla V_z\|{_{\mathcal{H}^{m-1}}^2}d\tau.
	\end{align}
	It follows from \eqref{4.16} that
	\begin{equation}\label{6.28}
		\int_{0}^{t}\|\nabla V_z\|{_{\mathcal{H}^{m-1}}^2}d\tau\leq\Lambda_\i(t)\int_{0}^{t}\|v\|{_{\mathcal{H}^{m-1}}^2}+\|\nabla v\|{_{\mathcal{H}^{m-1}}^2}+|\MZ^{m}h|{_\frac{1}{2}^2}d\tau.
	\end{equation}
	Combining \eqref{6.24}-\eqref{6.28} leads to that
	\begin{align}\label{6.29}
		\quad-\v\int_{0}^{t}\int_{\mathcal{S}}\MZ^\a\nabla^\varphi p\MZ^{\tilde{\beta}}(\frac{(1-z)V_z}{\g pz})Z_3\MZ^{\tilde{\nu}}\nabla^\varphi p\leq\Lambda_\i(t)\int_{0}^{t}Y_m(\tau)+\v|\MZ^{m}h|{_\frac{1}{2}^2}d\tau.
	\end{align}
Substituting \eqref{6.19}, \eqref{6.21}-\eqref{6.23} and \eqref{6.29} into \eqref{6.18} yields that
	\begin{align}\label{6.30}
		&\v\int_{0}^{t}\int_{\mathcal{S}}\MZ^\a\nabla^\varphi p\nabla^\varphi\MZ^\alpha\mbox{div}^\varphi vd\mathcal{V}_\tau d\tau\nonumber\\
		&\leq-\v\int_{\mathcal{S}}\frac{1}{2\g p}|\MZ^\a\nabla^\varphi p|^2d\mathcal{V}_t+\v\int_{\mathcal{S}}\frac{1}{2\g p_0}|\MZ^\a\nabla^{\varphi_0} p_0|^2d\mathcal{V}_0+C_\d\Lambda_\i(t)\v^2\int_{0}^{t}\|\nabla^2v\|^2_{\mathcal{H}^{m-2}}d\tau
		\nonumber\\
		&~~~~~~~~~~~~~+\d\int_0^t\|\nabla^\varphi p\|^2_{\mathcal{H}^{m-1}}d\tau+\Lambda_\i(t)(1+\|\nabla\mbox{div}^\varphi v\|_{\i,t}^2)\int_{0}^{t}Y_m(\tau)d\tau.
	\end{align}
	Now, we  estimate the RHS of \eqref{6.11}. It follows from  \eqref{2.7} \eqref{4.23} and \eqref{4.26} that
	\begin{align}\label{6.31}
		 J_1&=(2\m+\l)\v^2\int_{0}^{t}\int_{\mathcal{S}}|\nabla^\varphi\MZ^\a\mbox{div}^\varphi v|^2d\mathcal{V}_\tau d\tau+(2\m+\l)\v^2\int_{0}^{t}\int_{\mathcal{S}}[\MZ^\a,\nabla^\varphi]\mbox{div}^\varphi v\cdot\nabla^\varphi\MZ^\a\mbox{div}^\varphi vd\mathcal{V}_\tau d\tau\nonumber\\
		 &\geq\frac{2\m+\l}{2}\v^2\int_{0}^{t}\int_{\mathcal{S}}|\nabla^\varphi\MZ^\a\mbox{div}^\varphi v|^2d\mathcal{V}_\tau d\tau-\Lambda_\i(t)(1+\|\nabla\mbox{div}^\varphi v\|_{\i,t}^2)\int_{0}^{t}|\MZ^{m-1}h|{_\frac{1}{2}^2}d\tau\nonumber\\
		 &\qquad\qquad\qquad\qquad\qquad\qquad\qquad\qquad-\Lambda_\i(t)\v^2\int_{0}^{t}\|\nabla^2 v\|^2_{\mathcal{H}^{m-2}}d\tau,
	\end{align}
and 
    \begin{align}\label{6.33}
    |J_3|&\leq\f18(2\m+\l)\v^2\int_{0}^{t}\|\nabla^\varphi\MZ^{m-1}\mbox{div}^\varphi v\|^2d\tau+C\int_{0}^{t}\|\mathcal{C}_M^\a(p)\|^2d\tau\nonumber\\
    &\leq\f18(2\m+\l)\v^2\int_{0}^{t}\|\nabla\MZ^{m-1}\mbox{div}^\varphi v\|^2d\tau+\Lambda_\i(t)\int_{0}^{t}Y_m(\tau)d\tau.
    \end{align}	
Using integration by parts, one can obtain from \eqref{2.24} and the trace estimates \eqref{2.10} that
	\begin{align}\label{6.50-3}
	&|J_2|=\Big|\m\v^2\int_{0}^{t}\int_{\mathcal{S}}\MZ^\a(\nabla^\varphi\times(\nabla^\varphi\times v))\cdot\nabla^\varphi\MZ^\a\mbox{div}^\varphi vd\mathcal{V}_\tau d\tau\Big|\nonumber\\
	&\leq\f18(2\m+\l)\v^2\int_{0}^{t}\|\nabla^\varphi\MZ^\a\text{div}^\varphi v\|^2d\tau+\Lambda_0\int_0^t\v^2\|\nabla^2 v\|^2_{\mathcal{H}^{m-2}}+\Lambda_\i(t)\int_{0}^{t}Y_m(\tau)d\tau\nonumber\\
	&~~~~~~+\Big|\m\v^2\int_{0}^{t}\int_{\mathcal{S}}(\nabla^\varphi\times\MZ^\a(\nabla^\varphi\times v))\cdot\nabla^\varphi\MZ^\a\mbox{div}^\varphi vd\mathcal{V}_\tau d\tau\Big|\nonumber\\
	&\leq\f18(2\m+\l)\v^2\int_{0}^{t}\|\nabla^\varphi\MZ^\a\text{div}^\varphi v\|^2d\tau+
	\Lambda_0\v^2\int_0^t\|\nabla^2 v\|^2_{\mathcal{H}^{m-2}}+\Lambda_\i(t)\int_{0}^{t}Y_m(\tau)d\tau\nonumber\\
	&~~~~~~~+\m\v^2\Big|\int_{0}^{t}\int_{z=0}(\mathbf{N}\times\MZ^\a(\nabla^\varphi\times v))\cdot(\nabla_y,0)\MZ^\a\mbox{div}^\varphi vdyd\tau\Big|\nonumber\\
	&\leq\f18(2\m+\l)\v^2\int_{0}^{t}\|\nabla^\varphi\MZ^\a\text{div}^\varphi v\|^2d\tau+\Lambda_0\int_0^t\v^2\|\nabla^2 v\|^2_{\mathcal{H}^{m-2}}+
	\Lambda_{\infty}(t)\int_0^tY_m(\tau)d\tau\nonumber\\
	&~~~~~~~~+\m\v^2\int_{0}^{t}|(\mathbf{N}\times\MZ^\a(\nabla^\varphi\times v))|_{\f12}\cdot|\MZ^\a\mbox{div}^\varphi v|_{\f12}d\tau\nonumber\\
	&\leq \f14(2\m+\l)\v^2\int_{0}^{t}\|\nabla^\varphi\MZ^{m-1}\text{div}^\varphi v\|^2d\tau+\d_1\v^2\int_0^t\|\nabla^2v\|^2_{\mathcal{H}^{m-1}}d\tau+\f{\Lambda_0}{\d_1}\v^2\int_0^t\|\nabla V^m\|^2d\tau\nonumber\\
	&~~~~~~~~~~~+C_{\d_1}\Lambda_{\infty}(t)\int_0^tY_m(\tau)+\v^2\|\nabla^2 v\|_{\mathcal{H}^{m-2}}^2+\v|\MZ^m h|^2_{\f12}d\tau.
	\end{align}
Combining \eqref{6.11}, \eqref{6.17}, \eqref{6.30}-\eqref{6.50-3} and summing over $\a$($|\a|\leq m-1$) complete the proof of this lemma. $\hfill\Box$

\

In order to close the estimate in Lemma \ref{lem6.5}, it remains to bound $\v^2\int_{0}^{t}\|\nabla^2 v\|{_{\mathcal{H}^{m-2}}^2}d\tau$. We first estimate $\int_0^t\|\nabla\MZ^{m-2}\mbox{div}^{\varphi}u\|^2d\tau$ by using of the mass equation \eqref{6.1}.
It should be mentioned that, in general, it is difficult to derive an uniform estimate for $\int_0^t\| \MZ^{m-2}\partial_{zz}u\|^2d\tau$ due to the possible appearance of boundary layers. However, one can expect the estimate of  $\int_0^t\|\nabla\MZ^{m-2}\mbox{div}^{\varphi}u\|^2d\tau$ due to weaker boundary layer for  $\mbox{div}u$. Moreover, this estimate is very helpful to bound $\v\int_{0}^{t}\|\nabla^2 v\|{_{\mathcal{H}^{m-2}}^2}d\tau$, which will be used later. Precisely, we can prove the following lemma.
	\begin{lemma}\label{l6.3}
		For $m\geq 2$, it holds that
		\begin{equation}\label{6.37-1}
			\int_{0}^{t}\|\nabla\MZ^{m-2}\mbox{div}^\varphi v\|^2d\tau\leq\Lambda_\i(t)\int_{0}^{t}\|\nabla p(\tau)\|{_{\mathcal{H}^{m-1}}^2}+Y_m(\tau)d\tau,
		\end{equation}
		and
		\begin{equation}\label{6.37-2}
			 \v^2\int_0^t\|\nabla^2v\|^2_{\mathcal{H}^{m-2}}d\tau+\v^2\int_0^t\|\nabla^\varphi\nabla^\varphi v\|^2_{\mathcal{H}^{m-2}}d\tau\leq \Lambda_\i(t)\int_{0}^{t}Y_m(\tau)d\tau.
		\end{equation}

	\end{lemma}
	\noindent\textbf{Proof}. 
For $|\a|\leq m-2$, $m\geq 2$, it follows from \eqref{6.1} that
	\begin{align}\label{6.36}
		\nabla\MZ^\a\mbox{div}^\varphi v&=-\nabla(\frac{1}{\g p}(\partial_t+v_y\cdot\nabla_y+V_z\partial_z)\mathcal{Z}^\alpha p)\nonumber\\
		&\quad-\nabla(\frac{1}{\g p}([\mathcal{Z}^\alpha, v_y]\nabla_yp)-\nabla([\mathcal{Z}^\alpha,\frac{1}{\g p}](\partial_tp+v_y\cdot\nabla_y p))\\
		&\quad-\nabla(\frac{1}{\g p}([\mathcal{Z}^\alpha,V_z]\partial_zp+V_z[\mathcal{Z}^\alpha,\partial_z]p)+[\mathcal{Z}^\alpha,\frac{1}{\g p}](V_z\partial_z p))=:\sum_{i=1}^{3}J_i.\nonumber
	\end{align}
	It follows immediately from \eqref{2.6}-\eqref{2.7} and the fact $V_z|_{z=0}=0$ that
	\begin{equation}\label{6.37}
		\int_{0}^{t}\int_{\mathcal{S}}|J_1|^2d\mathcal{V}_\tau d\tau\leq\Lambda_\i(t)\int_{0}^{t}\|(p,\nabla p)\|{_{\mathcal{H}^{m-1}}^2},
	\end{equation}
	and
	\begin{equation}\label{6.38}
		 \int_{0}^{t}\int_{\mathcal{S}}|J_2|^2d\mathcal{V}_\tau d\tau\leq\Lambda_\i(t)\int_{0}^{t}\|(p,v)\|{_{\mathcal{H}^{m-2}}^2}+\|(\nabla p,\nabla v)\|{_{\mathcal{H}^{m-2}}^2}d\tau.
	\end{equation}
	To estimate $J_3$, one first notes that
	\begin{equation*}
		\nabla([\MZ^\a, V_z]\partial_z p)=\sum_{\substack{\beta+\nu=\a\\|\beta|\geq 1}}c_{\beta,\nu}(\nabla\MZ^\beta V_z\MZ^\nu\partial_zp+\MZ^\beta V_z\nabla\MZ^\nu\partial_zp)
	\end{equation*}
	Due to \eqref{2.6} and \eqref{6.28}, one can obtain that
	\begin{align}
		\int_{0}^{t}\|\nabla\MZ^\beta V_z\MZ^\nu\partial_zp\|^2d\tau&\leq\sum_{|\tilde{\b}|\leq m-2}\int_{0}^{t}\|\MZ^{\tilde{\beta}}\nabla V_z\MZ^\nu\partial_zp\|^2d\tau\nonumber\\
		&\leq\Lambda_\i(t)\int_{0}^{t}\|\nabla V_z\|{_{\mathcal{H}^{m-2}}^2}+\|\partial_zp\|{_{\mathcal{H}^{m-2}}^2}d\tau\nonumber\\
		 &\leq\Lambda_\i(t)\int_{0}^{t}\|v\|{_{\mathcal{H}^{m-2}}^2}+\|\nabla( v ,p)\|{_{\mathcal{H}^{m-2}}^2}+|h|{_{\mathcal{H}^{m}}^2}d\tau.
	\end{align}
 On the other hand, it follows from \eqref{2.6}, Hardy inequality, \eqref{4.26} and \eqref{6.28} that
	\begin{align}
		&\int_{0}^{t}\|\MZ^\beta V_z\nabla\MZ^\nu\partial_zp\|^2d\tau\leq\int_{0}^{t}\|\MZ^\beta (\frac{1-z}{z}V_z)\MZ\MZ^\nu\partial_zp\|^2d\tau\nonumber\\
		&\leq\Lambda_\i(t) \int_{0}^{t}\|\frac{1-z}{z}V_z\|{_{\mathcal{H}^{m-2}}^2}+\|\partial_zp\|{_{\mathcal{H}^{m-1}}^2}d\tau\leq\Lambda_\i(t)\int_{0}^{t}\|(V_z,\nabla V_z)\|{_{\mathcal{H}^{m-2}}^2}+\|\partial_zp\|{_{\mathcal{H}^{m-1}}^2}d\tau\nonumber\\
		 &\leq\Lambda_\i(t)\int_{0}^{t}\|(v,\nabla v)\|{_{\mathcal{H}^{m-2}}^2}+\|\nabla p\|{_{\mathcal{H}^{m-1}}^2}+|h|{_{\mathcal{H}^{m}}^2}d\tau.
	\end{align}
	Thus, 
	\begin{align}\label{6.39}
		&\int_{0}^{t}\int_{\mathcal{S}}|\nabla(\frac{1}{\g p}[\MZ^\a,V_z]\partial_zp)|^2d\mathcal{V}_\tau d\tau\leq\Lambda_\i(t)\int_{0}^{t}\int_{\mathcal{S}}|[\MZ^\a,V_z]\partial_zp|^2+|\nabla([\MZ^\a,V_z]\partial_zp)|^2d\mathcal{V}_\tau d\tau\nonumber\\
		&\leq\Lambda_\i(t)\int_{0}^{t}\|V_z\|{_{\mathcal{H}^{m-2}}^2}+\|\partial_zp\|{_{\mathcal{H}^{m-3}}^2}+\|\nabla\MZ^\beta V_z\MZ^\nu\partial_zp\|^2+\|\MZ^\beta V_z\nabla\MZ^\nu\partial_zp\|^2d\tau\nonumber\\
	    &\leq\Lambda_\i(t)\int_{0}^{t}\|(v,\nabla v)\|{_{\mathcal{H}^{m-2}}^2}+\|\nabla p\|{_{\mathcal{H}^{m-1}}^2}+|h|{_{\mathcal{H}^{m}}^2}d\tau.
	\end{align}
	Similarly, 
	\begin{align}\label{6.40}
		\int_{0}^{t}\int_{\mathcal{S}}|\nabla([\MZ^\a,\frac{1}{\g p}])(V_z\partial_zp)|^2d\mathcal{V}_\tau d\tau\leq\Lambda_\i(t)\int_{0}^{t}\|(v,p,\nabla v,\nabla p)\|{_{\mathcal{H}^{m-2}}^2}+|h|{_{\mathcal{H}^{m}}^2}d\tau.
	\end{align}
	Finally, it follows from \eqref{2.7}, \eqref{4.24} and \eqref{6.14} that
	\begin{equation}\label{6.41}
		\int_{0}^{t}\int_{\mathcal{S}}|\nabla(\frac{1}{\g p}V_z[\MZ^\a,\partial_z]p)|^2d\mathcal{V}_\tau d\tau
		\leq\Lambda_\i(t)\int_{0}^{t}\|\nabla p\|{_{\mathcal{H}^{m-2}}^2}d\tau.
	\end{equation}
	Combining \eqref{6.39}, \eqref{6.40} and \eqref{6.41}, one gets that
	\begin{equation}\label{6.42}
		\int_{0}^{t}\int_{\mathcal{S}}|J_3|^2d\mathcal{V}_\tau d\tau\leq\Lambda_\i(t)\int_{0}^{t}\|(v,\nabla v)\|{_{\mathcal{H}^{m-2}}^2}+\|(p,\nabla p)\|{_{\mathcal{H}^{m-1}}^2}+|h|{_{\mathcal{H}^{m}}^2}d\tau.
	\end{equation}
Then, \eqref{6.37-1} follows from \eqref{6.36}-\eqref{6.38} and \eqref{6.42}. 

By $\eqref{1.22}_2$, \eqref{6.37-1} and Corollary \ref{cor2.8}, one  obtains  that 
\begin{align}\label{6.37-10}
\v^2\int_0^t\|\partial_z^2v\|^2_{\mathcal{H}^{m-2}}d\tau &\leq\v^2\int_{0}^{t}\|\nabla\MZ^{m-2}\mbox{div}^\varphi v\|d\tau+\Lambda_\i(t)\int_{0}^{t}Y_m(\tau)d\tau\nonumber\\
&\leq \Lambda_\i(t)\int_{0}^{t}Y_m(\tau)d\tau.
\end{align}
Then, \eqref{6.37-2} follows immediately.   Therefore, the lemma is proved. $\hfill\Box$

Now we can estimate $\int_0^t\|\nabla p\|^2_{\mathcal{H}^{m-1}}d\tau$  in the following lemma.
	\begin{lemma}\label{lem6.4-11}
		For every $m\geq 1$, it holds that
		\begin{equation}\label{6.46-1}
			\int_{0}^{t}\|\nabla p\|{_{\mathcal{H}^{m-1}}^2}d\tau\leq \Lambda_0\v^2\int_0^t\|\nabla^2v\|^2_{\mathcal{H}^{m-1}}d\tau+\Lambda_\i(t)\int_{0}^{t} \Lambda(Y_m(\tau))+\v|\MZ^{m-1}\nabla_yh|_\frac{1}{2}^2d\tau.
		\end{equation}
 
	\end{lemma}
	\noindent\textbf{Proof}.   For any $|\alpha|\leq m-1$, it follows from $\eqref{6.2}_2$, \eqref{2.6}, \eqref{2.7}, and \eqref{2.24} that
	\begin{align}\label{6.43}
		&\int_{0}^{t}\|\MZ^\a \nabla^\varphi p\|^2d\tau\lesssim\int_{0}^{t}\|\varrho(\partial_t+v_y\cdot\nabla_y+V_z\partial_z)\mathcal{Z}^\alpha v\|^2d\tau\nonumber\\
		&~~~~~~~+\v^2\int_{0}^{t}\|\mbox{div}^\varphi S^\varphi v\|{_{\mathcal{H}^{m-1}}^2}+\|\nabla^\varphi\mbox{div}^\varphi v\|{_{\mathcal{H}^{m-1}}^2}d\tau+\int_{0}^{t}\|[\mathcal{Z}^\alpha,\varrho](\partial_tv+v_y\cdot\nabla_yv+V_z\partial_zv)\|^2d\tau\nonumber\\
		&~~~~~~~+\int_{0}^{t}\|\varrho([\mathcal{Z}^\alpha,v_y]\nabla_yv+[\mathcal{Z}^\alpha,V_z]\partial_zv+V_z[\mathcal{Z}^\alpha,\partial_z]v)\|^2d\tau\nonumber\\
		&\leq\Lambda_\i(t)\int_{0}^{t}\|v\|{_{\mathcal{H}^m}^2}+\|\nabla v\|{_{\mathcal{H}^{m-2}}^2}+|h|{_{\mathcal{H}^m}^2}+\|p\|{_{\mathcal{H}^m}^2}d\tau
		+C\v^2\int_0^t\|\nabla^\varphi\nabla^\varphi v\|{_{\mathcal{H}^{m-1}}^2}d\tau\nonumber\\
		&\leq \Lambda_0\v^2\int_0^t\|\nabla^2v\|^2_{\mathcal{H}^{m-1}}d\tau+\Lambda_\i(t)\int_{0}^{t} Y_m(\tau)+\v|\MZ^{m-1}\nabla_yh|_\frac{1}{2}^2d\tau,
	\end{align}
	which combining \eqref{2.23-1} yields the proof of this lemma.
	$\hfill\Box$
	
	\
	
	Since $\mbox{div}^\varphi v$ is not free, the following estimate is the key to estimate $\int_{0}^{t}\|\partial_zv\|_{\mathcal{H}^{m-1}}^2d\tau$:
		\begin{lemma}\label{div}
			For $m\geq 6$, it holds that
			\begin{equation}\label{9.4-0}
				\|\mbox{div}^\varphi v\|_{\mathcal{H}^{m-1}}^2\leq\Lambda(\frac{1}{c_0},\mathcal{Q}_m(t)+\mathcal{Q}(t)),
			\end{equation}
	where $\mathcal{Q}_m(t)\doteq\|(v,p)(t)\|_{\mathcal{H}^m}^2+\|(\nabla v,\nabla p)(t)\|{_{\mathcal{H}^{m-2}}^2}+|h(t)|{_{\mathcal{H}^m}^2}$.
		\end{lemma}
		\noindent\textbf{Proof}. Applying $\MZ^\alpha$ with $|\alpha|\leq m-1$ to \eqref{6.3} gives that
		\begin{equation}\label{9.5}
			\MZ^\alpha\mbox{div}^\varphi v=-\MZ^\alpha(\frac{\partial_tp}{\g p}+\frac{v_y\cdot\nabla_y p}{\g p})-\MZ^\alpha(\frac{V_z\partial_zp}{\g p}).
		\end{equation}
		Note that
		\begin{equation}
			 \MZ^\alpha(\frac{\partial_tp}{p})=\sum_{|\beta|=0}^{|\alpha|}\MZ^\beta(\frac{1}{p})\MZ^\nu\partial_tp,
		\end{equation}
		then it follows that
		\begin{equation}
			 \|\frac{\partial_tp}{p}\|_{\mathcal{H}^{m-1}}^2\leq\sum_{\substack{|\beta|=0\\|\beta|+|\nu|\leq m-1}}^{|\beta|=[\frac{m}{2}]}\|\MZ^\beta(\frac{1}{p})\|{_\i^2}\|\MZ^\nu \partial_tp\|^2+\sum_{\substack{|\beta|=1+[\frac{m}{2}]\\|\beta|+|\nu|\leq m-1}}^{|\beta|=m-1}\|\MZ^\nu \partial_tp\|{_\i^2}\|\MZ^\beta(\frac{1}{p})\|^2.
		\end{equation}
		For the $L^\i$ norm, it holds, for $|\beta|\leq[\frac{m}{2}]$,  that
		\begin{equation}
			\sum_{|\beta|\leq[\frac{m}{2}]}\|\MZ^\beta p\|{_{L^\i}^2}\lesssim\sum_{|\beta|\leq[\frac{m}{2}]}\|\nabla\MZ^\beta p\|_1\|\MZ^\beta p\|_2\lesssim\|p\|{_{\mathcal{H}^{[\frac{m}{2}]+2}}^2}+\|\nabla p\|{_{\mathcal{H}^{[\frac{m}{2}]+1}}^2}.
		\end{equation}
		Thus, for $m\geq 6$, that is $[\frac{m}{2}]+2\leq m$ and $[\frac{m}{2}]+1\leq m-2$, one gets that
		\begin{equation}
			\sum_{|\beta|\leq[\frac{m}{2}]}\|\MZ^\beta p\|{_\i^2}\leq\mathcal{Q}_m(t),
		\end{equation}
		which implies that
		\begin{equation}\label{9.10}
			\|\frac{\partial_tp}{p}\|_{\mathcal{H}^{m-1}}^2\leq \Lambda(\frac{1}{c_0},\mathcal{Q}_m(t)+\mathcal{Q}(t)),\quad m\geq 6.
		\end{equation}
		Similarly, one has that
		\begin{equation}
			\sum_{|\beta|\leq[\frac{m}{2}]}\|\MZ^\beta v\|{_\i^2}\leq \Lambda(\frac{1}{c_0},\mathcal{Q}_m(t)+\mathcal{Q}(t)),\quad m\geq 6,
		\end{equation}
		which yields that, for $m\geq 6$,
		\begin{align}\label{9.12}
			\|\frac{v_y\cdot\nabla_y p}{p}(t)\|_{\mathcal{H}^{m-1}}^2&\leq\sum_{\substack{|\beta|=0\\|\beta|+|\nu|\leq m-1}}^{|\beta|=[\frac{m}{2}]}\|\MZ^\beta(\frac{v}{p})\|{_\i^2}\|\MZ^\alpha\nabla_yp\|^2+\sum_{\substack{|\beta|=1+[\frac{m}{2}]\\|\beta|+|\nu|\leq m-1}}^{|\beta|=m-1}\|\MZ^\nu\nabla_yp\|{_\i^2}\|\MZ^\beta(\frac{v}{p})\|^2\nonumber\\
			&\leq \Lambda(\frac{1}{c_0},\mathcal{Q}_m(t)+\mathcal{Q}(t)).
		\end{align}
		Additional care is needed to estimate $\|\frac{V_z\partial_zp}{p}(t)\|_{\mathcal{H}^{m-1}}$ since it involves $\partial_zp$. Rewrite this term as
		\begin{equation}\label{9.14-1}
			 \MZ^\alpha(\frac{V_z\partial_zp}{p})=\frac{V_z}{p}\MZ^\alpha\partial_zp+\partial_zp\MZ^\alpha(\frac{V_z}{p})+\sum_{\substack{|\beta|=1\\\beta+\nu=\alpha}}^{m-2}c_k\MZ^\beta(\frac{V_z}{p})\MZ^\nu\partial_zp.
		\end{equation}
	Then, it follows from \eqref{4.26}  that
		\begin{equation}\label{9.14}
			 \|\frac{V_z}{p}\MZ^\alpha\partial_zp\|^2\leq\|\frac{(1-z)V_z}{zp}\|{_\i^2}\|p\|{_{\mathcal{H}^m}^2}\leq\Lambda(\frac{1}{c_0},\mathcal{Q}(t)+\mathcal{Q}_m(t)),
		\end{equation}
		\begin{equation}\label{9.15}
			 \|\partial_zp\MZ^\alpha(\frac{V_z}{p})(t)\|^2\leq\Lambda(\frac{1}{c_0},\mathcal{Q}(t)+\mathcal{Q}_m(t)),
		\end{equation}
		and
		\begin{align}\label{9.16}
			 &\sum_{\substack{|\beta|=1\\\beta+\nu=\alpha}}^{|\beta|=m-2}\|\MZ^\beta(\frac{V_z}{p})\MZ^\nu\partial_zp\|^2\nonumber\\
			 &\leq\|\MZ(\frac{V_z}{p})\|{_\i^2}\|\partial_zp\|_{\mathcal{H}^{m-2}}^2+\|\MZ\partial_zp\|{_{L^\i}^2}\|\frac{V_z}{p}\|_{\mathcal{H}^{m-2}}^2+\sum_{\substack{|\beta|=2\\\beta+\nu=\alpha}}^{|\beta|=m-3}\|\MZ^\beta(\frac{V_z}{p})\|{_\i^2}\|\MZ^\nu\partial_zp\|^2\nonumber\\
			 &\leq\Lambda(\frac{1}{c_0},\mathcal{Q}_m(t)+\mathcal{Q}(t))+\|\frac{V_z}{p}\|{_{\mathcal{H}^{m-3,\i}}^2}\|\partial_zp\|{_{\mathcal{H}^{m-3}}^2}\nonumber\\
			 &\leq\Lambda(\frac{1}{c_0},\mathcal{Q}_m(t)+\mathcal{Q}(t)).
		\end{align}
		Collecting the inequalities \eqref{9.5}, \eqref{9.10}, \eqref{9.12}-\eqref{9.16} shows \eqref{9.4-0}.  Therefore,  the proof of this lemma is completed. $\hfill\Box$

	\subsection{Normal derivative estimates: Part I}
	
	\setcounter{equation}{0}
	
	In this subsection, we will focus on the estimates of normal derivative of $v$. To close the a priori estimates, one needs to bound
	$\int_{0}^{t}\|\partial_z v\|{_{\mathcal{H}^{m-1}}^2}d\tau$, $\v\int_{0}^{t}\|\partial_z^2 v\|{_{\mathcal{H}^{m-2}}^2}d\tau$ and $\v^2\int_{0}^{t}\|\partial_z^2 v\|{_{\mathcal{H}^{m-1}}^2}d\tau$ in this and next subsection. Since
	  \begin{equation}\label{7.0}
		\partial_zv=(\partial_zv\cdot\mathbf{n})\mathbf{n}+\Pi(\partial_zv),
	\end{equation}
it suffices  to estimate $\partial_zv\cdot\mathbf{n}$ and $\Pi(\partial_zv)$.  First, one can bound the normal component of $\partial_zv$ in terms of $\text{div}^\varphi v$ as follows:
	\begin{lemma}\label{lem7.1}
		For any $k\in \mathbb{N}$, it holds that
		\begin{equation}\label{4.4.2}
			 \int_{0}^{t}\|\partial_z v\cdot\mathbf{n}\|_{\mathcal{H}^k}^2\leq\Lambda_0\int_{0}^{t}\|\text{div}^\varphi v\|_{\mathcal{H}^k}^2+\|V^{k+1}\|^2d\tau+\Lambda_{\i}(t)\int_{0}^{t}|\MZ^kh|_{\frac{1}{2}}^2d\tau
		\end{equation}
	\end{lemma}
	\noindent\textbf{Proof}. It follows from the definition of $\partial{_i^\varphi}$ that
	\begin{eqnarray}\label{7.2}
		\mbox{div}^\varphi v=\partial{_1^\varphi}v_1+\partial{_2^\varphi}v_2+\partial{_3^\varphi}v_3
		 =\frac{1}{\partial_z\varphi}\partial_zv\cdot\mathbf{N}+\partial_1v_1+\partial_2v_2.\nonumber
	\end{eqnarray}
	which yields  that
	\begin{equation}\label{7.3}
		\partial_zv\cdot\mathbf{n}=\frac{\partial_z\varphi}{\sqrt{1+|\nabla_y\varphi|^2}}(\mbox{div}^\varphi v-\partial_1v_1-\partial_2v_2).
	\end{equation}
	Hence, \eqref{4.4.2} follows from \eqref{2.6}, Lemma \ref{l2.4} and \eqref{5.2}. $\hfill\Box$

\
	
Next, we will estimate the tangential components of $\partial_zv$. Note that
	\begin{equation}\label{7.6}
		\nabla^\varphi v\mathbf{N}=-\partial_1\varphi\partial{_1^\varphi}v-\partial_2\varphi\partial{_2^\varphi}v+\partial{_z^\varphi}v=\frac{1+|\nabla_y\varphi|^2}{\partial_z\varphi}\partial_zv-\partial_1\varphi\partial_1v-\partial_2\varphi\partial_2v,
	\end{equation}
and 
    \begin{equation}\label{4.4.5}
    (\nabla^\varphi v)^t\mathbf{N}=(\partial_1^\varphi v\cdot\mathbf{N},\partial_2^\varphi v\cdot\mathbf{N},\partial_z^\varphi v\cdot\mathbf{N})^t=(\partial_1v\cdot\mathbf{N},\partial_2v\cdot\mathbf{N},0)^t+\frac{1}{\partial_z\varphi}(\partial_zv\cdot\mathbf{N})\mathbf{N}.
    \end{equation}
One can obtain immediately that
	\begin{align}\label{7.7}
		 \Pi(\partial_zv)&=\frac{\partial_z\varphi}{1+|\nabla_y\varphi|^2}(\Pi(\nabla^\varphi v\mathbf{N})+\Pi(\partial_1\varphi\partial_1v+\partial_2\varphi\partial_2v))\nonumber\\
		 &=\frac{\partial_z\varphi}{1+|\nabla_y\varphi|^2}(2\Pi(S^\varphi v\mathbf{N})-\Pi((\nabla^\varphi v)^t\mathbf{N})+\Pi(\partial_1\varphi\partial_1v+\partial_2\varphi\partial_2v))\nonumber\\
		 &=\frac{\partial_z\varphi}{1+|\nabla_y\varphi|^2}(2\Pi(S^\varphi v\mathbf{N})-\Pi\{(\partial_1v\cdot\mathbf{N},\partial_2v\cdot\mathbf{N},0)^t\}+\Pi(\partial_1\varphi\partial_1v+\partial_2\varphi\partial_2v)).
	\end{align}
As a consequence of $\partial_zv$ given in \eqref{7.3} and \eqref{7.7},  one have the following estimates for $\partial_z v$:
\begin{lemma}\label{lem7.2}
		For any $k\in \mathbb{N}$, it holds that
		\begin{equation}
			 \int_{0}^{t}\|\partial_z v\|_{\mathcal{H}^k}^2\leq\Lambda_0\int_{0}^{t}\|(\text{div}^\varphi v,S_\mathbf{n})\|_{\mathcal{H}^k}^2+\|V^{k+1}\|^2d\tau+\Lambda_{\i,t}\int_{0}^{t}|\MZ^kh|_{\frac{1}{2}}^2d\tau
		\end{equation}
	where $S_\mathbf{n}:=\Pi(S^\varphi v\mathbf{N})$.
	\end{lemma}
	Furthermore, it follows from the expression of $\partial_{zz}v$ and \eqref{6.37-1} that
	\begin{lemma}\label{lem7.3}
		For any $m\geq 1$, it holds that
		\begin{align}\label{7.3-10}
		&\v\int_0^t\|\nabla^2v(\tau)\|{_{\mathcal{H}^{m-2}}^2}d\tau\leq \Lambda_0\v\int_0^t \|\nabla^\varphi\MZ^{m-2} S_\mathbf{n}(\tau)\|^2d\tau+\Lambda_{\infty}(t)\int_0^tY_m(\tau)d\tau,
		\end{align}
		and 
		\begin{align}
	    &\v^2\int_0^t\|\nabla^2v(\tau)\|{_{\mathcal{H}^{m-1}}^2}d\tau\leq \Lambda_0\v^2\int_0^t \|(\nabla^\varphi\MZ^{m-1}\text{div}^\varphi v,\nabla^\varphi\MZ^{m-1} S_\mathbf{n},\nabla V^{m})(\tau)\|^2d\tau\nonumber\\
		&\qquad\qquad\qquad\qquad\qquad\qquad+\Lambda_{\infty}(t)\int_0^tY_m(\tau)
			+\v|\MZ^{m-1}\nabla_yh(\tau)|^2_{\frac{1}{2}}d\tau.\label{7.3-11}
		\end{align}
	\end{lemma}
	
	\

The key is to estimate $S_\mathbf{n}$. Now, we follow the argument in \cite{Masmoudi-R-1} for the convention-diffusion equation solved by $S_\mathbf{n}$ to derive the estimates of $\v\int_{0}^{t}\|\nabla^\varphi\MZ^{m-2} S_\mathbf{n}(\tau)\|^2d\tau$ and $\v^2\int_{0}^{t}\|\nabla^\varphi\MZ^{m-1} S_\mathbf{n}(\tau)\|^2d\tau$. The only difference is that $\text{div}^\varphi v$ will be involved but has been estimated in the last subsection. It should be noted that $S_\mathbf{n}$ satisfies the homogeneous Dirichlet boundary condition
	\begin{equation}\label{7.8}
		S_\mathbf{n}=0,\quad\quad on\quad z=0.
	\end{equation}
which follows from \eqref{1.23-1} and will be important for the following analysis. 

Now, we start with the estimate of $S_\mathbf{n}$ in $\mathcal{H}^{m-2}$ which bound $\v\int_{0}^{t}\|\nabla^\varphi\MZ^{m-2} S_\mathbf{n}(\tau)\|^2d\tau$. 
	\begin{lemma}\label{lem7.4}
		For any $t\in[0,T^\v]$ and $m\geq 2$, it holds that
		\begin{align}\label{7.4-1}
			 &\|S_\mathbf{n}(t)\|{_{\mathcal{H}^{m-2}}^2}+\v\int_{0}^{t}\|\nabla^\varphi\MZ^{m-2} S_\mathbf{n}(\tau)\|^2d\tau\nonumber\\
			 &\leq \Lambda_0\|\MZ^\a S_\mathbf{n}(0)\|^2+\d\Lambda_0\v^2\int_{0}^{t}\|(\nabla^\varphi\MZ^{m-1}\mbox{div}^\varphi v,\nabla^\varphi\MZ^{m-1} S_\mathbf{n},\nabla^\varphi V^m)(\tau)\|^2d\tau\nonumber\\
			 &\quad+C_\d\Lambda(\frac{1}{c_0},\mathcal{Q}(t)+\v\|\nabla^\varphi\mbox{div}^\varphi v\|{_{\mathcal{W}^{1,\i}_t}^2})\int_{0}^{t}Y_m(\tau)+\v|\MZ^m h|{_\frac{1}{2}^2}d\tau.
		\end{align}
	\end{lemma}
\noindent\textbf{Proof}. We first derive the equations satisfied by $S_\mathbf{n}$. Applying $\nabla^\varphi$ to $\eqref{1.22}_2$ gives that
	\begin{equation*}
		\varrho\partial{_t^\varphi}\nabla^\varphi v+\varrho (v\cdot\nabla^\varphi)\cdot \nabla^\varphi v-\m\v \Delta^\varphi(\nabla^\varphi v)=(\m+\l)\v(\nabla^\varphi)^2\mbox{div}^\varphi v-(\nabla^\varphi)^2p-\nabla^\varphi\varrho\otimes\partial_t^\varphi v-\nabla^\varphi v\nabla^\varphi(\varrho v),
	\end{equation*}
	where $(\nabla^\varphi)^2f$ denotes the Hessian matrix of $f$. Thus, 
	\begin{align*}
		&\quad\varrho\partial{_t^\varphi}(S^\varphi v)+\varrho(v\cdot\nabla^\varphi)(S^\varphi v)-\m\v\Delta^\varphi(S^\varphi v)=(\m+\l)\v(\nabla^\varphi)^2\mbox{div}^\varphi v-(\nabla^\varphi)^2p\nonumber\\
		&\quad-\frac{1}{2}(\nabla^\varphi\varrho\otimes\partial_t^\varphi v+\partial_t^\varphi v\otimes\nabla^\varphi\varrho)-\frac{1}{2}(\nabla^\varphi v\nabla^\varphi(\varrho v)+(\nabla^\varphi v\nabla^\varphi(\varrho v))^t).\nonumber
	\end{align*}
Therefore, $S_\mathbf{n}$ satisfies
	\begin{equation}\label{7.10}
		 \varrho\partial{_t^\varphi}S_\mathbf{n}+\varrho(v\cdot\nabla^\varphi)S_\mathbf{n}-\m\v\Delta^\varphi S_\mathbf{n}=F,
	\end{equation}
where $F$ is given by
	\begin{equation}
		F=F_1+F_2+F_3,
	\end{equation}
with
	\begin{align}
		 F_1&=\varrho(\partial{_t^\varphi}\Pi+v\cdot\nabla^\varphi\Pi)(S^\varphi v\mathbf{N})+\varrho\Pi(S^\varphi v(\partial{_t^\varphi\mathbf{N}+v\cdot\nabla^\varphi\mathbf{N}}))\\
		 &\quad-\frac{1}{2}\Pi\Big\{(\nabla^\varphi\varrho\otimes\partial_t^\varphi v+\partial_t^\varphi v\otimes\nabla^\varphi\varrho)\mathbf{N}\Big\}-\frac{1}{2}\Pi\Big\{(\nabla^\varphi v\nabla^\varphi(\varrho v)+(\nabla^\varphi v\nabla^\varphi(\varrho v))^t)\mathbf{N}\Big\},\nonumber\\\label{F2}
		F_2&=-\Pi\Big\{(\nabla_y,0)^t(\nabla^\varphi p\cdot \mathbf{N})+(\nabla^\varphi\mathbf{N})^t\nabla^\varphi p\Big\},\\\label{F3}
		 F_3&=(\m+\l)\v\Pi\Big\{(\nabla_y,0)^t(\nabla^\varphi\mbox{div}^\varphi v\cdot\mathbf{N})-(\nabla^\varphi N)^t\nabla^\varphi\mbox{div}^\varphi v\Big\}\nonumber\\
		&\quad-\m\v(\Delta^\varphi\Pi)(S^\varphi v\mathbf{N})-2\m\v\partial{_i^\varphi}\Pi\partial{_i^\varphi}(S^\varphi v\mathbf{N})-\m\v\Pi\Big\{S^\varphi v\Delta^\varphi\mathbf{N}+2\partial{_i^\varphi}S^\varphi v\partial{_i^\varphi}\mathbf{N}\Big\}.
	\end{align}
	where the summation convention has been used in the above expressions. It follows from  Lemma \ref{l2.2}, Lemma \ref{lem2.4} and \eqref{6.37-2}  that
	\begin{align}
		 &\int_{0}^{t}\|F_1\|{_{\mathcal{H}^{m-2}}^2}d\tau\leq\Lambda_\i(t)\int_{0}^{t}\|(\nabla v,\nabla p)(\tau)\|{_{\mathcal{H}^{m-2}}^2}+\|(p,v)\|{_{\mathcal{H}^{m-1}}^2}+|\MZ^{m-1}h|{_\frac{1}{2}^2}d\tau,\label{7.15}\\
		 &\int_{0}^{t}\|F_2\|{_{\mathcal{H}^{m-2}}^2}d\tau\leq\Lambda_\i(t)\int_{0}^{t}\|\nabla p\|{_{\mathcal{H}^{m-1}}^2}+|\MZ^{m-1}h|{_{\frac{1}{2}}^2}d\tau,\label{7.16}\\
		 &\int_{0}^{t}\|F_3\|{_{\mathcal{H}^{m-2}}^2}d\tau\leq\Lambda_0\int_{0}^{t}\|v\|_{\mathcal{H}^{m-2}}^2+\v\|\nabla^2 v\|_{\mathcal{H}^{m-2}}^2+\v^2\|\nabla^2v\|_{\mathcal{H}^{m-1}}^2d\tau\nonumber\\
		 &~~~~~~~~~~~~~~~~~~~~~~~~~~+\Lambda(\frac{1}{c_0},\mathcal{Q}(t)+\v\|\nabla\text{div}^\varphi v\|_{\mathcal{W}^{1,\i}_t}^2)\int_{0}^{t}|\MZ^{m-1}h|_\frac{1}{2}^2+\v|\MZ^mh|_\frac{1}{2}^2d\tau.\label{7.17}
	\end{align}
The lemma is proved by induction.  For the case $m=2$, multiplying \eqref{7.10} by $S_\mathbf{n}$, and then using Lemma \ref{l3.1} and the boundary condition \eqref{7.7}, one gets that
	\begin{equation}\label{7.20}
		 \int_{\mathcal{S}}\frac{1}{2}\varrho(t)|S_\mathbf{n}(t)|^2d\mathcal{V}_t+\m\v\int_{0}^{t}\int_{\mathcal{S}}|\nabla^\varphi S_\mathbf{n}|^2d\mathcal{V}_\tau d\tau=\int_{\mathcal{S}}\frac{1}{2}\varrho_0|S_\mathbf{n}(0)|^2d\mathcal{V}_0+\int_{0}^{t}\int_{\mathcal{S}}F\cdot S_\mathbf{n}d\mathcal{V}_\tau d\tau
	\end{equation}
It follows from \eqref{7.15}-\eqref{7.17} with $m=2$ that
	\begin{align}\label{7.21}
		&\int_{0}^{t}\int_{\mathcal{S}}F\cdot S_\mathbf{n}d\mathcal{V}_\tau d\tau\leq\Big(\int_{0}^{t}\|(F_1,F_2,F_3)\|^2d\tau\Big)^\frac{1}{2}\Big(\int_{0}^{t}\|S_\mathbf{n}\|^2d\tau\Big)^\frac{1}{2}\nonumber\\
		&\leq\int_{0}^{t}\|\nabla p\|{_{\mathcal{H}^{1}}^2}+\v^2\|\nabla\mbox{div}^\varphi v\|{_{\mathcal{H}^{1}}^2}+\v^2\|\nabla^2 v\|_{\mathcal{H}^1}^2d\tau\nonumber\\
		&\quad+\Lambda(\frac{1}{c_0},\mathcal{Q}(t)+\|\nabla^\varphi\mbox{div}^\varphi v\|{_{\i,t}^2})\int_{0}^{t}\|v\|{_{\mathcal{H}^2}^2}+\|(p,\nabla v)\|{_{\mathcal{H}^{1}}^2}+|\MZ h|{_\frac{1}{2}^2}+\v|\MZ^2 h|{_\frac{1}{2}^2}\nonumber\\
		&\leq \Lambda(\frac{1}{c_0},\mathcal{Q}(t)+\|\nabla^\varphi\mbox{div}^\varphi v\|{_{\i,t}^2})
		\int_0^t\Lambda(Y_m(\tau))+\v|\MZ^2h|^2_{\f12}d\tau,
	\end{align}
	where \eqref{6.37-2} has been used in the last inequality. Substituting \eqref{7.21} into \eqref{7.20}, one gets from \eqref{6.37-1} and \eqref{7.3-11} that
	\begin{align}
		 &\int_{\mathcal{S}}\frac{1}{2}\varrho(t)|S_\mathbf{n}(t)|^2d\mathcal{V}_t+\m\v\int_{0}^{t}\int|\nabla^\varphi S_\mathbf{n}|^2d\mathcal{V}_\tau d\tau\nonumber\\
		 &\leq\int_{\mathcal{S}}\frac{1}{2}\varrho_0|S_\mathbf{n}(0)|^2d\mathcal{V}_0+\Lambda_0\v^2\int_{0}^{t}\|\nabla^\varphi\MZ\text{div}^\varphi v,\nabla^\varphi\MZ S_\mathbf{n},\nabla V^1\|d\tau\nonumber\\
		 &\qquad\qquad+\Lambda(\frac{1}{c_0},\mathcal{Q}(t)+\|\nabla^\varphi\mbox{div}^\varphi v\|{_{\i,t}^2})\int_{0}^{t}\Lambda(Y_m(\tau))+\v|\MZ^2 h|{_\frac{1}{2}^2}d\tau.
	\end{align}
	Assume that \eqref{7.4-1} has been proved for  $k\leq m-3$, we shall prove   it   for $k=m-2$. First, applying $\MZ^\a$, with $|\a|=m-2$, to \eqref{7.10} yields that
	\begin{equation}\label{7.23}
		\varrho\partial{_t^\varphi}\MZ^\a S_\mathbf{n}+\varrho v\cdot\nabla^\varphi\MZ^\a S_\mathbf{n}-\m\v\Delta^\varphi\MZ^\a S_\mathbf{n}=\MZ^\a F+C_S,
	\end{equation}
	with $C_S=C{_S^1}+C{_S^2}$ given by
	\begin{equation}
	C{_S^1}\doteq-[\MZ^\a,\varrho]\partial_tS_\mathbf{n}-[\MZ^\a,\varrho v_y]\nabla_yS_\mathbf{n}-[\MZ^\a,\varrho V_z]\partial_zS_\mathbf{n}-\varrho V_z[\MZ^\a,\partial_z]S_\mathbf{n}\doteq\sum_{i=1}^4C_{S}^{1,i},
	\end{equation}
	and
	\begin{equation}
	C{_S^2}\doteq\m\v[\MZ^\a,\Delta^\varphi]S_\mathbf{n}.
	\end{equation}
	Multiplying \eqref{7.23} by $\MZ^\a S_\mathbf{n}$ with $|\a|=m-2$ and integrating give that
	\begin{eqnarray}\label{7.25}
		&&\int_{\mathcal{S}}\frac{1}{2}\varrho(t)|\MZ^\a S_\mathbf{n}(t)|^2d\mathcal{V}_t+\m\v\int_{0}^{t}\int_\mathcal{S}|\nabla^\varphi\MZ^\a S_\mathbf{n}|^2d\mathcal{V}_\tau d\tau\nonumber\\
		&&=\int_{\mathcal{S}}\frac{1}{2}\varrho_0|\MZ^\a S_\mathbf{n}(0)|^2d\mathcal{V}_0+\int_{0}^{t}\int_\mathcal{S}|\MZ^\a F\cdot\MZ^\a S_\mathbf{n}|d\mathcal{V}_\tau d\tau+\int_{0}^{t}\int_\mathcal{S}C_S\cdot\MZ^\a S_\mathbf{n}d\mathcal{V}_\tau d\tau.
	\end{eqnarray}
	As a consequence of \eqref{7.15}-\eqref{7.17}, \eqref{6.46-1}, \eqref{7.3-10} and \eqref{7.3-11}, we obtain that
	\begin{align}\label{7.26}
		&\int_{0}^{t}\int_\mathcal{S}|\MZ^\a F\cdot\MZ^\a S_\mathbf{n}|d\mathcal{V}_\tau d\tau\leq\Big(\int_{0}^{t}\|(F_1,F_2,F_3)\|^2_{\mathcal{H}^{m-2}}d\tau\Big)^\frac{1}{2}\Big(\int_{0}^{t}\|S_\mathbf{n}\|{_{\mathcal{H}^{m-2}}^2}d\tau\Big)^\frac{1}{2}\nonumber\\
		&\leq\d\int_{0}^{t}\|\nabla p\|{_{\mathcal{H}^{m-1}}^2}+\v\|\nabla^2 v\|_{\mathcal{H}^{m-2}}^2+\v^2\|\nabla^2v\|_{\mathcal{H}^{m-1}}^2d\tau\nonumber\\
		&~~~~~+\Lambda(\frac{1}{c_0},\mathcal{Q}(t)+\v\|\nabla\text{div}^\varphi v\|_{\mathcal{W}^{1,\i}_t}^2)\int_{0}^{t}Y^m(\tau)+\v|\MZ^mh|_\frac{1}{2}^2d\tau.\nonumber\\
		&\leq \frac{\mu\v}{4}\int_{0}^{t}\|\nabla^\varphi\MZ^{m-2}S_\mathbf{n}\|^2d\tau+\d\Lambda_0\v^2\int_0^t \|(\nabla^\varphi\MZ^{m-1}\text{div}^\varphi v,\nabla^\varphi\MZ^{m-1} S_\mathbf{n},\nabla^\varphi V^{m})(\tau)\|^2d\tau\nonumber\\
		&~~~~~~~~~~~+C_\d\Lambda(\frac{1}{c_0},\mathcal{Q}(t)+\v\|\nabla\text{div}^\varphi v\|_{\mathcal{W}^{1,\i}_t}^2)\int_0^tY_m(\tau)
		+\v|\MZ^{m}h(t)|^2_{\frac{1}{2}}d\tau.
	\end{align}
	It remains to estimate the terms involving $C_{S}^1,~C_{S}^2$. Similar to \cite{Masmoudi-R-1}, integrating by parts and using the Hardy inequality, one gets that
	\begin{equation}\label{7.32}
		\int_{0}^{t}\int_{\mathcal{S}}|C{_S^1}\MZ^\a S_\mathbf{n}|d\mathcal{V}_td\tau\leq\int_{0}^{t}\Lambda(Y_m(\tau))d\tau.
	\end{equation}
	For the term involving $C{_S^2}$, we   notice that $\Delta^\varphi$ can be rewritten as
	\begin{equation}
		\Delta^\varphi f=\frac{1}{\partial_z\varphi}\nabla\cdot(E\nabla f),
	\end{equation}
	with the matrix $E$ defined by
	\begin{equation}
		E=\left(
		\begin{array}{ccc}
			\partial_z\varphi & 0 & -\partial_1\varphi\\
			0 & \partial_z\varphi & -\partial_2\varphi\\
			-\partial_1\varphi & -\partial_2\varphi & \frac{1+|\nabla_y\varphi|^2}{\partial_z\varphi}\\
		\end{array}
		\right).
	\end{equation}
	This   yields immediately that
	\begin{equation*}
		C{_S^2}=C{_S^{2,1}}+C{_S^{2,2}}+C{_S^{2,3}},
	\end{equation*}
	with
	\begin{equation*}
		 C{_S^{2,1}}=\m\v[\MZ^\a,\frac{1}{\partial_z\varphi}]\nabla\cdot(E\nabla S_\mathbf{n}),~~ C{_S^{2,2}}=\m\v\frac{1}{\partial_z\varphi}[\MZ^\a,\nabla\cdot](E\nabla S_\mathbf{n}),~~ C{_S^{2,3}}=\m\v\frac{1}{\partial_z\varphi}\nabla\cdot[\MZ^\a,E\nabla]S_\mathbf{n}.
	\end{equation*}
By similar arguments as \cite{Masmoudi-R-1}, one obtains  that
	\begin{eqnarray}\label{7.40}
		&&|\int_{0}^{t}\int_\mathcal{S}C{_S^2}\MZ^\a S_\mathbf{n}d\mathcal{V}_\tau d\tau|\leq\f{\m\v}{16}\int_{0}^{t}\int|\nabla\MZ^{m-2} S_\mathbf{n}|^2d\mathcal{V}_\tau
		+\Lambda_0\v\int_0^t\int|\nabla\MZ^{m-3} S_\mathbf{n}|^2d\mathcal{V}_\tau d\tau\nonumber\\
		 &&~~~~~~~~~~~~~~~~~~~~~~~~~~~~~~~~~~~~+\Lambda_\i(t)\int_{0}^{t}Y_m(\tau)d\tau.
	\end{eqnarray}
	Plugging \eqref{7.26}, \eqref{7.32} and \eqref{7.40} into \eqref{7.25}, we obtain that
	\begin{align*}
		&\int_{\mathcal{S}}\frac{1}{2}\varrho(t)|\MZ^\a S_\mathbf{n}(t)|^2d\mathcal{V}_t+\m\v\int_{0}^{t}\int_\mathcal{S}|\nabla^\varphi\MZ^\a S_\mathbf{n}|^2d\mathcal{V}_\tau d\tau\nonumber\\
		&\leq\int_{\mathcal{S}}\frac{1}{2}\varrho_0|\MZ^\a S_\mathbf{n}(0)|^2d\mathcal{V}_0+\Lambda_0\v\int_0^t\int|\nabla\MZ^{m-3} S_\mathbf{n}|^2d\mathcal{V}_\tau d\tau +\Lambda_\i(t)\int_{0}^{t}\Lambda(Y_m(\tau))+\v|\MZ^mh|^2_{\f12}d\tau.
	\end{align*}
	Thus, the lemma is proved by \eqref{2.23} and the induction assumption to control $\Lambda_0\v\int_{0}^{t}\|\nabla^\varphi\MZ^{m-3}S_\mathbf{n}\|^2d\tau$. $\hfill\Box$

\

 The following lemma give the estimate of $\v^2\int_{0}^{t}\|\nabla^\varphi\MZ^{m-1} S_\mathbf{n}\|^2d\tau$.
	\begin{lemma}\label{lem7.5}
		For any $t\in[0,T^\v]$ and $m\geq1$, it holds that
		\begin{align}\label{7.47-1}
			 &\v\|S_\mathbf{n}(t)\|^2_{\mathcal{H}^{m-1}}+\v^2\int_{0}^{t}\|\nabla^\varphi\MZ^{m-1} S_\mathbf{n}\|^2d\tau\nonumber\\
			 &\leq\f{1}{\d_1}\Lambda(\f1{c_0}, Y_m(0))+\Lambda_0\d_1\v^2\int_{0}^{t}\int|\nabla^\varphi\MZ^{m-1}\mbox{div}^{\varphi}v|^2d\mathcal{V}_\tau d\tau\nonumber\\
			 &~~~~~~~~~~~+C_{\d_1}\Lambda_\i(t)\int_{0}^{t}\Lambda(Y_m(\tau))+\v|\MZ^mh|{_\frac{1}{2}^2}+\|\nabla v\|^2_{\mathcal{H}^{m-1}}d\tau,
		\end{align}
		where $\d_1>0$ and $\d>0$   will be chosen later.
	\end{lemma}
	\noindent\textbf{Proof}. Multiplying \eqref{7.23} by $\v\MZ^\a S_\mathbf{n}$ with $|\a|= m-1$ and integrating give that
		\begin{align}\label{7.48-1}
			&\int_{\mathcal{S}}\frac{\v}{2}\varrho(t)|\MZ^\a S_\mathbf{n}(t)|^2d\mathcal{V}_t+\m\v^2\int_{0}^{t}\int_\mathcal{S}|\nabla^\varphi\MZ^\a S_\mathbf{n}|^2d\mathcal{V}_\tau d\tau\nonumber\\
			&=\int_{\mathcal{S}}\frac{\v}{2}\varrho_0|\MZ^\a S_\mathbf{n}(0)|^2d\mathcal{V}_0+\v\Big|\int_{0}^{t}\int_\mathcal{S}\MZ^\a F\cdot\MZ^\a S_\mathbf{n}d\mathcal{V}_\tau d\tau\Big|+\v\Big|\int_{0}^{t}\int_\mathcal{S} C_S\cdot\MZ^\a S_\mathbf{n} d\mathcal{V}_\tau d\tau\Big|.
		\end{align}
		Similar to \eqref{7.15}, it is easy to obtain that
		\begin{equation*}
			 \int_{0}^{t}\|F_1\|{_{\mathcal{H}^{m-1}}^2}d\tau\leq\Lambda_\i(t)\int_{0}^{t}\|(\nabla v,\nabla p)\|{_{\mathcal{H}^{m-1}}^2}+\|(p,v)\|{_{\mathcal{H}^{m}}^2}+|\MZ^{m}h|{_\frac{1}{2}^2}d\tau,
		\end{equation*}
		which implies that
		\begin{equation}\label{7.50-1}
			\v\Big|\int_{0}^{t}\int_\mathcal{S}\MZ^\a F_1\cdot\MZ^\a S_\mathbf{n}d\mathcal{V}_\tau d\tau\Big|\leq\Lambda_\i(t)\int_{0}^{t}Y_m(\tau)+\v|\MZ^{m}h|^2_\frac{1}{2}d\tau.
		\end{equation}
	Using integration by parts, one can get that
		\begin{equation}\label{7.51-2}
			\v\Big|\int_{0}^{t}\int_\mathcal{S}\MZ^\a F_2\cdot\MZ^\a S_\mathbf{n}d\mathcal{V}_\tau d\tau\Big|\leq\v\int_{0}^{t}\|\nabla V^m\|^2d\tau+\Lambda_\i(t)\int_{0}^{t}Y_m(\tau)+\v|\MZ^{m}h|^2_\frac{1}{2}d\tau,
		\end{equation}
		\begin{align}\label{7.56}
			&\v\Big|\int_{0}^{t}\int_\mathcal{S}\MZ^\a F_3\cdot\MZ^\a S_\mathbf{n}d\mathcal{V}_\tau d\tau\Big|\nonumber\\
			&\leq \f{\mu\v^2}{16}  \int_{0}^{t}\int |\nabla^\varphi\MZ^\a S_\mathbf{n}|^2d\mathcal{V}_\tau d\tau+\Lambda_0\d_1\v^2\int_{0}^{t}\int|\nabla^\varphi\MZ^{m-1}\mbox{div}^{\varphi}v|^2d\mathcal{V}_\tau d\tau\nonumber\\
			&~~~~~~~+\f{\Lambda_0\v^2}{\d_1}\int_0^t\|\nabla V^m\|^2d\tau+C_{\d_1}\Lambda_\i(t)\int_{0}^{t}\Lambda(Y_m(\tau))+\v|\MZ^mh|{_\frac{1}{2}^2}d\tau
		\end{align}	
	and
		\begin{align}\label{7.57}
			&\v\Big|\int_{0}^{t}\int_{\mathcal{S}}C{_S^1}\MZ^\a S_\mathbf{n}d\mathcal{V}_\tau d\tau\Big|+\v\Big|\int_{0}^{t}\int_\mathcal{S}C{_S^2}\MZ^\a S_\mathbf{n}d\mathcal{V}_\tau d\tau\Big|\nonumber\\
			&\leq \f{\mu\v^2}{16}  \int_{0}^{t}\int |\nabla^\varphi\MZ^\a S_\mathbf{n}|^2d\mathcal{V}_\tau d\tau+ \Lambda_\i(t)\int_{0}^{t}\Lambda(Y_m(\tau))+\v|\MZ^mh|{_\frac{1}{2}^2}d\tau.
		\end{align}
	Substituting \eqref{7.50-1}-\eqref{7.57} into \eqref{7.48-1}, we prove \eqref{7.48-1}. Therefore, the proof of this lemma is completed.   $\hfill\Box$
	
\
	
Then, suitable choice of $\d$ and $\d_1$ yield the following proposition
	\begin{proposition}\label{prop7.6}
		For any $t\in [0,T^\v]$, $m\geq 5$, it holds that
		\begin{align}\label{7.60-1}
			 &\|(V^m,Q^m)(t)\|^2+|(h,\sqrt\sigma\nabla_yh)(t)|^2_{\mathcal{H}^m}+\|(\nabla p,\nabla v)(t)\|{_{\mathcal{H}^{m-2}}^2}+\v\|(\text{div}^\varphi v,\nabla^\varphi p,S_\mathbf{n})(t)\|{_{\mathcal{H}^{m-1}}^2}\nonumber\\
			&\qquad+\v|\MZ^mh|_\frac{1}{2}^2+\int_0^t\|\nabla p\|_{\mathcal{H}^{m-1}}^2d\tau+\v\int_{0}^{t}\|\nabla v\|_{\mathcal{H}^m}^2d\tau+\int_{0}^{t}\v\|\nabla^2v\|^2_{\mathcal{H}^{m-2}}+\v^2\|\nabla^2v\|^2_{\mathcal{H}^{m-1}}d\tau\nonumber\\
			&\leq  \Lambda(\f1{c_0},Y_m(0))+\Lambda_\i(t)(1+\|\nabla^\varphi\mbox{div}^\varphi v\|{_{\i,t}^2})\int_0^t\Lambda(Y_m(\tau))+\v|\MZ^m h|^2_{\f12}+\|\nabla v\|^2_{\mathcal{H}^{m-1}}d\tau.
		\end{align}
	\end{proposition}
	\noindent\textbf{Proof}. Combining Propositions \ref{p5.1}, \ref{p6.6},  Lemmas \ref{lem7.3}, \ref{lem7.4} and \ref{lem7.5},  one can obtain that	 
	\begin{align}\label{7.60-2}
	 &\|(V^m,Q^m)(t)\|^2+|(h,\sqrt\sigma\nabla_yh)(t)|^2_{\mathcal{H}^m}+\|\nabla v(t)\|{_{\mathcal{H}^{m-2}}^2}+\v\|(\text{div}^\varphi v,\nabla^\varphi p,S_\mathbf{n})(t)\|{_{\mathcal{H}^{m-1}}^2}\nonumber\\
	&\qquad+\int_{0}^{t}\v\|\nabla V^m(\tau)\|^2 +\v\|\nabla^2v(\tau)\|^2_{\mathcal{H}^{m-2}}+\v^2\|\nabla^2v(\tau)\|^2_{\mathcal{H}^{m-1}}d\tau\nonumber\\
	&\leq \f{1}{\d_1}\Lambda(\f1{c_0}, Y_m(0))+\Lambda_0(\d_1+\Lambda_0\f{\d}{\d_1})\v^2\int_{0}^{t}\|\nabla^2v(\tau)\|^2_{\mathcal{H}^{m-1}}d\tau+(\d_1+\Lambda_0\f{\d}{\d_1})\v\int_{0}^{t}\|\nabla V^m\|d\tau\nonumber\\
	&\qquad+C_\d\Lambda_\i(t)(1+\|\nabla^\varphi\mbox{div}^\varphi v\|{_{L{^\i_{t,x}}}^2})\int_{0}^{t}\Lambda(Y_m(\tau))+\v|\MZ^m h|{_\frac{1}{2}^2}+\|\nabla v\|^2_{\mathcal{H}^{m-1}}d\tau.
	\end{align}	
Choosing  $\d_1>0$ to be  small such that $\d_1\leq \f1{4\Lambda_0}$, and  $\d=\frac{\d_1^2}{\Lambda_0}$ yields that
	\begin{align}\label{7.60-3}
	 &\|(V^m,Q^m)(t)\|^2+|(h,\sqrt\sigma\nabla_yh)(t)|^2_{\mathcal{H}^m}+\|\nabla v(t)\|{_{\mathcal{H}^{m-2}}^2}+\v\|(\text{div}^\varphi v,\nabla^\varphi p,S_\mathbf{n})(t)\|{_{\mathcal{H}^{m-1}}^2}\nonumber\\
	 &+\int_{0}^{t}\v\|\nabla^\varphi V^m(\tau)\|^2 +\v\|\nabla^2v(\tau)\|^2_{\mathcal{H}^{m-2}}+\v^2\|\nabla^2v(\tau)\|^2_{\mathcal{H}^{m-1}}d\tau\\
	&\leq \Lambda(\f1{c_0}, Y_m(0))+\Lambda_\i(t)(1+\|\nabla^\varphi\mbox{div}^\varphi v\|{_{L{^\i_{t,x}}}^2})\int_0^t\Lambda(Y_m(\tau))+\v|\MZ^m h|^2_{\f12}+\|\nabla v\|^2_{\mathcal{H}^{m-1}}d\tau.\nonumber
	\end{align}
This, together with \eqref{10.4} and Lemma \ref{vh}, yields immediately \eqref{7.60-1}. Thus, the proof of this proposition is completed. $\hfill\Box$

	\subsection{Normal derivative estimates: Part II}
	\setcounter{equation}{0}
	
Now we start to estimate $\int_{0}^{t}\|\nabla v\|^2_{\mathcal{H}^{m-1}}d\tau$. It follows from Lemmas \ref{div} and \ref{lem7.2} that one needs only to bound $\|S_\mathbf{n}\|_{\mathcal{H}^{m-1}}$. However, this seems difficult since $\|\nabla^\varphi p\|_{\mathcal{H}^m}$ and $\|\nabla^\varphi\text{div}^\varphi v\|_{\mathcal{H}^m}$, which appear in $F_2$ and $F_3$ in the equation \eqref{7.10} for $S_\mathbf{n}$, cannot be estimated uniformly in $\v$(c.f. Lemma \ref{lem7.5}). As in \cite{Masmoudi-R-1}, we use the vorticity instead of $S_\mathbf{n}$ to perform this estimate. Let us set 
 \begin{equation}
 \omega_\mathbf{n}:=\omega\times\mathbf{N}=\Pi(\omega\times\mathbf{N})
 \end{equation}
 where $\omega$ is the vorticity defined by $\omega=\nabla^\varphi\times v$. Since \eqref{4.4.5} yields that
\begin{align}\label{8.17}
	2S_\mathbf{n}
	 &=\omega_\mathbf{n}+2\Pi\{(\partial_1v\cdot\mathbf{N},\partial_2v\cdot\mathbf{N},0)^t\},
\end{align}
one can obtain from \eqref{7.8} that
\begin{equation}\label{4.5.2}
(\omega_\mathbf{n})^b=-2\Pi\{(\partial_1v\cdot\mathbf{N},\partial_2v\cdot\mathbf{N},0)^t\}.
\end{equation}
which implies that
	\begin{equation}
		 (\MZ^\alpha\omega\times\mathbf{N})^b=-(\omega\times\MZ^\alpha\mathbf{N})^b-([\MZ^\a,\omega\times,\mathbf{N}])^b-2\MZ^\a(\Pi\{(\partial_1v\cdot\mathbf{N},\partial_2v\cdot\mathbf{N},0)^t\}).
	\end{equation}
Thus, one can obtain the sharp estimate that on the boundary, for $\a\leq m-1$, 
\begin{align}\label{9.27}
\sqrt{\v}\int_{0}^{t}|(\MZ^\alpha\omega\times\mathbf{N})^b|^2d\tau&\leq\Lambda_0\sqrt{\v}\int_{0}^{t}|\omega^b|_{\mathcal{H}^{m-2}}^2+|v^b|_{\mathcal{H}^{m}}^2d\tau+\Lambda_\i(t)\int_{0}^{t}\sqrt{\v}|\nabla_yh|_{\mathcal{H}^{m-1}}^2d\tau\nonumber\\
&\leq \Lambda_0\v\int_{0}^{t}\|\nabla^2v\|_{\mathcal{H}^{m-2}}+\|\nabla v\|_{\mathcal{H}^m}^2+\Lambda_\i(t)\int_{0}^{t}Y_m(\tau)d\tau.
\end{align}
Hence, one can only expect to bound $\int_{0}^{t}\|\MZ^{m-1}\omega\times\mathbf{N}\|^4d\tau$ by the similar argument in \cite{Masmoudi-R-1}. Here and henceforth, $\|\MZ^{m-1}\omega\times\mathbf{N}\|^2=\sum_{|\a|\leq m-1}\|\MZ^{\a}\omega\times\mathbf{N}\|^2$. It should be emphasized that, one may not expect to bound directly $\MZ^{m-1}\omega$ like \cite{Masmoudi-R-1}, since it involves $\sqrt{\v}\int_{0}^{t}|\omega^b|_{\mathcal{H}^{m-1}}^2d\tau$ which seems impossible to be bounded due to the invalid of \eqref{20.120} in the compressible flows. However, as shown in \eqref{8.17}, $\omega_\mathbf{n}$ involves only the tangential derivatives of $v$ so that one can follow the argument in \cite{Masmoudi-R-1} for $\omega_\mathbf{n}$ instead of $\omega$. Actually, we will bound $\MZ^{m-1}\omega\times\mathbf{N}$ to avoid too much regularity of $h$ involving. 

\ 

The main result of this subsection is the following proposition:
	\begin{proposition}\label{prop9.7} It holds that, for any $m\geq 6$ that
      \begin{align}
      \Big(\int_0^t\|\nabla v\|^4_{\mathcal{H}^{m-1}}\Big)^{\f12}&\leq \Lambda_0\|\MZ^\a\omega_0\|^2+ \v\int_{0}^{t}\|\nabla V^m\|^2+\|\nabla^{\varphi}\MZ^{m-2} \text{div}^\varphi v\|^2+\|\nabla^{\varphi}\MZ^{m-2} S_\mathbf{n}\|^2d\tau\nonumber\\
      &\quad~~+\Lambda(M(t)) \int_{0}^{t}\Lambda(Y_m(\tau))d\tau
      +\Lambda_{\infty}(t)t^{\f12}\int_0^t\|\nabla p\|^2_{\mathcal{H}^{m-1}}d\tau.\label{9.66}
      \end{align}
      where  
      \begin{equation}\label{4.5.35}
      M(t):=\sup_{0\leq\tau \leq t}\Big(\Lambda(\frac{1}{c_0},\mathcal{Q}(\tau)+\mathcal{Q}_m(\tau)+\v|\MZ^mh|^2_{\f12})
      +\v\int_{0}^{\tau}\|\nabla V^6\|^2+\|\nabla^2v\|_{\mathcal{H}^4}^2ds\Big).
      \end{equation} 	
	\end{proposition}
The proof of this proposition is a consequence of \eqref{4.5.7} and Proposition \ref{prop4.15} below. 

\

Since
    \begin{equation}
		 \MZ^\alpha\omega_\mathbf{n}=\MZ^\alpha(\omega\times\mathbf{N})=\MZ^\alpha\omega\times\mathbf{N}-[\MZ^\alpha,\mathbf{N}\times]\omega,\qquad|\alpha|\leq m-1,
	\end{equation}
and
	\begin{align}
		 \|-[\MZ^\alpha,\mathbf{N}\times]\omega\|^2&\leq\|\omega\|_\i^2\|\nabla_y\eta\|{_{\mathcal{H}^{m-1}}^2}
		 +\|\omega\|_{\mathcal{H} ^{1,\infty}}^2\|\nabla_y\eta\|{_{\mathcal{H}^{m-2}}^2}+\sum_{\substack{|\beta|=1\\\beta+\nu=\alpha}}^{|\beta|=m-3}\|\MZ^\beta\mathbf{N}\times\MZ^\nu\omega\|^2\nonumber\\
		 &\leq\Lambda(\frac{1}{c_0},\mathcal{Q}_m(t)+\mathcal{Q}(t)),\label{9.16-1}
	\end{align}
we can obtain from \eqref{7.0} \eqref{7.3}, \eqref{7.7} \eqref{9.4-0} and \eqref{8.17} that
	\begin{align}\label{4.5.7}
		 \int_{0}^{t}\|\partial_zv\|^4_{\mathcal{H}^{m-1}}d\tau&\leq\Lambda_0\int_0^t\|\MZ^{m-1}\omega\times\mathbf{N}\|^4+\|\text{div}^\varphi v\|_{\mathcal{H}^{m-1}}^4d\tau+\Lambda_\i(t)\int_{0}^{t}\Lambda(Y_m(\tau))d\tau\nonumber\\
		 &\leq\Lambda_0\int_0^t\|\MZ^{m-1}\omega\times\mathbf{N}\|^4d\tau+\Lambda_\i(t)\int_{0}^{t}\Lambda(Y_m(\tau))d\tau.
	\end{align}
So, it remains  to control $\MZ^\alpha\omega\times\mathbf{N}$, for $|\a|\leq m-1$. 

Now we first derive the equation solved by $\MZ^\alpha\omega\times\mathbf{N}$ as follows. Applying $\nabla^\varphi\times$ to $\eqref{1.22}_2$ yields that
	\begin{equation}\label{9.19}
		\varrho\partial{_t^\varphi}\omega+\varrho v\cdot\nabla^\varphi\omega-\m\v\Delta^\varphi\omega=F_1,
	\end{equation}
	where $F_1$
	\begin{equation}\label{9.20}
		 F_1=-\nabla^\varphi\varrho\times\partial{_t^\varphi}v-\nabla^\varphi\varrho\times((v\cdot\nabla^\varphi)v)+\varrho\omega\cdot\nabla^\varphi v-\varrho\mbox{div}^\varphi v\omega.
	\end{equation}
	Applying $\MZ^\a$, with $|\a|= m-1$, to \eqref{9.19} gives that
	\begin{equation}\label{9.21}
		\varrho\partial_t\MZ^\a\omega+\varrho v_y\cdot\nabla_y\MZ^\a\omega+\varrho V_z\partial_z\MZ^\a\omega-\m\v\Delta^\varphi\MZ^\a\omega=\MZ^\a F_1+F_2,
	\end{equation}
	where
	\begin{equation}\label{9.22}
		F_2=[\MZ^\a,\varrho]\partial_t\omega+[\MZ^\a,\varrho v_y]\nabla_y\omega+[\MZ^\a,\varrho V_z]\partial_z\omega+\varrho V_z[\MZ^\a,\partial_z]\omega+\m\v[\MZ^\a,\Delta^\varphi]\omega.
	\end{equation}
	Multiply \eqref{9.21} to get that
	\begin{equation}\label{9.23}
		\varrho\partial_t(\MZ^\a\omega\times\mathbf{N})+\varrho v_y\cdot\nabla_y(\MZ^\a\omega\times\mathbf{N})+\varrho V_z\partial_z(\MZ^\a\omega\times\mathbf{N})-\m\v\Delta^\varphi(\MZ^\a\omega\times\mathbf{N})=F,
	\end{equation}
	where the source term $F$ is given by
	\begin{equation}\label{9.24}
		F=\MZ^\a F_1\times\mathbf{N}+F_2\times\mathbf{N}+F_3,
	\end{equation}
	with $F_3$ given by
	\begin{equation}\label{9.25}
		F_3=\varrho\MZ^\a\omega\times\partial_t\mathbf{N}+\varrho v_y\cdot\MZ^\a\omega\times\nabla_y\mathbf{N}+\varrho V_z\MZ^\a\omega\times\partial_z\mathbf{N}+\m\v\MZ^\a\omega\times\Delta^\varphi\mathbf{N}+2\m\v\nabla^\varphi\MZ^\a\omega\times\nabla^\varphi\mathbf{N}.
	\end{equation}

Next, we follow the argument in \cite{Masmoudi-R-1}, where $\MZ^{m-1}\omega$ is bounded, to perform the estimate of $\MZ^{m-1}\omega\times\mathbf{N}$ as follows. 
\begin{proposition}\label{prop4.15}
		For any $t\in[0,T^\v]$ and $|\a|= m-1$,  there exists $\Lambda(M(t))$ such that
		\begin{align}\label{9.36-1}
			 \|\MZ^\alpha\omega\times\mathbf{N}\|_{L^4(0,t;L^2(\mathcal{S}))}^2&\leq\Lambda_0\|\MZ^\a\omega_0\|^2+\Lambda_0\v\int_{0}^{t}\|\nabla V^m\|^2+\|\nabla^2 v\|_{\mathcal{H}^{m-2}}^2d\tau\nonumber\\
			&~~~~ +\Lambda(M(t))\int_{0}^{t}\Lambda(Y_m(\tau))d\tau
			 +\Lambda_{\i}(t)t^{\f12}\int_{0}^{t}\|\nabla(p,v)(\tau)\|{_{\mathcal{H}^{m-1}}^2}d\tau.
		\end{align}
		where $M(t)$ is defined in \eqref{4.5.35}.
	\end{proposition}
The proof of this Proposition will be a consequence  of the following lemmas. As in \cite{Masmoudi-R-1}, we split $\MZ^\a\omega\times\mathbf{N}$ into two parts of the form
	\begin{equation}\label{spl}
		\MZ^\a\omega\times\mathbf{N}=\omega{_h^\a}+\omega{_{nh}^\a},
	\end{equation}
	where $\omega{_{nh}^\a}$ solves  the non-homogeneous problem:
	\begin{equation}\label{nhw}
	\begin{cases}
	\varrho\partial_t\omega{_{nh}^\a}+\varrho v_y\cdot\nabla_y\omega{_{nh}^\a}+\varrho V_z\partial_z\omega{_{nh}^\a}-\m\v\Delta^\varphi\omega{_{nh}^\a}=F,\\
	\omega{_{nh}^\a}|_{t=0}=\MZ^\a\omega_0\times\mathbf{N}_0,\\
	 \omega{_{nh}^\a}|_{z=0}=0,
	\end{cases}
	\end{equation}
while $\omega{_h^\a}$ solves  the homogeneous one:
	\begin{equation}\label{hw}
	\begin{cases}
	\varrho\partial_t\omega{_{h}^\a}+\varrho v_y\cdot\nabla_y\omega{_{h}^\a}+\varrho V_z\partial_z\omega{_{h}^\a}-\m\v\Delta^\varphi\omega{_{h}^\a}=0,\\
	\omega{_n^\a}|_{t=0}=0,\\
	\omega{_n^\a}|_{z=0}=(\MZ^\alpha\omega\times\mathbf{N})^b.
	\end{cases}	
	\end{equation}
	
The solution $\omega{_{nh}^\a}$ to \eqref{nhw} can be estimated by standard energy estimates.
	\begin{lemma}\label{p9.2}
		For any $|\a|\leq m-1$, it holds that
		\begin{align}\label{9.28}
			 &\|\omega{_{nh}^\a}\|^2+\m\v\int_{0}^{t}\|\nabla^\varphi\omega{_{nh}^\a}\|^2d\tau
			 \leq\Lambda_0\|\MZ^\a\omega_0\|^2+\Lambda_0\v\int_{0}^{t}\|\nabla^2 v\|_{\mathcal{H}^{m-2}}^2d\tau\nonumber\\
			 &~~~~~~~~~~~~~~~~~~~~~~~~~~~~~~~~~~~~~+\Lambda_\i(t)\int_{0}^{t}\v|\MZ^mh|^2_{\f12}+\|\nabla(p,v)\|{_{\mathcal{H}^{m-1}}^2}+\Lambda(Y_m(\tau))d\tau.
		\end{align}
	\end{lemma}
	\noindent\textbf{Proof}. Noting the homogeneous Dirichlet boundary condition for $\omega{_{nh}^\a}$, one deduces from \eqref{nhw} by a standard energy estimate that
	\begin{equation}\label{9.29}
		 \frac{1}{2}\int_\mathcal{S}\varrho|\omega{_{nh}^\a}|^2d\mathcal{V}_\tau+\m\v\int_{0}^{t}\int_\mathcal{S}|\nabla^\varphi\omega{_{nh}^\a}|d\mathcal{V}_\tau d\tau\leq\Lambda_0\int_\mathcal{S}\varrho_0|\MZ^\a\omega_0|^2d\mathcal{V}_\tau+|\int_{0}^{t}\int_\mathcal{S}F\cdot\omega{_{nh}^\a}d\mathcal{V}_\tau d\tau|,
	\end{equation}
	where $F$ is given by \eqref{9.24}. Then \eqref{2.6} implies that
	\begin{align}\label{9.30}
		&|\int_{0}^{t}\int_\mathcal{S}\MZ^\a F_1\times\mathbf{N}\cdot\omega{_{nh}^\a}d\mathcal{V}_\tau d\tau|+\int_{0}^{t}\int_{\mathcal{S}}|F_3\cdot\omega^\a_{nh}|^2d\mathcal{V}_\tau d\tau\\
		&\leq\f1{16}\sup_{0\leq\tau\leq t}\int|\omega^\a_{nh}(\tau)|^2d\mathcal{V}_\tau+\Lambda_\i(t)\int_{0}^{t}\|(p,v)\|{_{\mathcal{H}^m}^2}+\|(\nabla p,\nabla v)\|{_{\mathcal{H}^{m-1}}^2}+|\MZ^{m-1}h|{_{\frac{1}{2}}^2}d\tau.\nonumber
	\end{align}
Using similar arguments as \cite{Masmoudi-R-1},  we can obtain that 
	\begin{align}\label{9.54}
		 &\int_{0}^{t}\int_{\mathcal{S}}|F_2\times\mathbf{N}\cdot\omega^\a_{nh}|d\mathcal{V}_\tau d\tau\nonumber\\
		&\leq\frac{\m}{8}\v\int_{0}^{t}\int|\nabla^{\varphi} \omega{_{nh}^\a}|^2d\mathcal{V}_{\tau}d\tau+\Lambda_0\v\int_{0}^{t}\|\nabla^2 v\|^2_{\mathcal{H}^{m-1}}d\tau+\f1{16}\sup_{0\leq\tau\leq t}\int|\omega^\a_{nh}(\tau)|^2d\mathcal{V}_\tau\nonumber\\
		&~~~~~+\Lambda_\i(t)\int_{0}^{t}\v|\MZ^mh|^2_{\f12}+\|\nabla(p,v)\|{_{\mathcal{H}^{m-1}}^2}+\Lambda(Y_m(\tau))d\tau.
	\end{align}
Substituting \eqref{9.30} and \eqref{9.54}  into \eqref{9.29} proves  Proposition \ref{p9.2}. Thus,  the proof  is completed.$\hfill\Box$

	\

	To treat $\omega_h$,  we will modify the approach of  microlocal analysis  used in \cite{Masmoudi-R-1} to obtain that
	\begin{lemma}\label{p9.3}
		For any $|\a|= m-1$,  it holds that
		\begin{align}\label{9.48}
			 \Big(\int_{0}^{t}\|\omega{_{h}^\a}(\tau)\|^4d\tau\Big)^\frac{1}{2}\leq \v\int_{0}^{t}\|\nabla V^m\|^2+\|\nabla^2v\|^2_{\mathcal{H}^{m-2}}d\tau+\Lambda(M(t))\int_{0}^{t}\Lambda(Y_m(\tau))d\tau.
		\end{align}
		where $M(t)$ is defined in \eqref{4.5.35}.
	\end{lemma}
	
	In order to prove the above lemma,  we follow the argument in \cite{Masmoudi-R-1} to use Lagrangian coordinates to eliminate the convection term so that it is convenient to perform the microlocal symmetrizer method. Let us define a parametrization by
	\begin{equation}
		\partial_tX(t,y,z)=u(t,X(t,y,z))=v(t,\Phi(t,\cdot)^{-1}\circ X),~~ X(0,y,z)=\Phi(0,y,z),
	\end{equation}
	where $\Phi(t,\cdot)^{-1}$ stands for the inverse of the map $\Phi(t,\cdot)$ defined by \eqref{1.14}. Define $J(t, y,z)=|det \nabla X(t,y,z)|$ to be the Jacobian of the change of variable.
	The following estimates for $X$ are proved in Lemma 10.5 in \cite{Masmoudi-R-1}:
	\begin{lemma}\label{l9.4}
		It holds that for $t\in[0,T]$,
		\begin{align}\label{9.50}
			&\|\nabla X(t)\|_{L^\i}+\|\partial_t\nabla X(t)\|_{L^\i}\leq\Lambda_0e^{t\Lambda(M(t))},\\
			&\|\nabla X(t)\|_{1,\i}+\|\partial_t\nabla X(t)\|_{1,\i}\leq\Lambda(M(t))e^{t\Lambda(M(t))},\\
			&\sqrt{\v}\|\nabla^2X\|_{1,\i}+\sqrt{\v}\|\partial_t\nabla^2X\|_{L^\i}\leq\Lambda(M(t))(1+t^2)e^{t\Lambda(M(t))}.
		\end{align}
	\end{lemma}
	
	Now set
	\begin{eqnarray}
		\Om^\a=e^{-\g t}\omega{_h^\a}(t,\Phi^{-1}\circ X)
	\end{eqnarray}
	where $\g>0$ is a large parameter to be chosen. Then $\Om^\a$ solves in $\mathcal{S}$ the equation
	\begin{equation}\label{9.63}
		 a_0(\partial_t\Om^\a+\g\Om^\a)-\m\v\partial_i(a_{ij}\partial_j\Om^\a)=0,
	\end{equation}
	where
	\begin{equation}\label{9.56}
		a_0=\varrho(t,\Phi^{-1}\circ X)|J|^\frac{1}{2},\quad(a_{ij})=|J|^\frac{1}{2}P^{-1},\quad P_{ij}=\partial_iX\cdot\partial_jX.
	\end{equation}
	Thanks to Lemma \ref{l9.4}, the equation \eqref{9.63} is a parabolic equation. On the boundary, it holds that
	\begin{equation}\label{9.57}
		\Om^\a|_{z=0}=(\Om^\a)^b:=e^{-\g t}(\MZ^\a\omega\times\mathbf{N})(t,(\Phi^{-1}\circ X)(t,y,0).
	\end{equation}
	
	\
	
	The following theorem holds:
	\begin{theorem}[Masmoudi-Rousset \cite{Masmoudi-R-1}]\label{M-S}
		There exists $\g_0$ depending only on $M>0$ such that for $\g\geq\g_0$, the solution of \eqref{9.63} with the boundary condition \eqref{9.57} satisfies the estimate
		\begin{equation}
			 \|\Om^{m-1}\|{^2_{H^\frac{1}{4}(0,T;L^2(\mathcal{S}))}}\leq\Lambda(M(T))\sqrt{\v}\int_{0}^{T}|(\Om^{m-1})^b|{^2_{L^2(\mathbb{R}^2)}},
		\end{equation}
		where the norm $H^\frac{1}{4}(0,T;L^2)$ is defined by
		\begin{equation*}
			 \|f\|_{H^\frac{1}{4}(0,T;L^2)}=\inf\{\|Pf\|_{H^\frac{1}{4}(\mathbb{R},L^2(\mathcal{S}))}, Pf=f ~~\mbox{on}~~ [0,T]\times \mathcal{S}\},
		\end{equation*}
		with the norm on the whole space by Fourier transform in time.
	\end{theorem}
	
	\noindent\textbf{Proof of Theorem \ref{M-S}}. The proof of this Theorem is almost the same as the proof of Theorem 10.6 in \cite{Masmoudi-R-1} where the symmetrized method and paradifferential calculus are used. This is because the coefficients of \eqref{9.63} given in \eqref{9.56} is the same type (determined by the estimate in Lemma \ref{l9.4}, $\|(\varrho,\nabla\varrho)\|_{1,\i}\leq M$ and initial data which ensure $\Phi_0$) of symbols as the one in \cite{Masmoudi-R-1} which is essential in the progress of taking use of paradifferential calculus. We  refer to Theorem 10.6 in \cite{Masmoudi-R-1} for more details. $\hfill\Box$
	
	\
	
	\noindent\textbf{Proof of Lemma \ref{p9.3}}.  As a consequence of Theorem \ref{M-S} and the Sobolev embedding inequality, one has that
	\begin{equation*}
		\|\Om^\alpha\|_{L^4(0,T;L^2(S))}^2\leq C\|\Om^{m-1}\|_{H^\frac{1}{4}(0,T;L^2(S))}^2\leq \Lambda(M(T))\sqrt{\v}\int_{0}^{T}|(\Om^{m-1})^b|_{L^2(\mathbb{R}^2)}^2
	\end{equation*}
	Consequently, it follows from   changing of variable and \eqref{9.27} that
	\begin{align}\label{9.67}
		\|\omega{_h^\a}\|_{L^4(0,t;L^2(S))}^2&\leq \Lambda(M(t))\sqrt{\v}\int_{0}^{t}|(\MZ^\a\omega\times\mathbf{N})^b|_{L^2(\mathbb{R}^2)}^2\nonumber\\
		&\leq\v\int_{0}^{t}\|\nabla v\|_{\mathcal{H}^m}^2+\|\nabla^2v\|^2_{\mathcal{H}^{m-2}}d\tau+\Lambda(M(t)) \int_{0}^{t}Y_m(\tau)d\tau.
	\end{align}
	which yields Proposition \ref{p9.3}.  Therefore, the proof   is completed.  $\hfill\Box$

	\

\subsection{$L^\i$-estimates}

	\setcounter{equation}{0}

In order to close the estimates, one needs to bound the $L^\i$ norms of $v$, $h$ and $p$ contained in $\Lambda_\i(t)$.
We start with the following estimates through standard Sobolev embedding theorem and anisotropic Sobolev embedding theorem.
\begin{lemma}\label{l8.1}
	For any $t\in[0,T^\v]$, the following estimates hold:
	\begin{align}\label{8.1}
		&|h(t)|_{\mathcal{H}^{k,\i}}\lesssim|\MZ^{k}h(t)|_\frac{3}{2},~ k\in\mathbb{N},\\\label{8.2}
		&\|(p,v)(t)\|_{\mathcal{H}^{2,\i}}^2\leq\Lambda(\frac{1}{c_0},\mathcal{Q}_5(t)),\\\label{8.3}
		&\|\nabla p(t)\|{_{\mathcal{H}^{1,\i}}^2}\leq\|\Delta^\varphi p(t)\|{_{\mathcal{H}^1}^2}+\Lambda(\frac{1}{c_0},\mathcal{Q}_5(t))\|\nabla p(t)\|_{\mathcal{H}^4}^2,\\\label{8.4}
		&\|\mbox{div}^\varphi  v(t)\|{_{\mathcal{H}^{1,\i}}^2}\leq\|\Delta^\varphi p(t)\|{_{\mathcal{H}^1}^2}+\Lambda(\frac{1}{c_0},\mathcal{Q}_5(t))(1+\|\nabla p(t)\|_{\mathcal{H}^4}^2),\\\label{8.5}
		&\|\nabla\mbox{div}^\varphi v(t)\|{_{\i}^2}\leq\|\Delta^\varphi p\|_{\mathcal{H}^1}^2+\Lambda(\frac{1}{c_0},\mathcal{Q}_5(t)+\|\nabla v(t)\|_\i^2)(1+\|\nabla p\|_{\mathcal{H}^4}^2),\\\label{8.6}
		&\|\nabla\mbox{div}^\varphi v(t)\|{_{\mathcal{H}^{1,\i}}^2}\leq\delta\|\Delta^\varphi p(t)\|{_{\mathcal{H}^2}^2}+C_\d\Lambda(\frac{1}{c_0},\mathcal{Q}_6(t)+\|\nabla v(t)\|{_{\mathcal{W}^{1,\i}}^2})(1+\|\nabla p\|_{\mathcal{H}^5}^2).
	\end{align}
\end{lemma}
\noindent\textbf{Proof}. \eqref{8.1} follows from the two dimensional Sobolev embedding. The anisotropic Sobolev embedding \eqref{2.9} yields that
\begin{equation}\label{8.7}
	\|(p,v)(t)\|_{\mathcal{H}^{2,\i}}^2\lesssim\|(\nabla p,\nabla v)(t)\|_{\mathcal{H}^3}\|(p,v)(t)\|_{\mathcal{H}^4}\leq\Lambda(\frac{1}{c_0},\mathcal{Q}_5(t)).
\end{equation}
To prove \eqref{8.3}, one notes that
\begin{equation*}
	 \Delta^\varphi=\frac{1+|\nabla_y\varphi|^2}{|\partial_z\varphi|^2}\partial_{zz}+\sum_{i=1,2}\Big(\partial_{ii}-\frac{2\partial_i\varphi}{\partial_z\varphi}\partial_i\partial_z-2\partial_i(\frac{\partial_i\varphi}{\partial_z\varphi})\partial_z+\frac{\partial_i\varphi}{\partial_z\varphi}\partial_z(\frac{\partial_i\varphi}{\partial_z\varphi})\partial_z+\frac{1}{\partial_z\varphi}\partial_z(\frac{1}{\partial_z\varphi})\partial_z\Big),
\end{equation*}
which follows from the definition of $\partial{_i^\varphi}$ immediately. Thus, \eqref{2.9} implies that
\begin{align}\label{8.10}
	\|\nabla p(t)\|{_{\i}^2}&\lesssim\|\nabla p(t)\|_3\|\partial_z\nabla p(t)+\|\nabla p(t)\|^2_3\|\lesssim\|\nabla p(t)\|_3\Big(\|\partial_z p(t)\|_3+\|\partial_{zz}p(t)\|\Big)\nonumber\\
	&\lesssim\Lambda(\frac{1}{c_0},|h(t)|_{2,\i})\|\nabla p(t)\|_3\Big(\|\nabla p(t)\|_3+\|p(t)\|_2+\|\Delta^\varphi p(t)\|\Big)\\
	&\leq\|\Delta^\varphi p(t)\|^2+\Lambda(\frac{1}{c_0},\mathcal{Q}_5(t)).\nonumber
\end{align}
Similarly, 
\begin{align}
	\|\MZ\nabla p(t)\|{_{\i}^2}&\lesssim\|\MZ\nabla p\|_3\|\partial_z\MZ\nabla p\|+\|\nabla p\|^2_4\nonumber\\
	&\lesssim\Lambda(\frac{1}{c_0},|h(t)|_{3,\i})\|\nabla p(t)\|_{\mathcal{H}^4}\Big(\|\nabla p(t)\|_{\mathcal{H}^4}+\|p(t)\|_{\mathcal{H}^3}+\|\Delta^\varphi p(t)\|_{\mathcal{H}^1}\Big)\nonumber\\
	&\leq\|\Delta^\varphi p(t)\|_{\mathcal{H}^1}^2+\Lambda(\frac{1}{c_0},\mathcal{Q}_5(t))\|\nabla p\|_{\mathcal{H}^4}^2.
\end{align}
Therefore, \eqref{8.3} holds true. Next, it follows from \eqref{6.3},   \eqref{8.2} and \eqref{8.3} that
\begin{align}
	\|\text{div}^\varphi v(t)\|{_{\mathcal{H}^{1,\i}}^2}&\lesssim\|\frac{1}{p}(t)\|{_{\mathcal{H}^{1,\i}}^2}\Big(\|p(t)\|{_{\mathcal{H}^{2,\i}}^2}+\|v_y(t)\|{_{\mathcal{H}^{1,\i}}^2}\|p(t)\|{_{\mathcal{H}^{2,\i}}^2}+\|V_z(t)\|{_{\mathcal{H}^{1,\i}}^2}\|\partial_z p\|{_{\mathcal{H}^{1,\i}}^2}\Big)\nonumber\\
	&\lesssim\Lambda(\frac{1}{c_0},\mathcal{Q}_5(t))\Big\{1+\|\nabla p(t)\|_{\mathcal{H}^4}\Big[\|\nabla p(t)\|_{\mathcal{H}^4}+\|p(t)\|_{\mathcal{H}^3}+\|\Delta^\varphi p(t)\|_{\mathcal{H}^1}\Big]\Big\}\nonumber\\
	&\leq\|\Delta^\varphi p(t)\|{_{\mathcal{H}^1}^2}+\Lambda(\frac{1}{c_0},\mathcal{Q}_5(t))(1+\|\nabla p\|_{\mathcal{H}^4}^2).
\end{align}
Furthermore, it follows from the relation
$$\nabla\mbox{div}^\varphi v=-\frac{1}{\g p}(\partial_t\nabla p+v_y\cdot\nabla_y\nabla p+V_z\partial_z\nabla p)-\nabla(\frac{1}{\g p})\partial_tp-\nabla(\frac{v_y}{\g p})\cdot\nabla_yp-\nabla(\frac{V_z}{\g p})\partial_zp,$$
that
\begin{align}
	&\|\nabla\text{div}^\varphi v(t)\|_\i^2\nonumber\\
	&\leq\Lambda_0\Big\{\|\nabla p(t)\|_{\mathcal{H}^{1,\i}}^2+\|v_y(t)\|_\i^2\|\nabla p(t)\|_{\mathcal{H}^{1,\i}}^2+\|(V_z\partial_z\nabla p)(t)\|_\i^2+\|\nabla{p}(t)\|{_\i^2}\|p(t)\|{_{\mathcal{H}^{1,\i}}^2}\nonumber\\
	&\quad+(\|\nabla v_y(t)\|{_\i^2}+\|(v_y\nabla p)(t)\|_\i^2)\|p(t)\|{_{\mathcal{H}^{1,\i}}^2}+(\|\nabla V_z(t)\|{_\i^2}+\|(V_z\nabla p)(t)\|_\i^2)\|\partial_z p(t)\|{_\i^2}\Big\}\nonumber\\
	&\leq\Lambda(\frac{1}{c_0},\mathcal{Q}_5(t))(1+\|\nabla v\|_\i^2)(1+\|\nabla p\|_{\mathcal{H}^{1,\i}}^2)\nonumber\\
	&\leq\|\Delta^\varphi p\|_{\mathcal{H}^1}^2+\Lambda(\frac{1}{c_0},\mathcal{Q}_5(t)+\|\nabla v(t)\|_\i^2)(1+\|\nabla p\|_{\mathcal{H}^4}^2),
\end{align}
since $V_z$ vanishes on the boundary. Similarly, 
\begin{equation}
	\|\nabla\mbox{div}^\varphi v(t)\|{_{\mathcal{H}^{1,\i}}^2}\leq\delta\|\Delta^\varphi p(t)\|{_{\mathcal{H}^2}^2}+C_\d\Lambda(\frac{1}{c_0},\mathcal{Q}_6(t)+\|\nabla v(t)\|{_{\mathcal{H}^{1,\i}}^2})(1+\|\nabla p\|_{\mathcal{H}^5}^2).
\end{equation}
Therefore,  the proof of this lemma is completed.$\hfill\Box$

	\
	
We now estimate $\|\Delta^\varphi p\|{_{\mathcal{H}^1}^2}$. Applying $\mbox{div}^\varphi$ to the momentum equations yields that
	\begin{equation}\label{8.52}
	-(2\m+\l)\v\Delta^\varphi\mbox{div}^\varphi v+\Delta^\varphi p=\mbox{div}^\varphi(\varrho\partial{_t^\varphi v}+\varrho v\cdot\nabla^\varphi v)
	\end{equation}
	This, together with \eqref{6.3} and \eqref{4.15} shows that
	\begin{align}\label{8.53}
	&\v\partial_t\Delta^\varphi p+\frac{\g p}{2\m+\l}\Delta^\varphi p=\dfrac{\g p}{2\m+\l}\mbox{div}^\varphi(\varrho\partial{_t^\varphi v}+\varrho v\cdot\nabla^\varphi v)-\v(v_y\nabla_y+V_z\partial_z)\Delta^\varphi p\nonumber\\
	&\quad-\v\g p\Big(\partial{_t^\varphi}p\Delta^\varphi(\frac{1}{\g p})+\nabla^\varphi p\Delta^\varphi(\frac{v}{\g p})+2\nabla^\varphi\partial{_t^\varphi}p\nabla^\varphi(\frac{1}{\g p})+2\nabla^\varphi(\frac{v}{\g p})\nabla^\varphi\nabla^\varphi p\Big).
	\end{align}
	
	\begin{lemma}\label{prop4.26}
		For $m\geq 6$, it holds that
			\begin{align}\label{20.6}
			&\sup_{0\leq \tau\leq t}\Big(\|\Delta^\varphi p(\tau)\|^2_{\mathcal{H}^1}+\v\|\Delta^\varphi p(\tau)\|^2_{\mathcal{H}^2}\Big)+\int_{0}^{t}\|\Delta^\varphi p\|^2_{\mathcal{H}^2}d\tau\nonumber\\
			&\leq C\Lambda_0 \Big(\|(\Delta^\varphi p)(0)\|^2_{\mathcal{H}^1}+\v\|(\Delta^\varphi
			p)(0)\|^2_{\mathcal{H}^2}\Big)+\Lambda_\infty(t)\int_0^t\Lambda(Y_m(\tau))d\tau.
			\end{align}
	\end{lemma}
	\noindent\textbf{Proof}. Applying $\MZ^\a(|\a|= 2)$ to \eqref{8.53}, multiplying the resulting equation by $\MZ^\a \Delta^\varphi p$ and then integrating over $\mathcal{S}\times[0,t]$, one can get that
	\begin{align}\label{8.55}
	&\f12\v\|\MZ^\a\Delta^\varphi p\|^2+\dfrac{1}{2\m+\l}\int_{0}^{t}\int_{\mathcal{S}}\MZ^\a(\g p\Delta^\varphi p)\MZ^\a\Delta^\varphi pd\mathcal{V}_\tau d\tau\nonumber\\
	&=\frac{\v}{2}\|(\MZ^\a\Delta^\varphi p)(0)\|^2-\v\int_{0}^{t}\int_{\mathcal{S}}\MZ^\a(v_y\cdot\nabla_y\Delta^\varphi p+V_z\partial_z\Delta^\varphi p)\MZ^\a\Delta^\varphi pd\mathcal{V}_\tau d\tau\nonumber\\
	&\quad-\v\int_{0}^{t}\int_{\mathcal{S}}\MZ^\a(\g p\partial{_t^\varphi}p\Delta^\varphi(\frac{1}{\g p}))\MZ^\a\Delta^\varphi pd\mathcal{V}_\tau d\tau-\v\int_{0}^{t}\int_{\mathcal{S}}\MZ^\a(\g p\nabla^\varphi p\Delta^\varphi(\frac{v}{\g p}))\MZ^\a\Delta^\varphi pd\mathcal{V}_\tau d\tau\nonumber\\
	&\quad-2\v\int_{0}^{t}\int_{\mathcal{S}}\MZ^\a(\g p(\nabla^\varphi\partial{_t^\varphi}p\nabla^\varphi(\frac{1}{\g p})+\nabla^\varphi(\frac{v}{\g p})\nabla^\varphi\nabla^\varphi p))\MZ^\a\Delta^\varphi pd\mathcal{V}_\tau d\tau\\
	&\quad-\dfrac{1}{2\m+\l}\int_{0}^{t}\int_{\mathcal{S}}\MZ^\a\{\g p\mbox{div}^\varphi(\varrho(\partial{_t^\varphi}+v\cdot\nabla^\varphi)v)\}\MZ^\a\Delta^\varphi pd\mathcal{V}_\tau d\tau=:\frac{\v}{2}\|(\MZ^\a\Delta^\varphi p)(0)\|^2+\sum_{i=1}^{5}J_i\nonumber
	\end{align}
	The second term on the LHS above can be estimated directly as
	\begin{align}\label{8.56}
	&\dfrac{1}{2\m+\l}\int_{0}^{t}\int_{\mathcal{S}}\MZ^\a(\g p\Delta^\varphi p)\MZ^\a\Delta^\varphi pd\mathcal{V}_\tau d\tau\geq\dfrac{1}{2\m+\l}\int_{0}^{t}\int_{\mathcal{S}}\g p|\MZ^\a\Delta^\varphi p|^2d\mathcal{V}_\tau d\tau\nonumber\\
	&\quad+\dfrac{2}{2\m+\l}\int_{0}^{t}\int_\mathcal{S}\g\MZ p\MZ(\Delta^\varphi p)\MZ^\a(\Delta^\varphi p)d\mathcal{V}_\tau d\tau+\dfrac{1}{2\m+\l}\int_{0}^{t}\int_{\mathcal{S}}\g\MZ^\a p\Delta^\varphi p\MZ^\a(\Delta^\varphi p)d\mathcal{V}_\tau d\tau\nonumber\\
	&\geq\dfrac{3}{4(2\m+\l)}\int_{0}^{t}\int_{\mathcal{S}}\g p|\MZ^\a\Delta^\varphi p|^2d\mathcal{V}_\tau d\tau-\Lambda_0\int_{0}^{t}\|\MZ p\|{_\i^2}\|\Delta^\varphi p\|{_{\mathcal{H}^1}^2}+\|\MZ^\a p\|{_\i^2}\|\Delta^\varphi p\|^2d\tau\nonumber\\
	&\geq\frac{3}{4(2\m+\l)}\int_{0}^{t}\int_{\mathcal{S}}\g p|\MZ^\a\Delta^\varphi p|^2d\mathcal{V}_\tau d\tau-\int_{0}^{t}\Lambda(\frac{1}{c_0},\mathcal{Q}_m(t)+\|\Delta^\varphi p\|{_{\mathcal{H}^1}^2})d\tau,
	\end{align}
	where $m\geq 5$.
	Next, the terms on the RHS of  \eqref{8.55} can be estimated separately. First, it follows from integration by parts, \eqref{2.7} and \eqref{4.26} that
	\begin{align}\label{8.57}
	 |J_1|&\leq\v|\int_{0}^{t}\int_{\mathcal{S}}(v_y\nabla_y\MZ^\a\Delta^\varphi p+V_z\partial_z\MZ^\a\Delta^\varphi p)\MZ^\a\Delta^\varphi pd\mathcal{V}_\tau d\tau|\nonumber\\
	 &\quad+\v|\int_{0}^{t}\int_{\mathcal{S}}([\MZ^\a,v_y\nabla_y]\Delta^\varphi p+[\MZ^\a, V_z]\partial_z\Delta^\varphi p+V_z[\MZ^\a, \partial_z]\Delta^\varphi p)\MZ^\a\Delta^\varphi pd\mathcal{V}_\tau d\tau|\nonumber\\
	&\leq\Lambda_{\i}(t)\int_{0}^{t}\v\|\Delta^\varphi p\|{_{\mathcal{H}^2}^2}d\tau+\v\int_{0}^{t}\|\dfrac{1-z}{z}\MZ^\a V_z\|{_\i^2}\|\Delta^\varphi p\|{_{\mathcal{H}^1}^2}+\|\dfrac{1-z}{z}\MZ V_z\|{_\i^2}\|\Delta^\varphi p\|{_{\mathcal{H}^2}^2} d\tau\nonumber\\
	&\leq\v^2\int_{0}^{t}\|\nabla^2 v\|{_{\mathcal{H}^3}^2}d\tau+\Lambda_{\i}(t)\int_{0}^{t}\v\|\Delta^\varphi p\|{_{\mathcal{H}^2}^2}d\tau+\int_{0}^{t}\Lambda(\frac{1}{c_0}, \mathcal{Q}_m(t)+\|\Delta^\varphi p\|{_{\mathcal{H}^1}^2})d\tau\nonumber\\
	&\leq \Lambda_{\i}(t)\int_{0}^t\Lambda(Y_m(\tau))d\tau, ~~\mbox{for}~m\geq6,
	\end{align}
	where \eqref{6.37-2} has been used in the last inequality.
	
	Note that
	\begin{equation*}
	\Delta^\varphi(\frac{1}{p})=\frac{2}{p^3}|\nabla^\varphi p|^2-\frac{\Delta^\varphi p}{p^2},\quad\Delta^\varphi(\frac{v}{p})=\frac{\Delta^\varphi v}{p}-\frac{2\partial{_i^\varphi}v\partial{_i^\varphi}p}{p^2}+\frac{2v}{p^3}|\nabla^\varphi p|^2-\frac{v}{p^2}\Delta^\varphi p,
	\end{equation*}
	Thus, for $m\geq6$, it holds that
	\begin{align}\label{8.58}
	&|J_2|\leq\frac{1}{16(2\m+\l)}\int_{0}^{t}\int_{\mathcal{S}}\g p|\MZ^\a\Delta^\varphi p|^2d\mathcal{V}_\tau d\tau+C\v^2\int_{0}^{t} \|\MZ^\a(\g p\partial{_t^\varphi}p\Delta^\varphi(\frac{1}{\g p}))\|^2d\tau\nonumber\\
	&\leq \frac{1}{16(2\m+\l)}\int_{0}^{t}\int_{\mathcal{S}}\g p|\MZ^\a\Delta^\varphi p|^2d\mathcal{V}_\tau d\tau
	+\Lambda_{\i}(t)\int_0^t\Lambda(\f1{c_0}, \mathcal{Q}_m(t)+\|\Delta^\varphi p\|{_{\mathcal{H}^1}^2}+\v\|\Delta^\varphi p\|{_{\mathcal{H}^2}^2})d\tau\nonumber\\
	&\leq \frac{1}{16(2\m+\l)}\int_{0}^{t}\int_{\mathcal{S}}\g p|\MZ^\a\Delta^\varphi p|^2d\mathcal{V}_\tau d\tau
	+\Lambda_{\i}(t)\int_0^t\Lambda(Y_m(\tau))d\tau,
	\end{align}
	and
	\begin{align}\label{8.59}
	&|J_3|\leq \v\Big\{\int_{0}^{t}\|\MZ^\a\Delta^\varphi p\|^2d\tau\Big\}^\frac{1}{2}\nonumber\\
	&\qquad\times\Big\{\int_{0}^{t}\|\MZ^\a\nabla^\varphi p\cdot\Delta^\varphi(\frac{v}{p})\|^2+\|\MZ\nabla^\varphi p\MZ\Delta^\varphi(\frac{v}{p})\|^2+\|\nabla^\varphi p\MZ^\a\Delta^\varphi(\frac{v}{p})\|^2d\tau\Big\}^\frac{1}{2}\nonumber\\
	&\leq \frac{1}{16(2\m+\l)}\int_{0}^{t}\int_{\mathcal{S}}\g p|\MZ^\a\Delta^\varphi p|^2d\mathcal{V}_\tau d\tau+C\v^2\Lambda(\frac{1}{c_0},\mathcal{Q}(t)+\mathcal{Q}_m)\int_{0}^{t}\Big(\|\nabla^\varphi p\|{_{\mathcal{H}^2}^2}+\|\nabla^\varphi p\nabla^\varphi v\|{_{\mathcal{H}^2}^2}\nonumber\\
	&\quad+\|\MZ^\a\nabla^\varphi p\|{_2^2}\|\Delta^\varphi v\|{_\i^2}+\|\MZ^\a\nabla^\varphi p\|{_\i^2}\|\Delta^\varphi p\|{_2^2}+\|\Delta^\varphi (v,p)\|{_{\mathcal{H}^2}^2}+\|\nabla^\varphi(p,v)\|_{\mathcal{H}^1}^2\Big)d\tau\nonumber\\
	&\leq \frac{1}{16(2\m+\l)}\int_{0}^{t}\int_{\mathcal{S}}\g p|\MZ^\a\Delta^\varphi p|^2d\mathcal{V}_\tau d\tau+C\Lambda(\frac{1}{c_0},\mathcal{Q}(t)+\mathcal{Q}_m)\int_{0}^{t}\v^2\|\nabla^2v\|^2_{\mathcal{H}^2}+\Lambda(Y_m(\tau))d\tau\nonumber\\
	&\leq \frac{1}{16(2\m+\l)}\int_{0}^{t}\int_{\mathcal{S}}\g p|\MZ^\a\Delta^\varphi p|^2d\mathcal{V}_\tau d\tau+\Lambda_{\i}(t)\int_0^t\Lambda(Y_m(\tau))d\tau.
	\end{align}
	Similarly, for $m\geq6$, it holds that
	\begin{eqnarray}\label{8.60}
	|J_4|\leq \frac{1}{16(2\m+\l)}\int_{0}^{t}\int_{\mathcal{S}}\g p|\MZ^\a\Delta^\varphi p|^2d\mathcal{V}_\tau d\tau+\Lambda_{\i}(t)\int_0^t\Lambda(Y_m(\tau))d\tau.
	\end{eqnarray}
	Next the expression
	\begin{equation*}
	 \mbox{div}^\varphi(\varrho(\partial{_t^\varphi}+v\cdot\nabla^\varphi)v)=\varrho\partial_t\mbox{div}^\varphi v+\varrho v_y\nabla_y\mbox{div}^\varphi v+\frac{1-z}{z}\varrho V_zZ_3\mbox{div}^\varphi v+\nabla^\varphi\varrho\partial{_t^\varphi}v+\nabla^\varphi(\varrho v)^t\nabla^\varphi v,
	\end{equation*}
	implies immediately, that for $m\geq6$, 
	\begin{equation}
	\int_{0}^{t}\|p \mbox{div}^\varphi(\varrho(\partial{_t^\varphi}+v\cdot\nabla^\varphi)v)\|{_{\mathcal{H}^2}^2}d\tau\leq \Lambda_{\i}(t)\int_0^t\Lambda(Y_m(\tau))d\tau.
	\end{equation}
	Thus, 
	\begin{align}\label{8.61}
	&|J_5|\leq\frac{1}{16(2\m+\l)}\int_{0}^{t}\int_{\mathcal{S}}\g p|\MZ^\a\Delta^\varphi p|^2d\mathcal{V}_\tau d\tau+C\int_{0}^{t}\|p \mbox{div}^\varphi(\varrho(\partial{_t^\varphi}+v\cdot\nabla^\varphi)v)\|{_{\mathcal{H}^2}^2}d\tau\nonumber\\
	&~~~~~\leq \frac{1}{16(2\m+\l)}\int_{0}^{t}\int_{\mathcal{S}}\g p|\MZ^\a\Delta^\varphi p|^2d\mathcal{V}_\tau d\tau+\Lambda_{\i}(t)\int_0^t\Lambda(Y_m(\tau))d\tau.
	\end{align}
	Substituting  \eqref{8.56}-\eqref{8.61} into \eqref{8.55} shows that 
	\begin{align}\label{20.60}
\v\|\Delta^\varphi p(\tau)\|^2_{\mathcal{H}^2}+\int_{0}^{t}\|\Delta^\varphi p\|^2_{\mathcal{H}^2}d\tau 
\leq \Lambda_0 \v\|\Delta^\varphi
	p(0)\|^2_{\mathcal{H}^2}+\Lambda_0\int_0^t\Lambda(Y_m(\tau))d\tau.
	\end{align}	
On the other hand,  it is easy to obtain that
\begin{align}\label{20.7}
&\|\Delta^\varphi p\|{_{\mathcal{H}^1}^2}\leq \Lambda_0\|(\Delta^\varphi p)(0)\|{_{\mathcal{H}^1}^2}+\Lambda_0\int_{0}^{t}\|\partial_t(\Delta^\varphi p)\|{_{\mathcal{H}^1}^2}d\tau\nonumber\\
&~~~~~~~~~~~~\leq \Lambda_0\|(\Delta^\varphi p)(0)\|{_{\mathcal{H}^1}^2}+\Lambda_0\int_{0}^{t}\|\Delta^\varphi p\|{_{\mathcal{H}^2}^2}d\tau+C_\d\Lambda_{\i}(t)\int_0^t\Lambda(Y_m(\tau))d\tau.
\end{align}
Then, combining \eqref{20.7} and \eqref{20.60} proves \eqref{20.6}. Therefore, the proof of this lemma is completed.$\hfill\Box$
	
\

We now turn to the most difficult part of $L^\i$-estimates: the control of $\|\nabla v(t)\|_{\mathcal{H}^{1,\i}}^2$ and $\v\|\partial_{zz}v(t)\|_\i^2$.  It follows from \eqref{7.0},  \eqref{7.3}, \eqref{7.7}, \eqref{8.2}, \eqref{8.4} and  \eqref{20.6} that, for $m\geq6$
\begin{align}\label{8.15}
	\|\nabla v(t)\|_{\mathcal{H}^{1,\i}}^2&\leq\Lambda(\frac{1}{c_0},|h(t)|_{\mathcal{H}^{3,\i}}^2)(\|\text{div}^\varphi v(t)\|_{\mathcal{H}^{1,\i}}^2+\|S_\mathbf{n}(t)\|_{\mathcal{H}^{1,\i}}^2+\|v(t)\|_{\mathcal{H}^{2,\i}}^2)\nonumber\\
	&\leq\|\Delta^\varphi p\|_{\mathcal{H}^1}^2+\Lambda(\frac{1}{c_0},\mathcal{Q}_5(t))(1+\|\nabla p\|_{\mathcal{H}^4}^2)+\Lambda_0\|S_\mathbf{n}(t)\|_{\mathcal{H}^{1,\i}}^2\nonumber\\
	&\leq \Lambda(Y_6(0))+  \Lambda_0t\sup_{0\leq\tau\leq t}\Lambda(Y_m(\tau))(1+\int_0^t\|\nabla p(\tau)\|_{\mathcal{H}^5}^2d\tau)+\Lambda_0\|S_\mathbf{n}(t)\|_{\mathcal{H}^{1,\i}}^2,
\end{align}
where the following elementary estimate have been used 
\begin{align}
&\mathcal{Q}_5(t)\leq \mathcal{Q}_5(0)+C\int_0^tY_6(\tau)d\tau,\label{4.5.31}\\
&\|\nabla p(t)\|_{\mathcal{H}^4}^2\leq \|\nabla p(0)\|_{\mathcal{H}^4}^2+\int_0^t\|\nabla p(\tau)\|_{\mathcal{H}^5}^2d\tau.\label{4.5.32}
\end{align}
Similarly, one can get that
\begin{align}\label{8.16}
	 \v\|\partial_{zz}v\|_\i^2&\leq\Lambda(\frac{1}{c_0},|h(t)|_{\mathcal{H}^{3,\i}}^2)\v(\|\nabla\text{div}^\varphi v(t)\|_\i^2+\|\nabla S_\mathbf{n}(t)\|_\i^2+\|\nabla v(t)\|_{\mathcal{H}^{1,\i}}^2)\nonumber\\
	 &\leq\Lambda_0\|\Delta^\varphi p\|_{\mathcal{H}^1}^2+\Lambda(\frac{1}{c_0},\mathcal{Q}_5(t)+\|\nabla v(t)\|_\i^2)\v\|\nabla p\|_{\mathcal{H}^4}^2+\Lambda_0\v(\|S_\mathbf{n}(t)\|_{\mathcal{H}^{1,\i}}^2+\|\nabla S_\mathbf{n}(t)\|_\i^2)\nonumber\\
	 &\leq \Lambda(Y_6(0))+\Lambda_0t\sup_{0\leq\tau\leq t}\Lambda(Y_m(\tau))+\Lambda_0\v(\|S_\mathbf{n}(t)\|_{\mathcal{H}^{1,\i}}^2+\|\nabla S_\mathbf{n}(t)\|_\i^2),
\end{align}
Therefore, it suffices to bound $\|S_\mathbf{n}(t)\|_{\mathcal{H}^{1,\i}}^2$ and  $\v\|\nabla S_\mathbf{n}(t)\|_\i^2$.

However, it seems difficult to bound $\|S_\mathbf{n}(t)\|_{\mathcal{H}^{1,\i}}^2$ and $\v\|\nabla S_\mathbf{n}(t)\|_\i^2$ by using the approach in \cite{Masmoudi-R-1} directly. This is because $(\nabla_y,0)^t(\nabla^\varphi\text{div}^\varphi v\cdot\mathbf{N})$ must appear in \eqref{7.10} for compressible flows, but we can not expect to bound uniformly $\|\nabla \text{div}^\varphi v\|_{\mathcal{H}^{1,\i}}^2$ and $\v\|\nabla\text{div}^\varphi v\|_{\mathcal{H}^{2,\i}}^2$. Thus, we propose to estimate $\omega_\mathbf{n}$, which eliminates $\nabla^\varphi\text{div}^\varphi v$ in its equation, 
instead of $S_\mathbf{n}$. Indeed, it follows from \eqref{8.17} that, for $m\geq6$,
\begin{align}\label{8.18}
	\|S_\mathbf{n}(t)\|_{\mathcal{H}^{1,\i}}^2&\leq \Lambda_0\|\omega_\mathbf{n}(t)\|_{\mathcal{H}^{1,\i}}^2+\Lambda(\frac{1}{c_0},|h(t)|_{\mathcal{H}^{3,\i}}^2)\|v(t)\|_{\mathcal{H}^{2,\i}}^2\nonumber\\
	&\leq \Lambda(\f1{c_0}Y_m(0))+\Lambda_0\|\omega_\mathbf{n}(t)\|_{\mathcal{H}^{1,\i}}^2+ C\int_0^t\Lambda(Y_m(\tau))d\tau,
\end{align}
and
\begin{align}\label{8.19}
	\v\|\nabla S_\mathbf{n}(t)\|_\i^2&\leq\Lambda_0\v\|\nabla\omega_\mathbf{n}\|_\i^2+\Lambda(\frac{1}{c_0},|h(t)|_{\mathcal{H}^{3,\i}}^2)(\v\|\nabla v(t)\|_{\mathcal{H}^{1,\i}}^2+\v\|v(t)\|_{\mathcal{H}^{1,\i}}^2)\nonumber\\
	&\leq\Lambda_0\v\|\nabla\omega_\mathbf{n}\|_\i^2+\|\Delta^\varphi p\|_{\mathcal{H}^1}^2+\Lambda(\frac{1}{c_0},\mathcal{Q}_5(t))(1+\v\|\nabla p\|_{\mathcal{H}^4}^2)+\Lambda_0\|\omega_\mathbf{n}\|_{\mathcal{H}^{1,\i}}^2\nonumber\\
	&\leq \Lambda_0\v\|\nabla\omega_\mathbf{n}\|_\i^2+\Lambda_0\|\omega_\mathbf{n}\|_{\mathcal{H}^{1,\i}}^2
	+\Lambda(\f1{c_0},Y_m(0))+\Lambda_0\int_0^t\Lambda(Y_m(\tau))d\tau.
\end{align}

Thus, it suffices to bound $\|\omega_\mathbf{n}\|_{\mathcal{H}^{1,\i}}^2$ and $\v\|\nabla\omega_\mathbf{n}(t)\|_\i^2$. However, the disadvantage of the term $\omega_\mathbf{n}$ is that it does not vanish on the boundary. But fortunately, the expression of $\omega_\mathbf{n}$ in \eqref{8.17} and the dynamic boundary condition \eqref{1.22} yield that, on the boundary
\begin{equation}\nonumber
	 \omega_\mathbf{n}^b=-2\Pi\{(\nabla_y,0)^t\partial_th-(\nabla^\varphi\mathbf{N})^tv^b\}.
\end{equation}
Thus, we  introduce 
\begin{align}\label{4.6.32}
\zeta_\mathbf{n}:=\omega_\mathbf{n}+2\Pi\{(\nabla_y,0)^t\partial_t\eta-(\nabla^\varphi\mathbf{N})^tv\}
\end{align}
which vanishes on the boundary. Then, for $m\geq6$, it follows from \eqref{8.2}, \eqref{4.5.31} and \eqref{4.5.32} that
\begin{align}\label{8.21}
	 \|\omega_\mathbf{n}\|_{\mathcal{H}^{1,\i}}^2&\leq\|\zeta_\mathbf{n}\|_{\mathcal{H}^{1,\i}}^2+\Lambda(\frac{1}{c_0},|h(t)|_{\mathcal{H}^{3,\i}}^2)(1+\|v(t)\|_{\mathcal{H}^{1,\i}}^2)\leq\|\zeta_\mathbf{n}\|_{\mathcal{H}^{1,\i}}^2+\Lambda(\frac{1}{c_0},\mathcal{Q}_5(t))\nonumber\\
	&\leq \|\zeta_\mathbf{n}\|_{\mathcal{W}^{1,\i}}^2+\Lambda(\frac{1}{c_0},Y_m(0))+\Lambda_0\int_0^t\Lambda(Y_m(\tau))d\tau,
\end{align}
and
\begin{align}\label{8.22}
	 \v\|\nabla\omega_\mathbf{n}\|_\i^2&\leq\v\|\nabla\zeta_\mathbf{n}\|_\i^2+\Lambda(\frac{1}{c_0},|h(t)|_{\mathcal{H}^{3,\i}}^2)\v(1+\|\nabla v(t)\|_\i^2)\nonumber\\
	 &\leq\v\|\nabla\zeta_\mathbf{n}\|_\i^2+\Lambda(\frac{1}{c_0},Y_m(0))+\Lambda_0\int_0^tP(Y_m(\tau))d\tau.
\end{align}

As a consequence, it remains to estimate $\|\zeta_\mathbf{n}(t)\|_{\mathcal{H}^{1,\i}}^2$ and $\v\|\nabla\zeta_\mathbf{n}(t)\|_\i^2$. Moreover, by using the following lemma, one can see that actually it suffices to derive these $L^\i$-estimates near the boundary, that is $\|\chi \zeta_\mathbf{n}(t)\|{_{\mathcal{W}^{1,\i}}^2}$ and $\v\|\chi\nabla\zeta_\mathbf{n}(t)\|{_\i^2}$ with $\chi$ compact supported and equal to 1 in the vicinity of $z=0$.
\begin{lemma}\label{l8.2}
	For any smooth cut-off function $\chi$ such that $\chi=0$ in a vicinity of $z=0$, it holds that for $m\geq k+2$,
	\begin{equation}\label{8.13}
		\|\chi f\|_{\mathcal{H}^{k,\i}}\lesssim\|f\|_{\mathcal{H}^{m}}.
	\end{equation}
\end{lemma}
This Lemma follows from the Sobolev embedding and the fact that the standard Sobolev norms are equivalent to the conormal ones away from the boundary. In fact, Lemma \ref{l8.2} and \eqref{8.2} imply, for $m\geq6$, that
\begin{equation}\label{cut}
	 \|\zeta_\mathbf{n}\|_{\mathcal{H}^{1,\i}}^2\lesssim\|\chi\zeta_\mathbf{n}\|_{\mathcal{H}^{1,\i}}^2+\|v\|_{\mathcal{H}^{2,\i}}^2\lesssim\|\chi\zeta_\mathbf{n}\|_{\mathcal{H}^{1,\i}}^2+ \Lambda(\f1{c_0},Y_m(0))+\int_0^t\Lambda(\f1{c_0},Y_m(\tau))d\tau.
\end{equation}

\

One can obtains the following lemma by the same argument in \cite{Masmoudi-R-1} but using $\|\chi\zeta_\mathbf{n}\|_{\mathcal{H}^{1,\i}}^2$ instead of $\|\chi S_\mathbf{n}\|_{\mathcal{H}^{1,\i}}^2$:
\begin{lemma}\label{lem4.29}
	For $m\geq 6$, it holds that
	\begin{equation}\label{4.5.43}
		\|\nabla v(t)\|{_{\mathcal{H}^{1,\i}}^2}\leq \Lambda(\f1{c_0},Y_m(0))+t\sup_{0\leq\tau\leq t}\Lambda(\f1{c_0},Y_m(\tau))\cdot\Big(1+\v^2|\MZ^mh(t)|_{\f12}^2+\int_0^t\|\nabla p(\tau)\|_{\mathcal{H}^5}^2d\tau\Big).
	\end{equation}
\end{lemma}
\noindent\textbf{Proof}. As indicated by \eqref{8.15}, \eqref{8.18}, \eqref{8.21} and \eqref{cut}, we  need only to estimate $\|\chi\zeta_\mathbf{n}\|_{\mathcal{H}^{1,\i}}^2$.
 The main difficulty of this estimate is to handle the commutator between $\MZ_i$ and $\Delta^\varphi$. We shall use a normal geodesic coordinate system in the vicinity of the boundary which gives a simpler expression of Laplacian\cite{Masmoudi-R-1}. Here, we can use this coordinate system because $L^\i$-estimate  does not require the  highest regularity of the boundary. Let us define a new parametrization at the vicinity of the boundary
\begin{align}
	 \Psi(t,\cdot):\mathcal{S}=&\mathbb{R}^2\times(-\i,0)\rightarrow D(t),\nonumber\\
	& (y,z)\mapsto(y,h(t,y))+z\mathbf{n}^b(t,y),
\end{align}
where $D(t)$ is defined in \eqref{1.27} and  $\mathbf{n}^b$ is the unit exterior normal $\mathbf{n}^b(t,y)=(-\partial_1h,-\partial_2h,1)/|\mathbf{N}|$. Note that
\begin{equation}\label{dpsi}
	D\Psi(t,\cdot)=\left(
	\begin{array}{ccc}
		1-z\partial_1(\frac{\partial_1h}{|\mathbf{N}|}) & -z\partial_2(\frac{\partial_1 h}{|\mathbf{N}}|) & -\frac{\partial_1 h}{|\mathbf{N}|}\\
		-z\partial_1(\frac{\partial_2 h}{|\mathbf{N}}|) & 1-z\partial_2(\frac{\partial_2 h}{|\mathbf{N}}|) & -\frac{\partial_2 h}{|\mathbf{N}|}\\
		\partial_1h+z\partial_1(\frac{1}{|\mathbf{N}|}) & \partial_2h+z\partial_2(\frac{1}{|\mathbf{N}|}|) & \frac{1}{|\mathbf{N}|})
	\end{array}\right),
\end{equation}
which is of the form $M_0+R$ with $|R|_\i\lesssim z|h|_{2,\i}$ and
\begin{equation*}
	M_0=\left(
	\begin{array}{ccc}
		1 & 0 & -\frac{\partial_1 h}{|\mathbf{N}|}\\
		0 & 1 & -\frac{\partial_2 h}{|\mathbf{N}|}\\
		\partial_1h & \partial_2h & \frac{1}{|\mathbf{N}|})
	\end{array}\right),
\end{equation*}
is invertible. This yields that $\Psi(t,\cdot)$ is a diffeomorphism from $\mathbb{R}^2\times(-\d, 0)$ to a vicinity of $\partial\Om_t$ for some $\d>0$ which depends only on $c_0>0$. By this parametrization, the scalar product in $\Om_t$ induces a Riemman metric which is given by
\begin{equation*}
	g(y,z)=
	\left(
	\begin{array}{cc}
		\tilde{g}(y,z) & 0\\
		0 & 1
	\end{array}\right),
\end{equation*}
and the Laplacian in this coordinate system is of the form:
\begin{equation*}
	 \Delta_gf=\partial_{zz}f+\frac{1}{2}\partial_z(\ln|g|)\partial_zf+\Delta_{\tilde{g}}f,
\end{equation*}
where $|g|$ denotes the determinant of the matrix $g$ and $\Delta_{\tilde{g}}f$ is given by
\begin{equation*}
	\Delta_{\tilde{g}}f=\frac{1}{|\tilde{g}|^\frac{1}{2}}\sum_{1\leq i,j\leq 2}\partial_{y_i}(\tilde{g}^{ij}|g|^\frac{1}{2}\partial_{y_j}f),
\end{equation*}
which involves only the tangential derivatives.

In order to  use this normal geodesic coordinate system, one first localizes the equation for $\zeta_\mathbf{n}$. Set
\begin{equation}\label{4.6.40}
	\omega^\chi=\chi(z)\omega
\end{equation}
where $\chi(z)$ is given as $\chi(z)=\tilde\k(\frac{z}{\d(c_0)})$ with $\tilde\k(z)$ be a smooth function of compact support such that $\tilde\k(z)\in[0,1]$ taking the value 1 in the vicinity of $z=0$, and $\d(c_0)$ to be determined later. Note that this choice implies that $|\chi^k(z)|\leq\Lambda_0$.

Applying $\nabla^\varphi\times$ to $\eqref{1.22}_2$ yields that in $\mathcal{S}$,
\begin{equation}\label{8.34}
	\varrho\partial{_t^\varphi}\omega+\varrho v\cdot\nabla^\varphi\omega-\m\v\Delta^\varphi\omega=F_\omega,
\end{equation}
where
\begin{equation}\label{8.35}
	 F_\omega=-\nabla^\varphi\varrho\times\partial{_t^\varphi}v-\nabla^\varphi\varrho\times((v\cdot\nabla^\varphi)v)+\varrho\omega\cdot\nabla^\varphi v-\varrho\mbox{div}^\varphi v\omega.
\end{equation}
Furthermore, it follows from \eqref{4.6.40} and \eqref{8.34} that
\begin{equation}
	\varrho\partial{_t^\varphi}\omega^\chi+\varrho v\cdot\nabla^\varphi\omega^\chi-\m\v\Delta^\varphi\omega^\chi=F_{\omega^\chi}
\end{equation}
where $F_{\omega^\chi}:=F^\chi_1+\chi F_\omega$ with
\begin{equation}\label{8.37}
	 F^\chi_1=\varrho(V_z\partial_z\chi)\omega-\m\v\nabla^\varphi\chi\cdot\nabla^\varphi\omega-\m\v\Delta^\varphi\chi\omega.
\end{equation}
Since $F^\chi_1$ is supported away from the boundary, it follows  from \eqref{4.23}, \eqref{4.24} and \eqref{8.13} that
\begin{equation}\label{8.38}
	 \|F^\chi_1\|_{\mathcal{H}^{1,\i}}^2\lesssim\Lambda(\frac{1}{c_0},|h(t)|_{\mathcal{H}^{1,\i}}^2+\|(v,p)(t)\|_{\mathcal{H}^{1,\i}}^2)\|v(t)\|_{\mathcal{H}^5}^2\lesssim\Lambda(\frac{1}{c_0},\mathcal{Q}_5(t)).
\end{equation}
Next, we define $\omega^\Psi:=\omega^\chi(t,\Phi^{-1}(t,\cdot)\circ\Psi)$. This change of variable is well-defined if we choose $\d(c_0)$ small enough such that $D\Psi$ is invertible. Then, $\omega^\Psi$ solves in $\mathcal{S}$ the convection diffusion equations
\begin{equation}\label{8.39}
	\varrho\partial_t\omega^\Psi+\varrho \mathbf{b}\cdot\nabla\omega^\Psi-\m\v(\partial_{zz}+\frac{1}{2}\partial_z(\ln|g|)\partial_z)\omega^\Psi=F_{\omega^\chi}(t,\Phi^{-1}\circ\Psi)+\m\v\Delta_{\tilde{g}}\omega^\Psi,
\end{equation}
where
\begin{equation}\label{b}
	\mathbf{b}=(D\Psi)^{-1}(v(t,\Phi^{-1}\circ\Psi)-\partial_t\Psi).
\end{equation}
Set
\begin{equation}
	 \omega_\mathbf{n}^\Psi(t,y,z)=\Pi^b(t,y)(\omega^\Psi(t,y,z)\times\mathbf{N}^b(t,y)),
\end{equation}
with $\mathbf{N}^b=(-\partial_1h,-\partial_2h,1)$ and $\Pi^b=\text{Id}-\mathbf{n}^b\otimes\mathbf{n}^b$.
Note that $\Pi^b$ and $\mathbf{n}^b$ are independent of $z$. This yields that $\omega_\mathbf{n}^\Psi$ solves
\begin{equation}\label{omen}
	\varrho\partial_t\omega_\mathbf{n}^\Psi+\varrho \mathbf{b}\cdot\nabla\omega_\mathbf{n}^\Psi-\m\v(\partial_{zz}+\frac{1}{2}\partial_z(\ln|g|)\partial_z)\omega_\mathbf{n}^\Psi=F_\omega^\Psi,
\end{equation}
where $F_\omega^\Psi$ is given by
\begin{equation}\label{8.43}
	 F_\omega^\Psi=\Pi^bF_{\omega^\chi}\times\mathbf{N}^b+F_\omega^{\Psi,1}+F_\omega^{\Psi,2},
\end{equation}
with
\begin{align}\label{8.44}
	& F_\omega^{\Psi,1}=\varrho(\partial_t\Pi^b+\mathbf{b}_y\cdot\nabla\Pi^b)\omega^\Psi\times\mathbf{N}^b+\varrho\Pi\{\omega^\Psi\cdot(\partial_t+\mathbf{b}_y\cdot\nabla_y)\mathbf{N}^b\},\\\label{8.45}
	& F_\omega^{\Psi,2}=-\m\v\Pi^b(\Delta_{\tilde{g}}\omega^\Psi\times\mathbf{N}^b).
\end{align}

Define
\begin{equation}
	v^\Psi=(\chi v)(t,\Phi^{-1}\circ\Psi).
\end{equation}
Note that $\nabla^\varphi\mathbf{N}^b=\nabla \mathbf{N}^b$ is independent of $z$. Thus $\Pi^b\{(\nabla^\varphi\mathbf{N}^b)^tv^\Psi\}$ solves
\begin{equation}\label{vn}
	\varrho\partial_t\Pi^b\{(\nabla^\varphi\mathbf{N}^b)^tv^\Psi\}+\varrho \mathbf{b}\cdot\nabla\Pi^b\{(\nabla^\varphi\mathbf{N}^b)^tv^\Psi\}-\m\v(\partial_{zz}+\frac{1}{2}\partial_z(\ln|g|)\partial_z)\Pi^b\{(\nabla^\varphi\mathbf{N}^b)^tv^\Psi\}=F_v^\Psi,
\end{equation}
where
\begin{equation}\label{fv1}
	F_v^\Psi=\Pi^b\{(\nabla^\varphi\mathbf{N}^b)^t(\chi F_v+F^\chi_2)\}+F_v^{\Psi,1}+F_v^{\Psi,2},
\end{equation}
with
\begin{align}\label{fv2}
	& F_v=-\nabla^\varphi p+(\m+\l)\v\nabla^\varphi\text{div}^\varphi v,\\\label{fv3}
	& F^\chi_2=\varrho(V_z\partial_z\chi)v-\m\v\nabla\chi\cdot\nabla^\varphi v-\m\v\Delta^\varphi\chi v,\\\label{fv4}
	& F_v^{\Psi,1}=\varrho(\partial_t\Pi^b+\mathbf{b}_y\cdot\nabla\Pi^b)(\nabla^\varphi\mathbf{N}^b)^tv^\Psi+\varrho\Pi\{(\partial_t+\mathbf{b}_y\cdot\nabla_y)(\nabla^\varphi\mathbf{N})^tv^\Psi\},\\\label{fv5}
	& F_v^{\Psi,2}=-\m\v\Pi^b((\nabla^\varphi\mathbf{N}^b)^t\Delta_{\tilde{g}}v^\Psi).
\end{align}

\

Set
\begin{equation}\label{4.6.60}
	 \zeta_\mathbf{n}^\Psi=\omega_\mathbf{n}^\Psi-2\Pi^b\{(\nabla^\varphi\mathbf{N}^b)^tv^\Psi\}+2\Pi^b((\nabla_y,0)^t\partial_t\eta).
\end{equation}
It follows from \eqref{omen} and \eqref{vn} that $\zeta_\mathbf{n}^\Psi$ solves in $\mathcal{S}$
\begin{equation}\label{8.54}
	\varrho\partial_t\zeta_\mathbf{n}^\Psi+\varrho \mathbf{b}\cdot\nabla\zeta_\mathbf{n}^\Psi-\m\v(\partial_{zz}+\frac{1}{2}\partial_z(\ln|g|)\partial_z)\zeta_\mathbf{n}^\Psi=F_\omega^\Psi+F_v^\Psi+F_\eta^\Psi,
\end{equation}
where
\begin{equation}\label{fe}
	 F_\eta^\Psi=\varrho\partial_t(\Pi^b(\nabla_y,0)^t\partial_t\eta^\Psi)+\varrho \mathbf{b}\cdot\nabla(\Pi^b(\nabla_y,0)^t\partial_t\eta^\Psi)-\m\v(\partial_{zz}+\frac{1}{2}\partial_z(\ln|g|)\partial_z)(\Pi^b(\nabla_y,0)^t\partial_t\eta^\Psi),
\end{equation}
with $\eta^\Psi$ defined by $\chi\eta(t,\Phi^{-1}\circ\Psi)$.
In order to eliminate the term $\frac{1}{2}\partial_z(\ln|g|)\partial_z\zeta_\mathbf{n}^\Psi$, we define
\begin{equation*}
	 \zeta_\mathbf{n}^\Psi=\dfrac{1}{|g|^\frac{1}{4}}\tilde{\zeta}_\mathbf{n}^\Psi.
\end{equation*}
Plugging $\tilde{\zeta}_\mathbf{n}^\Psi$ into \eqref{8.54}, one obtains that
\begin{align}\label{zeta2}
	 &\varrho\partial_t\tilde{\zeta}_\mathbf{n}^\Psi+\varrho\mathbf{b}\cdot\nabla\tilde{\zeta}_\mathbf{n}^\Psi-\m\v\partial_{zz}\tilde{\zeta}_\mathbf{n}^\Psi=|g|^\frac{1}{4}(F_\omega^\Psi+F_v^\Psi+F_\eta^\Psi+F_g^\Psi)=:S,
\end{align}
where $F_g^\Psi$
\begin{equation}\label{fg}
	 F_g^\Psi=-(\varrho\partial_t(|g|^{-\frac{1}{4}})+\varrho\cdot\mathbf{b}\nabla(|g|^{-\frac{1}{4}})-\m\v\partial_{zz}(|g|^{-\frac{1}{4}})-2\m\v\frac{(\partial_z|g|^{-\frac{1}{4}})^2}{|g|^{-\frac{1}{4}}})\tilde{\zeta}_\mathbf{n}^\Psi.
\end{equation}
Moreover, the boundary condition reads as
\begin{equation}\label{8.42}
	\tilde{\zeta}_\mathbf{n}^\Psi=\zeta_\mathbf{n}^\Psi=\zeta_\mathbf{n}=0.
\end{equation}
Since Sobolev conormal spaces are invariant by diffeomorphisms which preserve the boundary (see Lemma 9.5 in \cite{Masmoudi-R-1}) and $|\Pi-\Pi^b|+|\mathbf{N}-\mathbf{N}^b|=\mathcal{O}(z)$ in the vicinity of the boundary, one  obtains that
\begin{align}\label{zeta}
	 \|\zeta_\mathbf{n}\|_{\mathcal{H}^{1,\i}}^2&\leq\Lambda_0(\|\zeta_\mathbf{n}^\Psi\|_{\mathcal{H}^{1,\i}}^2+|h(t)|_{\mathcal{H}^{3,\i}}^2+\|v\|_{\mathcal{H}^{2,\i}}^2)\leq\Lambda_0\|\tilde{\zeta}_\mathbf{n}^\Psi\|_{\mathcal{H}^{1,\i}}^2+\Lambda(\frac{1}{c_0},\mathcal{Q}_5(t))\nonumber\\
	 &\leq\Lambda_0\|\tilde{\zeta}_\mathbf{n}^\Psi\|_{\mathcal{H}^{1,\i}}^2+\Lambda(Y_m(0))+\int_0^t\Lambda(Y_m(\tau))d\tau
\end{align}
which implies that it suffices to estimate $\|\tilde{\zeta}_\mathbf{n}^\Psi\|_{\mathcal{H}^{1,\i}}^2$. Now, we will estimate $\|\tilde{\zeta}_\mathbf{n}^\Psi\|_{\mathcal{H}^{1,\i}}^2$ by taking advantage of \eqref{zeta2} and \eqref{8.42}. First, since
\begin{equation}
	 b_3(t,y,0)=\sum_{j=1}^{3}(D\Psi)_{3j}^{-1}(v_j-\partial_t\Psi_j)(t,y,0)=0,
\end{equation}
then  \eqref{zeta2} can be rewritten as
\begin{align}
	 &\quad\varrho(t,y,0)(\partial_t+b_1(t,y,0)\partial_1+b_2(t,y,0)\partial_2+z\partial_zb_3(t,y,0))\tilde{\zeta}_\mathbf{n}^\Psi-\m\v\partial_{zz}\tilde{\zeta}_\mathbf{n}^\Psi\nonumber
	\\
	 &=S-[\varrho(t,y,z)-\varrho(t,y,0)]\partial_t\tilde{\zeta}_\mathbf{n}^\Psi-\sum_{i=1,2}\varrho[b_i(t,y,z)-b_i(t,y,0)]\partial_i\tilde{\zeta}_\mathbf{n}^\Psi\\
	 &~-\varrho(t,y,0)[b(t,y,z)-z\partial_zb_3(t,y,0)]\partial_z\tilde{\zeta}_\mathbf{n}^\Psi-[\varrho(t,y,z)-\varrho(t,y,0)]z\partial_zb_3(t,y,0)\partial_z\tilde{\zeta}_\mathbf{n}^\Psi=:G.\nonumber
\end{align}
It follows from Lemma \ref{lemA.2} in the Appendix that, for $m\geq6$,
\begin{align}\label{8.46}
	 \|\tilde{\zeta}_\mathbf{n}^\Psi(t)\|_{\mathcal{H}^{1,\i}}&\lesssim\|\tilde{\zeta}_\mathbf{n}^\Psi(0)\|_{\mathcal{H}^{1,\i}}+\int_{0}^{t}\|\varrho^{-1}\|_\i\|G\|_{\mathcal{H}^{1,\i}}d\tau\nonumber\\
	 &\quad+\int_{0}^{t}(1+\|\varrho^{-1}\|_\i)(1+\sum_{i=0}^{2}\|\MZ_i(\varrho,b,\partial_zb_3)\|{_\i^2})\|\tilde{\zeta}_\mathbf{n}^\Psi\|_{\mathcal{H}^{1,\i}}d\tau\\
	 &\lesssim\|\tilde{\zeta}_\mathbf{n}^\Psi(0)\|_{\mathcal{H}^{1,\i}}+\Lambda_0\int_{0}^{t}\|G\|_{\mathcal{H}^{1,\i}}d\tau+\int_{0}^{t}\Lambda(\frac{1}{c_0},Y_m(\tau))d\tau.\nonumber
\end{align}
It remains to estimate $\int_{0}^{t}\|G\|_{\mathcal{H}^{1,\i}}d\tau$. First, it follows from \eqref{8.35}, \eqref{8.38} and \eqref{8.43}-\eqref{8.45}  that
\begin{align}\label{foo}
	 \||g|^\frac{1}{4}F_\omega^\Psi\|_{\mathcal{H}^{1,\i}}^2&\leq\Lambda(\frac{1}{c_0},\mathcal{Q}_5(t))((1+\|\nabla p\|_{\mathcal{H}^{1,\i}}^2)\|\nabla v\|_{\mathcal{H}^{1,\i}}^2+\|\nabla v\|_{\mathcal{H}^{1,\i}}^4+\v^2\|\nabla v\|_{\mathcal{H}^{3,\i}}^2)\nonumber\\
	&\leq\Lambda(\frac{1}{c_0},Y_m(t))(1+\v^2\|\nabla^2 v\|_{\mathcal{H}^{4}}\|\nabla v\|_{\mathcal{H}^{5}}),~m\geq6,
\end{align}
where Lemma \ref{l8.1} has been used several times. Next, in a similar way, one can get from \eqref{fv1}-\eqref{fv5} and Lemma  \ref{l8.1} that
\begin{align}\label{fvv}
	 &\||g|^\frac{1}{4}F_v^\Psi\|_{\mathcal{H}^{1,\i}}^2\leq\Lambda(\frac{1}{c_0},\mathcal{Q}_5(t))[1+\|\nabla p\|_{\mathcal{H}^{1,\i}}^2+\|\nabla v\|_{\mathcal{H}^{1,\i}}^2+\v^2\|\nabla \text{div}^\varphi v\|_{\mathcal{H}^{1,\i}}^2+\v^2\|\nabla v\|_{\mathcal{H}^{3,\i}}^2]\nonumber\\
	&\leq\v^2\|\Delta^\varphi p\|_{\mathcal{H}^2}^2+\Lambda(\frac{1}{c_0},\mathcal{Q}_6(t)+\|\Delta^\varphi p\|_{\mathcal{H}^1}^2+\|\nabla v\|_{\mathcal{H}^{1,\i}}^2)(\v^2\|\nabla p\|_{\mathcal{H}^5}^2+\v^2\|\nabla v\|_{\mathcal{H}^{3,\i}}^2)\nonumber\\
		&\leq\Lambda(\frac{1}{c_0},Y_m(t))(1+\v^2\|\nabla^2 v\|_{\mathcal{H}^{4}}\|\nabla v\|_{\mathcal{H}^{5}}),~m\geq6.
\end{align}
It follows from \eqref{fe}, \eqref{fg} and Lemma \ref{l8.1}  that
\begin{equation}\label{fgg}
	 \||g|^\frac{1}{4}F_\eta^\Psi\|_{\mathcal{H}^{1,\i}}^2+\||g|^\frac{1}{4}F_g^\Psi\|_{\mathcal{H}^{1,\i}}^2\leq\Lambda(\frac{1}{c_0},Y_m(t))(1+\v^2|\MZ^mh(t)|_{\f12}^2),~m\geq6.
\end{equation}
Therefore, combining \eqref{foo}-\eqref{fgg} leads to that
\begin{align}\label{8.67}
	 &\|S\|_{\mathcal{H}^{1,\i}}^2\leq\Lambda(\frac{1}{c_0},Y_m(t))(1+\v^2|\MZ^mh(t)|_{\f12}^2+\v^2\|\nabla^2 v\|_{\mathcal{H}^{4}}\|\nabla v\|_{\mathcal{H}^{5}}),~m\geq6.
\end{align}
Next, using the Taylor's formula and the fact that $\tilde{\zeta}$ is compactly supported in $z$, one obtains that
\begin{align}\label{8.64}
	 &\quad\|(\varrho(t,y,z)-\varrho(t,y,0))\partial_t\tilde{\zeta}_\mathbf{n}^\Psi\|{_{\mathcal{H}^{1,\i}}^2}\nonumber\\
	 &\lesssim\Lambda_0\|\tilde{\zeta}_\mathbf{n}^\Psi\|{_{\mathcal{H}^{1,\i}}^2}+\|\MZ(\varrho(t,y,z)-\varrho(t,y,0))\|{_\i^2}\|\partial_t\tilde{\zeta}_\mathbf{n}^\Psi\|{_\i^2}+\|(\varrho(t,y,z)-\varrho(t,y,0))\MZ\partial_t\tilde{\zeta}_\mathbf{n}^\Psi\|{_\i^2}\nonumber\\
	&\leq \Lambda_0\|\tilde{\zeta}_\mathbf{n}^\Psi\|{_{\mathcal{H}^{1,\i}}^2}+\|\MZ\varrho\|{_\i^2}\|\tilde{\zeta}_\mathbf{n}^\Psi\|{_\i^2}+\|\partial_z\varrho\|{_\i^2}\|\frac{z}{1-z}\MZ\partial_t\tilde{\zeta}_\mathbf{n}^\Psi\|{_\i^2}.
\end{align}
Due to \eqref{2.9}, the following inequality holds for $|\a|\geq 0$
\begin{equation}\label{8.65}
	 \|\frac{z}{1-z}\MZ^\a\tilde{\zeta}_\mathbf{n}^\Psi\|{_\i^2}\lesssim\|\nabla(\frac{z}{1-z}\MZ^\a\tilde{\zeta}_\mathbf{n}^\Psi)\|_1\|\frac{z}{1-z}\MZ^\a\tilde{\zeta}_\mathbf{n}^\Psi\|_2\lesssim\|\MZ^\a\tilde{\zeta}_\mathbf{n}^\Psi\|{_2^2}.
\end{equation}
Therefore, plugging \eqref{8.65} with $|\a|=2$ into \eqref{8.64}, one gets from \eqref{8.10} that
\begin{align}\label{8.70}
	 \|(\varrho(t,y,z)-\varrho(t,y,0))\partial_t\tilde{\zeta}_\mathbf{n}^\Psi\|{_{\mathcal{H}^{1,\i}}} 
	&\leq\Lambda(\frac{1}{c_0},Y_m(t)),~m\geq6.
\end{align}
Similarly, 
\begin{equation}\label{8.71}
	 \sum_{i=1,2}\|\varrho(b_i(t,y,z)-b_i(t,y,0))\partial_i\tilde{\zeta}_\mathbf{n}^\Psi\|{_{\mathcal{H}^{1,\i}}^2}\leq\Lambda(\frac{1}{c_0},Y_m(t)),~m\geq6.
\end{equation}
Since $b_3(t,y,0)=0$, it holds that
\begin{align}\label{8.68}
	 &\quad\|\varrho(t,y,0)(b_3(t,y,z)-z\partial_zb_3(t,y,0))\partial_z\tilde{\zeta}_\mathbf{n}^\Psi\|{_{\mathcal{H}^{1,\i}}^2}\nonumber\\
	&\leq\|\varrho\|{_{\mathcal{W}^{1,\i}}^2}(\|\nabla b_3\|{_\i^2}\|\tilde{\zeta}_\mathbf{n}^\Psi\|{_{\mathcal{H}^{1,\i}}^2}+\|\MZ(b_3(t,y,z)-z\partial_zb_3(t,y,0))\partial_z\tilde{\zeta}_\mathbf{n}^\Psi\|{_\i^2}\nonumber\\
	 &\qquad+\|(b_3(t,y,z)-z\partial_zb_3(t,y,0))\MZ\partial_z\tilde{\zeta}_\mathbf{n}^\Psi\|{_\i^2})\\
	&\leq\Lambda(\frac{1}{c_0},\mathcal{Q}_5(t))(1+\|\nabla v\|^2_{\mathcal{H}^{1,\i}})\|\tilde{\zeta}_\mathbf{n}^\Psi\|{_{\mathcal{H}^{1,\i}}^2}+\Lambda(\frac{1}{c_0},\mathcal{Q}_5(t))\|\partial_{zz}b_3\|{_\i^2}\|\frac{z^2}{(1-z)^2}\partial_z\MZ\tilde{\zeta}_\mathbf{n}^\Psi\|_\i^2\nonumber\\
	&\leq\Lambda(\frac{1}{c_0},\mathcal{Q}_5(t))(1+\|\nabla v\|^2_{\mathcal{H}^{1,\i}})\|\tilde{\zeta}_\mathbf{n}^\Psi\|{_{\mathcal{H}^{1,\i}}^2}\nonumber\\
	 &\qquad+\Lambda(\frac{1}{c_0},\mathcal{Q}_5(t))\|\partial_{zz}b_3\|{_\i^2}\|\frac{z^2}{(1-z)^2}\partial_z\MZ\tilde{\zeta}_\mathbf{n}^\Psi\|_{\mathcal{H}^2}\|\partial_z(\frac{z^2}{(1-z)^2}\partial_z\MZ\tilde{\zeta}_\mathbf{n}^\Psi)\|_{\mathcal{H}^1}\nonumber\\
	&\leq\Lambda(\frac{1}{c_0},Y_m(t))+\Lambda(\frac{1}{c_0},Y_m(t))\|\partial_{zz}b_3\|{_\i^2}\|\tilde{\zeta}_\mathbf{n}^\Psi\|_{\mathcal{H}^4}^2,~\mbox{for}~m\geq6.\nonumber
\end{align}
It follows from \eqref{dpsi} and \eqref{b} that
\begin{align}
	&\|\partial_{zz}b_3\|_\i^2\leq \|\partial_{zz}v\cdot\mathbf{N}\|_\i^2+\Lambda(\frac{1}{c_0},|h|_{\mathcal{H}^{4,\i}}^2)\nonumber\\
	&\lesssim\Lambda_0\|\partial_z\mbox{div}^\varphi v\|_\i^2+\Lambda_0\|\nabla v(t)\|_\i^2+\Lambda(\frac{1}{c_0},|h|_{\mathcal{H}^{4,\i}}^2)\leq \Lambda(\f1{c_0},Y_m(t)),~m\geq6,
\end{align}
where one has used that
\begin{align}
	 \partial_{zz}(v\cdot\mathbf{N})&=\partial_z(\partial_zv\cdot\mathbf{N})+\partial_z(v\cdot\partial_z\mathbf{N})\nonumber\\
	 &=\partial_z(\partial_z\varphi(\mbox{div}^\varphi v-\partial_1v_1-\partial_2v_2))+\partial_z(v\cdot\partial_z\mathbf{N})\nonumber
\end{align}
Thus, 
\begin{equation}\label{8.74}
	 \|\varrho(t,y,0)(b_3(t,y,z)-z\partial_zb_3(t,y,0))\partial_z\tilde{\zeta}_\mathbf{n}^\Psi\|{_{\mathcal{H}^{1,\i}}^2}
	\leq\Lambda(\f1{c_0},Y_m(t)),~m\geq6.
\end{equation}
Similarly, note that $(\varrho(t,y,z)-\varrho(t,y,0))z=\mathcal{O}(z^2)$ in a vicinity of the boundary, one gets that
\begin{equation}\label{8.75}
	 \|(\varrho(t,y,z)-\varrho(t,y,0))z\partial_zb_3(t,y,0)\partial_z\tilde{\zeta}_\mathbf{n}^\Psi\|{_{\mathcal{H}^{1,\i}}^2}\leq\Lambda(\frac{1}{c_0},Y_m(t)),~m\geq6.
\end{equation}
Thanks to \eqref{8.67}, \eqref{8.70}, \eqref{8.71}, \eqref{8.74} and \eqref{8.75}, one has that
\begin{align}\label{8.76}
	 &\|G\|{_{\mathcal{H}^{1,\i}}^2}\leq\Lambda(\frac{1}{c_0},Y_m(t))(1+\v^2|\MZ^mh(t)|_{\f12}^2+\v^2\|\nabla^2 v\|^2_{\mathcal{H}^{4}}),~m\geq6.
\end{align}
Substituting \eqref{8.76} into \eqref{8.46}  and using \eqref{6.37-2},  one gets that
\begin{align}\label{4.5.93}
\|\tilde{\zeta}_\mathbf{n}^\Psi(t)\|^2_{\mathcal{H}^{1,\i}}\leq\|\tilde{\zeta}_\mathbf{n}^\Psi(0)\|^2_{\mathcal{H}^{1,\i}}
+\int_{0}^{t}\Lambda(\frac{1}{c_0},Y_m(t))(1+\v^2|\MZ^mh(t)|_{\f12}^2)d\tau,~m\geq6.
\end{align}
Then, \eqref{4.5.43} follows from \eqref{8.15}, \eqref{8.18}, \eqref{8.21}, \eqref{cut}, \eqref{zeta} and \eqref{4.5.93}.
Thus, the proof of this lemma is completed.$\hfill\Box$

\begin{lemma}\label{lem21.1}
	For $m\geq6$, it holds that
	\begin{align}\label{21.1}
		&\v\|\partial_{zz} v(t)\|_\i^2\leq t\sup_{0\leq\tau\leq t}\Lambda\left(\f1{c_0},Y_m(\tau)+\v^2|\MZ^mh(\tau)|_{\f12}^2\right)\cdot(1+\int_0^t\|\nabla p(\tau)\|_{\mathcal{H}^5}^2d\tau)\nonumber\\
		&~~~~~~~~~~~~~~~~~~~~~~~~+\Lambda(\frac{1}{c_0},Y_m(0)).
	\end{align}
\end{lemma}
\noindent\textbf{Proof}. We shall reduce the problem to the estimate of $\v\|\partial_z\tilde{\zeta}_\mathbf{n}^\Psi\|_{\i}$. First, it follows from \eqref{8.16}, \eqref{8.19}, \eqref{8.22} and \eqref{4.5.43} that
\begin{align}\label{zzv}
	\v\|\partial_{zz} v\|_\i^2&\leq t\sup_{0\leq\tau\leq t}\Lambda(\f1{c_0},Y_m(\tau))\Big(1+\v^2|\MZ^mh(t)|_{\f12}^2+\int_0^t\|\nabla p(\tau)\|_{\mathcal{H}^5}^2d\tau\Big)\nonumber\\
	 &~~~~~~~~~+\Lambda(\f1{c_0},Y_m(0))+\Lambda_0\v\|\nabla\zeta_\mathbf{n}\|_\i^2,m\geq6.
\end{align}
Note that $|\Pi-\Pi^b|+|\mathbf{N}-\mathbf{N}^b|=\mathcal{O}(z)$ in the vicinity of the boundary, one gets  from the definitions of $\tilde{\zeta}_\mathbf{n}^\Psi$ and $\zeta_\mathbf{n}$ that
\begin{align}\label{4.6.86}
	 &\v\|\partial_z\zeta_\mathbf{n}(t)\|_\i^2\leq\Lambda_0(\v\|\partial_z\tilde{\zeta}_\mathbf{n}^\Psi\|_\i^2+\|(v,\nabla v)\|_{\mathcal{H}^{1,\i}}^2)\nonumber\\
	&\leq \Lambda(\f1{c_0},Y_m(0))+\Lambda_0 t\sup_{0\leq\tau\leq t}\Lambda(\f1{c_0},Y_m(\tau))\Big(1+\v^2|\MZ^mh(t)|_{\f12}^2+\int_0^t\|\nabla p(\tau)\|_{\mathcal{H}^5}^2d\tau\Big)\nonumber\\
	 &~~~~~~~~~+\Lambda_0\v\|\partial_z\tilde{\zeta}_\mathbf{n}^\Psi\|_{\i},m\geq6.
\end{align}
It remains to estimate $\v\|\partial_z\tilde{\zeta}_\mathbf{n}^\Psi\|_\i^2$. Since $\tilde{\zeta}_\mathbf{n}^\Psi$ solves \eqref{zeta2} in $\mathcal{S}$ with the homogeneous Dirichlet boundary condition, one can use the one-dimensional heat kernel of $z>0$
$$G(t,z,z')=\dfrac{1}{\sqrt{4\pi\m\v t}}\Big(\exp({-\dfrac{(z-z')^2}{4\m\v t}})-\exp({-\dfrac{(z+z')^2}{4\m\v t}})\Big).$$
to obtain
\begin{align}
	 \sqrt{\v}\partial_z\tilde{\zeta}_\mathbf{n}^\Psi=&\sqrt{\v}\int_{0}^{+\i}\partial_zG(t,z,z')\tilde{\zeta}_\mathbf{n}^\Psi(0)dz'\nonumber\\
	 &+\int_{0}^{t}\int_{0}^{+\i}\sqrt{\v}\partial_zG(t-\tau,z,z')(S-\rho\mathbf{b}\cdot\nabla\tilde{\zeta}_\mathbf{n}^\Psi-(\rho-1)\partial_t\tilde{\zeta}_\mathbf{n}^\Psi)dz'd\tau.
\end{align}
It follows from $\tilde{\zeta}_\mathbf{n}^\Psi(0)=0$ and  integration by parts that
\begin{equation}\label{8.82}
	 \sqrt{\v}\|\partial_z\tilde{\zeta}_\mathbf{n}^\Psi\|_\i\leq\sqrt{\v}\|\partial_z\tilde{\zeta}_\mathbf{n}^\Psi(0)\|_\i+\int_{0}^{t}\dfrac{1}{\sqrt{t-\tau}}(\|S\|_\i+\|\varrho\mathbf{b}\cdot\nabla \tilde{\zeta}_\mathbf{n}^\Psi\|_\i+\|(\varrho-1)\partial_t\tilde{\zeta}_\mathbf{n}^\Psi\|_\i).
\end{equation}
It then follows from \eqref{8.67} that
\begin{align}\label{8.83}
	 &\left(\int_0^t\dfrac{1}{\sqrt{t-\tau}}\|S\|_{\mathcal{H}^{1,\i}}d\tau\right)^2\nonumber\\
	 &\leq\left(\int_0^t\f{1}{\sqrt{t-\tau}}\Lambda(\frac{1}{c_0},Y_m(t))\Big[1+\v|\MZ^mh(t)|_{\f12}+\sqrt\v\|\nabla^2 v\|^{\f12}_{\mathcal{H}^{4}}\sqrt{\v}\|\nabla v\|^{\f12}_{\mathcal{H}^{5}}\Big] d\tau\right)^2\nonumber\\
	&\leq t\Lambda\left(\frac{1}{c_0},Y_m(t)+\v^2|\MZ^mh(t)|^2_{\f12}\right)+\left(\int_0^t\v^2\|\nabla^2 v\|^2_{\mathcal{H}^{4}} \right)^{\f12}t^{\f12}\Lambda(\frac{1}{c_0},Y_m(t))\nonumber\\
	&\leq t\Lambda\left(\frac{1}{c_0},Y_m(t)+\v^2|\MZ^mh(t)|^2_{\f12}\right),
\end{align}
where one has used \eqref{6.37-2} in the last inequality.
Since $b_3(t,y,0)=0$, so 
\begin{equation}\label{8.84}
	 \|\varrho\mathbf{b}\cdot\nabla\tilde{\zeta}_\mathbf{n}^\Psi\|_\i^2\leq\Lambda_0(\|b\|_\i^2+\|\nabla b\|_\i^2)\|\tilde{\zeta}_\mathbf{n}^\Psi\|^2_{\mathcal{H}^{1,\i}}\leq\Lambda(\frac{1}{c_0},Y_m(t)).
\end{equation}
Finally, it is easy to obtain that
\begin{equation}\label{8.85}
	 \|(\varrho-1)\partial_t\tilde{\zeta}_\mathbf{n}^\Psi\|_\i^2\leq\Lambda(\frac{1}{c_0},Y_m(t).
\end{equation}
Substituting \eqref{8.83}-\eqref{8.85} into \eqref{8.82}, one gets that
\begin{align}\label{zze}
	\v\|\partial_z\tilde{\zeta}_\mathbf{n}^\Psi\|_\i^2 	 \leq\sqrt{\v}\|\partial_z\tilde{\zeta}_\mathbf{n}^\Psi(0)\|^2_\i
	+t\Lambda(\frac{1}{c_0},Y_m(t)+\v^2|\MZ^mh(t)|^2_{\f12}).
\end{align}
Now\eqref{21.1} follows from  \eqref{zze}, \eqref{zzv} and \eqref{4.6.86}.
Therefore, the proof of this lemma is completed.$\hfill\Box$

	\subsection{Proof of Theorem \ref{thm3.1} }
	
	\setcounter{equation}{0}

Based on the estimates obtained so far, we can complete the proof of Theorem \ref{thm3.1} in this subsection.First,   \eqref{3.0-3} follows easily from conservation of mass.

In order to prove \eqref{3.0-1}, one chooses $A>0$ large so that
	\begin{equation}\label{4.6.3}
	\partial_z\varphi(0,y,z)\geq1,
	\end{equation}
	where $A$ depends only on the initial data $|h(0)|_{H^3}$. For two parameters $R>0$ and $c_0$ to be chosen $0<\f{1}{c_0}\ll R$, define
	\begin{align}
	&T^\v_\star=\sup\Big\{T\in[0,\f{d_0}{8(1+R)}];  ~\varTheta_m(T)\leq R,~\partial_z\varphi(t)\geq c_0,~|h|_{\mathcal{H}^{3,\infty}}+|\nabla_y h|_{\mathcal{H}^{[\f{m}2]+1}}\leq\frac{1}{c_0},\nonumber\\
	 &~~~~~~~~~~~~~~~~~~~~~~~~~~~~~-\partial{_z^\varphi}p|_{z=0}\geq\frac{c_0}{2},~0<\f{1}{4C_0}\leq\varrho(t)\leq  4C_0,\forall t\in[0,T]\Big\}.
	\end{align}	
The restriction $T\leq \f{d_0}{8(1+R)}$ guarantees that the fluid moves a distance no more than $\f{1}{8}d_0$.  Then    $\tilde{\Om}_{0,1}$, $\tilde{\Om}_{k,\f18},~k=1,\cdots,n$ is still a valid covering of $\Om_t$ with $t\leq T$.
	
	 Now define
	\begin{align}\label{3.0-2}
	&\mathcal{N}_m(t)\triangleq\mathcal{N}_m(p,v,h)(t)\nonumber\\
	&=\sup_{0\leq \tau\leq t}\Big\{1+\|(V^m,Q^m)(\tau)\|^2+|(h,\sqrt\sigma\nabla_yh)(\tau)|^2_{\mathcal{H}^m}
	+\|(\nabla p, \nabla v)(\tau)\|{_{\mathcal{H}^{m-2}}^2}\nonumber\\
	&+\|\Delta^{\varphi}p(\tau)\|^2_{\mathcal{H}^1}+\|\nabla(p,v)(\tau)\|^2_{\mathcal{H}^{1,\infty}}+\v\|(\text{div}^\varphi v,\nabla^\varphi p,S_\mathbf{n})(\tau)\|{_{\mathcal{H}^{m-1}}^2}+\v|\MZ^mh|_\frac{1}{2}^2+\v\|\nabla^2 v(\tau)\|^2_{\mathcal{H}^{1,\infty}}\nonumber\\
	&+\v\|\nabla(p,v)(\tau)\|^2_{\mathcal{H}^{m-1}}+\v\|\Delta^{\varphi}p(\tau)\|^2_{\mathcal{H}^2}\Big\}+\int_0^t\|\nabla p(\tau)\|^2_{\mathcal{H}^{m-1}}d\tau
	+\Big(\int_0^t\|\nabla v(\tau)\|^4_{\mathcal{H}^{m-1}}d\tau\Big)^{\f12}\nonumber\\
	&+\v\int_{0}^{t}\|\nabla v(\tau)\|^2_{\mathcal{H}^{m}}d\tau+\v\int_{0}^{t}\|\nabla^2v(\tau)\|^2_{\mathcal{H}^{m-2}}d\tau
	+\v^2\int_{0}^{t}\|\nabla^2v(\tau)\|^2_{\mathcal{H}^{m-1}}d\tau,
	\end{align}
	which implies immediately that
	\begin{equation}
		 Y_m(t)+M(t)\leq\Lambda_0\Lambda(\mathcal{N}_m(t)),~\mathcal{N}_m(t)\leq \Lambda_0\Lambda(\varTheta_m(t))~
		\mbox{and}~\varTheta_m(t)\leq \Lambda_0\Lambda(\mathcal{N}_m(t)).
	\end{equation}
It follows from Propositions \ref{prop9.7}, \ref{prop7.6} that
\begin{align}\label{22.1}
&\|(V^m,Q^m)(t)\|^2+|(h,\sqrt\sigma\nabla_yh)(t)|^2_{\mathcal{H}^m}+\|(\text{div}^\varphi v,\nabla^\varphi p,S_\mathbf{n})(t)\|{_{\mathcal{H}^{m-2}}^2}+\v\|(\text{div}^\varphi v,\nabla^\varphi p,S_\mathbf{n})(t)\|{_{\mathcal{H}^{m-1}}^2}\nonumber\\
&+\v|\MZ^mh|_\frac{1}{2}^2+\int_0^t\|\nabla p\|_{\mathcal{H}^{m-1}}^2d\tau+\Big(\int_0^t\|\nabla v\|^4_{\mathcal{H}^{m-1}}d\tau\Big)^{\f12}+\v\int_{0}^{t}\|\nabla v(\tau)\|^2_{\mathcal{H}^{m}}+\|\nabla^2v(\tau)\|^2_{\mathcal{H}^{m-2}}
d\tau\nonumber\\
&+\v^2\int_{0}^{t}\|\nabla^2v(\tau)\|^2_{\mathcal{H}^{m-1}} d\tau\leq  \Lambda_0\Lambda(\mathcal{N}_m(0)+\|\nabla v(0)\|^2_{\mathcal{H}^{m-1}})+t^{\f12}\Lambda_0\Lambda(\mathcal{N}_m(t)).
\end{align}	
Then, as a consequence of \eqref{22.1},  Lemmas \ref{prop4.26}, \ref{lem4.29} and \ref{lem21.1}, one  obtains that
\begin{equation}\label{3.0-12}
\mathcal{N}_m(t)\leq  \Lambda(C_0, \mathcal{N}_m(0)+\|\nabla v(0)\|^2_{\mathcal{H}^{m-1}})+t^{\f12}\Lambda(C_0, \mathcal{N}_m(t)), ~~\forall t\in[0,T^\v].
\end{equation}
Therefore,  
\begin{equation}\label{3.0-11}
\varTheta_m(t)\leq  \Lambda(C_0, \varTheta_m(0)+\|\nabla v(0)\|^2_{\mathcal{H}^{m-1}})+t^{\f12}\Lambda(C_0, R), ~~\forall t\in[0,T^\v].
\end{equation}
Moreover,
\begin{align}\label{4.6.10}
\f{1}{C_0}\exp(T\Lambda(R))\leq \varrho(t)\leq C_0\exp(T\Lambda(R)),\forall t\in[0,T],
\end{align}
\begin{align}\label{4.6.11}
|h|^2_{\mathcal{H}^{3,\infty}}+|\nabla_y h|^2_{\mathcal{H}^{[\f{m}2]+1}}&\leq |h(0)|^2_{\mathcal{H}^{3,\infty}}+|\nabla_y h(0)|^2_{\mathcal{H}^{[\f{m}2]+1}}+T\Lambda(R),~\forall t\in[0,T],
\end{align}
\begin{align}\label{4.6.12}
&\partial_z\varphi(t)\geq 1-\int_0^t\|\partial_t\nabla\eta(t)\|_{L^\i}d\tau\geq 1-T \Lambda(R),~\forall t\in[0,T],
\end{align}
once the choice of $A>0$ so that \eqref{4.6.3} is satisfied. Finally, note that
\begin{align}\label{4.6.13}
-\partial{_z^\varphi}p(t)\geq-\partial{_z^\varphi}p_0-\int_0^t\|\partial_t\partial_z^\varphi p(\tau)\|d\tau
\geq -\partial{_z^\varphi}p_0-T\Lambda(R),~\forall t\in[0,T].
\end{align}
Taking $R=2\Lambda(C_0, \varTheta_m(0)+\|\nabla v(0)\|^2_{\mathcal{H}^{m-1}})$, in view of \eqref{3.0-11}-\eqref{4.6.13}, one gets that there exists $T_\ast>0$ depending only on $C_0$, $d_0$ and the initial data $\varTheta_m(0)+\|\nabla v(0)\|^2_{\mathcal{H}^{m-1}}$ (hence independent of $\v$ and $\s$) such that for $T\leq \min(T_\ast,T^\v)$
\begin{equation}
\begin{cases}
\varTheta_m(t)\leq 2\Lambda(C_0, \varTheta_m(0)+\|\nabla v(0)\|^2_{\mathcal{H}^{m-1}}),\\
\partial_z\varphi(t)\geq 2c_0,~|h(t)|_{\mathcal{H}^{3,\infty}}+|\nabla_y h(t)|_{\mathcal{H}^{[\f{m}2]+1}}\leq\frac{1}{2c_0},\\
-\partial{_z^\varphi}p|_{z=0}\geq\frac{3c_0}{4},~0<\f{1}{2C_0}\leq\varrho(t)\leq  2C_0,
\end{cases}
\forall t\in[0,T].
\end{equation}
Of course, it holds that $T_\ast\leq T^{\v}_{\star}$. Then taking $T_a=T_\ast$, one obtains \eqref{3.0-1} and  closes the a priori assumption \eqref{5.1}. Therefore, the proof Theorem \ref{thm3.1} is completed. $\hfill\Box$
	
	
\section{Proof of the Main Results }
	\renewcommand{\theequation}{\arabic{section}.\arabic{equation}}
	
\noindent\textbf{Proof of Theorem \ref{thm1.1}}:  We shall combine the a priori estimates obtained to prove the uniform existence result. Fix $m\geq6$, consider initial data $(p^\v_0, v^\v_0,h^\v_0)\in X^{\v,m}_{NS}$ such that
\begin{equation}\label{18.1}
\mathcal{I}_m(0):=\sup_{\v\in(0,1],\s\in[0,1]}\left(\|(p_0^\v,v_0^\v,h^\v_0)\|^2_{X^\v_m}+\|\nabla v_0^\v\|^2_{\mathcal{H}^{m-1}}\right)
\leq \tilde{C}_0,\nonumber
\end{equation}
and
\begin{equation}\label{18.3}
0<\f1{C_0}\leq \varrho_0^{\v}\leq C_0,
\end{equation}
For such initial data, we are not aware of a local well-posedness result for \eqref{1.1} and \eqref{1.4} and \eqref{1.5}. Thus, we shall prove the local existence by using energy estimates and a classical iteration scheme. Since $(p^\v_0, u^\v_0, F^\v_0)\in X^{\v,m}_{NS}$, there exists a sequence of smooth approximate  initial data $(p_0^{\v,\delta}, u_0^{\v, \delta},F_0^{\v,\d})\in X^{\v,m}_{NS,ap}$($\delta$ being a regularization parameter), which have enough spatial regularity so that the time derivatives at the initial time can be defined by Navier-Stokes equations and the boundary  compatibility conditions are satisfied. We construct  approximate solutions, inductively, as follows\\
(1) Define $u^0=u_0^{\v,\d}$, and \\
(2) Assuming that $u^{k-1}$ was defined for $k\geq1$, let $(\r^k,u^k)$ be the unique solution to the following linearized  initial boundary value problem:
\begin{eqnarray}\label{18.5}
\begin{cases}
\r^k_t+\mbox{div}(\r^k u^{k-1})=0 ~~\mbox{in}~(0,T)\times\Omega_t,\\
\r^k u^k_t+\r^k u^{k-1}\cdot\nabla{u}^k+\nabla{p}^k=\v\Delta u^k+\v\nabla\mbox{div}u^k~~\mbox{in}~(0,T)\times\Omega_t,\\
\partial_tF^k+u^{k-1}\cdot\nabla F^{k}=0,~\mbox{on}~\partial\Om_t,\\
(\r^k, u^k,F^k)|_{t=0}=(\r_0^{\v,\delta}, u_0^{\v, \delta},F_0^{\v,\d}),~~\mbox{with}~\f2{3C_0}\leq \r_0^{\v,\d}\leq \f32C_0,\\
\mbox{with boundary conditions }~~\eqref{1.5}.
\end{cases}
\end{eqnarray}

Since $\r^k$, $u^k$ and $F^k$ are decoupled, the existence of the global unique  smooth solution $(\r^k,u^k)$ of \eqref{18.5} with $0<\r^k<\infty$ can be obtained by using known results. For example, one can first obtain a local existence result for \eqref{18.5} by using local existence result in \cite{Secchi} without surface tension and  \cite{Zaja,Zaja-1} with surface tension. Since \eqref{18.5} is a linear problem, one can extend the local solution  to a global one.  Although the local existence result of \cite{Secchi} without surface tension and  \cite{Zaja,Zaja-1} with surface tension is proved in Sobolev-Slobodetskii spaces $W^{l,\f{l}2}$(see \cite{Zaja}) with $l\geq0$, but it is enough for us to proceed our {\it a priori}  estimates in Section 4 by taking $l=3m$.
	\
	
Similar to Section 1, let $(p^k, v^k,h^k)$ be the localized version $(p^k,u^k,F^k).$ 	Applying the a priori estimates given in Theorem \ref{thm3.1} and by an induction argument, we obtain a uniform time $T_1>0$ and constant $\tilde{C}_3>0$(independent of $\v,\s$ and  $\d$),  such that it holds for $(p^k, v^k,h^k),~k\geq 1$  that
	\begin{equation}\label{18.6}
	\varTheta_m(p^k,v^k,h^k)(t)\leq \tilde{C}_3,~\forall t\in[0,T_1],
	\end{equation}
	and
	\begin{equation}\label{18.7}
		\f1{2C_0}\leq\varrho^k(t)\leq 2C_0, ~~\forall t\in[0,T_1],
	\end{equation}
	where $T_1$ and $\tilde{C}_3$ depend only on $C_0$, $\tilde{C}_0$ and $d_0$. Based on the above uniform estimates for $(\r^k, u^k,F^k)$,   using Lagrangian coordinates and considering the difference of $(\r^{k+1}-\r^{k},u^{k+1}-u^{k},F^{k+1}-F^{k})$, one can prove, by using the energy methods,  that  there exists an uniform time $T_2(\leq T_1)$(independent of $\v,\s$ and  $\d$) such that $(\r^k, u^k,F^k)$ converges to a limit $(\r^{\v,\d}, u^{\v,\d},F^{\v,\d})$ as $k\rightarrow +\infty$ in the following strong sense:
\begin{equation}\nonumber
\begin{cases}
&(\r^k, u^k)\rightarrow (\r^{\v,\d}, u^{\v,\d})~~\mbox{in}~~L^\infty(0,T_2; L^2),\\
&\nabla u^k\rightarrow \nabla u^{\v,\d}~\mbox{in}~ L^2(0,T_2, L^2),\\
&F^k(t,x(\cdot,\cdot))\rightarrow F^{\v,\d}(t,x(\cdot,\cdot))~\mbox{in}~L^2(U).
\end{cases}
\end{equation}
Then, it is easy to check $(\r^{\v,\d}, u^{\v,\d}, F^{\v,\d})$ is a weak solution to the problem \eqref{1.1}, \eqref{1.4} and \eqref{1.5} with initial data $(\r_0^{\v,\delta}, u_0^{\v, \delta},F_0^{\v,\d})$. Then, by virtue of the lower semi-continuity of norms, one can deduce from the uniform bounds \eqref{18.6} and \eqref{18.7} that $(\r^{\v,\d}, u^{\v,\d},F^{\v,\d})$ satisfies the following regularity estimates
	\begin{equation}\label{18.8}
		\varTheta_m(p^{\v,\d}, v^{\v,\d},h^{\v,\d})(t)\leq \tilde{C}_3,~\forall t\in[0,T_2],
	\end{equation}
	and
	\begin{equation}\label{18.9}
		\f1{2C_0}\leq\varrho^{\v,\d}(t)\leq 2C_0, ~~\forall t\in[0,T_2],
	\end{equation}
	Based on the uniform estimates \eqref{18.8} and \eqref{18.9} for $(\r^{\v,\d}, u^{\v,\d}, F^{\v,\d})$, we pass the limit $\d\rightarrow0$ to get a strong solution $(\r^\v,u^\v,F^\v)$ to \eqref{1.1}, \eqref{1.4} and \eqref{1.5} with initial data $(\r_0^\v, u_0^\v,F^\v_0)$ satisfying  \eqref{18.1} by using a strong compactness arguments. Since the compactness arguments are almost the same as   the ones need for the proof of Theorem \ref{thm1.2} below, we shall not give more details here.  Moreover, the uniform bounds \eqref{1.18} and \eqref{1.19} are immediately results of the lower semi-continuity of norms. Changing the variable into Lagrangian coordinate and using energy method, it is easy to  prove the uniqueness of  $(\r^\v,u^\v,F^\v)$   since we work on functions with enough  regularity. Taking $T_0=T_2$ and $\tilde{C}_1=\tilde{C}_3$,  we complete the proof of Theorem \ref{thm1.1}. $\hfill\Box$
	
	\
	
	\

	\noindent\textbf{Proof of Theorem \ref{thm1.2}}: For any fixed $\s\geq 0$, it follows from Theorem \ref{thm1.1}, for any $\v\in(0,1]$,  that
	\begin{equation}\label{19.1}
		\varTheta_m(p^{\v}, v^{\v},h^\v)(t)\leq \tilde{C}_1,~\forall t\in[0,T_0],
	\end{equation}
	and
	\begin{equation}\label{19.2}
		\f1{2C_0}\leq\varrho^{\v}(x,t)\leq 2C_0, ~~\forall t\in[0,T_0].
	\end{equation}
	Thus, $(p^\v,v^\v)$ is uniform bounded in $L^\infty([0,T];H^m_{co})$, $\nabla(p^\v,v^\v)$ is uniform bounded in $L^\infty([0,T];H^{m-2}_{co})$, $h^\v$ is uniform bounded in $L^\infty(0,T;H^m)$,  $\partial_t(p^\v,v^\v)$ is uniform bounded in $L^\infty([0,T];L^2)$ and $\partial_th^\v$ is uniform bounded in $L^\infty(0,T;L^2)$. Then, it follows from the strong compactness arguments(see \cite{Simon}) that $(p^\v,v^\v)$ is compact in $C([0,T];H^{m-2}_{co})$ and $h^\v$ is compact in $C([0,T];H^{m-1})$. In particular, there exist a sequence $\v_n\rightarrow 0+$ and $(p,v)\in C([0,T];H^{m-2}_{co})$ and $h\in C([0,T];H^{m-1})$ such that
	\begin{equation}\label{19.3}
	\begin{cases}
	(p^{\v_n},v^{\v_n})\rightarrow (p,v)~~\mbox{in}~~C([0,T];H^{m-2}_{co}),\\
	h^\v\rightarrow h ~~\mbox{in}~~C([0,T];H^{m-1}),
	\end{cases}
	\mbox{as}~~\v_n\rightarrow 0+,
	\end{equation}
Denote $(\r,u,F)$ the global version of $(\varrho,v,h)$. Applying the lower semi-continuity of norms to   the bounds \eqref{19.1} and \eqref{19.2}, we obtain the uniform  estimates \eqref{1.18-1} and \eqref{1.19-1}. Then, based on \eqref{19.3} and  the uniform estimates  \eqref{1.18-1}, \eqref{1.19-1},  it is easy to check that $(\r, u, F)$ is a solution of the Euler equations with the free surface
\begin{eqnarray}\label{E1-1}
\begin{cases}
\partial_t\rho+\mbox{div}(\rho u)=0,\\
\rho \partial_tu+\rho (u\cdot\nabla)u+\nabla p=0,
\end{cases}
x\in\Om_t,~t\in(0,T_0]
\end{eqnarray}
with the boundary conditions
\begin{equation}\label{E2-1}
\partial_tF+u\cdot \nabla F=0,~~~\mbox{and}~~p=p_e-\sigma H,~~x\in\partial\Om_t,
\end{equation}
where $\Om_t$ is the domain occupied by the fluid on time  $t\geq0$ and the boundary of $\Om_t$ is given by
$$\partial\Om_t=\{x\in\mathbb{R}^3|~F(x,t)=0 \}.$$

	By using \eqref{1.18-1}, \eqref{19.3} and  the anisotropic Sobolev inequality \eqref{3.3},  we have that
	\begin{equation*}\label{9.4}
		\sup_{t\in[0,T]}\|(\varrho^{\v_n}-\r,v^{\v_n}-v)\|^2_{L^\infty}\leq \sup_{t\in[0,T]}\Big( \|\nabla(\varrho^{\v_n}-\varrho,v^{\v_n}-v)\|_{H^1_{co}}\cdot\|(\varrho^{\v_n}-\varrho,v^{\v_n}-v)\|_{H^2_{co}}\Big)\rightarrow 0.
	\end{equation*}
	and
	$$|h^\v-h|_{W^{1,\i}}\leq |h^\v-h|_{H^{3}}\rightarrow0,$$
	as $\v_n\rightarrow0+$. Hence, \eqref{1.20} is proved.

Consider $(\r^i,u^i,F^i),~i=1,2$ two solutions of \eqref{E1-1}, \eqref{E2-1} defined on $[0,T]$ with the same initial data and the regularity stated in \eqref{1.18-1} and \eqref{1.19-1}. We shall prove that $(\r^1,u^1,F^1)=(\r^2,u^2,F^2)$. Applying the standard energy method to the difference $(\r^1,u^1,F^1)-(\r^2,u^2,F^2)$ in the Lagrangian coordinates and using the uniform bounds in \eqref{1.18-1} and \eqref{1.19-1}. Therefore,  $(\r,u,F)$ is a unique solution to Euler system with free surface \eqref{E1-1}-\eqref{E2-1}.
And the uniqueness yields immediately that  the whole family $(\r^\v,u^\v,F^\v)$ converges to $(\r,u,F)$. Thus,  the proof Theorem \ref{thm1.2} is completed.  $\hfill\Box$
	
	\

	\noindent\textbf{Proof of Theorem \ref{thm1.2-1}}: Based on the uniform in both $\v$ and $\s$ estimates in Theorem \ref{thm1.1}, one can prove Theorem \ref{thm1.2-1} by using the similar strong compactness argument in Theorem \ref{thm1.2}. The details are omitted here for brevity. $\hfill\Box$

	\section{Appendix}

	Let $S(t,\tau)$ be the $C^0$ evolution operator generated by the following equation
	\begin{equation}\label{A.1}
	[\partial_th+b_1(t,y)\partial_{y^1}h+b_2(t,y)\partial_{y^2}h+z b_3(t,y)\partial_zh]-\v d(t,y)\partial_{zz}h=0,~z>0,~t>\tau,
	\end{equation}
	with the boundary condition $h(t,y,0)=0$ and with the initial condition $h(\tau,y,z)=h_0(y,z)$. The coefficients are smooth and $d(t,y)$ satisfies
	\begin{equation}\label{A.1-1}
		c_1\leq d(t,y)\leq \f1{c_1}
	\end{equation}
	for some positive constant $c_1>0$.
	
	\
	
	Then we have the following Lemnas whose proof can be found in \cite{Wang-Xin-Yong}, the details are omitted here. 
	\begin{lemma}\label{lemA.1}
		It holds that, for $t\geq \tau\geq 0$
		\begin{align}
			&\|S(t,\tau)h_0\|_{L^\infty}\leq \|h_0\|_{L^\infty},\label{A.2}\\
			&\|z\partial_zS(t,\tau)h_0\|_{L^\infty}\leq C(\|h_0\|_{L^\infty}+\|z\partial_zh_0\|_{L^\infty}),\label{A.2-1}
		\end{align}
		where  $C>0$ is a uniform constant independent of the bound of  $d(t,y)$ and $b_j(t,y),~j=1,2,3$.
	\end{lemma}

\begin{lemma}\label{lemA.2}
		Consider $h$ a smooth solution of
		\begin{equation}\label{A.13}
			 a(t,y)[\partial_th+b_1(t,y)\partial_{y^1}h+b_2(t,y)\partial_{y^2}h+z b_3(t,y)\partial_zh]-\v\partial_{zz}h=G,~z>0,~h(t,y,0)=0,
		\end{equation}
		for some smooth function $a(t,y)$ satisfies \eqref{A.1-1} and vector fields $b=(b_1,b_2,b_3)^t(t,y)$. Assume that $h$ and $G$ are compactly supported in $z$. Then, one has the estimate:
		\begin{align}\label{A.14}
			\|h\|_{\mathcal{H}^{1,\infty}}&\lesssim \|h_0\|_{\mathcal{H}^{1,\infty}}
			+\int_0^t\|\f{1}{a}\|_{L^\infty} \|G\|_{\mathcal{H}^{1,\infty}}d\tau\nonumber\\
			 &~~~~+\int_0^t(1+\|\f{1}{a}\|_{L^\infty})(1+\sum_{i=0}^2\|Z_i(a,b)\|^2_{L^\infty}) \|h\|_{\mathcal{H}^{1,\infty}}d\tau.
		\end{align}
		
	\end{lemma}

	\noindent {\bf Acknowledgments:}
	
Yong Wang is partially supported by National Natural Sciences Foundation of China No. 11371064 and 11401565. This research is support in part by Zheng Ge Ru Funds, Hong Kong RGC Earmarked Research Grants CUHK4041/11P, and CUHK4048/13P, NSFC/RGC Joint Research Scheme Grant N. CUHK443/14, A grant from the Croucher Foundation, and a Focus Area Grant at The Chinese University of Hong Kong.

	\

\end{document}